\documentclass{amsart}
\usepackage[francais]{babel}
\usepackage{a4}
\usepackage{amsfonts}
\usepackage{amssymb}
\usepackage[latin1]{inputenc}
\DeclareMathOperator{\Gal}{Gal}
\DeclareMathOperator{\Ind}{Ind}
\DeclareMathOperator{\card}{card}
\DeclareMathOperator{\Diag}{Diag}
\DeclareMathOperator{\pgcd}{pgcd}
\author{François Courtès} 
\address{Université de Poitiers\\SP2MI - Département de Mathématiques\\Téléport 2 - Boulevard Marie et Pierre Curie\\86962 Futuroscope Chasseneuil Cedex}
\email{courtes@mathlabo.univ-poitiers.fr}
\title[Distributions invariantes]{Distributions invariantes sur les groupes réductifs quasi-déployés}
\begin{document}
\renewcommand{\baselinestretch}{1.5}
\normalsize
\begin{abstract}
Soit $F$ un corps local non archimédien, et $G$ le groupe des $F$-points d'un groupe réductif connexe quasi-déployé défini sur $F$; dans cet article, on s'intéresse aux distributions sur $G$ invariantes par conjugaison, et à l'espace de leurs restrictions à l'algèbre de Hecke $\mathcal{H}$ des fonctions sur $G$ à support compact biinvariantes par un sous-groupe d'Iwahori $I$ donné. On montre tout d'abord que les valeurs d'une telle distribution sur $\mathcal{H}$ sont entièrement déterminées par sa restriction au sous-espace de dimension finie des éléments de $\mathcal{H}$ à support dans la réunion des sous-groupes parahoriques de $G$ contenant $I$; on utilise ensuite cette propriété pour montrer, moyennant certaines conditions sur $G$, que cet espace est engendré d'une part par certaines intégrales orbitales semi-simples, d'autre part par les intégrales orbitales unipotentes, en montrant tout d'abord des résultats analogues sur les groupes finis.
\end{abstract}
\maketitle
\tableofcontents

\def \mth{\mathbb}
\newtheorem{theo}{Théorème}[section]
\newtheorem{prop}[theo]{Proposition}
\newtheorem{lemme}[theo]{Lemme}
\newtheorem{cor}[theo]{Corollaire}
\newcommand{\module}[1]{\left\vert#1\right\vert}
\newcommand{\norme}[1]{\left\Vert#1\right\Vert}

\section{Introduction}

Soit $F$ un corps local non archimédien, $\mathcal{O}$ son anneau des entiers, ${\mathfrak{p}}$ l'idéal maximal de $\mathcal{O}$, $\varpi$ une uniformisante de $F$; soit $p$ la caractéristique résiduelle et $q$ le cardinal du corps résiduel $\mathcal{O}/{\mathfrak{p}}$ de $F$.

Soit $\underline{G}$ un groupe réductif connexe défini sur $F$; on supposera $\underline{G}$ quasi-déployé sur $F$. Soit $G=\underline{G}\left(F\right)$ le groupe des $F$-points de $\underline{G}$.

Dans ce qui suit, on s'intéressera aux distributions sur $G$ invariantes par conjugaison et à support dans l'ensemble des éléments de $G$ compacts, c'est-à-dire fixant au moins un point de l'immeuble de $G$. Soit $\mathcal{D}_c$ l'espace de ces distributions, et $\mathcal{D}_{c,1}$ le sous-espace des éléments de $\mathcal{D}_c$ à support dans la réunion des sous-groupes parahoriques de $G$, c'est-à-dire des fixateurs connexes (au sens de Bruhat-Tits) des facettes de l'immeuble de $G$. Considérons d'autre part l'algèbre de Hecke $\mathcal{H}=\mathcal{H}_I$ des fonctions sur $G$ à support compact et biinvariantes par un sous-groupe d'Iwahori $I$ fixé. D'après la conjecture de Howe démontrée par Clozel, on sait que, dans le cas de la caractéristique 0, l'espace des restrictions des éléments de $\mathcal{D}_{c,1}$ à cette algèbre est de dimension finie, et c'est également vrai en caractéristique $p$ grâce à la théorie des corps locaux proches; on peut donc chercher à en exhiber des systèmes de générateurs finis.

Si $g$ est un élément compact de $G$, en fixant une mesure invariante par conjugaison sur l'orbite de $g$ dans $G$, on peut définir une distribution intégrale orbitale associée à $g$, qui est un élément de $\mathcal{D}_c$. On va considérer deux familles particulières d'éléments compacts de $G$:
\begin{itemize}
\item les éléments de $G$ semi-simples non ramifiés de réduction régulière, c'est-à-dire les éléments $g$ appartenant à un sous-groupe parahorique $K$ de $G$ dont l'image dans le groupe fini correspondant est un élément semi-simple régulier. De tels éléments existent pour $q$ assez grand;
\item les éléments unipotents de $G$.
\end{itemize}

Dans \cite[3.2]{wsp2}, J.L. Waldspurger a conjecturé que:
\begin{itemize}
\item la restriction à $\mathcal{H}$ de l'intégrale orbitale d'un élément semi-simple non ramifié de réduction régulière ne dépend que de la classe de conjugaison dans $G$ de l'unique tore non ramifié maximal le contenant;
\item l'espace engendré par ces restrictions coïncide avec l'espace engendré par les restrictions à $\mathcal{H}$ des intégrales orbitales unipotentes sur $G$.
\end{itemize}

Le but de ce qui suit est de montrer cette conjecture moyennant certaines conditions sur $G$; on va même montrer que d'une part les intégrales orbitales des éléments semi-simples non ramifiés de réduction régulière de $G$, d'autre part les intégrales orbitales unipotentes, engendrent l'espace des restrictions à $\mathcal{H}$ de $\mathcal{D}_{c,1}$ tout entier.

On va chercher à se ramener au cas des groupes finis. L'objet du chapitre $3$ est de montrer que les restrictions à $\mathcal{H}$ des éléments de $\mathcal{D}_c$ sont entièrement déterminées par des distributions sur un nombre fini de groupes finis (éventuel\-lement tordus par un groupe abélien libre de type fini).
Pour cela, on aura besoin de quelques préliminaires.

Soit $\Gamma=G/G^1$, où $G^1$ est le sous-groupe de $G$ engendré par les sous-groupes parahoriques; $\Gamma$ est canoniquement isomorphe à $N_G\left(I\right)/I$, et si $K$ est un sous-groupe parahorique de $K$, $\Gamma_K=N_G\left(K\right)/K$ est canoniquement isomorphe à un sous-groupe de $\Gamma$ (lemme \ref{normpar}). Soit $\mathcal{M}$ l'ensemble des classes de conjugaison dans $G$ de couples $\left(M,K_\mu\right)$, où $M$ est un sous-groupe de Levi de $G$ et $K_\mu$ un sous-groupe parahorique maximal de $M$; on peut associer de manière canonique à chaque sous-groupe parahorique de $G$ un élément de $\mathcal{M}$; si $\left(M,K_\mu\right)$ est un représentant de $\mu$, le lemme \ref{normpar} assure de plus que les groupes ${\mth{G}}=K_\mu/K_\mu^1$ et ${\mth{G}}^+=\gamma^{\mth{Z}}K_\mu/K_\mu^1$, où $K_\mu^1$ est le premier sous-groupe de congruence de $K_\mu^1$, ne dépendent pas de $K$. D'autre part, pour tout $\mu\in\mathcal{M}$, si $\left(M,K_\mu\right)$ est un représentant de $\mu$, le sous-groupe $\Gamma_\mu=N_G\left(K_\mu\right)/N_{G^1}\left(K_\mu\right)$ de $\Gamma$ ne dépend pas de $\left(M,K_\mu\right)$, et pour tout $K$ d'image $\mu$, $\Gamma_K$ est canoniquement isomorphe à un sous-groupe de $\Gamma_\mu$; si $\mathcal{N}$ est l'ensemble des couples $\left(\gamma,\mu\right)$, avec $\mu\in\mathcal{M}$ et $\gamma\in\Gamma_\mu$, on peut donc associer à chaque couple $\left(\gamma,K\right)$, où $K$ est un parahorique de $G$ et $\gamma\in \Gamma_K$, un élément de $\mathcal{N}$. 

Soit $D\in\mathcal{D}_c$. Pour tout $\nu=\left(\gamma,\mu\right)\in\mathcal{N}$ et tout sous-groupe parahorique $K$ de $G$ d'image $\nu$, la restriction de $D$ à l'espace des fonctions à support dans $\gamma K$ invariantes par le premier sous-groupe de congruence $K^1$ de $K$ s'identifie à une distribution invariante sur ${\mth{G}}^+$ à support dans $\gamma{\mth{G}}$, ou encore à une fonction $\phi$ sur ${\mth{G}}$ invariante par $\gamma$-conjugaison. Soit $C\left({\mth{G}}\right)^{{\mth{G}},\gamma}$ l'espace de ces fonctions et $C_{cusp}\left({\mth{G}}\right)^{{\mth{G}},\gamma}$ le sous-espace des fonction $\gamma$-cuspidales; si $R$ est un système de représentants des sous-groupes de Levi de ${\mth{G}}$ stables par $\gamma$, on a le résultat de décomposition suivant (proposition \ref{dcc}):

\begin{prop}
\[C\left({\mth{G}}\right)^{{\mth{G}},\gamma}=\bigoplus_{{\mth{M}}\in R}\Ind_{{\mth{M}},\gamma}^{\mth{G}}C_{cusp}\left({\mth{M}}\right)^{{\mth{M}},\gamma}\]
\end{prop}

Dans l'égalité ci-dessus, si ${\mth{P}}={\mth{M}}{\mth{U}}$ est un sous-groupe parabolique de ${\mth{G}}$ stable par $\gamma$, $\Ind_{{\mth{M}},\gamma}^{\mth{G}}$ désigne l'induite tordue par $\gamma$, définie par:
\[\Ind_{{\mth{M}},\gamma}^{\mth{G}}\left(f\right):g\in{\mth{G}}\longmapsto\dfrac 1{\card\left({\mth{P}}\right)}\sum_{x\in{\mth{G}},m\in{\mth{M}},u\in{\mth{U}},g=\gamma^{-1}x\gamma mux}f\left(m\right);\]
on déduit de la proposition \ref{resind} que cette induite ne dépend que de ${\mth{M}}$ et pas de ${\mth{P}}$.

Soit $\phi_\nu\left(D\right)$ la composante dans $C_{cusp}\left({\mth{M}}\right)^{{\mth{M}},\gamma}$ de $\phi$; elle ne dépend pas du choix de $K$ (lemme \ref{reges}). Soit $T_0$ un tore déployé maximal de $G$ dont l'unique sous-groupe parahorique $K_{T_0}$ est contenu dans $I$, et soit $M_\gamma$ le plus petit sous-groupe de Levi de $G$ semi-standard relativement à $T_0$ et contenant à la fois $M$ et $\gamma$, et soit l'intégrale:
\[D_{\phi_\nu\left(D\right),\gamma}\left(f\right)=\int_{Z_{M_\gamma}^0\backslash M_\gamma}\left(\int_{M_\gamma}f\left(m^{-1}ym\right)\phi_\nu\left(D\right)\left(\gamma ^{-1}y\right)dy\right)dm.\]
Cette intégrale converge pour tout $f$ (corollaire \ref{cspd2}), et définit une distribution invariante sur $M_\gamma$ à support dans l'ensemble des normalisateurs de sous-groupes parahoriques de $M_\gamma$; on en déduit la distribution invariante sur $G$ suivante:
\[D_{\phi_\nu\left(D\right),\gamma}^G:f\longmapsto D_{\phi_\nu\left(D\right),\gamma}\left(f^{P_\gamma}\right),\]
où $P_\gamma$ est un sous-groupe parabolique de $G$ de Levi $M_\gamma$, et $f^{P_\gamma}$ est le terme constant de $f$ suivant $P_\gamma$.
Posons maintenant:
\[\tilde{D}=\sum_{\nu\in\mathcal{N}}D_{\overline{\phi_\nu\left(D\right)}}^G;\]
pour les $\nu\in\mathcal{N}$ qui ne sont associés à aucun parahorique de $G$, on posera $D_{\overline{\phi_\nu\left(D\right)}}^G=0$.
Soit $\tilde{\mathcal{D}_c}$ le sous-espace de $\mathcal{D}_c$ constitué des distributions de la forme ci-dessus, et soit $\mathcal{C}_0$ l'espace engendré par les fonctions sur $G$ à support dans le normalisateur d'un sous-groupe parahorique et biinvariantes par le premier sous-groupe de congruence de ce même parahorique; on a (proposition \ref{dprime}):

\begin{prop}
Soit $D\in\mathcal{D}_c$. Alors il existe un unique $D'\in\tilde{\mathcal{D}_c}$ tel que les restrictions de $D$ et $D'$ à $\mathcal{C}_0$ sont identiques; de plus, on a $D'=\tilde{D}$.
\end{prop}

Fixons maintenant un sous-groupe parahorique $I$ de $G$ et considérons le sous-espace $\mathcal{C}$ de $C_c^\infty\left(G\right)$ engendré par les fonctions à support compact biinvariantes par le premier sous-groupe de congruence d'un sous-groupe parahorique de $G$ contenant $I$. On a (théorème \ref{dct}):

\begin{theo}
Pour tout $D\in\mathcal{D}_c$, $\tilde{D}$ coïncide avec $D$ sur $\mathcal{C}$. En particulier, si $\tilde{D}=0$, $D$ est nulle sur $\mathcal{C}$.
\end{theo}

Soit $\mathcal{C}_0$ le sous-espace de $\mathcal{C}$ engendré par les fonctions à support dans le normalisateur d'un sous-groupe parahorique de $G$ contenant $I$ et biinvariantes par le premier sous-groupe de congruence de ce même sous-groupe parahorique, $\tilde{D}$ ne dépend que de la restriction de $D$ à $\mathcal{C}_0$, et on déduit alors immédiatement du théorème précédent le résultat suivant (\ref{dc}):

\begin{theo}
Si $D$ est nulle sur $\mathcal{C}_0$, alors elle est nulle sur $\mathcal{C}$.
\end{theo}

En particulier, si le support de $D$ est contenu dans une classe de $G$ modulo $G^1$, ses valeurs sur $\mathcal{C}$ sont entièrement déterminées par ses valeurs sur un espace de dimension finie. (Si $G$ est semi-simple, c'est même vrai pour $D$ quelconque.)

(Remarque: Il devrait être possible de généraliser ce résultat aux espaces de fonctions de niveau $n$ quelconque.)

On va utiliser ce résultat pour montrer son analogue en remplaçant $\mathcal{C}$ par $\mathcal{H}$.
Considérons donc le sous-espace $\mathcal{H}_0$ des éléments de $\mathcal{H}$ à support dans la réunion des normalisateurs des sous-groupes parahoriques de $G$ contenant $I$. On a le résultat suivant (théorème \ref{dh}):

\begin{theo}
Si $D$ est nulle sur $\mathcal{H}_0$, alors elle est nulle sur $\mathcal{H}$.
\end{theo}

Pour montrer les résultats cherchés (chapitres $4$ et $5$), on se ramène donc au cas des groupes finis. Soit ${\mth{G}}$ le groupe des ${\mth{F}}_q$-points d'un groupe réductif connexe défini sur ${\mth{F}}_q$, ${\mth{P}}_0$ un sous-groupe parabolique minimal de ${\mth{G}}$ et $\mathcal{H}_{\mth{G}}$ l'algèbre de Hecke des fonctions sur ${\mth{G}}$ biinvariantes par ${\mth{P}}_0$. Pour toute fonction $f$ sur ${\mth{G}}$ et tout $g\in{\mth{G}}$, posons:
\[1_g\left(f\right)=\sum_{h\in{\mth{G}}}f\left(h^{-1}gh\right);\]
c'est la somme de $f$ sur l'orbite de $g$. On a les résultats suivants (proposition \ref{intorbf}, lemme \ref{ssuf} et corollaire \ref{ssufc}):

\begin{prop}
\begin{itemize}
\item Si $g$ est un élément régulier d'un tore maximal ${\mth{T}}$ de ${\mth{G}}$, la restriction à $\mathcal{H}_{\mth{G}}$ de la distribution $1_g$ ne dépend que de la classe de conjugaison de ${\mth{T}}$ dans ${\mth{G}}$.
\item Soit $R_{\mth{G}}$ un système de représentants des classes de conjugaison de tore maximaux de ${\mth{G}}$; supposons que pour tout ${\mth{T}}\in R_{\mth{G}}$, il existe un élément régulier $g_{\mth{T}}$ dans ${\mth{T}}$. Si le diagramme de Dynkin de ${\mth{G}}$ ne contient aucune composante connexe de type $E_7$ ou $E_8$, alors les distributions $1_{g_{\mth{T}}}|_{\mathcal{H}_{\mth{G}}}$, $T\in R_G$, engendrent l'espace des distributions invariantes sur $\mathcal{H}_{\mth{G}}$, et en constituent une base si ${\mth{G}}$ est déployé.
\end{itemize}
\end{prop}

\begin{lemme}
Soit ${\mth{T}}$ un tore maximal de ${\mth{G}}$, $g$ un élément régulier de $T$; alors la restriction à $\mathcal{H}_{\mth{G}}$ de l'intégrale orbitale $1_g$ est combinaison linéaire des restrictions des intégrales orbitales unipotentes sur ${\mth{G}}$.
\end{lemme}

\begin{cor}
Si ${\mth{G}}$ est tel que tout tore maximal ${\mth{T}}$ de ${\mth{G}}$ admet au moins un élément régulier, et si aucune composante connexe du diagramme de Dynkin de ${\mth{G}}$ n'est de type $E_7$ ou $E_8$, alors les intégrales orbitales unipotentes engendrent l'espace des restrictions à $\mathcal{H}_{\mth{G}}$ des distributions sur ${\mth{G}}$.
\end{cor}

On va ensuite utiliser ces résultats pour montrer leur analogue sur $G$. Si $K$ est un sous-groupe parahorique de $G$ contenant $I$, $K/K^1$ s'identifie au groupe des ${\mth{F}}_q$-points d'un groupe réductif connexe défini sur ${\mth{F}}_q$, et $I/K^1$ est un sous-groupe parabolique minimal de ce groupe; on peut donc appliquer les résultats précédents à chacun des $K$ contenant $I$. En recollant les distributions ainsi obtenues et en utilisant le théorème \ref{dh}, on démontre (théorèmes \ref{intorb} et \ref{intun}, corollaire \ref{intuncor}): 

\begin{theo}
\begin{itemize}
\item Si $g$ est un élément non ramifié de réduction régu\-lière d'un tore non ramifié maximal $T$ de $G$, la restriction à $\mathcal{H}$ de la distribution $J\left(g,\cdot\right)$ ne dépend que de la classe de conjugaison de $T$ dans $G$.
\item Supposons que $G$ vérifie:

{\em (C1): Tout tore maximal non ramifié $T$ de $G$ admet au moins un élément non ramifié de réduction régulière.}

Soit $R$ un système de représentants des classes de conjugaison de tores non ramifiés maximaux de $G$, et pour tout $T\in R$, soit $g_T$ un élément non ramifié de réduction régulière de $T$. Si le diagramme de Dynkin de $G$ ne contient aucune composante connexe de type $E_7$ ou $E_8$, alors les distributions $J\left(g_T,\cdot\right)|_\mathcal{H}$ engendrent $\mathcal{D}_{c,1}|_\mathcal{H}$, et en constituent une base si $G$ est déployé sur $F$.
\end{itemize}
\end{theo}

Notons que (C1) est vraie pour $q$ assez grand (lemme \ref{srg}).

Soit $F_{nr}$ l'extension non ramifiée maximale de $F$, et $G_{nr}=\underline{G}\left(F_{nr}\right)$. On dira que $p$ est un bon nombre premier pour $G_{nr}$ si, en notant $\Phi_{nr}$ le système de racines de $G_{nr}$ relativement à un tore déployé maximal et $\Delta=\left\{\alpha_1,\dots,\alpha_n\right\}$ un système de racines simples de $\Phi_{nr}$, pour tout $\beta=\sum_ic_i\alpha_i\in\Phi_{nr}$, si $c_i\neq 0$, $c_i$ n'est pas multiple de $p$. On a:

\begin{theo}
Supposons que $G$ vérifie (C1), que son diagramme de Dynkin ne contient aucune composante connexe de type $E_7$ ou $E_8$,et que $p$ est un bon nombre premier pour $G_{nr}$. Alors les restrictions à $\mathcal{H}$ des distributions intégrales orbitales unipotentes sur $G$ engendrent $\mathcal{D}_{c,1}|_\mathcal{H}$.
\end{theo}

\begin{cor}
Si $G$ vérifie les conditions du théorème précédent, les espaces engendrés par les restrictions à $\mathcal{H}$ des distributions intégrales orbitales respectivement sur les éléments semi-simples réguliers non ramifiés de $G$ et sur les éléments unipotents de $G$ sont égaux.
\end{cor}

\section{Préliminaires}

Dans ce qui suit, si $H$ est un groupe et $H'$ un sous-groupe de $H$, on notera $N_H\left(H'\right)$ (resp. $Z_H\left(H'\right)$) le normalisateur (resp. le centralisateur) de $H'$ dans $H$. Si $g$ est un élément de $H$, on notera également $Z_H\left(g\right)$ le centralisateur de $g$ dans $H$. Enfin, on notera $Z_H=Z_H\left(H\right)$ le centre de $H$.

\subsection{Sous-groupes parahoriques}

Soit $\mathcal{B}=\mathcal{B}_G$ l'immeuble de Bruhat-Tits de $G$, $\mathcal{A}=\mathcal{A}_G$ un appartement de $\mathcal{B}$. Soit $N_0$ l'ensemble des points de $G$ conservant $\mathcal{A}$; c'est un sous-groupe fermé de $G$, dont la composante neutre $H_0$ est un tore maximal de $G$. Soit $T_0$ la composante déployée de $H_0$; c'est un tore déployé maximal de $G$, et on a $N_0=N_G\left(T_0\right)$ et $H_0=Z_G\left(T_0\right)$. De plus, $\mathcal{A}$ est un ${\mth{R}}$-espace affine dont la dimension est égale au rang semi-simple de $G$, et le ${\mth{R}}$-espace vectoriel $V$ associé à $\mathcal{A}$ est canoniquement isomorphe à $X_*\left(T_0\right)\otimes{\mth{R}}/X_*\left(Z_{G,d}\right)\otimes{\mth{R}}$, où $Z_{G,d}$ est la composante déployée du centre de $G$.

Soit $x$ un point de $\mathcal{A}$, et soit $K_x$ le sous-groupe parahorique de $G$ attaché à $x$; 
c'est le fixateur connexe (cf. \cite[II. 4.6.28]{bt}) de la facette $A_x$ de $\mathcal{B}$ contenant $x$.
Si $\mathcal{A}'$ est un autre appartement de $\mathcal{B}$ contenant $x$, d'après \cite[I. 2.5.8]{bt}, il existe un élément $g$ de $G$ tel que $g\left(\mathcal{A}\right)=\mathcal{A}'$ et que $gx=x$; le tore déployé maximal $T'_0$ de $G$ correspondant à $\mathcal{A}'$ est alors conjugué à $T_0$ par un élément de $N_G\left(K_x\right)$.

Soit $\Phi_{T_0}$ le système de racines de $G$ relativement à $T_0$, et soit, pour tout $\alpha\in\Phi_{T_0}$, le sous-groupe radiciel $U_\alpha$ de $G$ correspondant; on notera $U_{\alpha,x}=U_{\alpha,K_x}=U_\alpha\cap K_x$. On va définir le premier sous-groupe de congruence $K_x^1$ de $K_x$ comme suit: puisque $T_0$ est un tore déployé maximal et $G$ est quasi-déployé, $Z_G\left(T_0\right)$ est un tore; son immeuble de Bruhat-Tits est alors constitué d'un unique appartement, lui-même  constitué d'une unique facette, et $Z_G\left(T_0\right)$ possède donc un unique sous-groupe parahorique $K_{Z_G\left(T_0\right)}$.
On définit tout d'abord le premier sous-groupe de congruence $K_{Z_G\left(T_0\right)}^1$ de $K_{Z_G\left(T_0\right)}$ de la manière suivante: soit $\underline{T'}$ le tore maximal de $\underline{G}$ tel que $Z_G\left(T_0\right)=\underline{T'}\left(F\right)$. Soit $\overline{F}$ la clôture algébrique de $F$, $\overline{\mathcal{O}}$ son anneau des entiers, $\overline{\mathfrak{p}}$ l'idéal maximal de $\overline{\mathcal{O}}$; fixons un isomorphisme $\phi$ de $\left(\overline{F}^*\right)^d$ dans $\underline{T'}$, où $d$ est la dimension de $\underline{T'}$. Alors on a $K_{Z_G\left(T_0\right)}=Z_G\left(T_0\right)\cap\phi\left(\left(\overline{\mathcal{O}}^*\right)^d\right)$, qui ne dépend pas de $\phi$, et on peut poser:
\[K_{Z_G\left(T_0\right)}^1=Z_G\left(T_0\right)\cap\phi\left(\left(1+\overline{\mathfrak{p}}\right)^d\right);\]
ce groupe ne dépend pas non plus du choix de $\phi$.
D'autre part, soit $A$ une chambre de $\mathcal{A}$ dont l'adhérence contient $x$; le sous-groupe d'Iwahori $K_A$ de $G$ fixant $A$ est contenu dans $K_x$, donc pour tout $\alpha\in\Phi_{T_0}$, on a $U_{\alpha,K_A}\subset U_{\alpha_x}$. En particulier, posons:
\[U_{\alpha,x+}=\bigcap_AU_{\alpha,K_A},\]
l'intersection portant sur les chambres $A$ de $\mathcal{A}$ dont l'adhérence contient $x$; $U_{\alpha,x+}$ peut également être vu comme l'ensemble des éléments de $U_\alpha$ fixant point par point un voisinage de $x$ dans $\mathcal{A}$. Si $V$ est l'espace vectoriel associé à $\mathcal{A}$ et si $H_\alpha$ est l'hyperplan vectoriel noyau de la symétrie $s_\alpha$ sur $V$ associée à $\alpha$, on a $U_{\alpha,x+}=U_{\alpha,x}$ si on est dans un des deux cas suivants:
\begin{itemize}
\item $x+H_\alpha$ n'est pas un mur de $\mathcal{A}$;
\item la facette $A_x$ contenant $x$ n'est pas incluse dans $x+H_\alpha$.
\end{itemize}
Si $x+H_\alpha$ est un mur de $\mathcal{A}$ et si $A_x\subset x+H_\alpha$, $U_{\alpha,x+}$ est le plus grand sous-groupe strict de $U_{\alpha,x}$ de la forme $U_{\alpha,y}$, $y\in\mathcal{A}$; si on pose:
\[K_x^1=K_{Z_G\left(T_0\right)}^1\prod_{\alpha\in\Phi_{T_0}}U_{\alpha,x+},\]
d'après \cite[proposition I.2.2 et corollaire I.2.3]{scst}, ce produit ne dépend pas de l'ordre dans lequel on place les $\alpha$, et $K_x^1$ est un sous-groupe normal de $K_x$; de plus, le quotient ${\mth{G}}=K_x/K_x^1$ est le groupe des ${\mth{F}}_q$-points d'un groupe réductif connexe défini sur ${\mth{F}}_q$. De plus, si $K_{T_0}$ est l'unique sous-groupe parahorique de $T_0$ et $K_{T_0}^1$ son premier sous-groupe de congruence, le quotient ${\mth{T}}_0=K_{T_0}/K_{T_0}^1$ s'identifie canoniquement à un tore déployé maximal de ${\mth{G}}$, et l'ensemble des $\alpha\in\Phi_{T_0}$ tels que $U_{\alpha,x+}\subsetneq U_{\alpha,x}$ s'identifie canoniquement au système de racines de ${\mth{G}}$ relativement à ${\mth{T}}_0$.

On notera $C_0\left(G\right)$ le sous-espace des éléments de $C_c^\infty\left(G\right)$ annulés par toute distribution sur $G$ invariante par conjugaison. On notera également $C_0\left({\mth{G}}\right)$ le sous-espace des fonctions sur ${\mth{G}}$ annulées par toute distribution sur ${\mth{G}}$ invariante par conjugaison.

\subsection{Normalisateur d'un sous-groupe parahorique}

Soit $A$ une chambre de $\mathcal{B}$, $I$ le sous-groupe d'Iwahori de $G$ qui fixe $A$, $\mathcal{A}$ un appartement de $\mathcal{B}$ contenant $A$, $T_0$ le tore déployé maximal de $G$ associé à $\mathcal{A}$ et $K_{Z_G\left(T_0\right)}=Z_G\left(T_0\right)\cap I$ l'unique sous-groupe parahorique de $Z_G\left(T_0\right)$; le groupe $W'=W'\left(G/T_0\right)=N_G\left(T_0\right)/K_{Z_G\left(T_0\right)}$ est le {\em groupe de Weyl affine} de $G$ relativement à $T_0$. D'après \cite[3.5]{cai}, c'est le produit semi-direct du groupe $\Gamma_I=N_G\left(I\right)/I$ et du sous-groupe distingué $W'_0$ de $W'$ engendré par les symétries orthogonales $s'_1,\dots,s'_{n+1}$ par rapport aux murs de $A$ ($n$ désignant le rang semi-simple de $G$); $\Gamma_I$ agit sur $W'_0$ par permutation des $s'_i$. De plus, on a la décomposition:
\[G=\bigsqcup_{w'\in W'}Iw'I.\]
Par un léger abus de notation, si $H$ est un sous-groupe de $G$ contenant $K_{Z_G\left(T_0\right)}$, on notera $W'\cap H$ l'ensemble des éléments de $W'$ contenus dans $H$.

Le groupe $W'$ agit transitivement sur l'ensemble des chambres de $\mathcal{A}$; pour tout $w'\in W$, on appelle {\em longueur} de $w'$ (relativement à $A$) et on note $l\left(w'\right)=l_A\left(w'\right)$ la longueur d'une galerie minimale entre $A$ et $w'A$; c'est également la longueur $l$ d'une décomposition minimale de $w'$ sous la forme $w'=\gamma s'_{i_1}\dots s'_{i_l}$, $\gamma\in\Gamma_I$. On a pour toute décomposition de ce type:
\[Iw'I=\gamma Is'_{i_1}I\dots Is'_{i_l}I;\]
d'autre part, on a pour tout $i$:
\[Is'_iIs'_iI\subset I\cup Is'_iI.\]
On en déduit en particulier que, si on pose:
\[G^1=\bigsqcup_{w'\in W'_0}Iw'I,\]
$G^1$ est un sous-groupe normal de $G$, ouvert et fermé dans $G$, qui ne dépend pas du choix de $I$;
de plus, posons:
\[\Gamma=G/G^1;\]
d'après ce qui précède, $\Gamma$ est canoniquement isomorphe à $\Gamma_I$; le groupe $\Gamma_I$ ne dépend donc, à isomorphisme près, pas de $I$, et chaque élément de $\Gamma$ possède un unique représentant dans $W'$ appartenant à $\Gamma_I$.

D'autre part, pour toute partie $J\subsetneq\left\{1,\dots,n+1\right\}$, il existe au moins un point de $\overline{A}$ fixé par les $s'_i$, $i\in J$; le sous-groupe de $G$ engendré par $I$ et les $s'_i$, $i\in J$ est alors un sous-groupe parahorique $K_J$ de $G$ contenant $I$, et tous les sous-groupes parahoriques de $G$ contenant $I$ peuvent être obtenus par ce moyen. $G^1$ contient alors tous ces sous-groupes, donc également tous les autres sous-groupes parahoriques de $G$ puisque ceux-ci leurs sont conjugués. En particulier, $G^1$ contient tous les sous-groupes compacts des sous-groupes radiciels de $G$ relativement à $T_0$; $G^1$ contient alors tous ces sous-groupes radiciels, donc également le groupe qu'ils engendrent, qui est le groupe dérivé de $G$; on en déduit que $\Gamma$ est abélien.

De plus, on a le lemme suivant:

\begin{lemme}\label{normpar}
Soit $J$ une partie de $\left\{1,\dots,n+1\right\}$ distincte de $\left\{1,\dots,n+1\right\}$ lui-même. On a:
\[N_G\left(K_J\right)=\Gamma_JK_J,\]
où $\Gamma_J$ est le sous-groupe des éléments de $\Gamma_I$ qui permutent les $\alpha_i$, $i\in J$.
\end{lemme}

\begin{proof}
Puisque $I\subset K_J$, $N_G\left(K_J\right)$ est une réunion de doubles classes modulo $I$, et il suffit de caractériser les éléments $w'\in W'$ qui normalisent $K_J$. Soit donc un tel $w'\in W'$; écrivons-le $w'=w'_0\gamma$, avec $\gamma\in\Gamma$ et $w'_0\in W'_0$. Puisque $\gamma$ normalise $I$, il stabilise $A$ et permute les faces de $A$; si $A_J$ est la face de $A$ dont $K_J$ est le fixateur connexe, il existe alors $J'$ telle que $\gamma A_J=A_{J'}$; pour que $w'$ normalise $K_J$, on doit donc avoir $w'_0A_{J'}=A_J$. Or puisque $A$ possède une unique face de chaque type, d'après \cite[I. 1.3.5]{bt}, il existe une unique face de $A$ transformée de $A_{J'}$ par un élément de $W'_0$, et cette face est forcément $A_{J'}$; on a donc $J'=J$, ce qui implique à la fois que $\gamma\in\Gamma_J$ et que $w'_0A_J=A_J$, d'où l'on déduit grâce à \cite[I. 1.3.5]{bt} que $w'_0$ est un élément du sous-groupe de $W'_0$ engendré par les $s_{\alpha_i}$, $i\in J$, ce qui démontre le lemme.
\end{proof}

On en déduit immédiatement le corollaire suivant:

\begin{cor}\label{normp2}
Pour tout $J$, $N_G\left(K_J\right)\cap G^1=K_J$.
\end{cor}

Soit maintenant $x$ un point quelconque de $G$, $K_x$ son fixateur connexe dans $G$; on posera $\Gamma_x=\Gamma_{K_x}=N_G\left(K_x\right)/K_x$.

Soit $K$ un sous-groupe parahorique de $G$. Comme son premier sous-groupe de congruence est déterminé de manière unique, il est également normalisé par tout élément de $N_G\left(K\right)$.
Supposons que $K$ contient $I$, et soit $\gamma\in N_G\left(I\right)$. Considérons le groupe ${\mth{G}}^+={\mth{G}}_\gamma^+=\gamma^{\mth{Z}}K/K^1$; c'est le produit semi-direct de ${\mth{G}}=K/K^1$ par le groupe $\gamma^{\mth{Z}}$. Puisque $\gamma$ peut être d'ordre infini, ${\mth{G}}^+$ n'est pas toujours le groupe des ${\mth{F}}_q$-points d'un groupe réductif, mais on peut faire la remarque suivante: puisque ${\mth{G}}$ est fini, son groupe d'automorphismes l'est également, et il existe donc un entier $b$ tel que $\gamma^b$ agit trivialement par conjugaison sur ${\mth{G}}$; $\gamma^b$ est alors central dans ${\mth{G}}^+$, et ${\mth{G}}^+/\gamma^{b{\mth{Z}}}$ est le groupe des ${\mth{F}}_q$-points d'un groupe réductif défini sur ${\mth{F}}_q$. On va donc pouvoir appliquer à ${\mth{G}}^+$ certains des résultats de \cite{dm}.

On appellera {\em sous-groupe parabolique} de ${\mth{G}}^+$ un sous-groupe de la forme ${\mth{P}}^+=N_{{\mth{G}}^+}\left({\mth{P}}\right)$, où ${\mth{P}}={\mth{M}}{\mth{U}}$ est un sous-groupe parabolique de ${\mth{G}}$. De même que dans \cite[proposition 1.5]{dm}, on a une décomposition de la forme ${\mth{P}}^+={\mth{M}}^+{\mth{U}}$, où ${\mth{M}}^+=N_{{\mth{P}}^+}\left({\mth{M}}\right)$; ${\mth{M}}^+$ est un {\em sous-groupe de Levi} de ${\mth{G}}^+$. De plus, on a ${\mth{P}}={\mth{P}}^+\cap{\mth{G}}$ et ${\mth{M}}={\mth{M}}^+\cap{\mth{G}}$; il y a donc bijection canonique entre les sous-groupes paraboliques (resp. de Levi) de ${\mth{G}}^+$ et ceux de ${\mth{G}}$.

D'autre part, il existe un représentant de $\gamma$ dans ${\mth{G}}^+$ qui est quasi-central (au sens de \cite[1.15]{dm}): c'est une conséquence immédiate de \cite[1.16 et 1.34]{dm}. En assimilant $\gamma$ à un tel représentant, et en posant ${\mth{G}}^\gamma={\mth{G}}\cap Z_{{\mth{G}}^+}\left(\gamma\right)$, on en déduit le résultat suivant:

\begin{lemme}\label{centgamma}
Il existe une bijection canonique entre les classes de conjugaison de sous-groupes de Levi de ${\mth{G}}^+$ contenant $\gamma$ et les classes de conjugaison de sous-groupes de Levi de ${\mth{G}}^\gamma$, induite par l'application canonique ${\mth{M}}^+\mapsto{\mth{M}}^+\cap{\mth{G}}^\gamma$.
\end{lemme}

\begin{proof}
Ce lemme est une conséquence immédiate de \cite[proposition 1.40]{dm}.
\end{proof}

Soit ${\mth{P}}^+={\mth{M}}^+{\mth{U}}$ un sous-groupe parabolique de ${\mth{G}}^+$. Si $f$ est un élément de l'espace $C\left({\mth{G}}^+\right)^{{\mth{G}}^+}$ des fonctions sur ${\mth{G}}^+$ invariante par conjugaison, on définit sa restriction à ${\mth{M}}^+$ par:
\[r_{{\mth{P}}^+}^{{\mth{G}}^+}f: m\in{\mth{M}}^+\longmapsto\sum_{u\in{\mth{U}}}f\left(mu\right).\]
Si maintenant $\phi$ est une fonction sur ${\mth{M}}^+$ invariante par conjugaison, on définit son induite à ${\mth{G}}^+$ par:
\[\Ind_{{\mth{P}}^+}^{{\mth{G}}^+}\phi:g\in{\mth{G}}^+\longmapsto\sum_{m\in{\mth{M}}^+,u\in{\mth{U}}^+,h\in{\mth{G}}^+,g=h^{-1}muh}\phi\left(m\right).\]

Considérons le produit hermitien défini positif sur l'espace des fonctions invariantes sur ${\mth{G}}^+$ à support dans $\gamma{\mth{G}}$, donné par:
\[\langle f,f'\rangle _{\mth{G}}=\dfrac 1{\card\left({\mth{G}}\right)}\sum_{g\in\gamma{\mth{G}}}\overline{f\left(g\right)}f'\left(g\right).\]
On a les résultats suivants:

\begin{lemme}\label{rpf}
Soit $f$ (resp. $\phi$) une fonction centrale sur ${\mth{G}}^+$ (resp. ${\mth{M}}^+$) à support dans $\gamma{\mth{G}}$ (resp. $\gamma{\mth{M}}$). On a:
\[\langle f,\Ind_{{\mth{P}}^+}^{{\mth{G}}^+}\phi\rangle _{\mth{G}}=\langle r_{{\mth{P}}^+}^{{\mth{G}}^+}f,\phi\rangle _{\mth{M}}.\]
\end{lemme}

\begin{proof}
En effet, on a:
\begin{equation*}\begin{split}\langle f,\Ind_{{\mth{M}}^+}^{{\mth{G}}^+}\phi\rangle _{\mth{G}}&=\dfrac 1{\card\left({\mth{G}}\right)}\sum_{g\in\gamma{\mth{G}}}\overline{f\left(g\right)}\left(\dfrac 1{\card\left({\mth{P}}\right)}\sum_{x\in{\mth{G}},m\in\gamma{\mth{M}},u\in{\mth{U}},g=x^{-1}mux}\phi\left(m\right)\right)\\&=\dfrac 1{\card\left({\mth{P}}\right)\card\left({\mth{G}}\right)}\sum_{m\in\gamma{\mth{M}}}\left(\sum_{x\in{\mth{G}},u\in{\mth{U}}}\overline{f\left(x^{-1}mux\right)}\phi\left(m\right)\right)\\&=\dfrac 1{\card\left({\mth{M}}\right)\card\left({\mth{U}}\right)}\sum_{m\in\gamma{\mth{M}},u\in{\mth{U}}}\overline{f\left(mu\right)}\phi\left(m\right)\\&=\langle r_{{\mth{P}}^+}^{{\mth{G}}^+}f,\phi\rangle _{\mth{M}}.\end{split}\end{equation*}
\end{proof}

\begin{prop}\label{resind}
Soit ${\mth{P}}^+={\mth{M}}^+{\mth{U}}$ un sous-groupe parabolique de ${\mth{G}}^+$ contenant $\gamma$.
\begin{itemize}
\item Soit $f\in C\left({\mth{G}}^+\right)^{{\mth{G}}^+}$ à support dans $\gamma{\mth{G}}$; la restriction $r_{{\mth{P}}^+}^{{\mth{G}}^+}f$ de $f$ à ${\mth{M}}^+$ ne dépend pas du choix de ${\mth{P}}^+$;
\item soit $\phi\in C\left({\mth{M}}^+\right)^{{\mth{M}}^+}$ à support dans $\gamma{\mth{M}}$; l'induite $\Ind_{{\mth{P}}^+}^{{\mth{G}}^+}\phi$ de $\phi$ à ${\mth{G}}^+$ ne dépend pas du choix de ${\mth{P}}^+$.
\end{itemize}
\end{prop}

\begin{proof}
Si ${\mth{P}}^+={\mth{G}}^+$, le résultat est trivial; supposons donc ${\mth{P}}^+$ propre. Soit ${\mth{P}}'{}^+={\mth{M}}^+{\mth{U}}'$ un autre sous-groupe parabolique de ${\mth{G}}^+$ contenant $\gamma$ et de Levi ${\mth{M}}^+$; on va montrer par récurrence sur le rang semi-simple de ${\mth{G}}^+$ que l'on a:
\[\langle \Ind_{{\mth{P}}^+}^{{\mth{G}}^+}\phi,\Ind_{{\mth{P}}^+}^{{\mth{G}}^+}\phi\rangle_{\mth{G}}=\langle \Ind_{{\mth{P}}^+}^{{\mth{G}}^+}\phi,\Ind_{{\mth{P}}'{}^+}^{{\mth{G}}^+}\phi\rangle_{\mth{G}}=\langle \Ind_{{\mth{P}}'{}^+}^{{\mth{G}}^+}\phi,\Ind_{{\mth{P}}'{}^+}^{{\mth{G}}^+}\phi\rangle_{\mth{G}},\]
ce qui montrera la seconde assertion, la première s'en déduisant au moyen du lemme \ref{rpf}.

On déduit de ce même lemme \ref{rpf} que ces égalités sont équivalentes à:
\[\langle \phi,r_{{\mth{P}}^+}^{{\mth{G}}^+}\Ind_{{\mth{P}}^+}^{{\mth{G}}^+}\phi\rangle=\langle \phi,r_{{\mth{P}}^+}^{{\mth{G}}^+}\Ind_{{\mth{P}}'{}^+}^{{\mth{G}}^+}\phi\rangle=\langle \phi,r_{{\mth{P}}'{}^+}^{{\mth{G}}^+}\Ind_{{\mth{P}}'{}^+}^{{\mth{G}}^+}\phi\rangle.\]
Ces dernières égalités se déduisent immédiatement de la formule de Mackey pour les paraboliques de ${\mth{G}}^+$ (cf. \cite[théorème 3.2]{dm}) et de l'hypothèse de récurrence.
\end{proof}

\subsection{Éléments compacts}

Le but de cette partie est de classifier les parties de $G$ de la forme $\gamma K$, où $K$ est un sous-groupe parahorique de $G$ et $\gamma\in N_G\left(K\right)$.

On dira qu'un élément de $G$ est {\em compact} dans $G$ s'il est contenu dans le normalisateur d'un sous-groupe parahorique de $G$; $g\in G$ est compact dans $G$ si et seulement si il appartient à un sous-groupe de $G$ compact modulo le centre $Z$ de $G$, ce qui est vrai si et seulement si toutes les valeurs propres de l'application linéaire $Ad\left(g\right)$ sur $Lie\left(G\right)$ sont de valuation $0$. La notion de compacité dépend du groupe considéré: en effet, si $H$ est un sous-groupe réductif fermé de $G$, tous les éléments compacts de $G$ sont compacts dans $H$, mais la réciproque est fausse.

\begin{lemme}\label{parferm}
L'ensemble $G_c$ des éléments compacts de $G$ est un ouvert fermé de $G$.
\end{lemme}

\begin{proof}
En effet, $G_c$ est réunion de sous-groupes ouverts de $G$; c'est donc un ouvert de $G$. D'autre part, considérons l'application $\phi$ de $G$ dans $F\left[X\right]$ telle que pour tout $g\in G$, $\phi\left(G\right)$ est le polynôme caractéristique de $Ad\left(g\right)$ sur $Lie\left(G\right)$; $\phi$ est continue, et $G_c$ est l'image réciproque par $\phi$ de l'ensemble des polynômes de degré $dim\left(G\right)$ dont toutes les valeurs propres sont de valuation $0$, qui est un fermé de $F\left[X\right]$; $G_c$ est donc un fermé de $G$ et le lemme est démontré.
\end{proof}

Soit $x\in\mathcal{B}$, $\mathcal{A}$ un appartement de $\mathcal{B}$ contenant $x$, $T_0$ le tore déployé maximal de $G$ associé à $\mathcal{A}$.
Soit $\Phi_x$ l'ensemble des racines de $G$ relativement à $T_0$ telles que $U_{\alpha,x+}$ est strictement contenu dans $U_{\alpha,x}$; d'après ce qui précède, c'est également l'ensemble des racines $\alpha$ telles que $H_\alpha$ est un mur de $\mathcal{A}$ et que $A_x\subset x+H_\alpha$; cet ensemble est un sous-système de racines de $\Phi$. Soit $M_x$ le groupe de $G$ engendré par $Z_G\left(T_0\right)$ et les $U_\alpha$, où $\alpha$ décrit l'ensemble des racines qui sont combinaisons linéaires à coefficients dans ${\mth{Q}}$ d'éléments de $\Phi_x$: c'est un sous-groupe de Levi de $G$, dont la classe de conjugaison ne dépend pas du choix de $T_0$. On pose $K_{M_x}=K_x\cap M_x$; $K_{M_x} $ est un sous-groupe parahorique maximal de $M_x$ dont la classe de conjugaison ne dépend pas non plus du choix de $T_0$; la classe de conjugaison $\mu_x$ (dans $G$) du couple $\left(M_x,K_{M_x}\right)$ est donc déterminée uniquement par $x$.

(Remarquons que le groupe $H_x$ engendré par $Z_G\left(T_0\right)$ et les $U_\alpha$, $\alpha\in\Phi_x$, est un sous-groupe réductif fermé de $G$ dont $T_0$ est un tore déployé maximal, mais n'est pas nécessairement un sous-groupe de Levi de $G$.)

Soit $x,y\in\mathcal{B}$; on dira que $K_x$ et $K_y$ sont {\em associés} si $\mu_x=\mu_y$, c'est-à-dire si deux couples $\left(M_x,K_{M_x}\right)$ et $\left(M_y,K_{M_y}\right)$ quelconques sont conjugués entre eux. Deux parahoriques de $G$ conjugués entre eux sont clairement associés; la réciproque est fausse.

On notera $\mathcal{M}_0$ l'ensemble des couples $\left(M,K_M\right)$, où $M$ est un sous-groupe de Levi de $G$ et $K_M$ un sous-groupe parahorique maximal de $M$, et $\mathcal{M}$ l'ensemble des classes de conjugaison d'éléments de $\mathcal{M}_0$.

Si $\mathcal{K}_0$ est l'ensemble des couples $\left(K_x,T_0\right)$, où $x$ est un élément de $\mathcal{B}$ et $T_0$ un tore déployé maximal de $G$ correspondant à un appartement de $\mathcal{B}$ contenant $x$, on a ainsi défini deux applications canoniques:
\begin{itemize}
\item une application $\kappa_0$ de $\mathcal{K}_0$ dans $\mathcal{M}_0$;
\item une application $\kappa$ de $\mathcal{B}$ dans $\mathcal{M}$.
\end{itemize}
Puisque $\kappa(x)$ ne dépend que de la facette contenant $x$, par un léger abus de notation, on notera parfois $\kappa\left(K_x\right)=\kappa\left(x\right)$.

Soit $\mu\in\mathcal{M}$ et $\left(M,K_\mu\right)$ un représentant de $\mu$; si $K_\mu^1$ est le premier sous-groupe de congruence de $K_\mu$, le groupe ${\mth{M}}=K_\mu/K_\mu^1$ est l'ensemble des ${\mth{F}}_q$-points d'un groupe réductif connexe $\underline{\mth{M}}$ défini sur ${\mth{F}}_q$, qui ne dépend pas du choix de $\left(M,K_\mu\right)$; de plus, si $x$ est un élément de $\mathcal{B}$ dont l'image par $\kappa$ est $\mu$, le groupe réductif fini $K_x/K_x^1$ est canoniquement isomorphe à ${\mth{M}}$.
Pour tout élément $k$ (resp. toute partie $X$) de $K_\mu$ ou de $K_x$, on notera $\overline{k}$ ou $\overline{k}_{\mth{M}}$ (resp. $\overline{X}$ ou $\overline{X}_{\mth{M}}$) son image dans ${\mth{M}}$.

Soit maintenant $\mathcal{A}$ un appartement de $\mathcal{B}$, $T_0$ le tore maximal de $G$ associé à $\mathcal{A}$, $x\in\mathcal{A}$, $\left(M,K_M\right)=\kappa_0\left(K_x,T_0\right)$, $\mathcal{A}_M$ l'appartement de l'immeuble de $\mathcal{B}_M$ de $M$ associé à $T_0$ et $x_M$ un point de la facette de $\mathcal{A}_M$ fixée par $K_M$; on va définir une surjection canonique entre $\mathcal{A}$ et $\mathcal{A}_M$. Puisque l'espace vectoriel $V_M$ associé à $\mathcal{A}_M$ est canoniquement isomorphe à $X_*\left(T_0\right)\otimes{\mth{R}}/X_*\left(Z_{M,d}\right)\otimes{\mth{R}}$, où $Z_{M,d}$ est la composante déployée du centre de $M$, et puisqu'on a clairement $Z_{G,d}\subset Z_{M,d}$, on a une surjection canonique entre l'espace vectoriel $V$ associé à $\mathcal{A}$ et $V_M$. Soit $\phi_x$ la bijection de $\mathcal{A}$ dans $\mathcal{A}_M$ induite par cette bijection canonique en posant $\phi_x\left(x\right)=x_M$; on va montrer que $\phi_x$ ne dépend pas du choix de $x$. Soit donc $y$ un autre élément de $\mathcal{A}$ tel que $\kappa_0\left(K_y,T_0\right)$ est de la forme $\left(M,K'_M\right)$;
on va vérifier que si $y_M$ est un point de $\mathcal{A}_M$ fixé par $K'_M$, on a $\phi_x\left(y\right)=y_M$.

Puisque l'adhérence de toute chambre de $\mathcal{A}$ contient au moins un $y$, par une récurrence évidente, on peut supposer que $x$ et $y$ sont séparés par un unique mur $H$ de $\mathcal{A}$; soit $\alpha_0'$ la racine affine de $G$ relativement à $T_0$ associée à $H$ et contenant $x$ et pas $y$, et $\alpha_0$ la racine de $G$ relativement à $T_0$ associée à $\alpha'_0$. Comme $\phi_x$ et $\phi_y$ sont deux applications affines de même partie linéaire, il suffit de considérer leurs images respectives en un point donné.

Pour toute racine $\alpha$ de $M$ relativement à $T_0$, si $U_\alpha$ est le sous-groupe radiciel de $M$ correspondant et $U_{\alpha,K}$ un sous-groupe ouvert compact de $U_\alpha$ maximal parmi ceux contenus dans $K'_M$, $U_{\alpha,K}$ fixe un demi-espace de $\mathcal{A}$ contenant $y$, donc également un demi-espace de $\mathcal{A}_M$ contenant $\phi_x\left(y\right)$ (c'est trivial si $\alpha\neq\pm\alpha_0$; si $\alpha=\alpha_0$ (resp. si $\alpha=-\alpha_0$), le demi-espace de $\mathcal{A}_M$ fixé par $U_{\alpha,K}$ est l'adhérence de la plus grande racine affine de $M$ relativement à $T_0$ strictement contenue dans $\phi_x\left(\alpha'_0\right)$ (resp. de la plus petite contenant strictement $\phi_x\left(\alpha'_0\right)$); or le seul sous-groupe parahorique de $M$ contenant tous les $U_{\alpha,K}$ est $K'_M$, ce qui impose $\phi_x\left(y\right)=y_M$. Par le même raisonnement, on voit que si maintenant $y$ est un point quelconque de $\mathcal{A}$, $\phi_x\left(y\right)$ est fixé par $K_y\cap M$; on a donc montré le résultat suivant:

\begin{lemme}\label{bjc}
Il existe une unique surjection affine de $\mathcal{A}$ dans $\mathcal{A}_M$ qui, à tout $x\in\mathcal{A}$, associe un point $x_M$ de $\mathcal{A}_M$ fixé par $K_x\cap M$.
\end{lemme}

On va maintenant utiliser ce lemme pour montrer le résultat suivant:

\begin{lemme}\label{musj}
$\kappa$ et $\kappa_0$ sont surjectives.
\end{lemme}

\begin{proof}
Montrons d'abord que $\kappa_0$ est surjective, c'est-à-dire que pour tout $\mu\in\mathcal{M}$ et tout représentant $\left(M,K_\mu\right)$ de $\mu$, il existe un point $x$ de $\mathcal{B}$ et un tore maximal $T_0$ de $G$ dont l'appartement $\mathcal{A}$ de $\mathcal{B}$ associé contient $x$, tels que $\left(M,K_\mu\right)=\kappa_0\left(K_x,T_0\right)$.

Soit donc $\mathcal{A}_M$ un appartement de l'immeuble de $M$ contenant un point $x_M$ fixé par $K_\mu$, $T_0$ le tore déployé maximal de $M$ associé à $\mathcal{A}_M$, $\mathcal{A}$ l'appartement de l'immeuble de $G$ associé à $T_0$,
et $\phi$ une isométrie de $\mathcal{A}$ dans $\mathcal{A}_M$ vérifiant les conditions du lemme précédent. Considérons l'image réciproque $E$ dans $\mathcal{A}$ par $\phi$ de la facette $A_M$ de $\mathcal{A_M}$ contenant $x_M$: c'est un sous-espace affine de $\mathcal{A}$, et pour tout $x\in E$, $K_x$ contient $K_\mu$; puisque $K_\mu$ est un sous-groupe parahorique maximal de $M$, on en déduit que l'on a $K_x\cap M=K_\mu$; réciproquement, pour tout $x\in\mathcal{A}$ tel que $K_x\cap M=K_\mu$, on a $\phi\left(x\right)\in A_M$, donc $x$ appartient à $E$. On en déduit que $E$ est réunion de facettes de $\mathcal{A}$.

Soit maintenant $A$ une facette de $\mathcal{A}$ de dimension maximale parmi celles contenues dans $E$, et soit $x\in A$.
Alors $K_x$ convient: en effet, soit $\left(M_x,K_{M_x}\right)=\kappa_0\left(K_x,T_0\right)$;
on déduit de la définition de $\kappa_0$ que $M\subset M_x$; d'autre part, le rang de $M_x$ est égal au rang de $G$ moins la dimension de $A$, qui est égale à celle de $E$, ce qui donne:
\[rg\left(M_x\right)=rg\left(G\right)-dim\left(E\right)=rg\left(M\right);\]
on en déduit que $M_x=M$. Enfin, $K_{M_x}=K_x\cap M=K_\mu$, ce qui achève la démonstration de la surjectivité de $\kappa_0$.

Montrons maintenant la surjectivité de $\kappa$: soit $\mu\in\mathcal{M}$, et $\left(M,K_\mu\right)$ un re\-pré\-sen\-tant de $\mu$ dans $\mathcal{M}_0$. Puisque $\kappa_0$ est surjective, il existe un couple $\left(K_x,T_0\right)\in\mathcal{K}$ dont l'image par $\kappa_0$ est $\left(M,K_\mu\right)$; mais alors, on a $\kappa\left(x\right)=\mu$. Donc $\kappa$ est surjective et le lemme est démontré.
\end{proof}

On va maintenant s'intéresser aux classes de la forme $\left(\gamma,K_x\right)$, avec $x\in\mathcal{B}$ et $\gamma\in N_G\left(K_x\right)$. Soit tout d'abord $\left(M,K\right)\in\mathcal{M}_0$, $T_0$ un tore déployé maximal de $M$, et $x\in\mathcal{B}$ tel que $\kappa_0\left(K_x,T_0\right)=\left(M,K\right)$; il est clair que $N_G\left(K_x\right)\subset N_G\left(M,K\right)$. Posons donc:
\[\Gamma_{\left(M,K\right)}=N_G\left(M,K\right)/N_G\left(M,K\right)\cap G^1.\]
Il est clair que si $\left(M',K'\right)$ est un élément de $\mathcal{M}_0$ conjugué à $\left(M,K\right)$, $\Gamma_{\left(M',K'\right)}=\Gamma_{\left(M,K\right)}$; $\Gamma_{\left(M,K\right)}$ ne dépend donc que de la classe de conjugaison $\mu$ de $\left(M,K\right)$, et on le notera alors $\Gamma_\mu$. Pour tout $\mu\in\mathcal{M}$, $\Gamma_\mu$ s'identifie donc à un sous-groupe de $\Gamma$ contenant tous les $\Gamma_x$, avec $x\in\mathcal{B}$ tel que $\kappa\left(x\right)=\mu$. Remarquons que les $\Gamma_x$ ne sont généralement pas égaux à $\Gamma_\mu$, ni même forcément égaux entre eux.

On a le résultat suivant:

\begin{prop}\label{gcj}
Soit $x,x'\in\mathcal{B}$ tels que $K_x$ et $K_{x'}$ sont associés, soit $T_0$ (resp. $T'_0$) le tore déployé maximal de $G$ associé à un appartement $\mathcal{A}$ (resp. $\mathcal{A}'$) contenant $x$ (resp. $x'$) et soit $\left(M,K\right)=\kappa_0\left(K_x,T_0\right)$ (resp. $\left(M',K'\right)=\kappa_0\left(K_{x'},T'_0\right)$). Soit $\gamma$ un élément de $W'\left(G/T_0\right)\cap N_G\left(K_x\right)$ et $\gamma'$ un élément de $W'\left(G,T'_0\right)\cap N_G\left(K_{x'}\right)$ représentant un même élément de $\Gamma$ et normalisant chacun un sous-groupe d'Iwahori de $G$; alors les triplets $\left(\gamma,M,K\right)$ et $\left(\gamma',M',K'\right)$ sont conjugués dans $G$.
\end{prop}

\begin{proof}
La démonstration de cette proposition fera l'objet du paragraphe \ref{ddlp}.
\end{proof}

On déduit de cette proposition que le groupe $\gamma^{\mth{Z}}K_x/K_x^1$ ne dépend que de $\gamma$ et de $\mu=\kappa\left(x\right)$.
Considérons donc l'ensemble $\mathcal{G}$ des couples de la forme $\left(\gamma,x\right)$, $x\in\mathcal{B}$ et $\gamma\in\Gamma_x$, et l'ensemble $\mathcal{N}$ des couples $\nu=\left(\gamma,\mu\right)$, $\mu\in\mathcal{M}$ et $\gamma\in\Gamma_\mu$, et soit l'application $\zeta$ de $\mathcal{G}$ dans $\mathcal{N}$ définie par:
\[\zeta:\left(\gamma,x\right)\mapsto\left(\gamma,\kappa\left(x\right)\right);\]
si $\nu\in\mathcal{N}$ possède au moins un antécédent $x$ par $\zeta$, on peut lui associer de manière canonique le groupe ${\mth{G}}^+={\mth{G}}^+_\nu=\gamma^{\mth{Z}}K_x/K_x^1$. Notons que $\zeta$ n'est généralement pas surjective.

Puisque $\zeta\left(\gamma,x\right)$ ne dépend que de $\gamma$ et de la facette contenant $x$, par un léger abus de notation, on notera également parfois $\zeta\left(\gamma,K_x\right)=\zeta\left(\gamma,x\right)$.

\subsection{Sous-groupes standard et semi-standard}

Fixons un tore déployé maximal $T_0$ de $G$ et un sous-groupe parabolique minimal $P_0=M_0U_0$ de $G$ tel que $T_0$ est contenu dans le centre de $M_0$. Un sous-groupe parabolique $P=MU$ de $G$ sera dit standard (resp. semi-standard) relativement à $P_0$ et $T_0$ (resp. $T_0$) si $P_0\subset P$ (resp. $M_0\subset P$); $M$ est un sous-groupe de Levi standard (resp. semi-standard) de $G$ s'il contient $M_0$ et s'il existe un $P=MU$ tel que $P$ est standard (resp. semi-standard).

Soit maintenant $M$ un sous-groupe de Levi semi-standard de $G$, $\mathcal{A}_M$ l'appartement de l'immeuble $\mathcal{B}_M$ de $M$ correspondant à $T_0$; soit $K$ un sous-groupe parahorique maximal de $M$ fixant un point de $\mathcal{A}_M$, et soit ${\mth{M}}=K/K^1$; il existe un unique sous-groupe d'Iwahori $I_K$ de $M$ contenu dans $K$ dont l'image $\overline{I}_K$ dans ${\mth{M}}$ est le sous-groupe parabolique minimal ${\mth{P}}_0=\overline{P_0\cap K}$ de ${\mth{M}}$; $I_K$ sera appelé {\em sous-groupe d'Iwahori standard} (relativement à $T_0$) de $K$. Si $K_x$ est un sous-groupe parahorique de $G$ fixant un point $x$ de $\mathcal{A}$ et tel que $\left(M,K\right)=\kappa_0\left(K_x,T_0\right)$, il existe un unique sous-groupe d'Iwahori $I$ de $G$ contenu dans $K_x$ et tel que $I\cap K_x=I_K$; on l'appellera le {\em sous-groupe d'Iwahori standard} de $K_x$ (relativement à $T_0$).

Si $M$ est un sous-groupe de Levi standard de $G$, il existe un unique sous-groupe parabolique standard $P$ de $G$ de Levi $M$; c'est $P=MP_0$. Soit $K_0$ un sous-groupe parahorique maximal spécial de $G$ fixant un sommet de l'appartement $\mathcal{A}$ de $\mathcal{B}$ associé à $T_0$; on a $G=K_0P_0=K_0P$, et $K_0\cap M$ est un sous-groupe parahorique maximal spécial de $M$. Fixons des mesures de Haar sur $K_0$ et sur $U$, et définissons, pour $f\in C_c^\infty\left(G\right)$, le terme constant de $f$ selon $P$, noté $f^P_{K_0}$ ou $f^P$, par:
\[f^P_{K_0}:m\in M\longmapsto\delta_P\left(m\right)^{\dfrac 12}\int_{K_0\times U}f\left(k^{-1}muk\right)dkdu,\]
où $\delta_P$ est la fonction module sur $P$; cette définition dépend du choix de $K_0$, mais si $K'_0$ est un sous-groupe de $G$ vérifiant les mêmes conditions, avec des normalisations convenables des mesures de Haar, $f^P_{K_0}$ et $f^P_{K'_0}$ ne diffèrent que d'un élément de $C_0\left(M\right)$; si $D$ est une distribution invariante sur $M$, $D\left(f^P_{K_0}\right)$ ne dépend pas de $K_0$ et on pourra donc écrire sans ambiguïté $D\left(f^P\right)$ sans préciser $K_0$.

\section{Distributions invariantes à support compact}

Le but de cette partie est de montrer que les restrictions à certains sous-espaces de $C_c^\infty\left(G\right)$ des distributions invariantes sur $G$ à support compact sont entièrement déterminées par leurs restrictions à des sous-espaces de dimension finie de fonctions à support dans des normalisateurs de sous-groupes parahoriques de $G$. Pour cela, on va utiliser certaines distributions particulières.

\subsection{Démonstration de la proposition}\label{ddlp}

On va maintenant démontrer la proposition \ref{gcj}.
Montrons d'abord que l'on peut supposer $\gamma=\gamma'$ et $T_0=T'_0$. Soit $A$ une chambre de $\mathcal{A}'$ dont l'adhérence contient $x'$, et $I$ le sous-groupe d'Iwahori de $G$ fixant $A$; on a $I\subset K_{x'}$. Quitte à remplacer $x$ par $gx$, avec $g$ convenable, et à conjuguer $T_0$ et $\gamma$ par $g$, on peut supposer que $A\subset\mathcal{A}$ et que $K_x$ contient également $I$; on peut alors supposer $\gamma\in W'\cap N_G\left(I\right)$ et $\gamma'=\gamma$.
De plus, d'après \cite[corollaire 2.2.6]{bt}, il existe $h\in G'$ fixant $A$ et tel que $h\mathcal{A}=\mathcal{A}'$; quitte à conjuguer $T_0$ et $\gamma$ par $h$, on peut donc supposer $T_0=T'_0$.

Supposons donc $\gamma=\gamma'$ et $T_0=T'_0$, et soit $g\in G$ tel que $g^{-1}\left(M',K'\right)g=\left(M,K\right)$, et soit $\mathcal{A}_M$ (resp. $\mathcal{A}'_M$) l'appartement de l'immeuble $\mathcal{B}_M$ de $M$ associé à $T_0$ (resp. $g^{-1}T_0g$). Alors $\mathcal{A}$ et $\mathcal{A}'$ contiennent le sommet $y$ de $\mathcal{B}_M$ fixé par $K$, donc d'après \cite[2.5.8]{bt}, il existe $h\in M$ fixant $y$ et tel que $h\mathcal{A}'=\mathcal{A}$. Mais alors on a $\left(gh\right)^{-1}\left(M',K'\right)gh=\left(M,K\right)$ et ${gh}^{-1}T_0gh=T_0$; $\left(M,K\right)$ et $\left(M',K'\right)$ sont donc conjugués par un élément $w$ de $W'$, et on a $w^{-1}\left(W'\cap K'\right)w=W'\cap K$.

Réciproquement, si $w$ est un élément de $W'$ vérifiant cette dernière égalité, posons $\left(M'',K''\right)=w^{-1}\left(M',K'\right)w$. Alors $W'\cap K''=W'\cap K$; les sous-espaces de $\mathcal{A}$ fixés par $K$ et $K''$ sont alors identiques, et on en déduit que l'on a $K=K''$, d'où, puisque $M$ (resp. $M''$) est la fermeture de Zariski de $K$ (resp. $K''$), $M=M''$. Pour montrer la proposition, il suffit donc de trouver $w\in W'$ qui commute avec $\gamma$ et tel que $w^{-1}\left(W'\cap K'\right)w=W'\cap K$.

Montrons maintenant que l'on peut supposer que le centre de $G$ est compact. Soit $Z_{G,d}=\underline{Z_{G,d}}\left(F\right)$ la composante déployée du centre de $G$; considérons le morphisme canonique:
\[\phi:G\hookrightarrow\underline{G}\longrightarrow\underline{G}/\underline{Z_{G,d}}.\]
Il est clair que le noyau de $\phi$ est $Z_{G,d}$; $\phi$ induit donc une injection de $G/Z_{G,d}$ dans $\underline{G}/\underline{Z_{G,d}}$. De plus, le groupe de cohomologie $H^1\left(\Gal\left(\overline{F}/F\right),\underline{Z_{G,d}}\right)$ est trivial (c'est vrai si $\underline{Z_{G,d}}$ est de dimension $1$ par \cite[II. 1, proposition 1]{ser}, et dans le cas général car $\underline{Z_{G,d}}$ est produit direct de tores de dimension $1$ déployés sur $F$), et l'on en déduit que l'image de $\phi$ est $\left(\underline{G}/\underline{Z_{G,d}}\right)\left(F\right)$, et donc que ce groupe est égal à $G/Z_{G,d}$. De plus, si $\left(\gamma,M,K\right)$ et $\left(\gamma',M',K'\right)$ sont conjugués dans $G$, leurs images le sont clairement dans $G/Z_{G,d}$; réciproquement, d'une part $Z_{G,d}$ est contenu à la fois dans $M$ et $M'$, et d'autre part $\gamma$ et $\gamma'$ (resp. $K$ et $K'$) sont contenus dans la même classe de $G$ modulo $G^1$, et sont invariants par multiplication par $K_{Z_G\left(T_0\right)}\supset Z_{G,d}\cap G^1$; on en déduit que si $\left(\gamma,M,K\right)$ et $\left(\gamma',M',K'\right)$ sont conjugués modulo $Z_{G,d}$, alors ils le sont dans $G$. Quitte à remplacer $G$ par $G/Z_{G,d}$, on peut donc supposer que le centre de $G$ est compact.

On va commencer par montrer la proposition dans le cas où le diagramme de Dynkin de $G$ est connexe.
Supposons d'abord $M=G$. On a alors également $M'=G$; de plus, il existe $w\in W'\left(G/T_0\right)$ tel que $K_x=w^{-1}K_{x'}w$. Soit $w'\in W'\cap K_x$ tel que $w'{}^{-1}w^{-1}Iww'=I$; $w'$ existe car les images de $I$ et de $w^{-1}Iw$ dans $K_x/K_x^1$ sont deux sous-groupes paraboliques minimaux de ce groupe. Alors $ww'$ normalise $I$, donc appartient à $\Gamma_I$; comme ce groupe est abélien et contient $\gamma$, $ww'$ et $\gamma$ commutent. On a donc $\left(ww'\right)^{-1}\left(\gamma,G,K_{x'}\right)ww'=\left(\gamma,G,K_x\right)$ et la proposition est démontrée dans ce cas.

Supposone maintenant $M$ et $M'$ propres.
Soit $\Delta=\left\{\alpha_1,\dots,\alpha_{n+1}\right\}$ l'ensemble des racines de $G$ relativement à $T_0$ correspondant aux murs $H_1,\dots,H_{n+1}$ de $A$ dans $\mathcal{A}$; $\Gamma$ agit sur $\Delta$ par l'intermédiaire de l'image canonique de $\Gamma_I$ dans le groupe de Weyl $W=N_G\left(T_0\right)/Z_G\left(T_0\right)$ de $G$ relativement à $T_0$. Pour toute partie $J$ de $\left\{1,\dots,n+1\right\}$, on notera $\Delta_J$ l'ensemble des $\alpha_j$, $j\in J$.

Soit $J$ (resp. $J'$) la partie de $\left\{1,\dots,n+1\right\}$ tel que $\Delta_J$ (resp. $\Delta_{J'}$) est un ensemble de racines simples de $M$ (resp. $M'$); $\Delta_J$ et $\Delta_{J'}$ sont stables par $\gamma$.

Soit $O_1,\dots,O_r$ les orbites de $\Delta$ pour l'action de $\gamma$. Pour tout $i$, on notera $K_i$ le sous-groupe parahorique de $G$ engendré par $I$ et par les éléments de $W'$ correspondant aux réflexions par rapport aux $H_j$, $j\not\in O_i$; $\gamma$ normalise $K_i$, et on a le lemme suivant:

\begin{lemme}\label{memepar}
Supposons qu'il existe $i\in\left\{1,\dots,r\right\}$ tel que $K$ et $K'$ sont tous deux inclus dans $K_i$, et conjugués par $w\in W'\cap K_i$; alors on peut supposer que $w$ commute avec $\gamma$.
\end{lemme}

\begin{proof}
En effet, considérons le groupe réductif fini non connexe ${\mth{G}}_i^+=\gamma^{\mth{Z}}K_i/K_i^1$; $W'\cap K_i$ s'identifie au groupe de Weyl de ${\mth{G}}_i=K_i/K_i^1$ relativement au tore déployé maximal ${\mth{T}}_0=T_0/T_0^1$. Les sous-groupes de Levi ${\mth{M}}^+$ et ${\mth{M}}'{}^+$ de ${\mth{G}}_i^+$, images respectivement de $\gamma^{\mth{Z}}K/K^1$ et $\gamma^{\mth{Z}}K'/K'{}^1$, sont conjugués, donc d'après le lemme \ref{centgamma}, en identifiant $\gamma$ à un de ses représentants quasi-centraux dans ${\mth{G}}_i^+$, si ${\mth{G}}_i^\gamma$ est le sous-groupe des éléments de ${\mth{G}}_i=K_i/K_i^1$ commutant avec $\gamma$,
les sous-groupes de Levi ${\mth{M}}^\gamma={\mth{M}}^+\cap{\mth{G}}_i^\gamma$ et ${\mth{M}}'{}^\gamma={\mth{M}}'{}^+\cap{\mth{G}}_i^\gamma$ de ${\mth{G}}_i^\gamma$ sont conjugués. Or ces sous-groupes sont semi-standard relativement au tore déployé maximal ${\mth{T}}_0^\gamma$ de ${\mth{G}}_i^\gamma$, donc il existe $w'\in W_\gamma=W\left({\mth{G}}_i^\gamma/{\mth{T}}_0^\gamma\right)$ tel que $w'{}^{-1}{\mth{M}}'{}^\gamma w'={\mth{M}}^\gamma$. On en déduit, toujours grâce au lemme \ref{centgamma}, que l'on a $w'{}^{-1}{\mth{M}}'w'={\mth{M}}$. Or $W_\gamma$ s'identifie à l'ensemble des éléments de $W'\cap K_i$ commutant avec $\gamma$, ce qui démontre le lemme.
\end{proof}

On va maintenant considérer les différents cas, suivant le type de $G$: puisque l'assertion du lemme ne dépend que de la structure de $W'$, il suffit de considérer le type du système de racines de $G$ relativement à $T_0$.

\begin{itemize}
\item Si $G$ est de type $A_n$, d'après \cite{bou}, il existe un entier $s$ divisant $n+1$ et tel que $\Gamma$ est cyclique d'ordre $\dfrac {n+1}s$, $r$ divise également $n+1$, est un multiple de $s$ et $\gamma$ est un élément d'ordre $\dfrac {n+1}r$ de $\Gamma$. Si $\alpha_1,\dots,\alpha_{n+1}$ sont numérotées dans l'ordre du diagramme de Dynkin, les $O_i$ sont les parties de $\Delta$ de la forme $\left\{\alpha_{i+rk}\mid 0\leq k<\dfrac{n+1}r\right\}$, et $J$ et $J'$ sont (modulo $n+1$) périodiques de période $r$.

Supposons d'abord qu'il existe un élément $j$ de $\left\{1,\dots,n+1\right\}$ n'appartenant ni à $J$ ni à $J'$; si $O_i$ est l'orbite contenant $\alpha_j$, $K$ et $K'$ sont contenus dans $K_i$, et $\left(M,K\right)$ et $\left(M',K'\right)$ sont conjugués par un élément de $W'\cap K_i$; c'est trivial si $K_i$ est un sous-groupe d'Iwahori de $G$, et si $K_i$ n'est pas un Iwahori, en posant ${\mth{G}}_i=K_i/K_i^1$, ${\mth{M}}=K/K^1$ et ${\mth{M}}'=K'/K'{}^1$, ${\mth{M}}$ et ${\mth{M}}'$ sont des sous-groupes de Levi de ${\mth{G}}_i$, et si $D_1,\dots,D_{\dfrac{n+1}r}$ sont les composantes connexes du diagramme de Dynkin de ${\mth{G}}_i$, puisque $\gamma$ agit transitivement sur les $D_k$, pour tout $k$, le diagramme de Dynkin de ${\mth{M}}$ (resp. ${\mth{M}}'$) est constitué de $\dfrac{n+1}r$ copies de son intersection avec $D_k$, et ces deux intersections sont isomorphes. De plus, $D_k$ est de type $A_{r-1}$, donc le groupe de Weyl de la composante correspondante de ${\mth{G}}_i$ est isomorphe au groupe symétrique $S_r$, et il est clair que deux sous-groupes de $S_r$ tous deux engendrés par des transpositions élémentaires et isomorphes entre eux sont conjugués. On en déduit que $K$ et $K'$ sont conjugués par un élément de $K_i$, et le lemme \ref{memepar} permet alors de conclure.

On va donc montrer que, quitte à conjuguer $\left(M',K'\right)$ par un élément de $W'$ commutant avec $\gamma$ et convenablement choisi, on peut toujours supposer qu'il existe un $j$ n'appartenant ni à $J$ ni à $J'$.

Le groupe $W'$ se plonge de manière canonique dans le groupe de Weyl affine $W'_n$ de $PGL_{n+1}\left(F\right)$ relativement à son tore diagonal, en identifiant, pour tout $j\in\left\{1,\dots,n+1\right\}$, $\alpha_j$ à la racine de ce tore donnée par:
\[\left(\lambda_1,\dots,\lambda_{n+1}\right)\longmapsto\lambda_j\lambda_{j+1}^{-1},\]
les indices étant pris modulo $n+1$; le groupe $W'$ s'identifie alors à l'ensemble des éléments de $W'_n$ dont la valuation du déterminant (qui est un élément de ${\mth{Z}}/\left(n+1\right){\mth{Z}}$) est multiple de $s$. Considérons l'ensemble:
\[L=\left\{j_1-j_2\mid j_1,j_2\not\in J\right\},\]
et soit $l$ le pgcd des éléments de $L$; on a le lemme suivant:

\begin{lemme}
Le groupe des valuations des déterminants des éléments de $N_{W'_n}\left(W'\cap K\right)$ est $l{\mth{Z}}/\left(n+1\right){\mth{Z}}$.
\end{lemme}

\begin{proof}
En effet, si $w\in N_{W'_n}\left(W'\cap K\right)$, la valuation du dé\-ter\-mi\-nant de $w$ est un multiple de $l$; réciproquement, il existe un élément de $N_{W'_n}\left(W'\cap K\right)$ dont la valuation du déterminant est $l$:
en effet, pour tout couple $\left(j_1,j_2\right)$, $j_1<j_2$, considérons l'élément:
\[\Diag\left(1,\dots,1,\varpi,\dots,\varpi,1,\dots,1\right)\]
du tore diagonal de $GL_{n+1}\left(F\right)$, où les termes valant $\varpi$ sont ceux d'indices $j_1+1$ à $j_2$ inclus; la valuation de son déterminant est $j_1-j_2$, et sa classe dans $PGL_{n+1}\left(F\right)$ normalise $W'\cap K$.
Le sous-groupe de ${\mth{Z}}/\left(n+1\right){\mth{Z}}$ des valuations des déterminants des éléments de $N_{W'_n}\left(W'\cap K\right)$ contient donc les classes de tous les $j_1-j_2$, donc contient celle de $l$ et le lemme est démontré.
\end{proof}

Revenons à la démonstration de la proposition.
Supposons d'abord qu'il existe un couple $\left(j,j'\right)$ et un entier $a$ tels que $j\not\in J$, $j'\not\in J'$ et $j-j'=la$. Puisque $l$ est le pgcd des éléments de $L$, on peut écrire (modulo $n+1$):
\[la=\sum_{k=1}^tl_i,\]
où les $l_i$ appartiennent à $L$; on supposera la décomposition choisie de façon à minimiser $t$. On va montrer l'assertion de la proposition par récurrence sur $t$, sachant que le cas $t=0$ correspond au cas $j=j'$ auquel on cherche à se ramener. Soit $j_1,j_2\not\in J$ tels que $j_1-j_2=l_1$, soit $O_i$ l'orbite de $\alpha_{j_1}$, et soit $w$ l'élément du groupe de Weyl de $K_i/K_i^1$ constitué de $r$ copies de l'élément:
\[\left(\begin{array}{ccccccc}&&&1\\&&&&\ddots\\&&&&&\ddots\\&&&&&&1\\1\\&\ddots\\&&1\end{array}\right)\left.\begin{array}{c}1\\\\\\l_1\\l_1+1\\\\n+1\end{array}\right.;\]
il est clair que $w$ commute avec $\gamma$. Posons $w'=\gamma_n^{j-j_1}w\gamma_n^{j_1-j}$; c'est un élément de $W'$ qui commute avec $\gamma$, et qui envoie $\left\{\alpha_j\mid j\in J\right\}$ sur une partie de $\Delta$ de la forme $\left\{\alpha_j\mid j\in J''\right\}$; $J''$ ne contient alors pas $j-l_1$, et on conclut en appliquant l'hypothèse de récurrence à $\left(M'',K''\right)=w'\left(M,K\right)w'{}^{-1}$ et $\left(M',K'\right)$.

Posons maintenant $l'=\pgcd\left(s,l\right)$ et supposons qu'il existe un couple $\left(j,j'\right)$ et un entier $b$ tels que $j\not\in J$, $j'\not\in J'$ et $j-j'=l'b$. Soit $c$ un entier tel que $cs+l'b$ est un multiple de $l$, et soit $\left(M'',K''\right)=\gamma ^{cs}\left(M,K\right)\gamma^{-cs}$, et $J''$ défini comme précédemment en remplaçant $w'$ par $\gamma^{cs}$; alors $J''$ ne contient pas $j+cs$, et on peut appliquer le cas précédent à $\left(M'',K''\right)$ et $\left(M',K'\right)$.

Supposons enfin qu'il n'existe aucun couple $\left(j,j'\right)$ tel que $j\not\in J$, $j'\not\in J'$ et $j-j'$ est multiple de $l'$. Alors d'après le lemme, la valuation du déterminant de tout élément $w'$ de $W'_n$ tel que $w'{}^{-1}\left(W'\cap K\right)w'=W'\cap K'$ est de la forme $j-j'+la$, donc non multiple de $l'$ sinon $j-j'$ le serait, et en particulier non multiple de $s$; un tel élément ne peut donc pas appartenir à $W'$, et on en déduit que $W'\cap K$ et $W'\cap K'$ ne sont pas conjugués dans $W'$. On arrive donc à une contradiction, ce qui achève la démonstration pour le cas $A_n$.

\item Supposons maintenant $G$ de type $B_n$, $n\geq 3$. Le diagramme de Dynkin étendu de $G$ est alors de la forme:

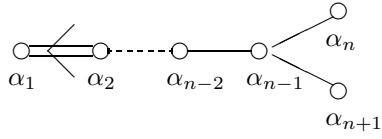
\begin{picture}(300,70)(0,0)
\put(145,20){\circle{6}}
\put(140,7){$\alpha_{n+1}$}
\put(145,50){\circle{6}}
\put(140,37){$\alpha_n$}
\put(143,22){\line(-2,1){23}}
\put(143,48){\line(-2,-1){23}}
\put(115,35){\circle{6}}
\put(110,22){$\alpha_{n-1}$}
\put(88,35){\line(1,0){24}}
\put(85,35){\circle{6}}
\put(80,22){$\alpha_{n-2}$}
\put(58,35){\dashbox{2}(24,0){}}
\put(55,35){\circle{6}}
\put(50,22){$\alpha_2$}
\put(28,33){\line(1,0){24}}
\put(28,37){\line(1,0){24}}
\put(35,35){\line(1,1){10}}
\put(35,35){\line(1,-1){10}}
\put(25,35){\circle{6}}
\put(20,22){$\alpha_1$}
\end{picture}

Si $\gamma$ est non trivial, d'après \cite{bou}, il est d'ordre $2$, permute $\alpha_n$ et $\alpha_{n+1}$ et fixe toutes les autres racines; $\alpha_n$ et $\alpha_{n+1}$ sont donc racines de $M$ si et seulement si elles sont racines de $M'$. De plus, $\alpha_1$, étant la seule racine simple courte, est racine de $M$ si et seulement si elle est racine de $M'$. Il y a donc trois cas possibles:
\begin{itemize}

\item si $\alpha_1$ n'est pas racine de $M$, $K$ et $K'$ sont contenus dans $K_1$; si $i$ est le plus grand élément de $\left\{1,\dots,n+1\right\}$ tel que $\alpha_i$ n'est pas racine de $M$, on a $K=K'$ si $i$ vaut $1$ ou $2$, et si $i>2$, $K$ et $K'$ sont conjugués par un élément de la composante de type $A_{i-2}$ du groupe de Weyl de $\left(K_1\cap K_i\right)/\left(K_1\cap K_i\right)^1$ (cette assertion se démontre de la même façon que l'assertion similaire qui se trouve dans la preuve du cas $A_n$);
le lemme \ref{memepar} permet alors de conclure;
\item si $\alpha_n$ et $\alpha_{n+1}$ ne sont pas racines de $M$, même chose en remplaçant $K_1$ par $K_n$ et en prenant pour $i$ le plus petit élément de $\left\{1,\dots,n+1\right\}$ tel que $\alpha_i$ n'est pas racine de $M$;
\item si $\alpha_1$, $\alpha_n$ et $\alpha_{n+1}$ sont racines de $M$, soit $j\in\left\{1,\dots,n+1\right\}$ maximal tel que $\alpha_j$ n'est pas racine de $M$; le diagramme de Dynkin de $M$ possède alors une composante de type $D_{n-j-1}$ (ou $A_1\times A_1$ avec $\alpha_n$ et $\alpha_{n+1}$ pour racines). Celui de $M'$ en possède alors également une, ce qui n'est possible que si $\alpha_j$ n'est pas non plus racine de $M'$; $K$ et $K'$ sont alors contenus dans $K_j$, et sont conjugués par un élément de la composante de type $B_{j-1}$ du groupe de Weyl de $K_j/K_j^1$; on termine comme dans le cas précédent.
\end{itemize}

\item Supposons maintenant $G$ de type $C_n$, $n\geq 2$. Le diagramme de Dynkin étendu de $G$ est alors de la forme:

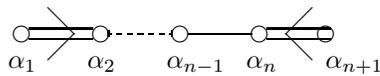
\begin{picture}(300,70)(0,0)
\put(140,35){\circle{6}}
\put(140,22){$\alpha_{n+1}$}
\put(118,33){\line(1,0){24}}
\put(118,37){\line(1,0){24}}
\put(125,35){\line(1,1){10}}
\put(125,35){\line(1,-1){10}}
\put(115,35){\circle{6}}
\put(110,22){$\alpha_{n}$}
\put(88,35){\line(1,0){24}}
\put(85,35){\circle{6}}
\put(80,22){$\alpha_{n-1}$}
\put(58,35){\dashbox{2}(24,0){}}
\put(55,35){\circle{6}}
\put(50,22){$\alpha_2$}
\put(28,33){\line(1,0){24}}
\put(28,37){\line(1,0){24}}
\put(45,35){\line(-1,1){10}}
\put(45,35){\line(-1,-1){10}}
\put(25,35){\circle{6}}
\put(20,22){$\alpha_1$}
\end{picture}

Si $\gamma$ est non trivial, d'après \cite{bou}, il est d'ordre $2$ et envoie, pour tout $j$, $\alpha_j$ sur $\alpha_{n+2-j}$; de plus, $\alpha_1$ et $\alpha_{n+1}$ étant les seules racines simples longues, $M$ les admet pour racines si et seulement si $M'$ les admet aussi. Il y a donc deux cas possibles:
\begin{itemize}
\item si $\alpha_1$ et $\alpha_{n+1}$ ne sont pas racines de $M$, $K$ et $K'$ sont contenus dans $K_1$, et sont conjugués par un élément de $W'\cap K_1$ qui est de type $A_{n-1}$; on conclut avec le lemme \ref{memepar};
\item si $\alpha_1$ et $\alpha_{n+1}$ sont racines de $M$, soit $j\in\left\{1,\dots,n+1\right\}$ minimal tel que $\alpha_j$ n'est pas racine de $M$; le diagramme de Dynkin de $M$ possède alors deux composantes de type $B_{j-1}$ (ou $A_1$ avec respectivement $\alpha_1$ et $\alpha_{n+1}$ pour racines), échangées par $\gamma$. Celui de $M'$ possède alors également de telles composantes, ce qui n'est possible que si $\alpha_j$ n'est pas non plus racine de $M'$; $K$ et $K'$ sont alors contenus dans $K_j$, et sont conjugués par un élément de la composante de type $A_{n+1-2j}$ du groupe de Weyl de $K_j/K_j^1$ si $j<\left[\dfrac n2+1\right]$, sachant que le cas $j=\left[\dfrac n2+1\right]$ impose $K=K'$ et que l'assertion du lemme est alors triviale; ici encore, on conclut grâce au lemme \ref{memepar}.
\end{itemize}

\item Supposons maintenant $G$ de type $D_n$, $n\geq 4$. Le diagramme de Dynkin étendu de $G$ est alors de la forme:

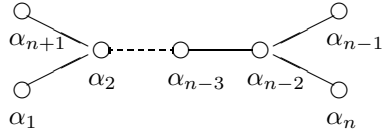
\begin{picture}(300,70)(0,0)
\put(145,20){\circle{6}}
\put(140,7){$\alpha_n$}
\put(145,50){\circle{6}}
\put(140,37){$\alpha_{n-1}$}
\put(143,22){\line(-2,1){23}}
\put(143,48){\line(-2,-1){23}}
\put(115,35){\circle{6}}
\put(110,22){$\alpha_{n-2}$}
\put(88,35){\line(1,0){24}}
\put(85,35){\circle{6}}
\put(80,22){$\alpha_{n-3}$}
\put(58,35){\dashbox{2}(24,0){}}
\put(55,35){\circle{6}}
\put(50,22){$\alpha_2$}
\put(27,22){\line(2,1){23}}
\put(27,48){\line(2,-1){23}}
\put(25,20){\circle{6}}
\put(20,7){$\alpha_1$}
\put(25,50){\circle{6}}
\put(20,37){$\alpha_{n+1}$}
\end{picture}

Supposons d'abord $n$ impair. Si $\gamma$ est non trivial, d'après \cite{bou}, il est d'ordre $2$ ou $4$, et $\Gamma$ est un groupe cyclique de cardinal $2$ ou $4$.

\begin{itemize}
\item Si $\gamma$ est d'ordre $4$, ses orbites sont $O_1=\left\{\alpha_1,\alpha_{n-1},\alpha_n,\alpha_{n+1}\right\}$, et pour tout $i\leq\dfrac{n-1}2$, $O_i=\left\{\alpha_i,\alpha_{n-i}\right\}$. Si $\Delta_J$ ne contient pas $O_1$, $\Delta_{J'}$ non plus; $K$ et $K'$ sont alors contenus dans $K_1$, et sont conjugués par un élément de $W'\cap K_1$ commutant avec $\gamma$ puisque $K_1/K_1^1$ est de type $A_{n-3}$ et grâce au lemme \ref{memepar}. Si $\Delta_J$ contient $O_1$, soit $i\in\left\{1,\dots,n-1\right\}$ minimal tel que $\Delta_i\not\in J$. Si $i=2$ (resp. si $i\geq 3$), le diagramme de Dynkin de $M$ comporte quatre composantes de type $A_1$ (resp. deux composantes de type $D_i$) permutées par $\gamma$ (resp. permutées par $\gamma$ et sur lesquelles l'action de $\gamma^2$ est non triviale); c'est donc également le cas pour $M'$, ce qui n'est possible que si $\Delta_{J'}$ contient tous les $O_{i'}$, $i'\leq i$, mais pas $O_i$. $K$ et $K'$ sont alors tous deux contenus dans $K_i$, et sont conjugués par un élément de la composante de type $A_{n-1-2i}$ du groupe de Weyl de $K_i/K_i^1$ si $i<\dfrac{n-1}2$, sachant que le cas $i=\dfrac{n-1}2$ impose $K=K'$ et que l'assertion du lemme est alors triviale. On conclut grâce au lemme \ref{memepar}.
\item Si maintenant $\gamma$ est d'ordre $2$, ses orbites sont $O_1=\left\{\alpha_1,\alpha_{n+1}\right\}$, $O_{n-1}=\left\{\alpha_{n-1},\alpha_n\right\}$, et pour tout $i\in\left\{2,\dots,n-2\right\}$, $O_i=\left\{\alpha_i\right\}$.
Si $\Delta_J$ ne contient ni $O_1$ ni $O_{n-1}$, $\Delta_{J'}$ non plus; $K$ et $K'$ sont alors contenus dans $K_1\cap K_{n-1}$ et on conclut comme d'habitude puisque $\left(K_1\cap K_{n-1}\right)/\left(K_1\cap K_{n-1}\right)^1$ est de type $A_{n-3}$.

Supposons ensuite que $\Delta_J$ contient soit $O_1$ soit $O_{n-1}$; c'est donc é\-ga\-le\-ment le cas de $\Delta_{J'}$.
Soit $i$ (resp. $i'$) minimal tel que $O_i$ n'est pas contenu dans $\Delta_J$ (resp. $\Delta_{J'}$). Si $i=i'$, alors $K$ et $K'$ sont contenus dans $K_i$, et sont égaux si $i=n-1$ ou $i=n-2$; si $i<n-2$, soit $i_1$ maximal tel que $O_{i_1}$ n'est pas contenu dans $\Delta_J$; $i_1$ est alors également l'indice maximal tel que $O_{i_1}$ n'est pas contenu dans $\Delta_{J'}$; de plus, si $i_1=i$ ou $i_1=i+1$, $K$ et $K'$ sont égaux, et sinon, ils sont conjugués par un élément de la composante de type $A_{i_1-i-1}$ du groupe de Weyl de $K_i/K_i^1$, et on conclut encore une fois grâce au lemme \ref{memepar}. Si $i\neq i'$, alors dans tous les cas, $\Delta_{J'}$ contient tous les $O_{n-j}$, $j>i$, et pas $O_{n-i}$. Si $\Gamma$ est de cardinal $4$ et si $\gamma'$ est un élément de $\Gamma$ distinct de $1$ et de $\gamma$, en remplaçant $\left(M',K'\right)$ par $\gamma'{}^{-1}\left(M',K'\right)\gamma$, on est ramené au cas $i=i'$.
\item Il reste à traiter le cas où $\Gamma$ est de cardinal $2$ et où $i\neq i'$; supposons par exemple $i>i'$.

Soit $L_1,\dots,L_t$ des sous-groupes parahoriques de $M$ contenant $I\cap M$ et contenus dans $K$ vérifiant les propriétés suivantes:
\begin{itemize}
\item on a $W'\cap K=\prod_{i=1}^t\left(W'\cap L_i\right)$;
\item pour tout $i\in\left\{2,\dots,t-1\right\}$, le diagramme de Dynkin de $L_i/L_i^1$ est non vide et connexe;
\item le diagramme de Dynkin de $L_1$ (resp. $L_t$) correspond à la réunion des $O_j$, $j<i$ (resp. $j>n-i'$).
\end{itemize}
De tels sous-groupes sont uniques à permutation des $L_i$, $2\leq i\leq t-1$, près. On définira de même $L'_1,\dots,L'_t$ en remplaçant $K$ par $K'$.

Soit $GSO'_{2n}$ la forme déployée de $GSO_{2n}$; ce groupe peut être vu comme le sous-groupe des éléments $g$ de $GL_{2n}$ tels que ${}^\tau g g$ est une homothétie, où ${}^\tau g$ désigne la transposée de $g$ suivant la diagonale non principale. Soit $PGSO'_{2n}$ le quotient de $GSO'_{2n}$ par son centre; c'est un groupe semi-simple, adjoint et déployé de type $D_n$. Soit $B_{2n}$ le sous-groupe de Borel des éléments triangulaires supérieurs de $PGSO'_{2n}$, $T_{2n}$ le tore diagonal de $B_{2n}$ et $\alpha'_1,\dots,\alpha'_n$ les racines simples de $PGSO'_{2n}$ relativement à $B_{2n}$ et $T_{2n}$, numérotées de façon similaire à $\alpha_1,\dots,\alpha_n$ sur le diagramme de Dynkin.

Si $p\neq 2$, on peut plonger $W'$ dans le groupe de Weyl affine $W'_n$ de $PGSO'_{2n}\left(F\right)$ relativement au tore diagonal (si $p=2$, cela marche encore en remplaçant $PGSO'_{2n}\left(F\right)$ par $PGSO'_{2n}\left(F'\right)$, où $F'$ est de caractéristique résiduelle différente de $2$; comme on s'intéresse uniquement à la structure de $W'$, ce qui suit reste valide);
un calcul matriciel montre que
si $w'\in W'_n$ normalise $W'\cap K$, on peut écrire $w'=w'_1w''w'_t$, où $w'_1\in W'\cap L_1$, $w'_t\in W'\cap L_t$ et $w''\in W'\cap M_0$, où $M_0$ est le sous-groupe de Levi standard de $PGSO'_{2n}\left(F\right)$ (relativement à $B_{2n}$ et $T_{2n}$) admettant $\left\{\alpha'_{i+1},\dots,\alpha'_{n-1-i'}\right\}$ comme ensemble de racines simples. De plus, si $I_n$ est le sous-groupe d'Iwahori standard de $PGSO'_{2n}\left(F\right)$ (c'est-à-dire le groupe des éléments de $PGSO'_{2n}\left(\mathcal{O}\right)$ triangulaires supérieurs modulo ${\mathfrak{p}}$) et si $\gamma'$ est un générateur de $N_G\left(I_n\right)/I_n$, $W'\cap K$ et $W'\cap K'$ sont conjugués par un élément de $\left(W'\cap M_0\right)\gamma'{}^{\mth{Z}}$; on en déduit que si $w\in W'_n$ est tel que $w^{-1}\left(W'\cap K\right)w=W'\cap K'$, il existe $w_1\in W'\cap L_1$, $w'_t\in W'\cap L_t$ et un entier $l$ tels que si $w''=ww'_1{}^{-1}w'_t{}^{-1}$, $w''{}^{-1}\left(W'\cap K'\right)w''=W'\cap K$ et $w''{}^{-1}\left(W'\cap M_0\right)w''=W'\cap\gamma'{}^{-l}M_0\gamma'{}^l$, d'où:
\[w''{}^{-1}\left(W'\cap L'_2L'_3\dots L'_{t-1}\right)w''=W'\cap L_2L_3\dots L_{t-1};\]
de plus, si $w$ appartient à $W'$, $w''$ aussi, et puisque $w'_1$ et $w'_t$ normalisent $W'\cap K$, on a $w''{}^{-1}\left(W'\cap K'\right)w''=W'\cap K$. Enfin, $L_2\dots L_{t-1}$ et $L'_2\dots L'_{t-1}$ sont contenus dans $K_1\cap K_{n-1}$; comme l'action de $\gamma$ est triviale sur $W'\cap K_1\cap K_{n-1}$, $w''$ commute avec $\gamma$ et l'assertion du lemme est démontrée.

\end{itemize}

Supposons maintenant $n$ pair. $\Gamma$ est un sous-groupe de $\left({\mth{Z}}/2{\mth{Z}}\right)^2$, et $\gamma$ est toujours d'ordre $2$ s'il est non trivial. On a les cas suivants:
\begin{itemize}
\item si $\gamma$ envoie $\alpha_1$ sur $\alpha_{n+1}$, ses orbites sont $O_1=\left\{\alpha_1,\alpha_{n+1}\right\}$, $O_{n-1}=\left\{\alpha_{n-1},\alpha_n\right\}$ et pour tout $i\in\left\{2,\dots,n-2\right\}$, $O_i=\left\{\alpha_i\right\}$; la suite de la démonstration est identique au cas $n$ impair et $\gamma$ d'ordre $2$;
\item si $\gamma$ envoie $\alpha_1$ sur $\alpha_{n-1}$, ses orbites sont $O_1=\left\{\alpha_1,\alpha_{n-1}\right\}$, $O_i=\left\{\alpha_i,\alpha_{n-i}\right\}$ pour tout $i\in\left\{2,\dots,\dfrac n2-1\right\}$, $O_{\dfrac n2}=\{\alpha_{\dfrac n2}\}$ et $O_{\dfrac n2+1}=\left\{\alpha_n,\alpha_{n+1}\right\}$. Si $\Delta_J$ contient $O_1$ et $O_{\dfrac n2+1}$, $\Delta_{J'}$ également et ce cas ce traite comme celui où $n$ est impair et $\gamma$ d'ordre $4$; même chose si $\Delta_J$ ne contient ni $O_1$ ni $O_{\dfrac n2+1}$. Supposons que $\Delta_J$ contient $O_1$ et pas $O_{\dfrac n2+1}$; s'il en est de même de $\Delta_J$, $K$ et $K'$ sont tous deux contenus dans $K_{\dfrac n2+1}$, dont le diagramme de Dynkin est de type $A_{n-1}$, et on conclut comme d'habitude. Si maintenant $\Delta_{J'}$ contient $O_{\dfrac n2+1}$ et pas $O_1$, si $\Gamma$ est de cardinal $4$, en conjugant $\left(M',K'\right)$ par un élément de $\Gamma$ non nul et distinct de $\gamma$, on est ramené au cas précédent. Supposons donc $\Gamma$ de cardinal $2$; en considérant des sous-groupes parahoriques $L_1,\dots,L_t$ (resp. $L'_1,\dots,L'_t$) de $G$ tels que l'on a:
\[W'\cap K=\prod_{i=1}^t\left(W'\cap L_i\right)\]
(resp.
\[W'\cap K'=\prod_{i=1}^t\left(W'\cap L'_1\right))\]
que pour tout $i$, le diagramme de Dynkin de $L_i/L_i^1$ (resp. $L'_i/L'_i{}^1$) est stable par $\gamma$ et est soit connexe, soit constitué de deux composantes connexes permutées par $\gamma$, et que celui de $L_1/L_1^1$ (resp. $L'_1/L'_1{}^1$) contient $O_1$ (resp. $O_{\dfrac n2+1}$), ce cas se traite de façon analogue au cas $n$ impair et $\gamma$ d'ordre $2$.
\item Enfin, le cas où $\gamma$ envoie $\alpha_1$ sur $\alpha_n$ est similaire au précédent.
\end{itemize}

\item Si $G$ est de type $E_6$, son diagramme de Dynkin étendu est de la forme suivante:

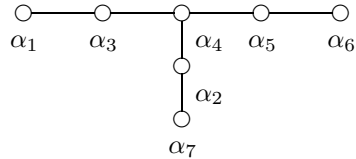
\begin{picture}(300,100)(0,0)
\put(25,65){\circle{6}}
\put(20,52){$\alpha_1$}
\put(28,65){\line(1,0){24}}
\put(55,65){\circle{6}}
\put(50,52){$\alpha_3$}
\put(58,65){\line(1,0){24}}
\put(85,65){\circle{6}}
\put(90,52){$\alpha_4$}
\put(85,62){\line(0,-1){14}}
\put(85,45){\circle{6}}
\put(90,32){$\alpha_2$}
\put(85,42){\line(0,-1){14}}
\put(85,25){\circle{6}}
\put(80,12){$\alpha_7$}
\put(88,65){\line(1,0){24}}
\put(115,65){\circle{6}}
\put(110,52){$\alpha_5$}
\put(118,65){\line(1,0){24}}
\put(145,65){\circle{6}}
\put(140,52){$\alpha_6$}
\end{picture}

Si $\gamma$ est non trivial, d'après \cite{bou}, il est d'ordre $3$, il engendre $\Gamma$ et on a $O_1=\left\{\alpha_1,\alpha_6,\alpha_7\right\}$, $O_2=\left\{\alpha_2,\alpha_3,\alpha_5\right\}$, $O_3=\left\{\alpha_4\right\}$. Le seul cas où $\left(M,K\right)$ est différent de $\left(M',K'\right)$ est celui où $\Delta_J=O_1$ et $\Delta_{J'}=O_2$ (ou l'inverse); $K$ et $K'$ sont alors contenus dans $K_4$ et sont conjugués par un élément de $w'\in W'\cap K_4$ commutant avec $\gamma$ ($w'$ envoie chaque élément de $O_1$ sur l'élément de $O_2$ situé dans la même composante connexe du diagramme de Dynkin de $K_4/K_4^1$).

\item Si $G$ est de type $E_7$, son diagramme de Dynkin étendu est de la forme suivante:

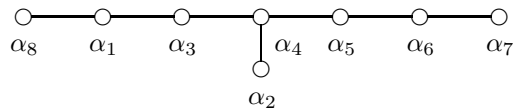
\begin{picture}(300,70)(0,0)
\put(25,45){\circle{6}}
\put(20,32){$\alpha_8$}
\put(28,45){\line(1,0){24}}
\put(55,45){\circle{6}}
\put(50,32){$\alpha_1$}
\put(58,45){\line(1,0){24}}
\put(85,45){\circle{6}}
\put(80,32){$\alpha_3$}
\put(88,45){\line(1,0){24}}
\put(115,45){\circle{6}}
\put(120,32){$\alpha_4$}
\put(115,42){\line(0,-1){14}}
\put(115,25){\circle{6}}
\put(110,12){$\alpha_2$}
\put(118,45){\line(1,0){24}}
\put(145,45){\circle{6}}
\put(140,32){$\alpha_5$}
\put(148,45){\line(1,0){24}}
\put(175,45){\circle{6}}
\put(170,32){$\alpha_6$}
\put(178,45){\line(1,0){24}}
\put(205,45){\circle{6}}
\put(200,32){$\alpha_7$}
\end{picture}

Si $\gamma$ est non trivial, d'après \cite{bou}, il est d'ordre $2$, il engendre $\Gamma$ et on a $O_1=\left\{\alpha_8,\alpha_7\right\}$, $O_2=\left\{\alpha_1,\alpha_6\right\}$, $O_3=\left\{\alpha_3,\alpha_5\right\}$, $O_4=\left\{\alpha_4\right\}$, $O_5=\left\{\alpha_2\right\}$. On a trois cas à considérer:
\begin{itemize}
\item si $\Delta_J$ et $\Delta_{J'}$ ne contiennent pas $O_3$, alors $K$ et $K'$ sont contenus dans $K_3$, et sont dans tous les cas conjugués par un élément de $W'\cap K_3$; on conclut grâce au lemme \ref{memepar};
\item si $\Delta_J$ contient à la fois $O_3$ et $O_4$, le diagramme de Dynkin de $M$, et donc également celui de $M'$, est de l'un des types suivants: $A_3$, $D_4$, $A_1\times A_3\times A_1$, $A_1\times D_4\times A_1$, $A_5$, $E_6$, $A_7$. Or on vérifie aisément qu'il existe un unique sous-diagramme stable par $\gamma$ du diagramme de Dynkin étendu de $E_7$ de chacun de ces types, et donc le seul choix possible pour $\left(M',K'\right)$ est $\left(M',K'\right)=\left(M,K\right)$;
\item si $\Delta_J$ contient $O_3$ et pas $O_4$, le diagramme de Dynkin de $M$, et donc également celui de $M'$, est de l'un des types suivants: $A_1\times A_1$, $A_1\times A_1\times A_1$, $A_2\times A_2$, $A_2\times A_1\times A_2$, $A_1\times A_1\times A_1\times A_1$, $A_1\times A_1\times A_1\times A_1\times A_1$, $A_3\times A_3$, $A_3\times A_1\times A_3$. Dans les quatre derniers cas, on conclut comme précédemment. Dans les cas $A_1\times A_1$ et $A_2\times A_2$, $\Delta_{J'}$ ne peut pas contenir $O_4$; $K$ et $K'$ sont donc tous deux contenus dans $K_4$ et on conclut comme d'habitude. Dans le cas $A_1\times A_1\times A_1$, on a $\Delta_J=O_3\cup O_5$, et $\Delta_{J'}$ ne peut pas contenir à la fois $O_1$ et $O_2$; $K$ et $K'$ sont alors tous deux contenus soit dans $K_1$, soit dans $K_2$, et on conclut toujours de la même manière. Dans le cas $A_2\times A_1\times A_2$, il est facile de voir que le seul cas où on ne peut pas faire de même est le cas où $\Delta_J=O_2\cup O_3\cup O_5$ et $\Delta_{J'}=O_1\cup O_2\cup O_4$; mais dans ce cas, si $x''$ est un point de la face de $A$ correspondant à $O_1\cup O_2\cup O_5$ et si $\left(M'',K''\right)=\kappa_0\left(K_{x''},T_0\right)$, il suffit d'appliquer le lemme successivement à $\left(\gamma,M,K\right)$ et $\left(\gamma,M'',K''\right)$, puis à $\left(\gamma,M'',K''\right)$ et $\left(\gamma,M',K'\right)$.
\end{itemize}

\item Si $G$ est de type $E_8$, $F_4$, $G_2$ ou $BC_n$, $n\geq 1$, d'après \cite{bou}, on a toujours $\Gamma=\left\{1\right\}$, d'où $\gamma=1$ et le résultat de la proposition est trivial.

\end{itemize}

Soit $W'_0=W'\cap G^1$: montrons maintenant le lemme suivant:

\begin{lemme}
Supposons que le diagramme de Dynkin de $G$ est connexe. Si $w'$ est un élément de $W'$ tel que $w'{}^{-1}\left(M',K'\right)w'=\left(M,K\right)$, il existe $w\in W'$ congru à $w'$ modulo $W'_0$ tel que $w^{-1}\left(\gamma,M',K'\right)w=\left(\gamma,M,K\right)$.
\end{lemme}

\begin{proof}
En effet, soit $\Gamma_K$ (resp. $\Gamma_{\gamma,K}$) l'image dans $\Gamma$ du sous-groupe $N_{W'}\left(W'\cap K\right)$ de $W'$ (resp. du sous-groupe $N_{W'}\left(\gamma,W'\cap K\right)$ des éléments de $W'$ qui commutent avec $\gamma$ et normalisent $W'\cap K$). Il est clair que $\Gamma_{\gamma,K}\subset\Gamma_K$; on va montrer que l'on a en fait l'égalité, ce qui suffira à démontrer le lemme puisque, en partant d'un $w$ quelconque vérifiant $w^{-1}\left(\gamma,M',K'\right)w=\left(\gamma,M,K\right)$, puisque $w\in w'N_{W'}\left(W'\cap K\right)$, on pourra alors toujours rendre $w$ congru à $w'$ modulo $W'_0$ en le multipliant à droite par un élément de $N_{W'}\left(\gamma,W'\cap K\right)$. Supposons donc qu'il existe $\gamma_0\in\Gamma_K$ n'appartenant pas à $\Gamma_{\gamma,K}$; $W'\cap K$ et $W'\cap\gamma_0^{-1}K\gamma_0$ sont alors conjugués par un élément de $W'_0$, mais tout élément $w$ de $W'$ tel que $w^{-1}\left(W'\cap\gamma_0^{-1}K\gamma_0\right))w=W'\cap K$ appartient à $\gamma_0^{-1}N_{W'}\left(\gamma,W'\cap K\right)$, qui ne rencontre pas $\Gamma_{\gamma,K}W'_0$. Or d'après ce qui précède, si $W'\cap K$ et $W'\cap\gamma_0^{-1}K\gamma$ sont conjugués dans $\Gamma_{\gamma,K}W'_0$, ils le sont par un élément commutant avec $\gamma$; on aboutit donc à une contradiction, ce qui démontre l'assertion cherchée.
\end{proof}

Supposons maintenant que $G$ est produit direct de groupes dont le diagramme de Dynkin est connexe; l'assertion de la proposition, ainsi que celle du lemme précédent, se déduisent trivialement du cas précédent.

Supposons enfin $G$ quelconque; en posant $T=Z_G\left(T_0\right)$, soit $\overline{\Phi}$ (resp. $\overline{\phi^\vee}$, le système de racines (resp. de poids radiciels) de $\underline{G}$ associé à $\underline{T}$, et $\overline{\Phi_1},\dots,\overline{\Phi_t}$ (resp. $\overline{\Phi_1^\vee},\dots,\overline{\Phi_t^\vee}$) la partition de $\overline{\Phi}$ (resp. $\overline{\Phi^\vee}$) correspondant aux composantes connexes du diagramme de Dynkin de $\underline{G}$. Pour tout $i\in\left\{1,\dots,t\right\}$, soit $\underline{G_i}$ le sous-groupe fermé de $\underline{G}$ engendré par $\underline{T}$ et par les sous-groupes radiciels $U_\alpha$, $\alpha\in\overline{\Phi_i}$,et soit $\underline{T_i}$ le sous-tore de $\underline{T}$ engendré par les images des $\xi\in\overline{\Phi^\vee}-\overline{\Phi_i^\vee}$; $\underline{G_i}$ et $\underline{T_i}$ sont clairement définis sur $F$, et $\underline{T_i}$ est central dans $\underline{G_i}$. Soit $G_i=\underline{G_i}\left(F\right)$, $T_i=\underline{T_i}\left(F\right)$ et $H_i=G_i/\left(T_i\cap T_0\right)$; $H_i$ est le groupe des $F$-points d'un groupe algébrique, et on a une injection canonique:
\[\phi: G\longmapsto H_1\times\dots\times H_t.\]
Soit, pour tout $i$, $K_{x,i}$, $K_{x',i}$, $M_i$, $M'_i$, $K_i$, $K'_i$ les images dans $H_i$ des intersections respectives de $K_x$, $K_{x'}$, $M$, $M'$, $K$, $K'$ avec $G_i$; on vérifie aisément que $K_{x,i}$ et $K_{x',i}$ sont des sous-groupes parahoriques de $H_i$ associés et que l'on a, avec un léger abus de notation, $\kappa_0\left(K_{x,i},T_0\right)=\left(M_i,K_i\right)$ et $\kappa_0\left(K_{x',i},T_0\right)=\left(M'_i,K'_i\right)$. De plus, si pour tout $i$, $W'_i$ est le groupe de Weyl affine de $H_i$ relqtivement à son tore déployé maximal $T_0/\left(T_i\cap T_0\right)$, on a une injection canonique de $W'$ dans le produit direct des $W'_i$, et si on pose:
\[\gamma=\gamma_1\dots\gamma_t,\]
avec pour tout $i$, $\gamma_i\in W'_i$, alors pour tout $i$, $\gamma_i$ normalise $K_{x,i}$, $K_{x',i}$, $M_i$, $M'_i$, $K_i$ et $K'_i$;
puisque le diagramme de Dynkin de $H_i$ est connexe, on peut lui appliquer la proposition, et obtenir $w_i\in W'_i$ tel que l'on a $w_i\left(\gamma_i,M_i,K_i\right)w_i^{-1}=\left(\gamma_i,M'_i,K'_i\right)$; en posant $w=\left(w_1,\dots,w_t\right)$, on obtient dans $H_1\times\dots\times H_t$:
\[w\left(\gamma,M,K\right)w^{-1}=\left(\gamma,M',K'\right).\]
Il reste donc à vérifier que l'on peut choisir $w$ dans $W'$. Or si $w'$ est un élément de $W'\subset W'_1\times\dots\times W'_t$ tel que $w'\left(M,K\right)w'{}^{-1}=\left(M',K'\right)$, d'après ce qui précède, on peut supposer que $w$ est dans la même classe que $w'$ modulo $\left(W'_1\times\dots\times W'_t\right)\cap\left(H_1\times\dots\times H_t\right)^1$; or ce dernier groupe contient $W'\cap G^1$, et est un groupe de Coxeter affine de même type que celui-ci; ces deux groupes sont donc égaux, et on en déduit que si $w'\in W'$, $w$ aussi, ce qui achève la démonstration. $\Box$

\subsection{Distributions associées à des fonctions $\gamma$-cuspidales}

Soit $\mu\in\mathcal{M}$, $\left(M,K_M\right)$ un représentant de $\mu$, et ${\mth{G}}=K_M/K_M^1$ le groupe réductif fini associé à $\mu$; si $\phi$ est une fonction sur ${\mth{G}}$ invariante par conjugaison, $\phi$ est {\em cuspidale} si pour tout sous-groupe parabolique propre ${\mth{P}}={\mth{M}}{\mth{U}}$ de ${\mth{G}}$  et tout $m\in{\mth{M}}$, on a:
\[\sum_{u\in{\mth{U}}}\phi\left(mu\right)=0.\]

La notion de fonction cuspidale peut se généraliser de la façon suivante:
soit $\nu=\left(\gamma,\mu\right)\in\mathcal{N}$, et ${\mth{G}}^+=\gamma^{\mth{Z}}{\mth{G}}$ le groupe associé.
Soit $\phi$ une fonction de ${\mth{G}}^+$ dans ${\mth{C}}$ invariante par conjugaison; on dira que $\phi$ est {\em cuspidale} si pour tout sous-groupe parabolique propre ${\mth{P}}^+={\mth{M}}^+{\mth{U}}$ de ${\mth{G}}^+$ et tout $m\in{\mth{M}}^+$, on a:
\[\sum_{u\in{\mth{U}}}\phi\left(mu\right)=0.\]

Soit $g,h\in{\mth{G}}$; $g$ et $h$ seront dit {\em $\gamma$-conjugués} si les éléments $\gamma g$ et $\gamma h$ de ${\mth{G}}^+$ sont conjugués.
Si $\phi$ est une fonction de ${\mth{G}}$ dans ${\mth{C}}$ invariante par $\gamma$-conjugaison, $\phi$ sera dite {\em $\gamma$-cuspidale} si la fonction $\phi'$ de ${\mth{G}}$ dans ${\mth{C}}$ définie par $\phi'\left(g\right)=\phi\left(\gamma^{-1}g\right)$ si $g\in\gamma{\mth{G}}$ et $\phi'\left(g\right)=0$ sinon est cuspidale.

Montrons d'abord qu'il suffit, dans la définition précédente, de considérer les sous-groupes paraboliques de ${\mth{G}}^+$ stables par $\gamma$; cela se déduit immédiatement du lemme suivant:

\begin{lemme}
Pour tout sous-groupe parabolique ${\mth{P}}^+$ de ${\mth{G}}^+$ tel que ${\mth{P}}^+\cap\gamma{\mth{G}}\neq\emptyset$, il existe un sous-groupe parabolique ${\mth{P}}'{}^+$ de ${\mth{G}}^+$ conjugué à ${\mth{P}}^+$ et stable par conjugaison par $\gamma$.
\end{lemme}

\begin{proof}
En effet, soit ${\mth{P}}_0$ un sous-groupe parabolique minimal de ${\mth{G}}$ stable par conjugaison par $\gamma$; il en existe d'après \cite[proposition 1.36]{dm}, et puisque tous les sous-groupes paraboliques minimaux de ${\mth{G}}$ sont conjugués entre eux, ${\mth{P}}^+$ contient un conjugué de ${\mth{P}}_0$. Soit donc $h\in{\mth{G}}$ tel que $h{\mth{P}}_0h^{-1}\subset{\mth{P}}^+$; posons:
\[{\mth{P}}'{}^+=h^{-1}{\mth{P}}^+h.\]
Alors $\gamma^{-1}{\mth{P}}'{}^+\gamma$ contient $\gamma^{-1}{\mth{P}}_0\gamma={\mth{P}}_0$; or puisque ${\mth{P}}^+$ rencontre $\gamma{\mth{G}}$, ${\mth{P}}'{}^+$ aussi, donc ${\mth{P}}'{}^+$ et $\gamma^{-1}{\mth{P}}'{}^+\gamma$ sont conjugués entre eux par un élément de ${\mth{G}}$; si ${\mth{P}}'={\mth{P}}'{}^+\cap{\mth{G}}$, ${\mth{P}}'$ et $\gamma^{-1}{\mth{P}}'\gamma$ sont donc deux sous-groupes paraboliques de ${\mth{G}}$ conjugués entre eux dans ${\mth{G}}$; puisqu'ils contiennent un même sous-groupe parabolique minimal de ${\mth{G}}$, ils sont égaux, donc leurs normalisateurs dans ${\mth{G}}^+$ le sont aussi et le lemme est démontré.
\end{proof}

Supposons qu'il existe $x\in\mathcal{B}$ tel que $\gamma\in\Gamma_x$ et $\zeta\left(\gamma,x\right)=\left(\gamma,\mu\right)$; on peut relever $\phi$ en une fonction sur $\gamma K_x$ invariante par conjugaison par $K_x$, que par un léger abus de notation on notera également $\phi$; on étendra $\phi$ à $G$ tout entier en posant $\phi\left(g\right)=0$ pour $g\not\in\gamma K_x$.

Soit $\mathcal{A}$ un appartement de $\mathcal{B}$ contenant $x$ et $T_0$ le tore déployé maximal de $G$ associé à $\mathcal{A}$.
Pour définir une distribution sur l'ensemble des conjugués de $\gamma K_x$ à partir de $\phi$, on va utiliser un sous-groupe de Levi de $G$ semi-standard relativement à $T_0$ un peu plus gros que $M$.
Soit $V$ l'espace vectoriel associé à $\mathcal{A}$; on a une décomposition de $V$ en facettes vectorielles pour l'action de $W$, et tout élément de $W$ qui fixe un point d'une facette donnée la fixe toute entière; on en déduit que si $V^\gamma$ est le sous-espace des éléments de $V$ fixés point par point par l'élément $w_\gamma$ de $W$ image de $\gamma$ par la surjection canonique de $W'$ dans $W$, $V^\gamma$ est réunion des facettes vectorielles qu'il rencontre.

Considérons l'immeuble sphérique $\mathcal{B}_s$ de $G$; $V$ s'identifie canoniquement à l'appartement de $\mathcal{B}_s$ associé à $T_0$. Soit ${\mathfrak{A}}$ une facette vectorielle de $V^\gamma\cap V_M$ de dimension maximale; l'ensemble des points de $G$ qui la fixent est un sous-groupe parabolique $P_\gamma$ de $G$ semi-standard relativement à $T_0$. Posons $P_\gamma=M_\gamma U_\gamma$, avec $M_\gamma$  également semi-standard relativement à $T_0$; le groupe de Weyl $W_{M_\gamma}$ de $M_\gamma$ fixe point par point le sous-espace de $V$ engendré par ${\mathfrak{A}}$, qui n'est autre que $V^\gamma\cap V_M$. Réciproquement, un élément de $W$ qui fixe cet espace appartient à $W\cap P_\gamma=W_{M_\gamma}$; on en déduit que $W_{M_\gamma}$, et donc $M_\gamma$, ne dépendent pas du choix de ${\mathfrak{A}}$, et que $M_\gamma$ est le fixateur de $V^\gamma\cap V_M$. En particulier, il contient $M$ et $\gamma$, et si $M'$ est un sous-groupe de Levi de $G$ semi-standard relativement à $T_0$ contenant $M$ et $\gamma$, l'ensemble des points de $\mathcal{B}_s$ qu'il fixe est contenu dans $V^\gamma\cap V_M$, donc $M_\gamma\subset M'$.
$M_\gamma$ est donc le plus petit sous-groupe de Levi semi-standard de $G$ (relativement à $T_0$) contenant à la fois $M$ et $\gamma$.

Soit $K_\gamma=K_x\cap M_\gamma$; $K_\gamma/K_\gamma^1$ est isomorphe à ${\mth{G}}$ et $\gamma^{\mth{Z}}K_\gamma/K_\gamma^1$ à ${\mth{G}}^+$; si $k$ est un élément de $\gamma^{\mth{Z}}K_\gamma$ (resp. si $X$ est une partie de $\gamma^{\mth{Z}}K_\gamma$), on notera $\overline{k}$ (resp. $\overline{X}$) son image dans ${\mth{G}}^+$. Pour pouvoir définir des distributions à support dans l'ensemble des conjugués de $\gamma K_\gamma$ dans $M_\gamma$, on va d'abord montrer la proposition suivante:

\begin{prop}\label{cspdp2}
Soit $H$ un sous-groupe ouvert de $M_\gamma$, $C$ une partie compacte de $M_\gamma$. Alors il existe une partie $C'$ de $M_\gamma$, compacte modulo le centre $Z_{M_\gamma}$ de $M_\gamma$, telle que pour tout $m\in M_\gamma-C'$, l'une des deux conditions suivantes soit vérifiée:
\begin{itemize}
\item $m^{-1}Cm\cap\gamma K_\gamma=\emptyset$;
\item il existe un sous-groupe parabolique propre ${\mth{P}}^+={\mth{M}}^+{\mth{U}}$ de ${\mth{G}}^+$ tel que l'on a $\overline{m^{-1}Cm\cap\gamma K_\gamma}\subset{\mth{P}}^+$ et $\overline{m^{-1}Hm\cap K_\gamma}\supset{\mth{U}}$.
\end{itemize}
\end{prop}

\begin{proof}
La démonstration de cette proposition fera l'objet du paragraphe suivant.
\end{proof}

\begin{cor}\label{cspd2}
Soit $\phi$ une fonction $\gamma$-cuspidale sur ${\mth{G}}$, relevée en une fonction sur $M_\gamma$ à support dans $\gamma K_\gamma$. Pour tout $f\in C_c^\infty\left(M_\gamma\right)$, l'intégrale:
\[\int_{Z_{M_\gamma}^0\backslash M_\gamma}\left(\int_{M_\gamma}f\left(m^{-1}ym\right)\phi\left(\gamma ^{-1}y\right)dy\right)dm\]
converge.
\end{cor}

\begin{proof}
Pour montrer ce corollaire, il suffit de montrer que pour tout $f\in C_c^{\infty}\left(M_\gamma\right)$, l'application:
\[\phi_f:m\longmapsto\int_{M_\gamma}f\left(m^{-1}ym\right)\phi\left(\gamma ^{-1}y\right)dy\]
est à support compact modulo $Z_{M_\gamma}^0$. Soit $C$ le support de $f$, et $H$ un sous-groupe ouvert compact de $M_\gamma$ tel que $f$ est biinvariante par $H$. Puisque $C$ est compact, on peut appliquer à $C$ et $H$ la proposition précédente, et puisque $\phi$ est $\gamma$-cuspidale, la partie $C'$ de $M_\gamma$ donnée par cette proposition contient le support de $\phi_f$, ce qui montre le corollaire.
\end{proof}

On peut donc, grâce à ce corollaire, définir une distribution $D_\phi^{M_\gamma}$ sur $M_\gamma$ par:
\[D_\phi^{M_\gamma}:f\in C_c^\infty\left(M_\gamma\right)\longmapsto\int_{Z_{M_\gamma}^0\backslash M_\gamma}\left(\int_{M_\gamma}f\left(m^{-1}ym\right)\phi\left(\gamma ^{-1}y\right)dy\right)dm.\]
Cette distribution est clairement invariante par conjugaison. On définira é\-ga\-le\-ment, si $P_\gamma=M_\gamma U_\gamma$ est un sous-groupe parabolique de $G$ dont $M_\gamma$ est un Levi:
\[D_\phi^G:f\in C_c^\infty\left(G\right)\longmapsto D_\phi^{M_\gamma}\left(f^{P_\gamma}\right):\]
cette distribution est également invariante. De plus, elle ne dépend pas du relèvement de $\phi$ à $M_\gamma$: en effet, soit $x'$ un autre élément de $\zeta^{-1}\left(\nu\right)$ tel que, en tant qu'élément de $\Gamma$, $\gamma\in\Gamma_x'$; soit $\mathcal{A}'$ un appartement de $\mathcal{B}$ contenant $x'$, $T'_0$ le tore déployé maximal de $G$ associé à $\mathcal{A}'$, $\left(M',K_{M'}\right)=\kappa\left(K_{x'},T'_0\right)$, $I'$ un sous-groupe d'Iwahori de $G$ contenu dans $K_{x'}$, $\gamma'$ l'unique élément de $\Gamma_{I'}\cap\gamma G^1$ et $M'_{\gamma'}$ le plus petit sous-groupe de Levi de $G$ semi-standard relativement à $T'_0$ et contenant $M'$ et $\gamma'$. D'après la proposition \ref{gcj}, quitte à remplacer $x'$ par $gx'$, $g\in G$, on peut supposer $T'_0=T_0$ et $\left(\gamma',M',K_{M'}\right)=\left(\gamma,M,K_M\right)$, d'où $M'_{\gamma'}=M_\gamma$; de plus, on a le lemme suivant:

\begin{lemme}\label{unicp}
Il existe un unique sous-groupe parahorique de $M_\gamma$ contenant $K_M$ et normalisé par $\gamma$.
\end{lemme}

\begin{proof}
Puisque l'assertion du lemme ne dépend que de $M_\gamma$ et pas de $G$, on aupposera $M_\gamma=G$. Soit $E$ l'ensemble des points de $\mathcal{A}$ fixés par $K_M$; c'est un sous-espace affine de $\mathcal{A}$ contenant $A_x$ et de même dimension que celle-ci, donc engendré par celle-ci. De plus, d'après le lemme \ref{uneorb}, $\gamma$ permute transitivement les murs de $A_x$; il ne stabilise donc aucune face de $A_x$ autre que $A_x$ elle-même. On en déduit que le sous-espace des éléments de $E$ fixés par $\gamma$ ne coupe aucune face de $A_x$ autre que $A_x$; il est donc entièrement contenu dans $A_x$, et on déduit du lemme du point fixe de Bruhat-Tits (\cite[I. 3.2.4]{bt}) que la seule facette de $E$ stabilisée par $\gamma$ est $A_x$, ce qui démontre l'assertion du lemme.
\end{proof}

D'après ce lemme, on a $K_x\cap M_\gamma=K_{x'}\cap M_\gamma$, et en notant $\phi_x$ (resp. $\phi_{x'}$) le relèvement de $\phi$ à $M_\gamma$ obtenu en considérant $x$ (resp. $x'$), il existe un élément $g$ de $N_{M_\gamma}\left(K_x\cap M_\gamma\right)$ tel que $\phi_{x'}=\phi_x\circ Ad\left(g\right)$; on en déduit que $D_{\phi_x}^{M_\gamma}=D_{\phi_{x'}}^{M_\gamma}$, donc que $D_{\phi_x}^G=D_{\phi_{x'}}^G$, ce qui démontre l'assertion cherchée.

Calculons maintenant certaines valeurs de ces distributions. Soit $\nu=\left(\mu,\gamma\right)\in\mathcal{N}$, $\left(M,K_\mu\right)$ un représentant quelconque de $\mu$, $T_0$ un tore déployé maximal de $M$, $M_\gamma$ le plus petit sous-groupe de Levi de $G$ semi-standard relativement à $T_0$ et contenant $M$ et $\gamma$, ${\mth{M}}=K_\mu/K_\mu^1$ et $\phi$ une fonction $\gamma$-cuspidale sur ${\mth{M}}$; supposons-les tels qu'il existe $x\in\mathcal{B}$ tel que $\left(M,K_\mu\right)=\kappa_0\left(K_x,T_0\right)$ et que $\gamma$ s'identifie à un élément de $W'$ normalisant $K_x$,
et considérons la distribution $D_\phi^G$: d'après ce qui précède, on peut supposer, en fixant un sous-groupe parabolique minimal $P_0$ de $G$ contenant $T_0$, que $M$ et $M_\gamma$ sont standards relativement à $P_0$ et que si $I$ est le sous-groupe d'Iwahori standard de $K_x$ relativement à $P_0$, $\gamma\in\Gamma_I$.

Définissons également $\nu,\mu',M',K_{\mu'},T'_0,\gamma',M'_{\gamma'},\phi',x',K_{x'}$ de manière similaire; quitte à les conjuguer, on peut supposer que $T'_0=T_0$ et que $K_{x'}$ admet également $I$ comme sous-groupe d'Iwahori standard relativement à $P_0$.

On a le lemme suivant:

\begin{lemme}\label{dscp}
On a:
\[D_\phi^G\left(\phi'\right)=0\]
si $\nu\neq\nu'$, et:
\[D_\phi^G\left(\phi'\right)=v_\nu\sum_{n\in N_G\left(\gamma,K_\mu\right)/N_{M_\gamma}\left(\gamma,K_\mu\right)}\langle \overline{\phi},Ad\left(n\right)\phi'\rangle _{\mth{M}},\]
où $v_\nu$ est une constante ne dépendant que de $\nu$, si $\nu=\nu'$.
\end{lemme}

\begin{proof}
Si $\gamma\neq\gamma'$, les supports de $D_\phi^G$ et de $\phi'$ sont disjoints, et on a $D_\phi^G\left(\phi'\right)=0$; supposons donc $\gamma=\gamma'$. Soit $K_0$ un sous-groupe parahorique maximal spécial de $G$; quitte à le conjuguer, on peut supposer qu'il contient $I$; l'ensemble des doubles classes de $K_0$ modulo $K_0\cap P_\gamma$ à gauche et $K_0\cap K_{x'}$ à droite admet alors pour système de représentants une partie de $W'\cap K_0$. Ecrivons, pour $m\in M$:
\[\phi'{}^{P_\gamma}\left(m\right)=\delta_{P_\gamma}^G\left(m\right)^{\dfrac 12}\int_{K_0}\int_{U_\gamma}\phi'\left(k^{-1}muk\right)dudk.\]
Il est clair que multiplier $k$ à gauche par un élément de $P_\gamma$ revient à conjuguer $m$ par un élément de $M_\gamma$; de plus, $\phi'$ étant invariante par conjugaison par $I$, $\phi'{}^{P_\gamma}$ est égale à un élément de $C_0\left(M_\gamma\right)$ près à:
\[f:m\longmapsto\delta_{P_\gamma}^G\left(m\right)^{\dfrac 12}\sum_{w\in\left(K_0\cap P_\gamma\right)\backslash K_0/I}vol\left(\left(K_0\cap P_\gamma\right)wI\right)\int_{U_\gamma}\phi'\left(w^{-1}muw\right)du.\]
Montrons d'abord les lemmes suivants:

\begin{lemme}\label{iwah}
Supposons $x'$ tel que $K_{x'}$ est un sous-groupe d'Iwahori de $G$; soit $w\in W'$. Alors $K_{x'}\cap w^{-1}M_\gamma w$ est un sous-groupe d'Iwahori de $w^{-1}M_\gamma w$.
\end{lemme}

\begin{proof}
En effet, puisque $K_{x'}$ est un sous-groupe d'Iwahori de $G$, $x'$ est contenu dans une chambre $A$ de $\mathcal{B}$; de plus, puisque $x'\in\mathcal{A}$, $A\subset\mathcal{A}$, donc puisque $w\in W'$, $wA\subset\mathcal{A}$. Soit $\mathcal{A}_{w^{-1}M_\gamma w}$ l'appartement de $w^{-1}M_\gamma w$ associé à $T_0$, et soit $A'$ la facette contenant l'image de $wA$ par la surjection canonique de $\mathcal{A}$ dans $\mathcal{A}_{w^{-1}M_\gamma w}$; c'est une chambre de $\mathcal{A}_{w^{-1}M_\gamma w}$, et tout élément de $K_x\cap w^{-1}Mw$ fixe $A'$, donc $K_x\cap w^{-1}M_\gamma w$ est contenu dans le sous-groupe d'Iwahori $K_{A'}$ de $w^{-1}M_\gamma w$ qui fixe $A'$.
Réciproquement, tout élément de $w^{-1}M_\gamma w$ qui fixe $A'$ fixe son image réciproque dans $\mathcal{A}$, et en particulier $wA$; on a donc $K_x\cap w^{-1}M_\gamma w=K_{A'}$ et le lemme est démontré.
\end{proof}

\begin{cor}\label{iw2}
Supposons $x'$ quelconque. Alors $K_{x'}\cap w^{-1}M_\gamma w$ est un sous-groupe parahorique de $w^{-1}M_\gamma w$.
\end{cor}

\begin{proof}
En effet, si $A$ est une chambre de $\mathcal{A}$ dans l'adhérence con\-tient $x'$, $K_{x'}$ contient le sous-groupe d'Iwahori de $G$ qui fixe $A$, dont l'intersection avec $w^{-1}M_\gamma w$ est un sous-groupe d'Iwahori de $w^{-1}M_\gamma w$ d'après le lemme précédent. De plus, $K_{x'}\cap w^{-1}M_\gamma w$ est un sous-groupe compact de $w^{-1}M_\gamma w$ contenu dans $\left(w^{-1}M_\gamma w\right)^1$; c'est donc un sous-groupe parahorique de $w^{-1}M_\gamma w$.
\end{proof}

\begin{lemme}\label{wmuw}
Soit $w\in W'\cap K_0$, $m\in M_\gamma$, $u\in U_\gamma$ tels que $w^{-1}muw\in\gamma K_{x'}$. Alors on a $w^{-1}mw\in\gamma K_{x'}$ et $w^{-1}uw\in K_{x'}$; de plus, $\gamma$ appartient à l'ensemble $\left(W'\cap w^{-1}M_\gamma w\right)\left(W'\cap K_{x'}\right)$.
\end{lemme}

\begin{proof}
En effet, en posant $h=mu$, $w^{-1}hw$ normalise le groupe $K_{x'}\cap w^{-1}P_\gamma w$, donc par projection sur $w^{-1}M_\gamma w$, on voit que $w^{-1}mw$ normalise le sous-groupe parahorique $K_{x'}\cap w^{-1}M_\gamma w$ de $w^{-1}M_\gamma w$; d'après le lemme \ref{normpar}, il existe alors $\gamma'\in W'\cap w^{-1}M_\gamma w$ tel que $w^{-1}mw\in\gamma'\left(K_{x'}\cap w^{-1}M_\gamma w\right)$. De plus, on peut supposer que $\gamma'$ normalise le sous-groupe d'Iwahori $I\cap w^{-1}M_\gamma w$ de $w^{-1}M_\gamma w$; comme $K_{x'}\cap w^{-1}M_\gamma w\subset\left(w^{-1}M_\gamma w\right)^1$, d'après le lemme \ref{normpar}, $\gamma'$ ne dépend alors pas de $h$. On en déduit:
\[\gamma K_{x'}\cap w^{-1}P_\gamma w\subset \gamma'\left(\left(K_{x'}\cap w^{-1}M_\gamma w\right)w^{-1}U_\gamma w\right),\]
d'où, puisque $W=N_G\left(T_0\right)/Z_G\left(T_0\right)$ et $Z_G\left(T_0\right)\subset K_{x'}\cap w^{-1}M_\gamma w$:
\[W'\cap\left(\gamma K_{x'}\cap w^{-1}P_\gamma w\right)\subset W'\cap\left(\gamma'\left(K_{x'}\cap w^{-1}M_\gamma w\right)\right).\]
Or $\gamma K_{x'}\cap w^{-1}P_\gamma w$ étant non vide, il contient au moins une double classe modulo le sous-groupe d'Iwahori $I\cap w^{-1}P_\gamma w$ de $w^{-1}P_\gamma w$, donc le membre de gauche de l'égalité ci-dessus n'est pas vide; on en déduit que $\gamma'\in W'\cap\gamma K_{x'}$, donc d'une part que $w^{-1}mw\in\gamma K_{x'}$ et $w^{-1}uw\in K_{x'}$, ce qui démontre la première assertion, et d'autre part que $\gamma\in \gamma'\left(W\cap K_{x'}\right)$, ce qui démontre également la deuxième assertion.
\end{proof}

Revenons à la démonstration du lemme \ref{dscp}. D'après le lemme \ref{wmuw}, on peut écrire, pour tout $m\in M_\gamma$:
\begin{multline*}f\left(m\right)=\delta_{P_\gamma}^G\left(m\right)^{\dfrac 12}\sum_{w\in\left(K_0\cap P_\gamma\right)\backslash K_0/I}vol\left(\left(K_0\cap P_\gamma\right)wI\right)\\\cdot\int_{w^{-1}U_\gamma w\cap K_{x'}}\phi'\left(w^{-1}mwu\right)du.\end{multline*}

Pour tout $w$, si l'image de $w^{-1}U_\gamma w\cap K_{x'}$ dans $K_{x'}/K_{x'}^1$ n'est pas réduite à $\left\{1\right\}$, le terme correspondant de la somme ci-dessus est nul puisque $\phi'$ est $\gamma$-cuspidale. Or on a $w^{-1}U_\gamma w\cap K_{x'}\subset K_{x'}^1$ si et seulement si $M'\subset w^{-1}M_\gamma w$. Supposons donc cette condition vérifiée: comme d'après le lemme \ref{wmuw}, $\phi_{\mu'}\left(w^{-1}mwu\right)$ est toujours nul si $\gamma\not\in\left(W'\cap w^{-1}M_\gamma w\right)\left(W'\cap K_{x'}\right)=W'\cap w^{-1}M_\gamma w$, on en déduit que les seuls termes éventuellement non nuls de la somme ci-dessus sont ceux correspondant aux $w$ tels que $M'_\gamma\subset w^{-1}M_\gamma w$. En particulier $f$ est nulle si $M'_\gamma$ n'est contenu dans aucun conjugué de $M_\gamma$.

Supposons donc que $M'_\gamma$ est contenu dans un conjugué de $M_\gamma$; quitte à conjuguer $M_\gamma$ et à remplacer $K_x$ par un sous-groupe parahorique associé, ce qui ne change pas $D_\phi^G$ d'après ce qu'on a déjà vu, on peut supposer que $M'_\gamma\subset M_\gamma$.
Soit, pour tout $m$, $W'_m$ l'ensemble des éléments $w$ de $\left(K_0\cap P_\gamma\right)\backslash K_0/I$ tels que $w^{-1}mw\in\gamma K_{x'}$; on a alors:
\[f\left(m\right)=\delta_{P_\gamma}^G\left(m\right)^{\dfrac 12}\sum_{w\in W'_m}vol\left(\left(K_0\cap P_\gamma\right)wI\right)vol\left(w^{-1}U_\gamma w\cap K_{x'}\right)\phi'\left(w^{-1}mw\right),\]
d'où:
\begin{equation*}\begin{split}D_\phi^G\left(\phi'\right)&=D_\phi^G\left(f\right)\\&=\int_{Z_{M_\gamma}^0\backslash M_\gamma}\int_{M_\gamma}\delta_{P_\gamma}^G\left(y\right)^{\dfrac 12}\phi\left(m^{-1}ym\right)\\&\quad\cdot\left(\sum_{w\in W'_y}vol\left(\left(K_0\cap P_\gamma\right)wI\right)vol\left(w^{-1}U_\gamma w\cap K_{x'}\right)\phi'\left(w^{-1}yw\right)\right)dydm\\&=\sum_{\begin{array}{c}w\in\left(K_0\cap P_\gamma\right)\backslash K_0/I\\M'_\gamma\subset w^{-1}M_\gamma w\end{array}}vol\left(\left(K_0\cap P_\gamma\right)wI\right)vol\left(w^{-1}U_\gamma w\cap K_{x'}\right)\\&\quad\cdot\int_{Z_{M_\gamma}^0\backslash M_\gamma}\int_{M_\gamma\cap m\gamma K_xm^{-1}\cap w\gamma K_{x'}w^{-1}}\phi\left(m^{-1}ym\right)\phi'\left(w^{-1}yw\right)dydm.\end{split}\end{equation*}
Considérons, pour tout $w$ et tout $m\in G$ tel que $M_\gamma\cap m\gamma K_xm^{-1}\cap w\gamma K_{x'}w^{-1}$ n'est pas vide, l'intégrale:
\[\int_{M_\gamma\cap m\gamma K_xm^{-1}\cap w\gamma K_{x'}w^{-1}}\phi\left(m^{-1}ym\right)\phi'\left(w^{-1}yw\right)dy.\]
On supposera que $w$ est un élément de $W'\cap K_0$ de longueur minimale dans sa classe à gauche modulo $W'\cap\left(K_0\cap P_\gamma\right)$; on a alors $w\left(I\cap M_\gamma\right)w^{-1}\subset I$.
Quitte à multiplier $m$ par un élément de $I\cap M_\gamma$ à gauche et par un autre à droite, ce qui ne change pas la valeur de l'intégrale puisque $\phi$ et $\phi'$ sont stables par conjugaison par $I$, on peut supposer $m\in N_G\left(T_0\right)$.
D'autre part, quitte à remplacer $m$ par $mw^{-1}$, on peut supposer $w=1$.

Considérons le sous-groupe compact $m^{-1}K_{x'}^1m\cap w^{-1}K_x w$ de $G$; pour que l'in\-té\-gra\-le soit non nulle, il est nécessaire que l'image de ce sous-groupe dans le groupe $\left(w^{-1}K_x w\right)/\left(w^{-1}K_x^1 w\right)$, identifié à ${\mth{M}}$, ne contienne aucun sous-groupe qui soit le radical unipotent d'un sous-groupe parabolique de $\left(w^{-1}\gamma^{\mth{Z}}K_x w\right)/\left(w^{-1}K_x^1w\right)$, identifié à ${\mth{M}}^+$, rencontrant $\gamma{\mth{M}}$, ou, de manière équivalente, d'un sous-groupe parabolique de ${\mth{M}}$ stable par conjugaison par au moins un élément de $\gamma{\mth{M}}$. On va donc montrer que cette condition implique $\nu=\nu'$.

On va d'abord vérifier que pour tout $m\in N_G\left(T_0\right)$, si $m\gamma K_xm^{-1}\cap\gamma K_{x'}$ est non vide, elle contient un élément de $N_G\left(T_0\right)$. Soit donc $k\in m\gamma K_xm^{-1}\cap\gamma K_{x'}$; $k$ stabilise à la fois $mA_x$ et $A_{x'}$, et ces deux facettes sont alors contenues dans l'appartement $k\mathcal{A}$ de $\mathcal{B}$. Or d'après \cite[I, proposition 2.5.8]{bt}, il existe $g\in G^1$ tel que $gk\mathcal{A}=\mathcal{A}$ et que $g$ fixe $\mathcal{A}\cap k\mathcal{A}$, donc en particulier $mA_x$ et $A_{x'}$, point par point. Autrement dit, $gk\in N_G\left(T_0\right)$ et $g\in mK_xm^{-1}\cap K_{x'}$, donc $gk\in m\gamma K_xm^{-1}\cap\gamma K_{x'}$.

Puisque $gk\in N_G\left(T_0\right)$, $W'\cap\left(m\gamma K_xm^{-1}\cap\gamma K_{x'}\right)$ est non vide; il existe alors $w_x\in W'\cap K_x$, $w_{x'}\in W'\cap K_{x'}$ tels que l'on a:
\[m\gamma w_xm^{-1}=\gamma w_{x'}.\]
On en déduit que si le groupe $\left(\overline{m^{-1}K_{x'}^1m\cap K_x}\right)_{\mth{M}}$ n'est pas réduit à l'identité, il contient le radical unipotent d'un sous-groupe parabolique de ${\mth{M}}$ stable par conjugaison par $\gamma w_x=m^{-1}\gamma w_{x'}m$, et l'intégrale est alors nulle. Dans le cas contraire, on a $m^{-1}K_{x'}^1m\cap K_x\subset K_x^1$, ce qui n'est possible que si on a $m^{-1}Mm=M'$ et $m^{-1}K_\mu m=K_{\mu'}$, c'est-à-dire $\mu=\mu'$. 
Donc la première assertion du lemme est démontrée.

Supposons maintenant $M_\gamma=M'_\gamma$, $\nu=\nu'$ et $K_\mu=K_{\mu'}$; on a alors $K_x\cap M_\gamma=K_{x'}\cap M_\gamma$, d'où $K_x=K_{x'}$ puisque $K_x$ et $K_{x'}$ contiennent tous deux $I$, et pour tous $m\in N_G\left(T_0\right)$ et $w\in N_G\left(T_0\right)\cap K_0$:
\begin{multline*}\int_{M_\gamma\cap m\gamma K_xm^{-1}\cap w\gamma K_xw^{-1}}\phi_\mu\left(m^{-1}ym\right)\phi_{\mu'}\left(w^{-1}yw\right)dy\\=vol\left(K_x\right)\langle \overline{\phi'},Ad\left(wm^{-1}\right)\phi\rangle _{\mth{M}}\end{multline*}
si $wm^{-1}\in N_G\left(\gamma,K_\mu\right)$, et $0$ sinon;
de plus, si $w\in W'\cap K_0$ normalise $M_\gamma$ et $I\cap M_\gamma$, d'après le lemme \ref{unicp}, il normalise aussi $\left(\gamma,K_\mu\right)$ et on a alors $m\in N_{M_\gamma}\left(\gamma,K_\mu\right)$.
D'autre part, le groupe $N_G\left(\gamma,K_\mu\right)/N_{M_\gamma}\left(\gamma,K_\mu\right)$ admet un système de représentants dans $N_G\left(T_0\right)$, et puisque $K_0$ est spécial, on a:
\[N_G\left(T_0\right)=\left(N_G\left(T_0\right)\cap K_0\right)\left(N_G\left(T_0\right)\cap M_\gamma\right);\]
enfin, pour tout $w\in W'\cap K_0$ tel que $M_\gamma=w^{-1}M_\gamma w$, on a $w^{-1}U_\gamma w\cap K_x=w^{-1}U_\gamma w\cap I$ et:
\begin{equation*}\begin{split}vol\left(\left(K_0\cap P_\gamma\right)wI\right)&vol\left(w^{-1}U_\gamma w\cap K_x\right)\\&=\dfrac{vol\left(K_0\cap P_\gamma\right)vol\left(w^{-1}U_\gamma w\cap I\right)vol\left(I\right)}{vol\left(w^{-1}P_\gamma w\cap I\right)}\\&=\dfrac{vol\left(K_0\cap P_\gamma\right)vol\left(I\right)}{vol\left(M_\gamma\cap I\right)}.\end{split}\end{equation*}
Donc en posant:
\[v_\nu=\dfrac{vol\left(K_0\cap P_\gamma\right)vol\left(K_x\right)vol\left(Z_{M_\gamma}^0\backslash N_{M_\gamma}\left(\gamma,K_\mu\right)\right)vol\left(I\right)}{vol\left(M_\gamma\cap I\right)},\]
l'assertion du lemme est vérifiée.
\end{proof}

\subsection{Démonstration de la proposition}

On va maintenant démontrer la proposition \ref{cspdp2}. L'assertion ne dépendant que de $M_\gamma$ et pas de $G$, on supposera $G=M_\gamma$, d'où $K_x=K_\gamma$. De plus, quitte à remplacer $G$ par $G/Z_{G,d}$, on pourra supposer que le centre de $G$ est compact.
D'autre part, il est clair que la proposition est vraie pour $H$ et $C$ si elle est vraie pour un groupe plus petit que $H$ et une partie compacte de $G$ plus grande que $C$; on pourra donc faire les hypothèses suivantes:
\begin{itemize}
\item il existe un sous-groupe parahorique maximal spécial $K_0$ de $G$ contenant $K_{Z_G\left(T_0\right)}$ et tel que $C$ est réunion de doubles classes modulo $K_0$;
\item $H$ est compact, et stable par conjugaison par $K_0$.
\end{itemize}
En effet, si $K_0$ est un bon sous-groupe parahorique maximal de $G$, l'ensemble des doubles classes de $G$ modulo $K$ rencontrant $C$ est fini puisque $C$ est compact, donc leur réunion est également une partie compacte de $C$;
d'autre part, si $H'$ est un sous-groupe ouvert compact de $G$, il existe un entier $s$ tel que $H'$ contient le $s$-ème sous-groupe de congruence $K_0^s$ de $K_0$, qui est stable par conjugaison par $K_0$.

Si $m$ est un élément de $G$ vérifiant l'une des deux conditions de l'énoncé, alors pour tout $m'\in K_0mK_x$, $m'$ vérifie la même condition; c'est clair pour la première, et pour la seconde, en posant $m'=kmk'$, $k\in K_0$, $k'\in K_x$, si ${\mth{P}}^+$ est un sous-groupe parabolique propre de ${\mth{G}}^+$ vérifiant la condition voulue pour $m$, $\overline{k'}^{-1}{\mth{P}}^+\overline{k'}$ convient pour $m'$; pour montrer la proposition, il suffit donc de considérer un ensemble de représentants de $K_0\backslash G/K_x$.

Soit $W$ le groupe de Weyl de $G$ relativement à $T_0$.
On a $W=N_G\left(T_0\right)/Z_G\left(T_0\right)$ et $W'=N_G\left(T_0\right)/K_{Z_G\left(T_0\right)}$, où $K_{Z_G\left(T_0\right)}$ est l'unique sous-groupe parahorique de $Z_G\left(T_0\right)$; il existe donc une surjection canonique de $W'$ dans $W$, dont le noyau est $Z_G\left(T_0\right)/K_{Z_G\left(T_0\right)}$.
De plus, soit $I$ un sous-groupe d'Iwahori de $G$ contenant $K_{Z_G\left(T_0\right)}$ et contenu dans $K_x$; considérons la décomposition de Bruhat:
\[G=\bigsqcup_{w\in W'}IwI.\]
Soit $I'$ un sous-groupe d'Iwhori de $G$ contenant également $K_{Z_G\left(T_0\right)}$ et contenu dans $K_0$; il existe $t_0\in T_0$ tel que $I'=t_0It_0^{-1}$, et en multipliant tous les termes de la décomposition ci-dessus par $t_0$ à gauche, on obtient:
\[G=\bigsqcup_{w\in W'}I'wI.\]
D'autre part, puisque $K_0$ est spécial, pour tout $w\in W$, il existe un unique $w'\in W'\cap K_0$ dont l'image dans $W$ par la surjection canonique est $w$; on en déduit que l'on a:
\[G=\bigsqcup_{t\in X}K_0tI=\bigcup_{t\in X}K_0tK_x,\]
où $X$ est un système de représentants de $Z_G\left(T_0\right)/K_{Z_G\left(T_0\right)}$.
Ce dernier ensemble étant un groupe abélien libre de type fini, on pourra supposer pour alléger les notations que $X$ est lui-même un groupe. De plus, soit $X_0=X\cap T_0K_{Z_G\left(T_0\right)}$; $X_0$ est canoniquement isomorphe à $T_0/K_{T_0}$, donc à $X_*\left(T_0\right)$; de plus, il est d'indice fini dans $X$; $X$ est donc canoniquement isomorphe à un sous-groupe de type fini de $X_*\left(T_0\right)\otimes{\mth{Q}}$.

D'après ce qui précède, il suffit de montrer que tous les éléments de $X$ sauf un nombre fini d'entre eux vérifient au moins une des deux conditions de la proposition. Pour cela, montrons d'abord le lemme suivant:

\begin{lemme}\label{mder}
Il existe un compact $T_1\subset Z_G\left(T_0\right)$ tel que pour tout $t\in Z_G\left(T_0\right)-T_1\left(T_0\cap M_{der}\right)$, on a $t^{-1}Ct\cap\gamma K_x=\emptyset$.
\end{lemme}

\begin{proof}
Pour montrer ce lemme, on aura besoin du résultat suivant:

\begin{lemme}\label{ilexpas}
Soit $t\in X$. Si $\gamma^{-1}t\gamma\not\in w^{-1}twK_{Z_G\left(T_0\right)}$ pour tout $w\in W'\cap K_x$, alors $t^{-1}\gamma K_x t$ et $\gamma K_x$ sont disjoints.
\end{lemme}

\begin{proof}
En effet, supposons qu'il existe $g\in\gamma^{-1}t^{-1}\gamma K_x t\cap K_x$; alors si $A_x$ est la facette de $\mathcal{A}$ contenant $x$, $g$ fixe $A_x$ et envoie $t^{-1}A_x$ sur $\gamma^{-1}t^{-1}A_x$. D'après \cite[2.5.8]{bt}, il existe $h\in G^1$ tel que $hg\mathcal{A}=\mathcal{A}$ et que $h$ fixe $\mathcal{A}\cap g\mathcal{A}$, donc en particulier $A_x$ et $\gamma^{-1}t^{-1}A_x$, ou encore $h\in K_x\cap\gamma^{-1}t^{-1}K_xt\gamma$; on en déduit:
\[hg\in Z_G\left(T_0\right)\cap K_x\cap\gamma^{-1}t^{-1}\gamma K_x t.\]
Il existe donc $w,w'\in W'\cap K_x$ tels que l'on a:
\[w=\gamma^{-1}t^{-1}\gamma w't;\]
de plus, $w$ et $w'$ ont même image par la surjection canonique de $W'$ dans $W$, ce qui n'est possible que s'ils sont égaux. L'assertion du lemme s'en déduit immédiatement. 
\end{proof}

Revenons à la démonstration du lemme \ref{mder}. Considérons le sous-groupe $T_2$ de $Z_G\left(T_0\right)$ engendré par les $t$ tels qu'il existe $w\in W'\cap K_x$ tel que $\gamma^{-1}t\gamma\in w^{-1}twK_{Z_G\left(T_0\right)}$; on va montrer d'une part qu'il existe un compact $T'_1$ de $Z_G\left(T_0\right)$ tel que $T_2\subset T'_1\left(T_0\cap M_{der}\right)$, et d'autre part qu'il existe une partie finie $T''_1$ de $Z_G\left(T_0\right)$ telle que pour tout $t\in Z_G\left(T_0\right)-T''_1T_2$, $t^{-1}Ct\cap\gamma K_x=\emptyset$; $T''_1T'_1$ est alors un compact de $Z_G\left(T_0\right)$ vérifiant la condition du lemme.

Montrons d'abord la première assertion.
Soit $s$ le rang semi-simple de $G$, $s'$ le nombre de composantes connexes de son diagramme de Dynkin, et soit $\beta_1,\dots,\beta_{s+s'}$ l'ensemble des racines de $G$ relativement à $T_0$ correspondant aux murs de $A$; on supposera les indices choisis de telle sorte que:
\begin{itemize}
\item $\left\{\beta_1,\dots,\beta_s\right\}$ est un ensemble de racines simples de $G$ relativement à $T_0$;
\item si $r$ est le rang semi-simple de $M$, $\left\{\beta_1,\dots,\beta_r\right\}$ est un ensemble de racines simples de $M$.
\end{itemize}
Alors $w_\gamma$ permute les $\beta_j$, $1\leq j\leq s+s'$; on notera pour tout $j$, $\beta_{\gamma\left(j\right)}=w_\gamma\beta_j$.

Soit $\langle \cdot,\cdot\rangle $ le produit de dualité canonique entre $X^*\left(T_0\right)$ et $X_*\left(T_0\right)$, c'est-à-dire le produit tel que pour tous $\chi\in X^*\left(T_0\right)$ et $\xi\in X_*\left(T_0\right)$, on a:
\[\chi\left(\xi\left(\varpi\right)\right)=\varpi^{\langle \chi,\xi\rangle }.\]
Ce produit s'étend de façon canonique en un produit de dualité entre $X^*\left(T_0\right)\otimes{\mth{Q}}$ et $X_*\left(T_0\right)\otimes{\mth{Q}}$ à valeurs dans ${\mth{Q}}$, qui se relève en un produit de dualité entre $X^*\left(T_0\right)\otimes{\mth{Q}}$ et $Z_G\left(T_0\right)$.

On va d'abord montrer que si $t\in T_2$, les valeurs de $\langle \beta_{r+1},t\rangle ,\dots,\langle \beta_{s+s'},t\rangle $ sont déterminées par celles de $\langle \beta_1,t\rangle ,\dots,\langle \beta_r,t\rangle $. L'application qui à $t\in Z_G\left(T_0\right)$ associe $\left(\langle \beta_1,t\rangle ,\dots,\langle \beta_{s+s'},t\rangle \right)$ étant un morphisme de groupes de $Z_G\left(T_0\right)$ dans ${\mth{Q}}^{s+s'}$, il suffit de le montrer pour une partie génératrice de $T_2$; on peut donc supposer qu'il existe $w\in W'\cap K_x$ tel que $\gamma^{-1}t\gamma\in w^{-1}twK_{Z_G\left(T_0\right)}$.

Supposons d'abord $w=1$.
Puisque par hypothèse, $\gamma^{-1}t\gamma\in tK_{Z_G\left(T_0\right)}$, on a pour tout $i$:
\[\langle \beta_{\gamma\left(j\right)},t\rangle =\langle \beta_j,t\rangle .\]
De plus, si $J_1,\dots,J_{s'}$ sont les parties de $\left\{1,\dots,s+s'\right\}$ correspondant aux composantes connexes du diagramme de Dynkin, on a pour tout $i$ une relation de la forme:
\[\sum_{j\in J_i}c_i\langle \beta_j,t\rangle =0,\]
les $c_i$ étant des entiers strictement positifs. Enfin, soit, pour tout $i$, $\Delta_i$ l'ensemble des $\beta_j$, $j\in J_i$; on a le lemme suivant:

\begin{lemme}\label{uneorb}
Pour tout $i$, les éléments de $\Delta_i$ qui ne sont pas racines de $M$, s'il en existe, sont tous situés dans une même orbite $O_i$ pour l'action de $w_\gamma$ sur $\Phi$.
\end{lemme}

\begin{proof}
Il est clair que l'on peut supposer $s'=1$.
De plus, puisque $G=M_\gamma$, en définissant $V$, $V^\gamma$ et $V_M$ comme dans le paragraphe précédent, on a $V^\gamma\cap V^M=0$.

Considérons l'ensemble $\Delta_1=\left\{\beta_1,\dots,\beta_{s+1}\right\}$. On sait que $w_\gamma$ agit sur $\Delta_1$; puisque $\gamma$ normalise $M$, $w_\gamma$ aussi, donc il permute les $\beta_1,\dots,\beta_r$. Supposons l'assertion du lemme fausse; il existe alors au moins une orbite $O$ de $\Phi$ pour l'action de $w_\gamma$ incluse dans $\left\{\beta_{r+1},\dots,\beta_s\right\}$. Si de plus cette orbite est $\left\{\beta_{r+1},\dots,\beta_s\right\}$ tout entier, l'ensemble $\left\{\beta_1,\dots,\beta_s\right\}$ est stable par $w_\gamma$, ce qui n'est possible que si $w_\gamma=1$, c'est-à-dire si $\gamma\in T_0\subset M$; mais alors $G=M$ et l'assertion du lemme est trivialement vraie, d'où contradiction. On en déduit que le cardinal de $O$ est strictement inférieur à $s-r$.

D'autre part, soit ${\mathfrak{a}}=Hom\left(X^*\left(T_0\right),{\mth{R}}\right)$ l'algèbre de Lie réelle de $T_0$; ${\mathfrak{a}}$ est canoniquement isomorphe à $X_*\left(T_0\right)\otimes{\mth{R}}$, donc à $V$ puisqu'on a supposé le centre de $G$ compact, par l'application $\phi$ définie, pour tous $\lambda_i\in{\mth{R}}$ et tous $\chi_i\in X_*\left(T_0\right)$, par:
\[\phi\left(\sum_i\lambda_i\xi_i\right)\left(\chi\right)=\sum_{i\in I}\lambda_i\langle \chi,\xi_i\rangle .\]
Soit ${\mathfrak{a}}^M$ le sous-espace de ${\mathfrak{a}}$ engendré par les coracines de $M$ relativement à $T_0$;
puisque $\gamma$ normalise $M$, son image $w_\gamma$ par l'application canonique de $W'$ dans $W$ laisse stable ${\mathfrak{a}}^M$, donc également son orthogonal ${\mathfrak{a}}_M$ par un produit scalaire $W$-invariant sur ${\mathfrak{a}}$; si $p_M$ est la projection orthogonale de ${\mathfrak{a}}$ sur ${\mathfrak{a}}_M$, $p_M$ commute avec l'action de $w_\gamma$. Soit, maintenant, pour tout $\beta_j\in O$, $e_j$ le vecteur unité de $V$ normal au mur de $A$ correspondant à $\beta_j$ et dirigé vers $A$; considérons le vecteur:
\[e_O=\sum_{j,\beta_j\in O}e_j.\]
Puisque $w_\gamma$ permute les $\beta_j\in O$, il permute également les $e_j$, donc fixe $e_O$. D'autre part, on a, pour tout $\xi\in V$, en normalisant convenablement le produit scalaire $W$-invariant sur ${\mathfrak{a}}$:
\[\left(e_j,\xi\right)=\langle \beta_j,\xi\rangle \]
pour tout $j$, d'où:
\[\left(e_O,\xi\right)=\sum_{\beta_j\in O}\langle \beta_j,\xi\rangle .\]
Puisque $p_M\left(e_O\right)-e_O$ est combinaison linéaire des coracines de $M$, on a une égalité de la forme:
\[\left(p_M\left(e_O\right),\xi\right)=\sum_{\beta_j\in O}\langle \beta_j,\xi\rangle -\sum_{j=1}^r\lambda_j\langle \beta_j,\xi\rangle \]
pour tout $\xi\in V$, donc $p_M\left(e_O\right)\neq 0$. Puisque $p_M\left(e_O\right)$ est fixé par $w_\gamma$, le sous-espace des éléments de ${\mathfrak{a}}_M$ fixés par $w_\gamma$ n'est pas réduit à $0$; or on a ${\mathfrak{a}}_M=V_M$, d'où contradiction et le lemme est démontré.
\end{proof}

On déduit de ce lemme et de ce qui précède que pour tout $i$ et tout $j\in I_i$ tel que $j>r$, on a:
\[\left(\sum_{k\in I_i;k>r}c_k\right)\langle \beta_j,t\rangle =-\sum_{k\in I_i;k\leq r}c_k\langle \beta_k,t\rangle ;\]
comme la somme des $c_k$ du membre de gauche est non nulle, $\langle \beta_j,t\rangle $ ne dépend que des $\langle \beta_k,t\rangle $, avec $k\in I_i$ et $k\leq r$.

Supposons maintenant $w$ quelconque; on a le lemme suivant:

\begin{lemme}\label{wprim}
Soit $t\in X$ tel que $\gamma^{-1}t\gamma\in w^{-1}twK_{Z_G\left(T_0\right)}$; il existe $w'\in W'\cap K_x$ tel que $\gamma^{-1}w'{}^{-1}tw'\gamma\in w'{}^{-1}tw'K_{Z_G\left(T_0\right)}$.
\end{lemme}

\begin{proof}
Supposons d'abord $s'=1$. Si $K_x$ est un sous-groupe d'Iwahori de $G$, alors $w=1$, et $w'=1$ convient. Supposons donc que $K_x$ n'est pas un Iwahori: $A_x$ n'est alors pas de dimension maximale, et est contenue dans exactement $r$ murs de $\mathcal{A}$. Soit donc $H_1,\dots,H_r$ les murs de $\mathcal{A}$ la contenant, et ${\mathfrak{D}}_1,\dots,{\mathfrak{D}}_k$ les adhérences des composantes connexes du complémentaire dans $\mathcal{A}$ de la réunion des $H_i$: $\gamma$ permute les ${\mathfrak{D}}_j$, et si ${\mathfrak{D}}_1$ est la composante contenant les $t'A_x$, avec $t'$ tel que $\langle \beta_k,t'\rangle >0$ pour tout $k\leq r$, $\gamma$ laisse ${\mathfrak{D}}_1$ globalement stable. Or on a, pour tout $w''\in W'\cap K_x$:
\[\gamma^{-1}w''{}^{-1}tw''\gamma\in\left(w\gamma^{-1}w''\gamma\right)^{-1}t\left(w\gamma^{-1}w''\gamma\right)K_{Z_G\left(T_0\right)},\]
donc $\gamma$ permute les $w''{}^{-1}tw''K_{Z_G\left(T_0\right)}$, $w''\in W'\cap K_x$; on en déduit que $\gamma$ permute les facettes de la forme $w''{}^{-1}tw''A_x$. De plus, chaque ${\mathfrak{D}}_i$ contient exactement une des facettes $w''{}^{-1}tw''A_x$; donc si $w'$ est un élément de $W'\cap K_x$ tel que $w'{}^{-1}tw'A_x\in{\mathfrak{D}}_1$, on a $\gamma^{-1}w'{}^{-1}tw'\gamma A_x=w'{}^{-1}tw'A_x$, d'où l'on déduit que $\gamma^{-1}w'{}^{-1}tw'\gamma\in w'{}^{-1}tw'K_{Z_G\left(T_0\right)}$ (comme on a supposé la composante déployée du centre de $G$ triviale, le stabilisateur de $w'{}^{-1}tw'A_x$ ne contient aucun élément de $X$ autre que $1$).

Supposons maintenant $s'$ quelconque. Définissons, pour tout $i\in\left\{1,\dots,s'\right\}$, des groupes $T_i$, $H_i$ et $K_{x,i}$ de la même façon que dans la démonstration du lemme \ref{gcj}; on a une injection canonique de $G$ dans $H_1\times\dots\times H_{s'}$, et chacun des $H_i$ possède un diagramme de Dynkin connexe. De plus, si pour tout $i$, $W'_i$ est le groupe de Weyl affine de $H_i$ relativement à $T_0/\left(T_i\cap T_0\right)$, $W'\cap K_x$ est canoniquement isomorphe au produit direct des $W'_i\cap K_{x,i}$. Ecrivons donc $\gamma=\gamma_1\dots\gamma_{s'}$ et $w=w_1\dots w_{s'}$, où pour tout $i$, $\gamma_i$ et $w_i$ sont des éléments de $W'_i$. Pour tout $i$, $w_i\in W'_i\cap K_{x,i}$ et $\gamma_i$ normalise $K_{x,i}$, qui est un sous-groupe parahorique de $H_i$; on peut donc appliquer l'assertion du lemme à l'image de $t$ dans $H_i$, $\gamma_i$ et $w_i$ et obtenir $w'_i\in W'_i\cap K_{x,i}$ vérifiant la condition voulue; l'élément $w'=w'_1\dots w'_{s'}$ de $W'\cap K_x$ convient alors pour $t,\gamma,w$.
\end{proof}

D'après ce lemme, on peut appliquer le cas $w=1$ à $w'{}^{-1}tw'$; on en déduit que les valeurs de $\langle \beta_{r+1},w'{}^{-1}tw'\rangle ,\dots,\langle \beta_{s+s'},w'{}^{-1}tw'\rangle $, et donc toutes les valeurs de $\langle \alpha,w'{}^{-1}tw'\rangle $, $\alpha\in\Phi$, sont déterminées par celles des $\langle \beta_i,w'{}^{-1}tw'\rangle$, $i\leq r$. Or on a pour tout $i$:
\[\langle \beta_i,w'{}^{-1}tw'\rangle =\langle w'\left(\beta_i\right),t\rangle ;\]
on en déduit que les valeurs des $\langle \alpha,t\rangle $ sont déterminées par celles des $\langle w'\left(\beta_1\right),t\rangle$, $i\leq r$. Or $w'$ appartient à $W'\cap K_x$, donc son image dans $W$ est un élément du groupe de Weyl de $M$ relativement à $T_0$; les racines $w'\left(\beta_1\right),\dots,w'\left(\beta_r\right)$ constituent donc un ensemble de racines simples de $M$ et l'assertion cherchée s'en déduit.

Soit $X_2=X\cap T_2$. D'apràs ce qui précède, si $R$ est un système de représentants des orbites de $\gamma$ dans $\left\{\beta_1,\dots,\beta_r\right\}$, l'application:
\[t\in X_2\longmapsto\left(\langle \beta_j,t\rangle \right)_{\beta_j\in R}\]
est un isomorphisme entre $X_2$ et un réseau de ${\mth{Q}}^{\card\left(R\right)}$; $X_2$ est donc un groupe abélien libre de rang $\card\left(R\right)$.

Considérons maintenant le groupe $T_0\cap M_{der}$. Le groupe $T_0\cap M_{der}/K_{T_0\cap M_{der}}$ est en bijection canonique avec le sous-groupe de $X_*\left(T_0\right)$ engendré par les coracines de $M_{der}$ relativement à $T_0$; si $X_d$ est un sous-groupe de $T_0\cap M_{der}$ formé d'un système de représentants des classes modulo $K_{T_0\cap M_{der}}$, d'une part $X_d$ est en bijection canonique avec un sous-groupe de $X_*\left(T_0\right)$, et d'autre part un élément $t'$ de $X_d$ est entièrement déterminé par la valeur des $\langle \beta_k,t'\rangle $, $k\leq r$. Si maintenant $X_{d,2}=X_d\cap T_2$, d'une part on a une injection canonique de $X_{d,2}$ dans $X_2$, et d'autre part, par le même argument que pour $X_2$, $X_{d,2}$ est un groupe abélien libre de rang $\card\left(R\right)$; il est donc isomorphe à un sous-groupe d'indice fini de $X_2$. Soit $X_1$ un système de représentants de $X_2$ modulo ce sous-groupe; posons:
\[T'_1=X_1K_{Z_G\left(T_0\right)};\]
d'après ce qui précède, $T'_1\left(T_0\cap M_{der}\right)$ contient $T_2$.

Montrons maintenant l'existence de $T''_1$; pour cela, on va utiliser le lemme suivant:

\begin{lemme}\label{cjp}
La réunion $C_{x,\gamma,0}$ des conjugués de $\gamma K_x$ par un élément de $Z_G\left(T_0\right)$ est un ouvert fermé de $G$.
\end{lemme}

\begin{proof}
Puisque $C_{x,\gamma,0}$ est réunion de parties ouvertes de $G$, c'est un ouvert de $G$. Il reste à montrer que c'est aussi un fermé de $G$.

Soit $\left(g_n\right)_{n\geq 1}$ une suite d'éléments de $C_{x,\gamma,0}$ convergeant vers un élément $g$ de $G$; on va montrer que $g\in C_{x,\gamma,0}$. Puisque d'après le lemme \ref{parferm}, l'ensemble des éléments compacts de $\gamma G^1$ est fermé, $g$ est compact; soit donc $K'$ un sous-groupe parahorique de $G$ et $\gamma'\in N_G\left(K'\right)$ tels que $g\in\gamma'K'$, et $A'$ la facette de $\mathcal{B}$ stabilisée par $K'$. Puisque $\gamma'K'$ est ouvert, tous les $g_n$, pour $n$ assez grand, sont dans $\gamma'K'$, et stabilisent alors à la fois $A'$ et une facette de la forme $t_nA_x$, $t_n\in X$. On va montrer qu'il existe une partie finie $X'$ de $X$ telle que pour tout $n$, $g_n$ stabilise une facette de la forme $t'A_x$, $t'\in X'$;
il existe alors au moins un $t'\in X'$ tel que $t'\gamma K_\gamma t'{}^{-1}$ contient une infinité de $g_n$, et on en déduit que ce groupe contient également $g$.

Soit $\Omega$, $\Omega'$ deux parties de $\mathcal{B}$; on définit la {\em distance combinatoire} $d_c\left(\Omega,\Omega'\right)$ de $\Omega$ à $\Omega'$ comme la longueur de la plus petite galerie reliant $\Omega$ à $\Omega'$, c'est-à-dire le plus petit entier $l$ tel qu'il existe des chambres $C_0,\dots,C_l$ de $\mathcal{B}$ telles que:
\begin{itemize}
\item l'adhérence de $C_0$ rencontre $\Omega$;
\item l'adhérence de $C_l$ rencontre $\Omega'$;
\item pour tout $i$, $C_i$ et $C_{i+1}$ sont adjacentes.
\end{itemize}
Si $x,y$ sont deux points de $\mathcal{B}$, on définit de manière analogue $d_c\left(x,y\right)$ et $d_c\left(x,\Omega\right)$.

Soit $d_k$ la distance combinatoire maximale entre deux chambres possédant un sommet commun; on a le lemme suivant:

\begin{lemme}\label{intr}
Pour toutes facettes $B,B',B''$ de $\mathcal{B}$, on a:
\[d_c\left(B,B''\right)\leq d_c\left(B,B'\right)+d_c\left(B',B''\right)+d_k.\]
\end{lemme}

\begin{proof}
En effet, soit $\left(C_0,\dots,C_l\right)$ et $\left(C''_0,\dots,C''_{l''}\right)$ deux galeries tendues respectivement entre $B$ et $B'$ et entre $B'$ et $B''$; les chambres $C_l$ et $C''_0$ ont alors une face commune $B'$, donc au moins un sommet commun, et il existe alors une galerie tendue de la forme $\left(C_l=C'_0,C'_1,\dots,C'_{l'}=C''_0\right)$, avec $l'\leq d_k$; on en déduit l'asserton du lemme.
\end{proof}

Soit $Q$ un quartier fermé de $\mathcal{A}$ de sommet $x$ et rencontrant $A'$, et $\Omega$ l'ensemble des points de $Q$ n'appartenant à aucun $tQ$, avec $t\in X-\left\{0\right\}$ tel que $tx\in Q$; $Q$ rencontre un nombre fini de facettes, et l'entier $d_m=\sup_{B,B'}d_c\left(B,B'\right)$, où $B,B'$ décrivent l'ensemble des facettes rencontrant $\Omega$, est bien défini; de plus, l'ensemble des $tQ$, $t\in X$, constitue une partition de $\mathcal{A}$.

Posons également $d_0=d_c\left(A',\mathcal{A}\right)$; puisque pour tout entier $a$, l'ensemble des facettes de $\mathcal{A}$ situées à une distance combinatoire d'au plus $a$ de $A'$ est fini, l'assertion cherchée se déduit immédiatement du lemme suivant:

\begin{lemme}
Pour tout $t\in X$, il existe $t'\in X$ tel que $d_c\left(A',t'A_x\right)\leq d_0+d_m+2d_k$ et que $\gamma'K'\cap t\gamma K_xt^{-1}\subset t'\gamma K_xt'{}^{-1}$.
\end{lemme}

\begin{proof}
Si $\gamma'K'$ et $t\gamma K_xt^{-1}$ sont disjoints, l'assertion du lemme est triviale; supposons donc $\gamma'K'\cap t\gamma K_xt^{-1}$ non vide.
Soit $E$ l'enclos de $A'\cup tA_x$ dans $\mathcal{B}$. Supposons d'abord $\gamma=1$; on va montrer par récurrence sur $d_0$ qu'il existe $t'\in X$ tel que $t'A_x\subset E$ et que $d_c\left(A',t'A_x\right)\leq d_0+d_m$. Si $d_0=0$, $A'\subset\mathcal{A}$ et il existe un unique $t'\in X$ tel que $A'$ rencontre $t'Q$; $t'$ convient alors par définition de $d_m$.

Supposons maintenant $d_0\geq 1$; soit $A_0$ une chambre de $\mathcal{B}$ dont l'adhérence contient $A'$ et telle que $d_c\left(A_0,\mathcal{A}\right)=d_0$, $A_1$ une chambre adjacente à $A$ et telle que $d_c\left(A_1,\mathcal{A}\right)=d_0-1$, $\mathcal{A}'$ un appartement contenant à la fois $A_1$ et $A'$ et $H$ le mur de $\mathcal{A}'$ contenant la cloison séparant $A_0$ et $A_1$.
Par hypothèse de récurrence appliquée à $A_1$, il existe $t'\in X$, une galerie tendue de la forme $\left(A_1,A_2,\dots,A_l\right)$ et un entier $l'$ tels que $l'\leq d_0+d_m$, $tA_x\subset\overline{A_l}$ et $t'A_x\subset\overline{A_{l'}}$. De plus, pour tout $i\geq 1$, $A_i\subset\mathcal{A}$ et est située du même côté de $H$ que $A_1$ et $A_x$; on déduit alors de \cite[I, corollaire 2.8.4]{bt} qu'il existe un appartement $\mathcal{A}''$ de $\mathcal{B}$ contenant $A_0$ et le demi-plan de $\mathcal{A}'$ délimité par $H$ contenant les $A_i$, $i>0$; puisque $H$ sépare $A_0$ et $A_1$, $\left(A_0,\dots,A_l\right)$ est une galerie tendue. Donc si $E$ contient au moins une chambre, il contient tous les $\overline{A_i}$, et en particulier $t'A_x$.

Il reste donc à traiter le cas où $E$ ne contient aucune chambre de $\mathcal{B}$. Soit $H'$ le sous-espace de $\mathcal{A}''$ engendré par $A'$ et $tA_x$. Puisque $d_0\geq 1$, $H'$ ne peut pas être réduit à un point, et on peut supposer que $\overline{A_1}$ rencontre $H'$: en effet, tous les murs de $\mathcal{A}''$ qui séparent $A'$ et $tA_x$ coupent un segment reliant $tx$ et un point de $A'$, donc coupent $E$, donc contiennent une facette de $E$; on peut alors choisir $A_1,\dots,A_l$ de façon à ce que $\overline{A_1}$ contienne une facette de $E$. On peut alors déduire l'assertion du lemme de l'hypothèse de récurrence appliquée à une telle facette.

Supposons maintenant $\gamma$ quelconque. Soit $\mathcal{A}'$ un appartement de $\mathcal{B}$ contenant $E$, et $T'_0$ le tore maximal de $G$ associé; on peut supposer que $\gamma'\in N_G\left(T'_0\right)$; d'autre part, puisque $\gamma'K'\cap t\gamma K_xt^{-1}$ n'est pas vide, d'après les lemmes \ref{ilexpas} et \ref{wprim} et quitte à remplacer $\gamma'$ par $\gamma'w'$, où $w'$ est un élément de $N_G\left(T'_0\right)\cap K'$, on peut supposer que $\gamma'$ stabilise $tA_x$. L'ensemble $E$ est alors également stable par $\gamma'$, et l'assertion du lemme est vraie s'il existe $t'\in X$ tel que $t'A_x$ est contenue dans $E$ et stable par $\gamma'$ et que $d_c\left(A',t'A_x\right)\leq d_0+d_m+2d_k$. Soit $V$ l'espace vectoriel associé à $\mathcal{A}$, et $\mathcal{A}'{}^{\gamma'}$ le sous-espace des éléments de $\mathcal{A}'$ stables par $\gamma'$; il existe $v\in V$ tel que $\mathcal{A}'{}^{\gamma'}+v$ soit réunion des facettes qu'il contient; il est clair que c'est alors également le cas de $\mathcal{A}'{}^{\gamma'}+wv$ pour tout $w\in W'\left(G/T'_0\right)$ commutant avec $\gamma$.
Soit $A_0$ (resp. $A_l$; on supposera $v$ tel que $\mathcal{A}'{}^{\gamma'}+v$ contient un sommet $y'$ de $A_l$. Soit $y=t^{-1}y'$, et soit $A''$ une facette contenue dans $\overline{A_0}\cap\left(\mathcal{A}'{}^{\gamma'}+v\right)$ (il en existe: en effet, on a $A_0=wA_l$, avec $w\in W'\left(G/T'_0\right)$ commutant avec $\gamma'$, et $\overline{A_0}\cap\left(\mathcal{A}'{}^{\gamma'}+v\right)=w\left(\overline{A_l}\cap\left(\mathcal{A}'{}^{\gamma'}+w^{-1}v\right)\right)$ est donc non vide); d'après le cas $\gamma=1$, il existe $t'\in X$ tel que $d_c\left(A'',t'y\right)\leq d_0+d_m$ et que $t'y$ appartient à l'enclos de $A''\cup ty$. Or cet enclos est inclus dans $\left(\mathcal{A}'{}^{\gamma'}\cap E\right)+v$; on en déduit que $t't^{-1}$ stabilise $\mathcal{A}'{}^{\gamma'}$, donc que $t'A_x$ rencontre $\mathcal{A}'{}^{\gamma'}$; puisque $\gamma t'A_x$ est de même type que $t'A_x$, on en déduit qu'elles sont égales. D'autre part, $t'y-v\in E\cap t'A_x$, donc $t'A_x\subset E$. Enfin, on a, grâce au lemme \ref{intr}:
\begin{multline*}d_c\left(A',t'A_x\right)\leq d_c\left(A',A''\right)+d_c\left(A'',t'y\right)+d_c\left(t'y,t'A_x\right)+2d_k\\\leq 0+\left(d_0+d_m\right)+0+2d_k,\end{multline*}
ce qui achève la démonstration.
\end{proof}
\end{proof}

D'après ce lemme, l'ensemble:
\[C_0=C\cap C_{x,\gamma,0}=C\cap\bigcup_{t\in Z_G\left(T_0\right)}t^{-1}\gamma K_x t\]
est un fermé de $C$. Il est donc compact; puisque les $t^{-1}\gamma K_x t$ sont des ouverts, il existe alors une partie finie $T''_1$ de $Z_G\left(T_0\right)$ telle que l'on a:
\[C_0\subset\bigcup_{t\in T''_1}t\gamma K_x t^{-1}.\]
Pour tout $t'\in Z_G\left(T_0\right)$ tel que $t'{}^{-1}Ct'\cap\gamma K_x\neq\emptyset$, il existe donc $t\in T''_1$ tel que $t'\gamma K_x t'{}^{-1}$ rencontre $t\gamma K_x t^{-1}$; d'après le lemme \ref{ilexpas}, on a alors $t'\in tT_2$, ce qui achève la démonstration du lemme \ref{mder}.
\end{proof}

D'après ce lemme, il suffit de considérer les éléments de $T_1\left(T_0\cap M_{der}\right)$. Cet ensemble est réunion finie de classes à gauche modulo $K_{Z_G\left(T_0\right)}X_d$; quitte à remplacer, pour $t$ décrivant un système de représentants de ces classes, $C$ par $tCt^{-1}$ (et $K_0$ par $tK_0t^{-1}$), il suffit de montrer qu'il n'existe qu'un nombre fini d'éléments de $X_d$ ne vérifiant aucune des deux conditions de l'énoncé de la proposition. De plus, soit $X^\gamma$ le sous-groupe des éléments $t$ de $X_d$ tels que $\gamma^{-1}t\gamma\in t\left(K_{T_0}\cap M_{der}\right)$; d'après les lemmes \ref{ilexpas} et \ref{wprim}, les seuls éléments de $X_d$ susceptibles de ne pas les vérifier sont ceux appartenant à $w'{}^{-1}X^\gamma w'$ pour au moins un $w'\in W'\cap K_x$; puisque ce groupe est fini, quitte à conjuguer, pour tout $w'$, $C$ et $K_0$ par $w'$, il suffit de considérer le cas où $w'=1$.

Montrons d'abord que l'on peut se limiter à l'étude d'une chambre positive de $T_0\cap M_{der}$. Soit $\left\{\alpha_1,\dots,\alpha_r\right\}$ l'ensemble de racines simples de ${\mth{G}}$ relativement à ${\mth{T}}_0=K_{T_0}/K_{T_0}^1$ correspondant au sosu-groupe parabolique minimal ${\mth{P}}_0$ image de $I$; pour tout $i\leq r$, soit $s_i$ la réflexion de $W$ correspondant à $\alpha_i$, et soit $s'_i$ (resp. $s'_{i,0}$) l'unique élément de $W'$ possédant des représentants dans $K_x$ (resp. $K_0$) et dont l'image dans $W$ est $s_i$. Les $s'_i$ engendrent le sous-groupe $W'\cap K_x$ de $W'$; les $s'_{i,0}$ engendrent un sous-groupe $W'_0$ de $W'\cap K_0$ (qui est $W'\cap K_0$ tout entier si et seulement si $K_x$ est un sous-groupe parahorique maximal spécial de $G$). Les groupes $W'\cap K_x$ et $W'_0$ sont tous deux canoniquement isomorphes (par la surjection canonique de $W'$ dans $W$) au sous-groupe $W_x$ de $W$ engendré par les $s_i$. Pour tout $w\in W_x$, il existe donc $t_w\in X$ tel que pour tout $j\in\left\{1,\dots,r\right\}$, si $w'$ (resp. $w'_0$) est l'unique élément de $W'\cap K_x$ (resp. $W'\cap K_0$) d'image $w$ dans $W$, on a:
\[\langle \beta_j,w'_0{}^{-1}tw'\rangle =\langle \beta_j,w^{-1}tw\rangle +\langle \beta_j,t_w\rangle .\]
D'autre part, on a le lemme suivant:

\begin{lemme}
Pour tout $j\in\left\{1,\dots,r\right\}$, il existe un élément $w_j$ de $W_x$ tel que $\langle \beta_j,w_j^{-1}tw_j\rangle$ est positif ou nul.
\end{lemme}

\begin{proof}
En effet, considérons l'élément $t_x=\prod_{w\in W_x}w^{-1}tw$ de $T_0$. Il est invariant par conjugaison par $W_x$, donc en particulier, pour tout $k$, par $s_k$, d'où:
\[\langle \alpha_k,t_x\rangle =\langle s_{\alpha_k}\alpha_k,s_k,t_x\rangle =\langle -\alpha_k,t_x\rangle,\]
d'où $\langle \alpha_k,t_x\rangle =0$. Or $\beta_j$ est combinaison linéaire des $\alpha_k$, donc on a:
\[0=\langle \beta_j,t_x\rangle =\sum_{w\in W_x}\langle \beta_j,w^{-1}tw\rangle ;\]
on en déduit qu'au moins un terme de la somme est positif ou nul, ce qui démontre le lemme.
\end{proof}

On en déduit:
\[\langle \beta_j,w'_0{}^{-1}tw'\rangle \geq\langle \beta_j,t_{w_j}\rangle \geq\inf_{w\in W_x}\langle \beta_j,t_w\rangle .\]
Or les $w'_0{}^{-1}tw'$, $w'\in W_{\mu'}$, représentent tous la même double classe de $G$ modulo $K_0$ à gauche et $K_x$ à droite; on en déduit, en posant pour tout $j$:
\[c_j=\inf_{w\in W_x}\langle \beta_j,t_w\rangle ,\]
que ces doubles classes sont indexées par une partie de l'ensemble suivant:
\[E=\left(c_1+{\mth{N}}\right)\times\dots\times\left(c_r+{\mth{N}}\right);\]
il suffit donc de considérer l'ensemble $X^{\gamma,+}$ des $t\in X^\gamma$ tels que l'on a:
\[\left(\langle \beta_1,t\rangle ,\dots,\langle \beta_r,t\rangle \right)\in E.\]

Soit $O_1,\dots,O_k$ les orbites de $\left\{\beta_1,\dots,\beta_r\right\}$ pour l'action de $\gamma$. Pour tout $i\in\left\{1,\dots,k\right\}$, posons:
\[a_i=\sup_{j\mid \beta_j\in O_i}\sup_{t''\in T''_1}\langle \beta_j,t''\rangle .\]
Pour tout $i$, tout $t''\in T_1$ et tout $t'\in X^{\gamma,+}$ tel que si $\beta_j\in O_i$, $\langle \beta_j,t'\rangle >a_i$, l'image de $t'{}^{-1}t''\gamma K_x t''{}^{-1}t'\cap\gamma K_x$ dans ${\mth{G}}^+$ est contenue dans le sous-groupe parabolique propre ${\mth{P}}_i^+={\mth{M}}_i^+{\mth{U}}_i$ de ${\mth{G}}^+$ standard relativement à ${\mth{P}}_0$ et ${\mth{T}}_0$ et tel que ${\mth{M}}_i^+$ admet pour système de racines l'ensemble des racines de ${\mth{G}}^+$ qui sont combinaisons linéaires d'éléments de $O_{i'}$, $i'\neq i$; comme d'après la démonstration du lemme \ref{mder}, $t'{}^{-1}Ct'\cap\gamma K_x$ est réunion des $t'{}^{-1}t''\gamma K_x t''{}^{-1}t'\cap\gamma K_\gamma$, $t''\in T''_1$, son image est alors contenue dans ${\mth{P}}_i^+$. De plus, soit $\alpha$ une racine de ${\mth{G}}^+$ qui n'est pas racine de ${\mth{M}}_i^+$, et soit $U_{\alpha,H}=U_\alpha\cap H$; $U_{\alpha,H}$ est un sous-groupe ouvert de $U_\alpha$, et il existe donc $b_{i,\alpha}\in{\mth{Q}}$ possédant la propriété suivante: si $t\in X^\gamma$ est tel que $\langle \beta_j,t\rangle \geq b_{i,\alpha}$ pour au moins un $j\in O_i$, $t^{-1}U_{\alpha,H}t$ contient $U_\alpha\cap K_x$. Posons donc:
\[b_i=\sup\left(a_i,\sup_{\alpha}b_{i,\alpha}\right),\]
la somme portant sur les $\alpha$ qui sont racines de ${\mth{G}}^+$ et pas de ${\mth{M}}_i^+$. Si $\langle \beta_j,t'\rangle \geq b_i$, d'une part l'image de $t'{}^{-1}Ct'\cap\gamma K_\gamma$ dans ${\mth{G}}^+$ est contenue dans ${\mth{P}}_i^+$, et d'autre part $t^{-1}Ht$ contient $U_\alpha\cap K_x$ pour tout $\alpha$, donc son image dans ${\mth{G}}$ contient ${\mth{U}}_i$; $t'$ vérifie donc la deuxième condition de la proposition.

On en déduit que pour que $t'\in X^{\gamma,+}$ ne vérifie aucune des deux conditions de la proposition, il est nécessaire que pour tout $i$ et tout $\beta_j\in O_i$, $\langle \beta_j,t'\rangle \leq b_{i,\alpha}$. Or les éléments de $X^{\gamma,+}$ vérifiant ces inégalités sont en nombre fini, ce qui achève la démonstration de la proposition. $\Box$

\subsection{Distributions $\tilde{D}$}

Fixons à nouveau un sous-groupe d'Iwahori $I$ de $G$.
Pour tout sous-groupe parahorique $K$ de $G$ contenant $I$, soit $\mathcal{C}_K$ l'espace des fonctions sur $G$ à support compact et biinvariantes par $K^1$, et $\mathcal{C}_{K,0}$ le sous-espace des éléments de $\mathcal{C}_K$ à support dans $N_G\left(K\right)$; si ${\mth{G}}^+=N_G\left(K\right)/K^1$, $\mathcal{C}_{K,0}$ s'identifie canoniquement à l'espace des fonctions sur ${\mth{G}}^+$. Si $x$ est un point de $\mathcal{B}$ dont le fixateur connexe est $K$, on notera également $\mathcal{C}_x=\mathcal{C}_K$, $\mathcal{C}_{x,0}=\mathcal{C}_{K_0}$.

Posons également:
\[\mathcal{C}=\sum_{K\supset I}\mathcal{C}_K;\]
\[\mathcal{C}_0=\sum_{K\supset I}\mathcal{C}_{K,0}.\]

Soit $\mathcal{D}_c$ l'espace des distributions sur $G$ invariantes par conjugaison et à support dans l'ensemble $G_c$ des éléments compacts de $G$. Pour tout $\nu=\left(\gamma,\mu\right)\in\mathcal{N}$, tout couple $\left(\gamma,K\right)$ d'image $\left(\gamma,\mu\right)$ par $\zeta$ et toute fonction $\gamma$-cuspidale $\phi$ sur ${\mth{G}}=K/K^1$, $D_\phi^G$ est un élément de $\mathcal{D}_c$. Soit donc $\tilde{\mathcal{D}_c}=\tilde{\mathcal{D}_c}\left(G\right)$ le sous-espace de $\mathcal{D}_c$ engendré par les $D_\phi^G$;
le but de cette section est d'associer de manière canonique à tout $D\in\mathcal{D}_c$ un élément $\tilde{D}$ de $\tilde{\mathcal{D}_c}$, qui coïncide avec $D$ sur $\mathcal{C}_0$.

Soit $D\in\mathcal{D}_c$, $\nu\in\mathcal{N}$, $\left(\gamma,x\right)\in\mathcal{G}$ tels que $\zeta\left(\gamma,x\right)=\nu$ et que $K_x$ contient $I$; on identifie $\gamma$ à un élément du normalisateur de $I$. Soit $C_{cusp}\left(\gamma K_x\right)$ l'espace des fonctions de $G$ dans ${\mth{C}}$ à support dans $\gamma K_x$, invariantes par $K_x^1$, invariantes par conjugaison par $K_x$ et $\gamma$-cuspidales, c'est-à-dire qui sont le relèvement à $G$ tout entier d'une fonction cuspidale sur ${\mth{G}}_x^+=\gamma^{\mth{Z}}K_x/K_x^1$; $C_{cusp}\left(\gamma K_x\right)$ est canoniquement isomorphe à l'espace $C_{cusp}\left({\mth{G}}_x\right)^{{\mth{G}}_x,\gamma}$ des fonctions invariantes par $\gamma$-conjugaison et $\gamma$-cuspidales sur ${\mth{G}}_x=K_x/K_x^1$, et la restriction de $D$ à $C_{cusp}\left(\gamma K_x\right)$ s'identifie à une distribution sur $C_{cusp}\left({\mth{G}}_x\right)^{{\mth{G}}_x,\gamma}$ invariante par $\gamma$-conjugaison; il existe donc une unique fonction $\phi_{\gamma,x}\left(D\right)\in C_{cusp}\left({\mth{G}}_x\right)^{{\mth{G}}_x,\gamma}$ telle que l'on a, pour tout $f\in C_{cusp}\left(\gamma K_x\right)$:
\[D\left(f\right)=v'_\nu\langle \phi_{\gamma,x}\left(D\right),f\rangle _{{\mth{G}}_x},\]
où $v'_\nu$ est défini de la manière suivante:
soit $\mu=\kappa\left(x\right)$, $\left(M,K_\mu\right)$ un représentant de $\mu$, $T_0$ un tore déployé maximal de $M$, $M_\gamma$ le plus petit sous-groupe de Levi de $G$ semi-standard relativement à $T_0$ et contenant $M$ et $\gamma$ (assimilé à un élément de $W'\left(G/T_0\right)$), et $v_\nu$ la constante définie dans l'énoncé du lemme \ref{dscp}; on pose:
\[v'_\nu=\left[N_G\left(\gamma,K_\mu\right):N_{M_\gamma}\left(\gamma,K_\mu\right)\right]v_\nu.\]

De plus, on a le lemme suivant:

\begin{lemme}\label{reges}
$\phi_{\gamma,x}$ ne dépend que de $\nu$.
\end{lemme}

\begin{proof}
En effet, soit $x'$ un autre élément de $\mathcal{B}$ tel que $\kappa\left(x'\right)=\mu$ et $\gamma\in N_G\left(K_{x'}\right)$. Soit $\mathcal{A}$ un appartement de $\mathcal{B}$ contenant $x$ et $x'$, et $T_0$ le tore maximal associé à $\mathcal{A}$; grâce à la proposition \ref{gcj} et quitte à remplacer $x'$ par $wx'$, où $w$ est un élément de $W'\left(G/T_0\right)$, on peut supposer que $\left(K_x,T_0\right)$ et $\left(K_{x'},T_0\right)$ ont même image $\left(M,K_\mu\right)$ par $\kappa_0$, et que $\gamma$ normalise également $K_{x'}$. Si $A_x$ (resp. $A_{x'}$) est la facette de $\mathcal{A}$ contenant $x$ (resp. $x'$), $A_x$ et $A_{x'}$ sont stables par $\gamma$ et engendrent le même sous-espace affine $E$ de $\mathcal{A}$, qui est le sous-espace des points fixés par $K_\mu$. En considérant les facettes de $E$ stables par $\gamma$ maximales, on voit qu'il existe une suite $A_x=A_0,A_1,\dots,A_k=A_{x'}$ de telles facettes telles que pour tout $i$, $A_i$ et $A_{i+1}$ ont en commun au moins une face stable par $\gamma$; comme pour tout $i$, tout point de $A_i$ admet $\mu$ comme image par $\kappa$, il suffit de montrer le lemme dans le cas où $A_x$ et $A_{x'}$ ont en commun une telle face; mais alors il existe un sous-groupe parahorique $K$ de $G$ stable par $\gamma$ et contenant à la fois $K_x$ et $K_{x'}$. Soit ${\mth{G}}=K/K^1$; ${\mth{G}}_x$ et ${\mth{G}}_{x'}=K_{x'}/K_{x'}^1$ s'identifient à un même sous-groupe de Levi de ${\mth{G}}$. D'autre part, 
si ${\mth{P}}={\mth{M}}{\mth{U}}$ est un sous-groupe parabolique de ${\mth{G}}$ stable par $\gamma$ et de Levi ${\mth{M}}$, et si $f$ est une fonction de ${\mth{M}}$ dans ${\mth{C}}$, en définissant l'induite de $f$ à ${\mth{G}}$ tordue par $\gamma$ par:
\[\Ind_{{\mth{M}},\gamma}^{\mth{G}}\left(f\right):g\in{\mth{G}}\longmapsto\dfrac 1{\card\left({\mth{P}}\right)}\sum_{x\in{\mth{G}},m\in{\mth{M}},u\in{\mth{U}},g=\gamma^{-1}x\gamma mux}f\left(m\right),\]
cette induite est un élément de $C\left({\mth{G}}\right)^{{\mth{G}},\gamma}$; de plus, on déduit immédiatement de la proposition \ref{resind} appliquée à ${\mth{G}}^+=\gamma^{\mth{Z}}{\mth{G}}$ qu'elle ne dépend pas du choix de ${\mth{P}}$.
Si $\gamma=1$, on retrouve bien évidemment l'induite usuelle. On a:

\begin{lemme}\label{find}
Soit ${\mth{P}}={\mth{M}}{\mth{U}}$ un sous-groupe parabolique de ${\mth{G}}$ stable par $\gamma$, et $K_P$ son image réciproque dans $K$.
Soit $f$ une fonction sur ${\mth{M}}$ invariante par $\gamma$-conjugaison; on a, en relevant $f$ et son induite, ainsi que les éléments de $C_0\left({\mth{G}}\right)$, à des fonctions sur $G$ à support dans $K$:
\[\Ind_{{\mth{M}},\gamma}^{\mth{G}}\left(f\right)\in \dfrac{\card\left({\mth{G}}\right)}{\card\left({\mth{P}}\right)}f+C_0\left({\mth{G}}\right).\]
\end{lemme}

\begin{proof}
Puisque $K_x\subset K$ et $K_x^1\supset K^1$, $f$ (en tant que fonction sur $G$) est une fonction à support dans $K$ et biinvariante par $K^1$, et s'identifie donc à une fonction sur ${\mth{G}}$; plus précisément, c'est alors une fonction à support dans ${\mth{P}}$ constante sur les $m{\mth{U}}$, $m\in{\mth{M}}$. Or on a, pour tout $g\in{\mth{G}}$:
\begin{equation*}\begin{split}\Ind_{{\mth{M}},\gamma}^{\mth{G}}f\left(g\right)&=\dfrac 1{\card\left({\mth{P}}\right)}\sum_{h\in{\mth{G}},m\in{\mth{M}},u\in{\mth{U}},g=\gamma^{-1}h^{-1}\gamma muh}f\left(m\right)\\&=\dfrac 1{\card\left({\mth{P}}\right)}\sum_{h\in{\mth{G}},p\in{\mth{P}},g=\gamma^{-1}h^{-1}\gamma ph}f\left(p\right)\\&=\dfrac 1{\card\left({\mth{P}}\right)}\sum_{h\in{\mth{G}}}f\left(\gamma^{-1}h\gamma gh^{-1}\right).\end{split}\end{equation*}
On a donc:
\[\Ind_{{\mth{M}},\gamma}^{\mth{G}}f=\dfrac 1{\card\left({\mth{P}}\right)}\sum_{h\in{\mth{G}}}Ad_\gamma\left(h\right)f,\]
et l'assertion du lemme s'en déduit immédiatement.
\end{proof}

En appliquant ce lemme successivement à ${\mth{P}}_x=K_x/K^1$ et à ${\mth{P}}_{x'}=K_{x'}/K^1$, puisque ces groupes admettent tous deux ${\mth{G}}_x$ comme composante de Levi, on en déduit que pour $f\in C_{cusp}\left({\mth{G}}_x\right)^{{\mth{G}}_x,\gamma}$, les relèvements de $f$ à respectivement $\gamma K_x$ et $\gamma K_{x'}$ sont des éléments de $C_c^\infty\left(G\right)$ qui ne diffèrent que d'un élément de $C_0\left(G\right)$; les restrictions de $D$ à $C_{cusp}\left(\gamma K_x\right)$ et $C_{cusp}\left(\gamma K_{x'}\right)$ induisent donc la même distribution sur $C_{cusp}\left({\mth{G}}_x\right)^{{\mth{G}}_x,\gamma}$, ce qui démontre le lemme \ref{reges}.
\end{proof}

Dans ce qui suit, on écrira $\phi_\nu\left(D\right)$ à la place de $\phi_{\gamma,x}\left(D\right)$; pour tout $\nu\in\mathcal{N}$ n'admettant aucun antécédent par $\zeta$, on posera $\phi_\nu\left(D\right)=0$. Posons:
\[\tilde{D}=\sum_{\nu\in\mathcal{N}}D_{\overline{\phi_\nu\left(D\right)}},\]
C'est un élément de $\tilde{\mathcal{D}_c}$, et l'application $D\mapsto\tilde{D}$ est clairement linéaire.

D'après le lemme \ref{dscp}, on a, pour tout $\nu=\left(\gamma,\mu\right)$ et toute fonction $\phi$ $\gamma$-cuspidale sur le groupe réductif fini ${\mth{M}}$ associé à $\mu$, si $M$, $K_\mu$ et $M_\gamma$ sont définis comme dans l'énoncé de ce lemme:
\begin{equation*}\begin{split}\tilde{D}\left(\phi\right)=D_{\phi_\nu\left(D\right)}^G\left(\phi\right)&=v_\nu\sum_{n\in N_G\left(\gamma,K_\mu\right)/N_{M_\gamma}\left(\gamma,K_\mu\right)}\langle \phi_\nu,Ad\left(n\right)\phi\rangle _{\mth{M}}\\&=v_\nu\sum_{n\in N_G\left(\gamma,K_\mu\right)/N_{M_\gamma}\left(\gamma,K_\mu\right)}\langle Ad\left(n^{-1}\right)\phi_\nu,\phi\rangle _{\mth{M}}.\end{split}\end{equation*}
Or puisque $D$ est invariante, $\phi_\nu$ est invariante par conjugaison par $N_G\left(\gamma,K_\mu\right)$ et on en déduit:
\[\tilde{D}\left(\phi\right)=D_{\phi_\nu\left(D\right)}^G\left(\phi\right)=v'_\nu\langle \phi_\nu,\phi\rangle _{\mth{M}}=D\left(\phi\right).\]

Les distributions $D$ et $\tilde{D}$ coïncident donc sur l'espace engendré par les $\phi$; on va maintenant montrer que modulo $C_0\left(G\right)$, les $\phi$ engendrent $\mathcal{C}_0$. Pour cela, on va montrer un résultat de décomposition pour les fonctions invariantes par $\gamma$-conjugaison sur un groupe fini.

Soit ${\mth{G}}$ le groupe des ${\mth{F}}_q$-points d'un groupe réductif connexe défini sur ${\mth{F}}_q$, et ${\mth{G}}^+$ un groupe engendré par ${\mth{G}}$ et un élément $\gamma$ normalisant ${\mth{G}}$.
Soit ${\mth{M}}$ un sous-groupe de Levi de ${\mth{G}}$ stable par conjugaison par $\gamma$; son normalisateur $N_{{\mth{G}},\gamma}\left({\mth{M}}\right)$ pour l'action ci-dessus est également stable par $\gamma$. On notera $C_{cusp}\left({\mth{M}}\right)^{N_{{\mth{G}},\gamma}\left({\mth{M}}\right),\gamma}$ le sous-espace des éléments $\gamma$-cuspidaux de $C\left({\mth{M}}\right)^{N_{{\mth{G}},\gamma}\left({\mth{M}}\right),\gamma}$.

Soit $R$ un système de représentants des classes de conjugaison des sous-groupes de Levi de ${\mth{G}}$ stables par $\gamma$; on a le résultat suivant:

\begin{prop}\label{dcc}
On a les décompositions suivantes:
\begin{equation*}\begin{split}C\left({\mth{G}}\right)^{{\mth{G}},\gamma}&=\bigoplus_{{\mth{M}}\in R}\Ind_{{\mth{M}},\gamma}^{\mth{G}}C_{cusp}\left({\mth{M}}\right)^{{\mth{M}},\gamma}\\&=\bigoplus_{{\mth{M}}\in R}\Ind_{{\mth{M}},\gamma}^{\mth{G}}C_{cusp}\left({\mth{M}}\right)^{N_{{\mth{G}},\gamma}\left({\mth{M}}\right),\gamma}.\end{split}\end{equation*}
\end{prop}

\begin{proof}
Montrons d'abord l'égalité des deux sommes de l'énoncé: pour cela, il suffit de montrer que pour tout ${\mth{M}}\in R$, on a:
\[\Ind_{{\mth{M}},\gamma}^{\mth{G}}C_{cusp}\left({\mth{M}}\right)^{{\mth{M}},\gamma}\subset \Ind_{{\mth{M}},\gamma}^{\mth{G}}C_{cusp}\left({\mth{M}}\right)^{N_{{\mth{G}},\gamma}\left({\mth{M}}\right),\gamma},\]
l'autre inclusion étant évidente. Soit donc $f$ un é\-lé\-ment de $C_{cusp}\left({\mth{M}}\right)^{{\mth{M}},\gamma}$; pour tout $m\in N_{{\mth{G}},\gamma}\left({\mth{M}}\right)$, en posant:
\[Ad_\gamma\left(m\right)\left(f\right):g\longmapsto f\left(\gamma^{-1}m\gamma gm^{-1}\right),\]
$Ad_\gamma\left(m\right)f$ est également $\gamma$-cuspidale et on a:
\[\Ind_{{\mth{M}},\gamma}^{\mth{G}}\left(Ad_\gamma\left(m\right)f\right)=\Ind_{{\mth{M}},\gamma}^{\mth{G}}\left(f\right).\]
Posons donc:
\[f_0=\dfrac 1{\card\left(N_{{\mth{G}},\gamma}\left({\mth{M}}\right)\right)}\sum_{m\in N_{{\mth{G}},\gamma}\left({\mth{M}}\right)}Ad_\gamma\left(m\right)f;\]
on a $f_0\in C_{cusp}\left({\mth{M}}\right)^{N_{{\mth{G}},\gamma}\left({\mth{M}}\right),\gamma}$ et $\Ind_{{\mth{M}},\gamma}^{\mth{G}}\left(f_0\right)=\Ind_{{\mth{M}},\gamma}^{\mth{G}}\left(f\right)$, ce qui montre l'assertion cherchée. De plus, si la première somme est directe, la seconde l'est aussi.

On va maintenant montrer, par récurrence sur le rang semi-simple de ${\mth{G}}$, la première égalité, ainsi que le fait que la première somme est directe. Si $rg\left({\mth{G}}\right)=0$, $R=\left\{{\mth{G}}\right\}$, toute fonction invariante par $\gamma$-conjugaison sur ${\mth{G}}$ est $\gamma$-cuspidale et l'assertion est triviale. Supposons donc ${\mth{G}}$ de rang non nul; montrons d'abord les lemmes suivants:

\begin{lemme}\label{vin}
On a $C\left({\mth{G}}\right)^{{\mth{G}},\gamma}=C_{cusp}\left({\mth{G}}\right)^{{\mth{G}},\gamma}+\sum_{{\mth{M}}\subsetneq{\mth{G}}}\Ind_{{\mth{M}},\gamma}^{\mth{G}}C\left({\mth{M}}\right)^{{\mth{M}},\gamma}$.
\end{lemme}

\begin{proof}
Pour montrer ce lemme, il suffit de montrer que $C_{cusp}\left({\mth{G}}\right)^{{\mth{G}},\gamma}$ contient l'orthogonal de $\sum_{{\mth{M}}\subsetneq{\mth{G}}}\Ind_{{\mth{M}},\gamma}^{\mth{G}}C\left({\mth{G}}\right)^{{\mth{G}},\gamma}$ dans ${\mth{G}}$ pour $\langle \cdot,\cdot\rangle _{\mth{G}}$; puisque ce produit hermitien est défini positif, l'assertion du lemme s'en déduit im\-mé\-dia\-te\-ment. Soit donc $f$ un élément de cet orthogonal; pour tout ${\mth{M}}\subsetneq{\mth{G}}$ et tout $f'\in C\left({\mth{M}}\right)^{{\mth{M}},\gamma}$, on a:
\[\langle f,\Ind_{{\mth{M}},\gamma}^{\mth{G}}f'\rangle _{\mth{G}}=0.\]
D'après le lemme précedent, on a donc:
\[\langle r_{\mth{M}}^{\mth{G}}f,f'\rangle _{\mth{M}}=0;\]
comme ceci est vrai pour tout ${\mth{M}}$ et tout $f'$, on a $r_{\mth{M}}^{\mth{G}}f=0$ pour tout ${\mth{M}}$, donc $f$ est $\gamma$-cuspidale.
\end{proof}

Pour tout ${\mth{M}}\in\mathcal{F}$, soit $W_{{\mth{M}},\gamma}$ l'ensemble des éléments commutant avec $\gamma$ du groupe de Weyl de ${\mth{M}}$ relativement à ${\mth{T}}_0$.

\begin{lemme}\label{fdm}
Soit ${\mth{M}},{\mth{M}}'\in R$; il existe un système de représentants $W_0$ de $W_{{\mth{M}}',\gamma}\backslash W_{{\mth{G}},\gamma}/W_{{\mth{M}},\gamma}$ tel que l'on a:
\[r_{{\mth{M}}'}^{\mth{G}}\circ \Ind_{{\mth{M}},\gamma}^{\mth{G}}=\sum_{w\in W_0}\Ind_{{\mth{M}}'\cap w{\mth{M}}w^{-1},\gamma}^{{\mth{M}}'}\circ Ad\left(w\right)\circ r_{{\mth{M}}\cap w^{-1}{\mth{M}}'w,\gamma}^{\mth{M}}.\]
\end{lemme}

\begin{proof}
Ce lemme est une conséquence immédiate du théorème 3.2 de \cite{dm}.
\end{proof}

\begin{lemme}\label{vdx}
Soit ${\mth{M}},{\mth{M}}'\in R$, $\phi$ (resp. $\phi'$) une fonction $\gamma$-cuspidale sur ${\mth{M}}$ (resp. ${\mth{M}}'$) invariante par $N_{{\mth{G}},\gamma}\left({\mth{M}}\right)$ (resp. $N_{{\mth{G}},\gamma}\left({\mth{M}}'\right)$). Alors on a:
\[\langle \Ind_{{\mth{M}},\gamma}^{\mth{G}}\phi,\Ind_{{\mth{M}}',\gamma}^{\mth{G}}\phi'\rangle _{\mth{G}}=0\]
si ${\mth{M}}$ et ${\mth{M}}'$ sont distincts, et:
\[\langle \Ind_{{\mth{M}},\gamma}^{\mth{G}}\phi,\Ind_{{\mth{M}}',\gamma}^{\mth{G}}\phi'\rangle _{\mth{G}}=\card\left(N_{{\mth{G}},\gamma}\left({\mth{M}}\right)/{\mth{M}}\right)\langle \phi,\phi'\rangle _{{\mth{G}}'}\]
si ${\mth{M}}={\mth{M}}'$.
\end{lemme}

\begin{proof}
Supposons d'abors ${\mth{M}}\neq{\mth{M}}'$. Puisque ${\mth{M}}$ et ${\mth{M}}'$ sont deux é\-lé\-ments de $R$ distincts, au moins l'un d'eux n'est pas contenu dans un conjugué de l'autre; supposons par exemple que ${\mth{M}}$ n'est pas contenu dans un conjugué de ${\mth{M}}'$. On a par le lemme \ref{rpf}:
\[\langle \Ind_{{\mth{M}},\gamma}^{\mth{G}}\phi,\Ind_{{\mth{M}}',\gamma}^{\mth{G}}\phi'\rangle _{\mth{G}}=\langle r_{{\mth{M}}'}^{\mth{G}}\Ind_{{\mth{M}},\gamma}^{\mth{G}}\phi,\phi'\rangle _{{\mth{M}}'}.\]
Or le lemme \ref{fdm} donne:
\begin{equation}\label{mck}
r_{{\mth{M}}'}^{\mth{G}}\Ind_{{\mth{M}},\gamma}^{\mth{G}}\phi=\sum_{w\in W_0}\Ind_{{\mth{M}}'\cap w{\mth{M}}w^{-1},\gamma}^{{\mth{M}}'}\circ Ad\left(w\right)\circ r_{{\mth{M}}\cap w^{-1}{\mth{M}}'w,\gamma}^{\mth{M}}\left(\phi\right),
\end{equation}
où $W_0$ est défini comme dans l'énoncé de ce lemme; d'autre part, pour tout $w$, puisque ${\mth{M}}$ n'est contenu dans aucun conjugué de ${\mth{M}}'$, ${\mth{M}}\cap w^{-1}{\mth{M}}'w$ est un sous-groupe de Levi propre de ${\mth{M}}$; comme $\phi$ est $\gamma$-cuspidale, on en déduit $r_{{\mth{M}}\cap w^{-1}{\mth{M}}'w}^{\mth{M}}\left(\phi\right)=0$. On a donc:
\[r_{{\mth{M}}'}^{\mth{G}}\Ind_{{\mth{M}},\gamma}^{\mth{G}}\phi'=0,\]
ce qui démontre la première assertion du lemme \ref{vdx}.

Supposons maintenant ${\mth{M}}={\mth{M}}'$. Alors les seuls termes non nuls dans le membre de droite de (\ref{mck}) sont ceux correspondant aux $w$ tels que ${\mth{M}}\cap w^{-1}{\mth{M}}w={\mth{M}}$, c'est-à-dire $w\in N_{{\mth{G}},\gamma}\left({\mth{M}}\right)$; pour un tel $w$, on a:
\[\Ind_{{\mth{M}}\cap w{\mth{M}}w^{-1},\gamma}^{\mth{M}}\circ Ad\left(w\right)\circ r_{{\mth{M}}\cap w^{-1}{\mth{M}}w}^{\mth{M}}\left(\phi\right)=Ad\left(w\right)\phi=\phi\]
puisque $\phi$ est invariante par $N_{{\mth{G}},\gamma}\left({\mth{M}}\right)$. On en déduit la deuxième assertion du lemme \ref{vdx}.
\end{proof}

Revenons à la démonstration de la proposition \ref{dcc}. Par hypothèse de ré\-cur\-rence, on a pour tout ${\mth{M}}\subsetneq{\mth{G}}$:
\[C\left({\mth{M}}\right)^{{\mth{M}},\gamma}=\bigoplus_{{\mth{M}}'\in R_{\mth{M}}}\Ind_{{\mth{M}}',\gamma}^{\mth{M}}C_{cusp}\left({\mth{M}}'\right)^{{\mth{M}}',\gamma},\]
où $R_{\mth{M}}$ est défini pour ${\mth{M}}$ de façon similaire à $R$ pour ${\mth{G}}$ (notons que deux sous-groupes de Levi de ${\mth{M}}$ peuvent être conjugués dans ${\mth{G}}$ et pas dans ${\mth{M}}$). On en déduit, par transitivité des induites et grâce au lemme \ref{vin}:
\[C\left({\mth{G}}\right)^{{\mth{G}},\gamma}=\sum_{{\mth{M}}\in R}\Ind_{{\mth{M}},\gamma}^{\mth{G}}C_{cusp}\left({\mth{M}}\right)^{{\mth{M}},\gamma}.\]
Il reste à montrer que la somme est directe. Pour cela, soit, pour tout ${\mth{M}}\in R$, $\phi_{\mth{M}}\in C_{cusp}\left({\mth{M}}\right)^{{\mth{M}},\gamma}$; supposons les $\phi_{\mth{M}}$ tels que l'on a:
\[\sum_{{\mth{M}}\in R}\Ind_{{\mth{M}},\gamma}^{\mth{G}}\phi_{\mth{M}}=0.\]
Pour tout ${\mth{M}}'\in R$, on a alors, d'après le lemme \ref{vdx}:
\[0=\langle \sum_{{\mth{M}}\in R}\Ind_{{\mth{M}},\gamma}^{\mth{G}}\phi_{\mth{M}},\Ind_{{\mth{M}}',\gamma}^{\mth{G}}\phi_{{\mth{M}}'}\rangle _{\mth{G}}=\card\left(N_{{\mth{M}}',\gamma}\left({\mth{M}}'\right)/{\mth{M}}'\right)\langle \phi_{{\mth{M}}'},\phi_{{\mth{M}}'}\rangle _{{\mth{M}}'};\]
on en déduit que $\phi_{{\mth{M}}'}=0$ pour tout ${\mth{M}}'$ et la proposition est démontrée.
\end{proof}

D'après cette proposition, pour tout $x\in\overline{A}$ et tout $\gamma$ normalisant $x$, l'espace des éléments de $\mathcal{C}_{x,0}$ invariants par conjugaison par $K_x$ s'identifie à l'espace engendré par les induites des fonctions $\gamma$-cuspidales sur les sous-groupes de Levi stables par $\gamma$ de ${\mth{G}}_x=K_x/K_x^1$, relevées à des fonctions sur $G$; on va en déduire, en utilisant le lemme \ref{find}, que ces fonctions $\gamma$-cuspidales engendrent $\mathcal{C}_0$, et donc que $D$ et $\tilde{D}$ coïncident sur $\mathcal{C}_0$.
Plus précisément, on a la proposition suivante:

\begin{prop}\label{dprime}
Soit $D\in\mathcal{D}_c$. Alors il existe un unique $D'\in\tilde{\mathcal{D}_c}$ tel que les restrictions de $D$ et $D'$ à $\mathcal{C}_0$ sont identiques; de plus, on a $D'=\tilde{D}$.
\end{prop}

\begin{proof}
Montrons d'abord l'unicité de $D'$: supposons que $D'$ et $D''$ sont deux éléments de $\tilde{\mathcal{D}_c}$ tels que $D'|_{\mathcal{C}_0}=D''|_{\mathcal{C}_0}=D|_{\mathcal{C}_0}$. Alors $D_0=D'-D''$ est un élément de $\tilde{\mathcal{D}_c}$ dont la restriction à $\mathcal{C}_0$ est nulle; il suffit donc de montrer qu'une telle distribution est forcément nulle. Posons:
\[D_0=\sum_\nu D_{\phi_\nu},\]
où pour tout $\nu=\left(\mu,\gamma\right)\in\mathcal{N}$, $\phi_\nu$ est une fonction $\gamma$-cuspidale sur le groupe réductif fini ${\mth{G}}$ associé à $\mu$. Si $\left(M,K_\mu\right)$ est un représentant de $\mu$, en identifiant ${\mth{G}}$ à $K_\mu/K_\mu^1$, $\phi_\nu$ est invariante par conjugaison par $N_G\left(\gamma,K_\mu\right)$; pour tout $\nu'=\left(\mu',\gamma'\right)\in\mathcal{N}$ et toute fonction $\phi'$ $\gamma$-cuspidale sur le groupe réductif fini ${\mth{G}}'$ associé à $\mu'$, on a alors, en utilisant le lemme \ref{dscp}:
\[D_0\left(\phi'\right)=\sum_\nu D_{\phi_\nu}\left(\phi'\right)=v'_{\nu'}\langle \overline{\phi_{\nu'}},\phi'\rangle _{{\mth{G}}'}.\]
Or puisque $\phi'$ se relève en un élément de $\mathcal{C}_0$, on a $D_0\left(\phi'\right)=0$; on en déduit:
\[\langle \overline{\phi_{\nu'}},\phi'\rangle _{{\mth{G}}'}=0.\]
Comme ceci est vrai pour tout $\phi'$, on a $\phi_{\nu'}=0$ pour tout $\nu'$, donc $D_0=0$.

Il reste à montrer que $\tilde{D}$ vérifie $\tilde{D}|_{\mathcal{C}_0}=D|_{\mathcal{C}_0}$. Comme les sous-espaces de la forme $\mathcal{C}_{x,0}$, $x\in\overline{A}$, engendrent $\mathcal{C}_0$, il suffit de vérifier l'égalité des restrictions de $D$ et $\tilde{D}$ à chacun de ces sous-espaces. De plus, pour tout $x\in\overline{A}$ et tout $\gamma\in\Gamma_I$ qui normalise $K_x$, soit $\mathcal{C}_{x,\gamma,0}=\mathcal{C}_{K_x,\gamma,0}$ le sous-espace des éléments de $\mathcal{C}_{x,0}$ à support dans $\gamma K_x$; $\mathcal{C}_{x,0}$ est somme directe des $\mathcal{C}_{x,\gamma,0}$. On va donc vérifier l'égalité des restrictions de $D$ et $\tilde{D}$ à chacun des $\mathcal{C}_{x,\gamma,0}$.

Soit $x\in\mathcal{B}$, $K_x$ le sous-groupe parahorique de $G$ fixant $x$ et $\gamma$ un élément de $\Gamma_I$ qui normalise $K_x$; on va montrer par récurrence sur le rang de ${\mth{G}}_x=K_x/K_x^1$ que $D$ et $\tilde{D}$ coïncident sur $\mathcal{C}_{x,\gamma,0}$. Cet espace est canoniquement isomorphe à l'espace des fonctions sur le groupe fini ${\mth{G}}_x$; il existe donc une fonction $\phi_{D,x,\gamma}$ sur ${\mth{G}}_x$ telle que pour tout $f\in\mathcal{C}_{x,\gamma,0}$, on a:
\[D\left(f\right)=\langle \phi_{D,x,\gamma},f\rangle _{{\mth{G}}_x}.\]
Puisque $D$ est invariante, $\phi_{D,x,\gamma}$ est invariante par $\gamma$-conjugaison; grâce à la proposition \ref{dcc}, on peut alors écrire de manière unique:
\[\phi_{D,x,\gamma}=\sum_{{\mth{M}}\in R}\Ind_{{\mth{M}},\gamma}^{{\mth{G}}_x}\phi_{D,x,\gamma,{\mth{M}}},\]
où $R$ est un système de représentants des classes de conjugaison de sous-groupes de Levi ${\mth{M}}$ de ${\mth{G}}_x$ stables par $\gamma$, et où pour tout ${\mth{M}}$, $\phi_{D,x,\gamma,{\mth{M}}}$ est une fonction $\gamma$-cuspidale sur ${\mth{M}}$ invariante par $N_{{\mth{G}}_x,\gamma}\left({\mth{M}}\right)$. D'autre part, pour $f\in\mathcal{C}_{x,\gamma,0}$, posons:
\[f':g\in G\longmapsto\int_{K_x}f\left(k^{-1}gk\right)dk.\]
Puisque $D$ et $\tilde{D}$ sont invariantes, on a $D\left(f\right)=D\left(f'\right)$ et $\tilde{D}\left(f\right)=\tilde{D}\left(f'\right)$; de plus, $f'$ est un élément de $\mathcal{C}_{x,\gamma,0}$, invariant par $\gamma$-conjugaison en tant que fonction sur ${\mth{G}}_x$; on peut donc écrire, en utilisant à nouveau la proposition \ref{dcc}:
\[f'=\sum_{{\mth{M}}\in R}\Ind_{{\mth{M}},\gamma}^{{\mth{G}}_x}f_{\mth{M}},\]
où pour tout ${\mth{M}}$, $f_{\mth{M}}$ est une fonction $\gamma$-cuspidale sur ${\mth{M}}$.
Or soit, pour tout ${\mth{M}}\in R$, ${\mth{P}}$ un sous-groupe parabolique de ${\mth{G}}_x$ invariant par $\gamma$ et de Levi ${\mth{M}}$, et $K'$ le sous-groupe parahorique de $G$ contenu dans $K_x$ et dont l'image dans ${\mth{G}}_x$ est ${\mth{P}}$; $f_{\mth{M}}$ se relève en une fonction sur $G$ à support dans $\gamma K'$ biinvariante par $K'_1$, c'est-à-dire en un élément de $\mathcal{C}_{K',\gamma,0}$, et si ${\mth{M}}\subsetneq{\mth{G}}_x$, on a donc par hypothèse de récurrence:
\[D\left(f_{{\mth{M}}}\right)=\tilde{D}\left(f_{{\mth{M}}}\right).\]
De plus, d'après le lemme \ref{find}, pour tout ${\mth{M}}\subsetneq{\mth{G}}_x$, on a:
\begin{equation*}\begin{split}D\left(\Ind_{{\mth{M}},\gamma}^{{\mth{G}}_x}\left(f_{\mth{M}}\right)\right)&=\dfrac{\card\left({\mth{G}}_x\right)}{\card\left({\mth{P}}\right)}D\left(f_{{\mth{M}}}\right)\\&=\dfrac{\card\left({\mth{G}}_x\right)}{\card\left({\mth{P}}\right)}\tilde{D}\left(f_{{\mth{M}}}\right)\\&=\tilde{D}\left(\Ind_{{\mth{M}},\gamma}^{{\mth{G}}_x}\left(f_{{\mth{M}}}\right)\right).\end{split}\end{equation*}
Pour achever la démonstration de la proposition, il reste donc à vérifier que l'on a $D\left(f_{{\mth{G}}_x}\right)=\tilde{D}\left(f_{{\mth{G}}_x}\right)$. 
Or d'après le lemme \ref{vdx}, on a:
\[D\left(f_{{\mth{G}}_x}\right)=\langle \phi_{D,x,\gamma,{\mth{G}}_x},f_{{\mth{G}}_x}\rangle _{{\mth{G}}_x}.\]
Ceci est vrai pour toute $f_{{\mth{G}}_x}\in C_{cusp}\left({\mth{G}}_x\right)^{{\mth{G}}_x,\gamma}$, donc si $\nu=\zeta\left(\gamma,x\right)$, on a, par définition de $\phi_\nu\left(D\right)$:
\[\phi_{D,x,\gamma,{\mth{G}}_x}=v'_\nu\phi_\nu\left(D\right).\]
D'autre part, on a, en utilisant la définition de $\tilde{D}$ et le lemme \ref{dscp}:
\begin{equation*}\begin{split}\tilde{D}\left(f_{{\mth{G}}_x}\right)&=\sum_{\nu'}D_{\overline{\phi_{\nu'}\left(D\right)}}\left(f_{{\mth{G}}_x}\right)\\&=v'_\nu\langle \phi_\nu\left(D\right),f_{{\mth{G}}_x}\rangle _{{\mth{G}}_x}\\&=D\left(f_{{\mth{G}}_x}\right),\end{split}\end{equation*}
ce qui achève la démonstration.
\end{proof}

\subsection{Restriction des distributions à $\mathcal{C}$}

On va maintenant montrer le résultat suivant:

\begin{theo}\label{dc}
Soit $D\in\mathcal{D}_c$. Si $D$ est nulle sur $\mathcal{C}_0$, alors elle est nulle sur $\mathcal{C}$.
\end{theo}

\begin{proof}
Si $D$ est nulle sur $\mathcal{C}_0$, toutes les $\phi_\nu\left(D\right)$ sont nulles, donc $\tilde{D}=0$; le résultat se déduit alors immédiatement du théorème suivant:

\begin{theo}\label{dct}
Pour tout $D\in\mathcal{D}_c$, $\tilde{D}$ coïncide avec $D$ sur $\mathcal{C}$. En particulier, si $\tilde{D}=0$, $D$ est nulle sur $\mathcal{C}$.
\end{theo}

Avant de commencer la démonstration, introduisons quelques notations sup\-plé\-men\-tai\-res. Soit $\pi$ une représentation irréductible tempérée de $G$; on définit la distribution trace compacte de $\pi$ de la manière suivante:
\[tr_c\left(\pi\right):f\in C_c^\infty\left(G\right)\longmapsto\int_Gtr\left(\pi\left(g\right)\right)1_c\left(g\right)f\left(g\right)dg,\]
où $1_c$ est la fonction caractéristique de l'ensemble des éléments compacts de $G$.

Soit $\pi$ une représentation irréductible admissible de $G$, $V$ son espace, et $K$ un sous-groupe parahorique de $G$; soit $V^{K^1}$ le sous-espace des éléments de $V$ invariants par le premier sous-groupe de congruence $K^1$ de $K$. Puisque $\pi$ est admissible, $V^{K^1}$ est de dimension finie; de plus, puisque $K^1$ est un sous-groupe normal de $N_G\left(K\right)$, $V^{K^1}$ est stable par $N_G\left(K\right)$; le groupe $N_G\left(K\right)/K^1$ agit donc sur $V^{K^1}$, et la représentation de $N_G\left(K\right)/K^1$ ainsi obtenue sera notée $\pi^{N_G\left(K\right)}$. Si $H$ est un sous-groupe de $N_G\left(K\right)$ contenant $K$, on notera $\pi^H$ la restriction de $\pi^{N_G\left(K\right)}$ à $H$.

\begin{proof}
On va montrer le théorème par récurrence sur le rang semi-simple de $G$. Si ce rang est nul, $G$ possède un unique sous-groupe parahorique $K$ et on a $N_G\left(K\right)=G$; on a alors $\mathcal{C}_0=\mathcal{C}$ et le théorème se déduit immédiatement de la proposition \ref{dprime}. Supposons maintenant $G$ de rang non nul;
on va d'abord montrer la deuxième assertion. Soit $D\in\mathcal{D}_c$ telle que $\tilde{D}=0$; d'après la proposition \ref{dprime}, $D$ est alors nulle sur $\mathcal{C}_0$.

Supposons d'abord que $D$ est combinaison linéaire de traces compactes de re\-pré\-sen\-ta\-tions tempérées; posons:
\[D=\sum_{i=1}^t\lambda_itr_c\left(\pi_i\right),\]
où les $\lambda_i$ sont des constantes et les $\pi_i$ des représentations irréductibles tempérées de $G$.

\begin{prop}
Si $D$ est nulle sur $\mathcal{C}_0$, alors $D$ est nulle sur l'espace des fonctions à support dans l'ensemble $G_{ell}$ des éléments elliptiques de $G$.
\end{prop}

\begin{proof}
Dire que $D$ est nulle sur $\mathcal{C}_0$ revient à dire que pour tout sous-groupe parahorique $K$ de $G$, on a:
\[\sum_{i=1}^t\lambda_itr\left(\pi_i^{N_G\left(K\right)}\right)=0.\]
Pour montrer cette proposition, on va donc tout d'abord étendre les résultats de \cite[III.4]{scst} au cas $G$ réductif quasi-déployé quelconque. Fixons un sous-groupe d'Iwahori $I$ de $G$; soit $\mathcal{P}$ l'ensemble de tous les sous-groupes parahoriques de $G$, et $\mathcal{P}_I$ l'ensemble de ceux qui contiennent $I$. Pour toute représentation irréductible tempérée $\pi$ de $G$, posons comme dans \cite[III.4]{scst}:
\[f_{EP}\left(\pi\right)=\sum_{K\in\mathcal{P}_I}\left(-1\right)^{dim\left(A_K\right)}vol\left(N_G\left(K\right)/Z\right)^{-1}\varepsilon_{A_K}tr\left(\pi^{N_G\left(K\right)}\right),\]
où pour tout $K\in\mathcal{P}$, $A_K$ est la facette de l'immeuble de $G$ correspondant à $K$, où $tr\left(\pi_i^{N_G\left(K\right)}\right)$ est considérée comme une fonction sur $N_G\left(K\right)$ étendue à $G$ par zéro, et où $\varepsilon_{A_K}$ est défini comme dans \cite{scst}.

Soit $h\in G_{ell}$. On a les lemmes suivants:

\begin{lemme}
L'ensemble $X^h$ des points de l'immeuble de $G$ fixés par $h$ recouvre un nombre fini de facettes.
\end{lemme}

\begin{proof}
En effet, le lemme est vrai pour $G/Z_{G,d}$ et l'image de $h$ dans $G/Z_{G,d}$ d'après \cite[III.4.9]{scst}; or les immeubles de $G$ et de $G/Z_{G,d}$ sont isomorphes en tant que complexes simpliciaux, et toute facette fixée par $h$ dans l'immeuble de $G$ l'est aussi par son image dans $G/Z_{G,d}$, et réciproquement.
\end{proof}

\begin{lemme}
On a pour tout $\pi$:
\[\int_{G/Z_{G,d}}f_{EP}\left(\pi\right)\left(g^{-1}hg\right)dg=\sum_{K\in\mathcal{P}}\left(-1\right)^{dim\left(A_K\right)}\overline{tr\left(\pi^{N_G\left(K\right)}\right)}\left(h\right).\]
\end{lemme}

\begin{proof}
La démonstration est identique à celle de \cite[III.4.10]{scst}, en remplaçant $G$ par $G/Z_{G,d}$; la somme du membre de droite est finie grâce au lemme précédent.
\end{proof}

\begin{lemme}
Soit $\psi\in C_c^\infty\left(G\right)$ à support dans $G_{ell}$. On a pour tout $\pi$:
\[tr\left(\pi\right)\left(\psi\right)=\int_G\psi\left(h\right)\int_{G/Z_{G,d}}f_{EP}\left(\pi\right)\left(g^{-1}h^{-1}g\right)dgdh.\]
\end{lemme}

\begin{proof}
La démonstration est identique à celle de \cite[III.4.16]{scst}.
\end{proof}

\begin{cor}
Sous les mêmes hypothèses, on a pour tout $\pi$:
\[tr\left(\pi\right)\left(\psi\right)=\int_G\psi\left(h\right)\sum_{K\in\mathcal{P}}\left(-1\right)^{dim\left(A_K\right)}\overline{tr\left(\pi^{N_G\left(K\right)}\right)}\left(h^{-1}\right)dh.\]
\end{cor}

Revenons à la démonstration de la proposition. Pour tout $\psi\in C_c^\infty\left(G\right)$ à support dans $G_{ell}$, on a:
\begin{equation*}\begin{split}\left(\psi\right)&=\sum_i\lambda_itr\left(\pi_i\right)\left(\psi\right)\\&=\int_G\psi\left(h\right)\left(\sum_{K\in\mathcal{P}}\left(-1\right)^{dim\left(A_K\right)}\sum_i\lambda_i\overline{tr\left(\pi_i^{N_G\left(K\right)}\right)}\right)\left(h^{-1}\right)dh.\end{split}\end{equation*}
Or l'ensemble des $K\in\mathcal{P}$ intervenant de façon non triviale dans le membre de droite est fini: en effet, pour tout $h$ elliptique, il existe un voisinage de $h$ dans $G$ dont tout élément fixe les mêmes facettes que $h$ (par exemple $hK'$, où $K'$ est un sous-groupe ouvert compact de $G$ contenu dans les fixateurs de toutes les facettes de $X^h$ et de toutes les facettes adjacentes à $X^h$), et le support de $\psi$, étant compact, admet un recouvrement fini constitué de tels voisinages. On peut donc écrire:
\[D\left(\psi\right)=\sum_{K\in\mathcal{P}}\left(-1\right)^{dim\left(A_K\right)}\int_G\psi\left(h\right)\left(\sum_i\lambda_i\overline{tr\left(\pi_i^{N_G\left(K\right)}\right)}\right)\left(h^{-1}\right)dh.\]
Supposons d'abord $N_G\left(K\right)/K^1$ fini; on a alors, pour tout $i$ et tout $h$:
\[\overline{tr\left(\pi_i^{N_G\left(K\right)}\right)}\left(h^{-1}\right)=tr\left(\pi_i^{N_G\left(K\right)}\right)\left(h\right),\]
d'où l'on déduit:
\[D\left(\psi\right)=\sum_{K\in\mathcal{P}}\left(-1\right)^{dim\left(A_K\right)}\int_G\psi\left(h\right)\left(\sum_i\lambda_itr\left(\pi_i^{N_G\left(K\right)}\right)\left(h\right)\right)dh=0\]
grâce à ce qui précède, et la proposition est démontrée dans ce cas.

Supposons maintenant $N_G\left(K\right)/K^1$ quelconque, et considérons le groupe $\gamma^{\mth{Z}}K$. On a vu précédemment qu'il existe un entier $b>0$ tel que $\gamma^b$ est dans le centre de $G$; pour tout $i$, il existe donc un caractère $\chi_i$ de ${\mth{Z}}$ tel que pour tout entier $z$, $\pi_i\left(\gamma^{bz}\right)$ est une homothétie de rapport $\chi_i\left(bz\right)$; de plus, $\chi_i$ est unitaire puisque $\pi_i$ l'est. Considérons donc la représentation:
\[\chi_i^{-1}\otimes\pi_i^{\gamma^{\mth{Z}}K}:\gamma^{\mth{Z}}K/K^1\longmapsto V\left(\pi_i\right)^{K^1},\]
$\chi_i$ étant étendu à un caractère sur $\gamma^{\mth{Z}}K$ par $\chi_i|_K=1$. Cette représentation est triviale sur $\gamma^{b{\mth{Z}}}$, donc induit une représentation du groupe fini $\gamma^{b{\mth{Z}}}K/\gamma^{b{\mth{Z}}}K^1$ canoniquement isomorphe à $K/K^1$. On a donc, pour tout $h\in\gamma^{\mth{Z}}K$:
\[\overline{tr\left(\chi_i^{-1}\otimes\pi_i^{\gamma^{\mth{Z}}K}\right)}\left(h^{-1}\right)=tr\left(\chi_i^{-1}\otimes\pi_i^{\gamma^{\mth{Z}}K}\right)\left(h\right),\]
d'où, puisque $\chi_i$ est unitaire:
\[\overline{tr\left(\pi_i^{N_G\left(K\right)}\right)}\left(h^{-1}\right)=tr\left(\pi_i^{N_G\left(K\right)}\right)\left(h\right).\]
On conclut comme dans le cas où $N_G\left(K\right)/K^1$ est fini.
\end{proof}

Revenons à la démonstration du théorème. La distribution $D$ est intégrable, invariante et nulle sur $G_{ell}$; d'après \cite[corollaire 3 au théorème 3.1]{kzh}, elle est somme d'induites de distributions invariantes $D_M$ sur un ensemble de re\-pré\-sen\-tants des classes de conjugaison de Levi propres de $G$; de plus, fixons un sous-groupe parahorique maximal spécial $K_0$ de $G$ contenant $I$; on a le lemme suivant:

\begin{lemme}
Soit $f\in C_c^\infty\left(G\right)$ à support dans l'ensemble des éléments compacts (resp. non compacts) de $G$, et soit $P=MU$ un sous-groupe parabolique de $G$. Alors $f^P=f^P_{K_0}$ est à support dans l'ensemble des éléments compacts (resp. non compacts) de $G$ contenus dans $M$.
\end{lemme}

\begin{proof}
Si $P=G$, le résultat est trivial; supposons donc $P$ propre. Supposons que $f$ est à support dans l'ensemble des éléments compacts de $G$, et soit $l\in M$ non compact dans $G$; on a:
\[f^P\left(l\right)=\int_{K_0\times U}f\left(k^{-1}luk\right)dk.\]
Supposons qu'il existe $k\in K_0$ et $u\in U$ tels que $k^{-1}luk$ est compact dans $G$; alors $lu$ aussi est compact dans $G$. Soit $K_U$ un sous-groupe compact de $U$ contenant $u$; on déduit de \cite[1.4.3]{cass} qu'il existe un élément $z$ contenu dans le centre de $M$ tel que la suite des $z^{-n}K_Uz^n$, $n>0$, est décroissante et que l'intersection des $z^{-n}K_Uz^n$ est réduite à $\left\{1\right\}$. La suite des $z^{-n}luz^n=lz^{-n}uz^n$ converge donc vers $l$; or, tous ces éléments sont compacts dans $G$, donc par le lemme \ref{parferm}, $l$ aussi, ce qui est exclu. Donc pour tous $k,u$, $k^{-1}luk$ n'est pas compact dans $G$, donc $f\left(k^{-1}luk\right)=0$; on a donc $f^P\left(l\right)=0$. La démonstration dans le cas où $f$ est à support dans les éléments non compacts de $G$ est similaire.
\end{proof}

Soit $1_c$ la fonction caractéristique des éléments compacts de $G$. On a:
\[D=D\left(1_c\cdot\right)=\sum_M\left(\Ind_M^GD_M\right)\left(1_c\cdot\right)=\sum_M\Ind_M^G\left(D_M\left(1_c\cdot\right)\right)\]
grâce au lemme précédent; on peut donc supposer que pour tout $M$, on a $D_M=D_M\left(1_c\cdot\right)$, c'est-à-dire que $D_M$ est à support dans $M\cap G_c\subset M_c$.

Soit, pour tout $M$, $\tilde{D}_M$ la distribution sur $M$ définie pour $D_M$ de manière analogue à $\tilde{D}$ pour $D$; posons:
\[D^\#=\sum_M\Ind_M^G\tilde{D}_M.\]
Par hypothèse de récurrence, les restrictions de $\tilde{D}_M$ et de $D_M$ à $\mathcal{C}_M$ sont égales. De plus, si $K$ est un sous-groupe parahorique de $M$ et $\gamma$ un élément de $N_M\left(K\right)$ tel que $\gamma K\not\subset G_c\cap M$, puisque $D_M$ est nulle sur $\gamma K$, $\tilde{D}_M$ aussi; on en déduit que $\tilde{D}_M$ est à support dans $G_c\cap M$; comme ceci est vrai pour tout $M$, $D^\#\in\tilde{\mathcal{D}_c}\left(G\right)$ d'après le lemme précédent. Enfin, on a le lemme suivant:

\begin{lemme}\label{tc}
Soit $f\in\mathcal{C}$, et $P=MU$ un sous-groupe parabolique de $G$. Alors le terme constant $f^P$ de $f$ est somme d'un élément de $\mathcal{C}_M$ et d'un élément de $C_0\left(M\right)$.
\end{lemme}

\begin{proof}
Il suffit de montrer ce lemme pour $f$ appartenant à une partie génératrice de $\mathcal{C}$; on supposera donc qu'il existe $g\in G$ et $x\in\overline{A}$ tels que $f$ est la fonction caractéristique de $K_x^1gK_x^1$. On peut supposer $P$ et $M$ standard relativement à $P_0$ et $T_0$, où $P_0$ est le sous-groupe parabolique minimal de $G$ d'image $I$ dans $K_0/K_0^1$; on a alors, pour tout $l\in M$:
\begin{equation*}\begin{split}f^P\left(l\right)&=\delta_P^G\left(l\right)^{\dfrac 12}\int_{K_0\times U}f\left(k^{-1}luk\right)dudk\\&=vol\left(K_x^1\right)\delta_P^G\left(l\right)^{\dfrac 12}\sum_{k\in K_0/K_x^1}vol\left(K_x^1gK_x^1\cap k^{-1}lUk\right)\\&=vol\left(K_x^1\right)\delta_P^G\left(l\right)^{\dfrac 12}\sum_{k\in K_0/K_x^1}vol\left(kK_x^1gK_x^1k^{-1}\cap lU\right).\end{split}\end{equation*}
Il suffit de vérifier que pour tout $k\in K_0$, l'application:
\[l\in M\longmapsto \sum_{k'\in\left(K_0\cap M\right)/\left(K_x^1\cap M\right)}vol\left(k'kK_x^1gK_x^1k^{-1}k'{}^{-1}\cap lU\right)\]
est biinvariante par le premier sous-groupe de congruence d'un sous-groupe parahorique de $M$. Elle l'est déjà par le sous-groupe ouvert compact $kK_x^1k^{-1}\cap M$ de $M$, qui ne dépend que de la classe de $k$ à droite modulo $K_x$; de plus, elle-même ne dépend que de la classe de $k$ à gauche modulo $K_0\cap P$. Il suffit donc de montrer l'assertion pour $k$ appartenant à un système de représentants $X$ de $\left(K_0\cap P\right)\backslash K_0/\left(K_0\cap K_x\right)$; puisque $K_x$ contient $K_0^1$ et puisque les images de $K_0\cap P$ et de $K_0\cap K_x$ dans $K_0/K_0^1$ sont des sous-groupes paraboliques standard de $K_0/K_0^1$ (relativement au parabolique minimal $I/K_0^1$ et au tore déployé maximal $T_0/K_{T_0}$), on peut supposer que $X$ est également un système de représentants de $\left(W'\cap\left(K_0\cap P\right)\right)\backslash\left(W'\cap K_0\right)/\left(W'\cap\left(K_0\cap K_x\right)\right)$.

Soit donc $w\in W'\cap K_0$. D'après le corollaire \ref{iw2}, $wK_xw^{-1}\cap M$ est un sous-groupe parahorique de $M$; $wK_x^1w^{-1}\cap M$ est alors son premier sous-groupe de congruence, ce qui démontre l'assertion cherchée.
\end{proof}

On a donc pour tout $M$, tout $P=MU$ et pour tout $f\in\mathcal{C}$, puisque $f^P\in\mathcal{C}_M+C_0\left(M\right)$ d'après le lemme \ref{tc} et grâce à l'hypothèse de récurrence:
\[\Ind_M^G\tilde{D}_M\left(f\right)=\tilde{D}_M\left(f^P\right)=D_M\left(f^P\right)=\Ind_M^GD_M\left(f\right),\]
d'où l'on déduit que $D^\#\left(f\right)=D\left(f\right)$. En particulier, $D$ et $D^\#$ coïncident sur $\mathcal{C}_0$; comme $D^\#\in\tilde{\mathcal{D}_c}\left(G\right)$, la proposition \ref{dprime} impose $D^\#=\tilde{D}$. Comme par hypothèse $\tilde{D}=0$, $D$ est nulle sur $\mathcal{C}$ et la deuxième assertion du théorème est montrée pour $D$ combinaison linéaire de traces compactes.

Supposons maintenant $D$ quelconque.
Pour tout $\gamma\in\Gamma$, soit $\mathcal{C}_\gamma$ (resp. $\mathcal{C}_{0,\gamma}$) le sous-espace des éléments de $\mathcal{C}$ (resp. $\mathcal{C}_0$) à support dans $\gamma G^1$; on va montrer que si $D$ est nulle sur $\mathcal{C}_0$, elle est nulle sur $\mathcal{C}_\gamma$. Comme $\mathcal{C}$ est somme directe des $\mathcal{C}_\gamma$, la seconde assertion s'en déduit.

Fixons donc $\gamma\in\Gamma$.
On va montrer qu'il existe une distribution $D'\in\mathcal{D}_c$, combinaison linéaire de traces compactes, qui coïncide avec $D$ sur $\mathcal{C}_\gamma$; la deuxième assertion du théorème se déduit alors du cas précédent. Si $D=0$, c'est évident; supposons donc $D\neq 0$, et soit $\mathcal{D}_{tc}$ le sous-espace des éléments de $\mathcal{D}_c|_{\mathcal{C}_\gamma}$ qui sont combinaisons linéaires de traces compactes; ce sous-espace est de dimension finie car, d'après ce qui précède, un élément de $\mathcal{D}_{tc}$ est entièrement déterminé par ses valeurs sur les fonctions à support dans un $\gamma K$, avec $K\supset I$ et $\gamma$ normalisant $K$, biinvariantes par $K^1$, et l'espace de ces fonctions est de dimension finie.
Soit $\left(D_1,\dots,D_t\right)$ une base de $\mathcal{D}_{tc}$ et $e_1,\dots,e_t$ des éléments de $\mathcal{C}_\gamma$ vérifiant $D_i\left(e_j\right)=\delta_{ij}$ pour tous $i,j$; un élément de $\mathcal{D}_c|_{\mathcal{C}_\gamma}$ qui est combinaison linéaire de traces compactes est entièrement défini par ses valeurs sur les $e_i$. D'autre part, soit $\mathcal{C}_\gamma^0$ le sous-espace des éléments de $\mathcal{C}_\gamma$ annulés par toute distribution trace compacte de représentation tempérée; les éléments de $\mathcal{C}_\gamma^0$ sont alors annulés par toutes les distributions traces compactes, donc par tous les éléments de $\mathcal{D}_c$ d'après \cite[théorème 0]{kzh} (si la caractéristique de $F$ est nulle) et \cite[théorème B]{kz2} (si la caractéristique de $F$ est positive; la démonstration de Kazhdan marche également pour les groupes quasi-déployés) et en particulier par $D$. Posons donc:
\[D'=\sum_{i=1}^tD\left(e_i\right)D_i.\]
On a $D'=D=0$ sur $\mathcal{C}_\gamma^0$, et pour tout $i$, $D'\left(e_i\right)=D\left(e_i\right)$; or $\mathcal{C}_\gamma^0$ et les $e_i$ engendrent $\mathcal{C}_\gamma$, donc $D|_{\mathcal{C}_\gamma}=D'|_{\mathcal{C}_\gamma}$, ce qui achève la démonstration de l'assertion cherchée.

Montrons maintenant la première assertion du théorème. Soit $D\in\mathcal{D}_c$ quelconque; d'après la proposition \ref{dprime}, la distribution $D-\tilde{D}$ est nulle sur $\mathcal{C}_0$. Comme la distribution nulle est un élément de $\tilde{\mathcal{D}_c}$ qui coïncide avec $D-\tilde{D}$ sur $\mathcal{C}_0$, en appliquant à nouveau la proposition \ref{dprime}, on en déduit que $\widetilde{D-\tilde{D}}=0$; d'après la deuxième assertion du théorème, $D-\tilde{D}$ est alors nulle sur $\mathcal{C}$ et le théorème est démontré. 
\end{proof}
\end{proof}

\subsection{Restrictions des distributions à $\mathcal{H}$}

Conservons le sous-groupe d'Iwahori $I$ de $G$ fixé précedemment. Soit $\mathcal{H}=\mathcal{H}_I$ l'algèbre de Hecke des fonctions sur $G$ à support compact et biinvariantes par $I$, et $\mathcal{H}_0=\mathcal{H}_{0,I}$ le sous-espace des éléments de $\mathcal{H}_I$ à support dans la réunion des normalisateurs des sous-groupes parahoriques de $G$ contenant $I$; $\mathcal{H}$ (resp. $\mathcal{H}_0$) est clairement un sous-espace de $\mathcal{C}$ (resp. $\mathcal{C}_0$).
On va montrer un résultat similaire au théorème précédent pour $\mathcal{H}$ et $\mathcal{H}_0$:

\begin{theo}\label{dh}
Soit $D\in\mathcal{D}_c$ nulle sur $\mathcal{H}_0$. Alors $D$ est nulle sur $\mathcal{H}$.
\end{theo}

\begin{proof}
Montrons d'abord la proposition suivante: soit $\left(\pi,V\right)$ une re\-pré\-sen\-ta\-tion irréductible tempérée de $G$, $I$ un sous-groupe d'Iwahori de $G$, et $V^I=V_\pi^I$ le sous-espace des éléments de $V$ invariants par $I$. On a:

\begin{prop}
Si $V^I$ est réduit à $0$, la distribution $tr_c\left(\pi\right)$ est nulle sur $\mathcal{H}$.
\end{prop}

\begin{proof}
En effet, d'après \cite[corollaire à la proposition 2.1]{cl2}, on a, pour tout $f\in C_c^\infty\left(G\right)$:
\[tr_c\left(\pi\right)\left(f\right)=\sum_{P=MU}\left(-1\right)^{rg\left(G\right)-rg\left(M\right)}tr\left(\left(\delta_P^G\right)^{-\dfrac 12}\pi_P\right)\left(\hat{\chi}_Uf^P\right),\]
où la somme porte sur l'ensemble des sous-groupes paraboliques standard de $G$, et où pour tout $P=MU$, $\hat{\chi}_U$ est une fonction sur $M$ constante sur les classes modulo $M^1$, et en particulier à la fois invariante par conjugaison et biinvariante par tout sous-groupe d'Iwahori de $M$.

D'autre part, soit $K_0$ le sous-groupe parabolique maximal spécial de $G$ intervenant dans la définition du terme constant;
on supposera qu'il contient $I$.

Pour tout sous-groupe de Levi $M$ semi-standard de $G$, on notera $\mathcal{H}_M$ l'algèbre de Hecke des fonctions sur $M$ biinvariantes par $I_M=I\cap M$; le lemme \ref{iwah} montre que ce sous-groupe est bien un Iwahori de $M$.

On a les lemmes suivants:

\begin{lemme}
Soit $f\in\mathcal{H}$, $P=MU$ un sous-groupe parabolique standard de $G$. Alors $f^P\in\mathcal{H}_M+C_0\left(M\right)$.
\end{lemme}

\begin{proof}
En effet, puisque les sous-groupes d'Iwahori de $M$ sont tous conjugués entre eux et par le même raisonnement que dans le lemme \ref{tc}, il suffit de montrer que pour tout $k\in K_0$, l'application:
\[l\in M\longmapsto \sum_{k'\in\left(K_0\cap kMk^{-1}\right)/\left(K_x^1\cap kMk^{-1}\right)}vol\left(k'kIgIk^{-1}k'{}^{-1}\cap lU\right)\]
est biinvariante par un sous-groupe d'Iwahori de $M$. Or elle l'est par $kIk^{-1}$, et elle ne dépend que de la classe de $k$ à gauche modulo $K_0\cap P$; on peut donc supposer que $k\in N_G\left(T_0\right)$, et le lemme \ref{iwah} permet alors de conclure.
\end{proof}

\begin{cor}
Si $f\in\mathcal{H}$, on a $\hat{\chi}_Uf^P\in\mathcal{H}_M+C_0\left(M\right)$.
\end{cor}

\begin{proof}
En effet, écrivons grâce au lemme précédent $f^P=f'+f''$, $f'\in\mathcal{H}_M$, $f''\in C_0\left(M\right)$. Puisque $\hat{\chi}_U$ est biinvariante par tout sous-groupe d'Iwahori de $M$, on a $\hat{\chi}_Uf'\in\mathcal{H}_M$; d'autre part, si $\gamma\in M/M^1$, les restrictions de $f''$ et de $\hat{\chi}_Uf''$ à $\gamma M^1$ sont proportionnelles, donc si $D$ est une distribution invariante sur $M$ à support dans $\gamma M^1$, puisque $D\left(f''\right)=0$, $D\left(\hat{\chi}_Uf''\right)=0$. Si maintenant $D$ est quelconque, on peut écrire de manière unique:
\[D=\sum_{\gamma\in M/M^1}D_\gamma,\]
où pour tout $\gamma$, $D_\gamma$ est à support dans $\gamma M^1$; $M/M^1$ est en général infini mais pour tout $\phi\in C_c^\infty\left(M\right)$, tous les $D_\gamma\left(\phi\right)$ sont nuls sauf au plus un nombre fini d'entre eux, donc la somme converge. On a alors:
\[D\left(\hat{\chi}_Uf''\right)=\sum_{\gamma\in M/M^1}D_\gamma\left(\hat{\chi}_Uf''\right)=0,\]
ce qui démontre le corollaire.
\end{proof}

\begin{lemme}
Soit $P=MU$ un sous-groupe parabolique standard de $G$, $I_M=I\cap M$, $\left(\pi,V\right)$ une représentation irréductible admissible de $G$. Si $V^I$ est réduit à zéro, alors $\pi_P$ n'admet aucun vecteur non nul invariant par $I_M$.
\end{lemme}

\begin{proof}
Ce lemme est une conséquence immédiate de \cite[théorème 3.3.3]{cass}.
\end{proof}

Revenons à la démonstration de la proposition. Pour tout $P=MU$, on vérifie facilement que si $f'\in\mathcal{H}_M$, et si $\sigma$ est une représentation irréductible admissible de $M$ et $V_\sigma$ son espace, $\sigma\left(f'\right)$ est un élément de $V_\sigma^{I_M}$; si cet espace est réduit à $0$, on a donc $\sigma\left(f'\right)=0$. Puisque $V_\pi^I=0$, on déduit de ce qui précède que si $f\in\mathcal{H}_I$, tous les $tr\left(\left(\delta_P^G\right)^{-\dfrac 12}{\pi_i}_P\right)\left(\hat{\chi_U}f^P\right)$ sont nuls, ce qui implique $tr_c\left(\pi\right)\left(f\right)=0$. La distribution $tr_c\left(\pi\right)$ est donc nulle sur $\mathcal{H}$ et la proposition est démontrée.
\end{proof}

Montrons maintenant le théorème. Pour tout $\gamma\in\Gamma$, soit $\mathcal{H}_\gamma$ l'espace des éléments de $\mathcal{H}$ à support dans $\gamma G^1$; on va montrer que si $D$ est nulle sur $\mathcal{H}_0$, elle est nulle sur $\mathcal{H}_\gamma$; comme $\mathcal{H}$ est somme directe des $\mathcal{H}_\gamma$, le théorème s'en déduit immédiatement.

Puisque $\mathcal{H}_\gamma\subset\mathcal{C}_\gamma$, d'après la démonstration du théorème \ref{dc}, il existe un élément de $\mathcal{D}_c$ qui est combinaison linéaire de traces compactes et qui coïncide avec $D$ sur $\mathcal{H}_\gamma$. On pourra donc, ici encore, supposer que $D$ est combinaison linéaire de traces compactes. Ecrivons:
\[D=\sum_{i=1}^t\lambda_i tr_c\left(\pi_i\right),\]
où les $\lambda_i$ sont des constantes et les $\pi_i$ des représentations irréductibles tempérées de $G$.
Pour tout $i$, soit $V_i$ l'espace de $\pi_i$ et $V_i^I$ le sous-espace des vecteurs de $V_i$ invariants par $I$; d'après la proposition précédente, on peut supposer que tous les $V_i^I$ sont non triviaux.
De plus, soit, pour tout $i$, $\chi_{\pi_i}$ le caractère central de $\pi_i$; puisque $\pi_i$ admet des vecteurs invariants par $I$, $\chi_{\pi_i}$ est trivial sur $Z_G\cap I=Z_G\cap G^1$, et on peut donc l'étendre à un caractère de $G$ trivial sur $G^1$; la distribution:
\[D'=\sum_{i=1}^t\lambda_i\chi_{\pi_i}\left(\gamma\right)tr_c\left(\chi_{\pi_i}^{-1}\pi_i\right),\]
qui est combinaison linéaire de traces compactes de représentations irréductibles tempérées de $G$ de caractère central trivial, coïncide avec $D$ sur $\mathcal{H}_\gamma$; on pourra donc supposer $D=D'$.
On va montrer que dans ce cas, $D$ est nulle sur $\mathcal{C}_0$; d'après le théorème \ref{dc}, $D$ est alors nulle sur $\mathcal{C}$ tout entier, et en particulier sur $\mathcal{H}_\gamma$.

Soit, pour tout $x\in\overline{A}$, $K_x$ le fixateur connexe de $x$ dans $G$; c'est un sous-groupe parahorique de $G$ qui contient $I$. Posons ${\mth{G}}_x=K_x/K_x^1$, et pour tout $\gamma\in\Gamma_x$, ${\mth{G}}^+_{x,\gamma}=\gamma^{\mth{Z}}K_x/K_x^1$; pour tous $i,x,\gamma$, considérons la représentation $\pi_i^{\gamma^{\mth{Z}}K_x}$ de ${\mth{G}}^+_{x,\gamma}$, dont l'espace est le sous-espace $V_i^{K_x^1}$ des vecteurs de $V_i$ invariants par $K_x^1$.

\begin{lemme}
Toute composante irréductible de $\pi_i^{\gamma^{\mth{Z}}K_x}$ admet des vecteurs non nuls invariants par le sous-groupe parabolique minimal ${\mth{P}}_{x,0}={\mth{M}}_{x,0}{\mth{N}}_{x,0}$ image de $I$ dans ${\mth{G}}_x$.
 \end{lemme}

\begin{proof}
Puisque $\pi_i^{K_x}$ est, pour tout $\gamma$ normalisant $K_x$, la restriction à ${\mth{G}}_x$ de $\pi_i^{\gamma^{\mth{Z}}K_x}$, il suffit de montrer ce lemme lorsque $\gamma=1$.
Soit $\rho$ une sous-représentation irréductible de $\pi_i^{K_x}$; alors d'une part il existe un sous-groupe parabolique ${\mth{P}}={\mth{M}}{\mth{U}}$ de ${\mth{G}}_x$ et une représentation cuspidale $\sigma$ de ${\mth{M}}$ telle que $\sigma$ est contenue dans le module de Jacquet $\rho_{\mth{P}}$, et d'autre part, d'après le théorème 5.2 de \cite{mp} (dont la démonstration de l'assertion d'unicité est valable en caractéristique quelconque), si $\rho'$ est une autre sous-représentation irréductible de $\pi_i^{K_x}$ et si ${\mth{P}}'={\mth{M}}'{\mth{U}}'$ et $\sigma'$ sont définis de la même façon pour $\rho'$, alors $\left({\mth{M}},\sigma\right)$ et $\left({\mth{M}}',\sigma'\right)$ sont conjugués (dans $G$). En particulier, si ${\mth{M}}={\mth{M}}_{x,0}$ et $\sigma$ est la représentation triviale, $\rho_{\mth{P}}$ est somme directe de copies de cette représentation triviale, donc $\rho$ admet des vecteurs invariants par ${\mth{P}}_{x,0}$, et $\rho'$ admet des vecteurs non nuls invariants par le sous-groupe parabolique minimal ${\mth{P}}'$ de ${\mth{G}}_x$; or tous les sous-groupes paraboliques minimaux de ${\mth{G}}_x$ sont conjugués, donc $\rho'$ admet aussi des vecteurs non nuls invariants par ${\mth{P}}_{x,0}$.

Il reste donc à montrer que ${\mth{M}}={\mth{M}}_{x,0}$ et $\sigma=1$ conviennent pour au moins une sous-représentation irréductible de $\pi_i^{K_x}$. Soit $v$ un élément non nul de $V_i^I\subset V_i^{K_x^1}$, si $V_i^{K_x^1}=W_1\oplus\dots\oplus W_t$ est une décomposition de $V_i^{K_x^1}$ en sous-$K_x$-modules irréductibles, écrivons:
\[v=\sum_{j=1}^tv_j,\]
où pour tout $j$, $v_j\in W_j$. Puisque $\pi_i^{K_x}$ stabilise chacun des $W_j$, on a, pour tout $h\in I$ et tout $j$, $\pi_i\left(h\right)v_j=v_j$, donc $v_j\in V_i^I$; puisque $v$ est non nul, au moins un des $v_j$ est non nul, donc la sous-représentation $\rho$ de $\pi_i^{K_x}$ d'espace $W_j$ admet des vecteurs non nuls invariants par ${\mth{P}}_{x,0}$. Mais alors le module de Jacquet $\rho_{{\mth{P}}_{x_0}}$ admet des vecteurs non nuls invariants par ${\mth{M}}_{x,0}$, ce qui achève la démonstration du lemme.
\end{proof}

D'après ce lemme, il existe des représentations irréductibles $\rho_1,\dots,\rho_s$ de ${\mth{G}}^+_{x,\gamma}$, appartenant toutes à la série principale, et des constantes $\mu_1,\dots,\mu_s$ telles que pour tout $f\in\mathcal{C}_{x,0}$ à support dans $\gamma^{\mth{Z}}K_x$, on a:
\[D\left(f\right)=\sum_{j=0}^s\mu_jtr\left(\rho_j\right)\left(f\right).\]
Supposons d'abord $\gamma$ d'ordre fini dans ${\mth{G}}^+_{x,\gamma}$. Alors ${\mth{G}}^+_{x,\gamma}$ est fini, donc d'après \cite[proposition 11.25]{cr}, puisque $D$ est nulle sur $\mathcal{H}_0$, tous les $\mu_j$ sont nuls, donc $D$ est nulle sur l'ensemble des éléments de $\mathcal{C}_{x,0}$ à support dans $\gamma^{\mth{Z}}K_x$. Supposons maintenant $\gamma$ quelconque; on va en fait se ramener au cas d'un groupe fini.
Il existe un entier $b$ tel que $\gamma^b$ possède des représentants dans $Z_G$; le sous-groupe $\gamma^{b{\mth{Z}}}$ de ${\mth{G}}_x$ est alors contenu dans le centre de ${\mth{G}}^+_{x,\gamma}$, donc pour tout $j$, il existe un caractère $\chi_{\rho_j}$ de ${\mth{Z}}$ tel que pour tout $z\in{\mth{Z}}$, $\rho_j\left(\gamma^{bz}\right)$ est une homothétie de rapport $\chi_{\rho_j}\left(z\right)$. De plus, puisque tous les $\chi_{\pi_i}$ sont triviaux, tous les $\chi_{\rho_j}$ aussi; on peut alors passer au quotient ${\mth{G}}^+_{x,\gamma}/\gamma^{b{\mth{Z}}}$, qui est fini, et on est ramené au cas précédent.
On en déduit que $D$ est nulle sur l'ensemble des éléments de $\mathcal{C}_{x,0}$ à support dans $\gamma^{\mth{Z}}K_x$; comme ceci est vrai pour tout $x\in\overline{A}$ et tout $\gamma\in\Gamma$ normalisant $K_x$, $D$ est nulle sur $\mathcal{C}_0$, donc sur $\mathcal{H}_\gamma$ et le théorème est démontré.
\end{proof}

\section{Intégrales orbitales semi-simples}

\subsection{Intégrales orbitales semi-simples sur $G$}

On va maintenant s'intéresser aux distributions intégrales orbitales semi-simples sur $G$. Pour tout $g\in G$ semi-simple régulier et tout $f\in C_c^\infty\left(G\right)$, posons:
\[J\left(g,f\right)=\int_{h\in Z_G^0\left(g\right)\backslash G}f\left(h^{-1}gh\right)dg.\]
Cette intégrale converge et la distribution $J\left(g,\cdot\right)$ est invariante; de plus, si $g$ est compact, elle est à support dans l'ensemble des éléments compacts de $G$. On peut donc s'intéresser au sous-espace de $\mathcal{D}_c$ engendré par les $J\left(g,\cdot\right)$, avec $g$ compact.

Dans ce qui suit, on s'intéressera uniquement aux $J\left(g,\cdot\right)$, $g\in G^1$. Soit $\mathcal{D}_{c,1}$ le sous-espace des éléments de $\mathcal{D}_c$ à support dans $G^1$; on va décrire, moyennant quelques conditions sur $G$, un système de générateurs de $\mathcal{D}_{c,1}|_\mathcal{H}$ constitué d'in\-té\-gra\-les orbitales semi-simples.

Soit $T$ un tore de $G$. C'est le groupe des $F$-points d'un tore $\underline{T}$ de $\underline{G}$ défini sur $F$, et il existe une extension finie $E$ de $F$ telle que $\underline{T}$ est déployé sur $E$. On dit que $T$ est {\em non ramifié} si $E$ peut être choisie non ramifiée; $T$ sera dit {\em non ramifié maximal} si $T$ est maximal parmi les tores non ramifiés de $G$.

Supposons maintenant que $T$ est un tore non ramifié maximal de $G$. Soit $K_T$ l'unique sous-groupe parahorique de $T$, $K_T^1$ son premier sous-groupe de congruence; ${\mth{T}}=K_T/K_T^1$ est alors le groupe des ${\mth{F}}_q$-points d'un tore de même dimension $d$ que $T$ et défini sur ${\mth{F}}_q$; si $K$ est un sous-groupe parahorique de $G$ contenant $K_T$, en posant ${\mth{G}}=K/K^1$, ${\mth{T}}$ s'identifie à un tore maximal de ${\mth{G}}$.

Soit $g$ un élément de $K_T$. Si son image $\overline{g}$ dans ${\mth{T}}$ est un élément régulier de ${\mth{T}}$ dans ${\mth{G}}$, alors il est clair que $g$ est un élément régulier de $T$ dans $G$. La réciproque est fausse, et il existe même des cas où ${\mth{T}}$ ne contient aucun élément régulier dans ${\mth{G}}$ (par exemple si $G=GL_n\left(F\right)$, ${\mth{G}}=GL_n\left({\mth{F}}_q\right)$, avec $n\geq q$, et $T$ est le tore diagonal de $G$).
On dira que $g$ est {\em non ramifié de réduction régulière} si pour tout sous-groupe parahorique $K$ de $G$ contenant $K_T$, l'image de $g$ dans $K/K^1$ est un élément régulier de $K/K^1$.

On va supposer que $G$ vérifie la condition suivante:

{\em (C1): Tout tore maximal non ramifié $T$ de $G$ admet au moins un élément non ramifié de réduction régulière.}

Cette condition est vraie pour $q$ suffisamment grand.
Plus précisément, con\-si\-dé\-rons le groupe $G/Z_{G,d}$; on a une décomposition de la forme:
\[G/Z_{G,d}=G_0G_1G_2\ldots G_t,\]
où $G_0$ est un groupe fini et $G_1,\dots,G_t$ sont les groupes des $F$-points de groupes réductifs simples et connexes $\underline{G_1},\dots,\underline{G_t}$. Soit $F^{nr}$ l'extension non ramifiée maximale de $F$; pour tout $i$, soit $G_i^{nr}=\underline{G_i}\left(F^{nr}\right)$, $h_i$ le nombre de Coxeter de $G_i^{nr}$, $T_i$ un tore déployé maximal de $G_i^{nr}$, $X_i=X_*\left(T_i\right)$ et $P_i$ le groupe des poids du système des coracines de $G_i^{nr}$ relativement à $T_i$. On notera, pour tout $i$, $\mathcal{T}_i$ le couple composé du diagramme de Dynkin de $G_i^{nr}$ et de $P_i$, et $\mathcal{T}$ l'ensemble des $\mathcal{T}_i$.
On a le résultat suivant:

\begin{lemme}\label{srg}
Il existe un entier $q_\mathcal{T}$, ne dépendant que de $\mathcal{T}$, tel que si $q\geq q_\mathcal{T}$, tout tore non ramifié maximal de $G$ admet des éléments non ramifiés de réduction régulière.
\end{lemme}

\begin{proof}
Quitte à remplacer $G$ par $G/Z_{G,d}$ et à le décomposer en produit de groupes simples, on peut le supposer simple; on notera alors $h=h_1$, $l=l_1$, $X=X_1$, $P=P_1$. Soit $T=\underline{T}\left(F\right)$ un tore non ramifié maximal de $G$. Supposons d'abord $T$ déployé; soit $\Phi=\Phi\left(G/T\right)$ le système de racines réduit de $G$ relativement à $T$, qui est aussi le système de racines réduit de $G_{nr}$ relativement au tore déployé maximal $T_{nr}=\underline{T}\left(F_{nr}\right)$; si $K$ est un sous-groupe parahorique maximal spécial de $G$ contenant l'unique parahorique $K_T$ de $T$, $\phi$ est également le système de racines du groupe réductif $K$ relativement à son tore maximal $K_T/K_T^1$ (un tel système est toujours réduit car le groupe est déployé). On va montrer que pour $q$ assez grand, il existe $\xi\in X_*\left(T\right)$ tel que pour tout $\alpha\in\Phi$, $\langle \alpha,\xi\rangle $ n'est pas un multiple de $q-1$; si $x$ est un élément de $\mathcal{O}^*$ dont la réduction modulo ${\mathfrak{p}}$ est un générateur de ${\mth{F}}_q^*$, $\xi\left(x\right)$ est alors un élément non ramifié de réduction régulière de $T$.

Soit $\Delta=\left\{\alpha_1,\dots,\alpha_n\right\}$ un ensemble de racines simples de $\Phi$; considérons l'application $\phi$ de $X_*\left(T\right)$ dans ${\mth{N}}^n$ définie par:
\[\phi:\xi\longmapsto\left(\langle \alpha_1,\xi\rangle ,\dots,\langle \alpha_n,\xi\rangle \right);\]
cette application est injective et son image est un sous-groupe d'indice $\left[P:X\right]$ de ${\mth{Z}}^n$.
Soit $\alpha_M$ la plus grande racine de $\Phi$; écrivons-la:
\[\alpha_M=\sum_ic_i\alpha_i.\]
Pour tout $\xi\in X_*\left(T\right)$ d'image $\left(a_1,\dots,a_n\right)$ par $\phi$, posons:
\[q_\xi=\sum_{i=1}^nc_ia_i+1,\]
et soit $q_\mathcal{T}$ le minimum des $q_\xi$, avec $\xi$ tel que pour tout $i$, $a_i\geq 1$. Soit $\xi_\mathcal{T}$ un élément de $X_*\left(T\right)$ vérifiant cette condition et tel que $q_{\xi_\mathcal{T}}=q_\mathcal{T}$; 
si $q-1\geq q_\mathcal{T}$, pour toute racine positive $\alpha$, on a $\langle \alpha,\xi_\mathcal{T}\rangle \in\left\{1,\dots,q-2\right\}$, qui ne contient aucun multiple de $q-1$ par hypothèse, et l'assertion du lemme est démontrée.

Supposons maintenant $T=\underline{T}\left(F\right)$ quelconque, et toujours $q\geq q_\mathcal{T}$; soit $E$ une extension non ramifiée de $F$ telle que $\underline{T}\left(E\right)$ est déployé, $k=\left[E:F\right]$, et $\sigma$ le générateur de $\Gal\left(E/F\right)$ dont la réduction sur les corps résiduels est le Frobenius $x\mapsto x^q$. Soit $\Phi$ le système de racines de $\underline{G}\left(E\right)$ relativement à $\underline{T}\left(E\right)$, et $\Delta=\left\{\alpha_1,\dots,\alpha_n\right\}$ un ensemble de racines simples de $\Phi$, et soit $\xi_\mathcal{T}\in X_*\left(\underline{T}\left(E\right)\right)$ construit de la même façon que dans le cas déployé; considérons l'élément suivant:
\[\xi'=\sum_{j=0}^{k-1}\sigma^j\left(\xi_\mathcal{T}\right).\]
C'est une application de $E^*$ dans $\underline{T}\left(E\right)$, stable par $\Gal\left(E/F\right)$ donc à valeurs dans $T$; d'autre part, puisque pour tous $i,j$, on a:
\[\langle \alpha_i,\sigma^j\xi_\mathcal{T}\rangle =q^j\langle \sigma^{-j}\alpha_i,\xi_\mathcal{T}\rangle ,\]
on en déduit que l'on a, pour tout $j$ et toute $\alpha\in\Phi$ positive:
\[0<\langle \alpha,\sigma^j\xi_\mathcal{T}\rangle <q^j\left(q-1\right),\]
d'où:
\[0<\langle \alpha,\xi'\rangle <\sum_{j=0}^{k-1}q^j\left(q-1\right)=q^{j^k}-1.\]
Si $x$ est un élément de l'anneau des entiers $\mathcal{O}_E$ de $E$ dont la réduction modulo l'idéal maximal ${\mathfrak{p}}_E$ est un générateur de ${\mth{F}}_{q^k}^*$, l'élément $\xi'\left(x\right)$ de $T$ est alors non ramifié de réduction régulière dans $\underline{G}\left(E\right)$, donc a fortiori dans $G$ et le lemme est démontré.
\end{proof}

On calculera explicitement les valeurs de $q_\mathcal{T}$ en fonction de $\mathcal{T}$ dans le paragraphe \ref{cqt}.

On va maintenant montrer le théorème:

\begin{theo}\label{intorb}
\begin{itemize}
\item Si $g$ est un élément non ramifié de réduction ré\-gu\-liè\-re d'un tore non ramifié maximal $T$ de $G$, la restriction à $\mathcal{H}$ de la distribution $J\left(g,.\right)$ ne dépend que de la classe de conjugaison de $T$ dans $G$.
\item Supposons que $G$ vérifie (C1). Soit $R$ un système de représentants des classes de conjugaison de tores non ramifiés maximaux de $G$, et pour tout $T\in R$, soit $g_T$ un élément non ramifié de réduction régulière de $T$. Si le diagramme de Dynkin de $G$ ne contient aucune composante connexe de type $E_7$ ou $E_8$, alors les distributions $J\left(g_T,\cdot\right)|_\mathcal{H}$ engendrent $\mathcal{D}_{c,1}|_\mathcal{H}$, et en constituent une base si $G$ est déployé sur $F$.
\end{itemize}
\end{theo}

\begin{proof}
D'après le théorème \ref{dh}, on peut remplacer $\mathcal{H}$ par $\mathcal{H}_0$ dans l'énoncé ci-dessus; de plus, pour tout élément non ramifié de réduction régulière $g$ d'un tore non ramifié maximal $T$ de $G$, on a $g\in G^1$, donc la distribution $J\left(g,\cdot\right)$ est un élément de $\mathcal{D}_{c,1}$, dont la restriction à $\mathcal{H}_0$ ne dépend que de sa valeur sur le sous-espace $\mathcal{H}'_0$ des éléments de $\mathcal{H}_0$ à support dans $G^1$.

Soit $A$ la facette de $\mathcal{B}$ fixée par $I$; $\mathcal{H}'_0$ est somme, pour $x$ décrivant $\overline{A}$, des sous-espaces $\mathcal{H}'_{0,x}$ des éléments de $\mathcal{H}'_0$ à support dans le fixateur connexe $K_x$ de $x$. Si, comme dans la démonstration du théorème \ref{dh}, on note ${\mth{G}}_x=K_x/K_x^1$ et ${\mth{P}}_{0,x}=I/K_x^1$, alors ${\mth{P}}_{0,x}$ est un sous-groupe parabolique minimal de ${\mth{G}}_x$ et $\mathcal{H}'_{0,x}$ s'identifie à l'espace des fonctions sur ${\mth{G}}_x$ biinvariantes par ${\mth{P}}_{0,x}$; de plus, si $T$ est un tore non ramifié maximal de $G$ tel que $K_T\subset K_x$, $g$ un élément non ramifié de réduction régulière de $T$
et $G_g$ un système de représentants de $G$ modulo $Z_G^0\left(g\right)$ à gauche et $K_x$ à droite; on a, pour tout $f\in\mathcal{H}'_{0,x}$:
\[J\left(g,f\right)=\sum_{h\in G_g}vol\left(Z_G^0\left(g\right)\backslash Z_G^0\left(g\right)hK_x\right)\int_{K_x}f\left(k^{-1}h^{-1}ghk\right)dk.\]
Les seuls termes non nuls de la somme ci-dessus sont ceux correspondant aux $h$ tels que $h^{-1}gh\in K_x$; de plus, on a le lemme suivant:

\begin{lemme}\label{cpt}
Soit $h\in G$ tel que $h^{-1}gh\in K_x$; alors $h^{-1}K_Th\subset K_x$.
\end{lemme}

\begin{proof}
En effet, supposons d'abord $T$ déployé; $T$ est alors un tore déployé maximal de $G$, et puisque $K_T\subset K_x$, l'appartement $\mathcal{A}'$ de $\mathcal{B}$ associé à $T$ contient $x$. D'après \cite[lemme 2.5.8]{bt}, il existe $h'\in G^1$ tel que $h'\mathcal{A}=\mathcal{A}'$ et que $h'$ fixe $\mathcal{A}\cap\mathcal{A}'$ point par point; on a donc en particulier $h'x=x$. On en déduit que $h'\in K_x$ et que $h'{}^{-1}Th'=T_0$; on peut donc supposer $T=T_0$. De plus, puisque $g$ est régulier dans $T_0$ et fixe $hx$, d'après \cite[3.6.1]{t}, $hx\in\mathcal{A}$; on en déduit que $h^{-1}K_{T_0}h\subset K_x$, ce qui démontre l'assertion du lemme.

Supposons maintenant $T$ quelconque, et soit $E$ une extension finie non ramifiée de $F$ telle que $T_E=\overline{T}\left(E\right)$ est déployé. $T_E$ est un tore non ramifié maximal déployé, donc un tore déployé maximal, de $G_E=\overline{G}\left(E\right)$; de plus, si $\mathcal{B}_E$ est l'immeuble de $G_E$, $\mathcal{B}$ s'identifie canoniquement à l'ensemble des points de $\mathcal{B}_E$ fixés par $\Gal\left(E/F\right)$. Si $g$ est régulier dans $T$, il l'est aussi dans $T_E$, donc d'après ce qui précède, $h^{-1}K_{T_E}h$ est contenu dans le fixateur connexe $K_{x,E}$ de $x$ dans $G_E$. En considérant les intersections avec $G$ de ces deux groupes, on obtient $h^{-1}K_Th\subset K_x$ et le lemme est démontré.
\end{proof}

\begin{cor}
$\overline{h^{-1}gh}$ est un élément semi-simple régulier de ${\mth{G}}_x$.
\end{cor}

\begin{proof}
En effet, $h^{-1}K_Th\subset K_x$ d'après le lemme précédent, et $\overline{h^{-1}gh}$ est un élément régulier du tore maximal $\overline{h^{-1}K_Th}$ de ${\mth{G}}_x$.
\end{proof}

D'après ce corollaire, la restriction de $J\left(g,\cdot\right)$ à $\mathcal{H}'_{0,x}$ s'identifie à une somme de distributions intégrales orbitales d'éléments semi-simples réguliers de ${\mth{G}}_x$. On va donc considérer les intégrales orbitales d'éléments semi-simples réguliers sur des groupes réductifs finis.

On va ainsi montrer un résultat analogue à celui du théorème, mais concernant les groupes finis.
Soit $\underline{\mth{G}}$ un groupe réductif connexe défini sur un corps fini ${\mth{F}}_q$, et soit $\mathbf{F}$ une application de Frobenius de $\underline{\mth{G}}$ dans lui-même, c'est-à-dire un morphisme bijectif de $\underline{\mth{G}}$ dans lui-même tel qu'il existe des entiers $k,l>0$ tels que ${\mathbf F}^k$ est le morphisme induit par l'application $x\mapsto x^{q^l}$ de $\overline{{\mth{F}}_q}$ dans lui-même. Soit ${\mth{G}}=\underline{\mth{G}}^{\mathbf F}$ l'ensemble des points de $\underline{\mth{G}}$ stables par $\mathbf{F}$; ${\mth{G}}$ est un groupe réductif connexe et fini.

Soit $L$ l'application de Lang de $\underline{\mth{G}}$ dans lui-même définie par:
\[L:g\mapsto g^{-1}{\mathbf F}\left(g\right).\]
D'après \cite[corollaire au théorème 1]{lg}, pour tout sous-groupe fermé connexe $\underline{\mth{H}}$ de $\underline{\mth{G}}$, $L\left(\underline{\mth{H}}\right)=\underline{\mth{H}}$; en particulier, $L$ est surjective.

Soit $\underline{\mth{B}}=\underline{\mth{T}}\underline{\mth{U}}$ un sous-groupe de Borel de $\underline{\mth{G}}$, avec $\underline{\mth{T}}$ est stable par $\mathbf{F}$ ($\underline{\mth{B}}$ lui-même pouvant ne pas l'être). Posons $\tilde{X}=L^{-1}\left(\underline{\mth{U}}\right)$; c'est une variété algébrique affine, stable par multiplication à gauche par ${\mth{G}}$ et par multiplication à droite par ${\mth{T}}=\underline{\mth{T}}^{\mathbf F}$: en effet, pour $x\in\tilde{X}$, $g\in{\mth{G}}$ et $t\in{\mth{T}}$, on a:
\[L\left(gxt\right)=t^{-1}x^{-1}g^{-1}g{\mathbf F}\left(x\right)t=t^{-1}L\left(x\right)t\in\underline{\mth{U}}.\]
Le groupe ${\mth{G}}\times{\mth{T}}$ peut donc être assimilé à un groupe fini d'automorphismes de $\tilde{X}$; pour tout $\left(g,t\right)\in{\mth{G}}\times{\mth{T}}$, on peut alors définir le nombre de Lefschetz $\mathcal{L}\left(\left(g,t\right),\tilde{X}\right)$ comme dans \cite[7.1.4]{car}.

Pour tout caractère $\theta$ de ${\mth{T}}$, soit $R_{\underline{\mth{T}},\theta}$ le caractère de Deligne-Lusztig de ${\mth{G}}$ dans ${\mth{C}}$ associé à $\underline{\mth{T}}$, défini par:
\[R_{\underline{\mth{T}},\theta}:g\longmapsto\dfrac 1{\card\left({\mth{T}}\right)}\sum_{t\in{\mth{T}}}\theta\left(t\right)^{-1}\mathcal{L}\left(\left(g,t\right),\tilde{X}\right);\]
d'après \cite[proposition 7.3.6]{car}, cette application ne dépend que de $\underline{\mth{T}}$ et pas de $\underline{\mth{B}}$.

Soit $\mathcal{H}_{\mth{G}}$ l'algèbre de Hecke des fonctions sur ${\mth{G}}$ biinvariantes par un sous-groupe parabolique minimal ${\mth{P}}_0$ de ${\mth{G}}$ fixé; on a:

\begin{prop}\label{intorbf}
\begin{itemize}
\item Si $g$ est un élément régulier d'un tore maximal ${\mth{T}}$ de ${\mth{G}}$, la restriction à $\mathcal{H}_{\mth{G}}$ de la distribution $1_g$ ne dépend que de la classe de conjugaison de ${\mth{T}}$ dans ${\mth{G}}$.
\item Soit $R_{\mth{G}}$ un système de représentants des classes de conjugaison de tore maximaux de ${\mth{G}}$; supposons que pour tout ${\mth{T}}\in R_{\mth{G}}$, il existe un élément régulier $g_{\mth{T}}$ dans ${\mth{T}}$. Si le diagramme de Dynkin de ${\mth{G}}$ ne contient aucune composante connexe de type $E_7$ ou $E_8$, alors les distributions $1_{g_{\mth{T}}}|_{\mathcal{H}_{\mth{G}}}$ engendrent l'espace des distributions invariantes sur $\mathcal{H}_{\mth{G}}$, et en constituent une base si ${\mth{G}}$ est déployé.
\end{itemize}
\end{prop}

\begin{proof}
La démonstration de cette proposition fera l'objet de la section suivante.
\end{proof}

Revenons à la démonstration du théorème. Puisque pour $g\in G$ semi-simple régulier, $J\left(g,\cdot\right)$ ne dépend que de la classe de conjugaison de $g$, il suffit, pour montrer la premiére assertion, de vérifier que si $T$ est un tore non ramifié maximal de $G$ et si $g,g'$ sont deux éléments non ramifiés de réduction régulière de $T$, alors $J\left(g,\cdot\right)$ et $J\left(g',\cdot\right)$ coïncident sur $\mathcal{H}'_0$. Puisque le sous-groupe parahorique $K_T$ de $T$ est un sous-groupe compact de $G^1$, il est contenu dans un sous-groupe parahorique maximal de $G$, donc il fixe au moins un point de $\mathcal{B}$; soit donc $A_T$ une facette de $\mathcal{B}$ de dimension maximale parmi celles fixées par $K_T$, et $K$ le sous-groupe parahorique correspondant. Quitte à conjuguer $T$, on peut supposer que $K$ contient $I$. Pour tout $x\in\overline{A_T}$, si $K_x$ est le fixateur connexe de $x$, on a $K_T\subset K_x$; $K_T/K_T^1$ est alors un tore maximal de ${\mth{G}}_x=K_x/K_x^1$, et en appliquant la première assertion de la proposition \ref{intorbf} à ${\mth{G}}_x$, $\overline{g}$ et $\overline{g'}$, en posant, pour tout $f\in C_c^\infty\left(G\right)$:
\[J_x\left(g,f\right)=J_{K_x}\left(g,f\right)=\int_{K_x}f\left(k^{-1}gk\right)dk\]
et de même pour $g'$, les distributions $J_x\left(g,\cdot\right)$ et $J_x\left(g',\cdot\right)$ coïncident sur l'espace $\mathcal{H}'_{0,x}$ des fonctions sur $K_x$ biinvariantes par $I$. De plus, si $f$ est un élément de $\mathcal{H}'_0$ dont le support est disjoint de $K_x$, on a clairement $J_x\left(g,f\right)=J_x\left(g',f\right)=0$; $J_x\left(g,\cdot\right)$ et $J_x\left(g',\cdot\right)$ coïncident donc sur $\mathcal{H}'_0$. 

Soit maintenant $A$ la chambre de $\mathcal{B}$ fixée par $I$, $x\in\overline{A}$ et $G_g$ un système de représentants de $G$ modulo $Z_G^0\left(g\right)$ à gauche et $K_x$ à droite; on a, pour tout $f\in\mathcal{H}'_{0,x}$:
\begin{equation}\label{decj}\begin{split}J\left(g,f\right)&=\sum_{h\in G_g}vol\left(Z_G^0\left(g\right)\backslash Z_G^0\left(g\right)hK_x\right)\int_{K_x}f\left(k^{-1}h^{-1}ghk\right)dk\\&=\sum_{h\in G_g,h^{-1}gh\in K_x}vol\left(Z_G^0\left(g\right)\backslash Z_G^0\left(g\right)hK_x\right)\int_{K_x}f\left(k^{-1}h^{-1}ghk\right)dk\\&=\sum_{h\in G_g,x\in\overline{A_{h^{-1}Th}}}vol\left(Z_G^0\left(g\right)\backslash Z_G^0\left(g\right)hK_x\right)J_x\left(h^{-1}gh,f\right),\end{split}
\end{equation}
où $A_{h^{-1}Th}$ est la face de $A$ de dimension maximale parmi celles fixées par $K_{h^{-1}Th}$ s'il en existe, et $A_{h^{-1}Th}=\emptyset$ sinon; d'après le lemme \ref{cpt}, $x\in\overline{A_{h^{-1}Th}}$ si et seulement si $h^{-1}gh\in K_x$.
On a une décomposition similaire pour $J\left(g',f\right)$; de plus, on a le lemme suivant:

\begin{lemme}\label{eqggp}
On a $Z_G^0\left(g\right)=Z_G^0\left(g'\right)$, et pour tout $h\in G_g$, $h^{-1}gh\in K_x$ si et seulement si $h^{-1}g'h\in K_x$.
\end{lemme}

\begin{proof}
Puisque $g$ et $g'$ sont réguliers dans $T$, on a:
\[Z_G^0\left(g\right)=Z_G^0\left(T\right)=Z_G^0\left(g'\right).\]
De plus, soit $h\in G_g$ tel que $h^{-1}gh\in K_x$; d'après le lemme \ref{cpt}, $h^{-1}K_Th\subset K_x$, donc $h^{-1}g'h\in K_x$. Puisque d'après la première assertion, on peut supposer $G_{g'}=G_g$, l'autre implication est symétrique.
\end{proof}

D'après ce lemme et ce qui précède, on a:
\[\sum_{h\in G_g,x\in\overline{A_{h^{-1}Th}}}vol\left(Z_G^0\left(g\right)\backslash Z_G^0\left(g\right)hK_x\right)J_x\left(h^{-1}gh,f\right)\]
\[=\sum_{h\in G_{g'},x\in\overline{A_{h^{-1}Th}}}vol\left(Z_G^0\left(g'\right)\backslash Z_G^0\left(g'\right)hK_x\right)J_x\left(h^{-1}g'h,f\right),\]
autrement dit $J\left(g,f\right)=J\left(g',f\right)$, pour tout $f\in\mathcal{H}'_{0,x}$. Comme ceci est vrai pour tout $x\in\overline{A}$, $J\left(g,\cdot\right)$ et $J\left(g',\cdot\right)$ coïncident sur $\mathcal{H}'_{0,I}$ et la première assertion du théorème est démontrée.

Montrons maintenant la seconde. Supposons donc que $G$ vérifie (C1), et que son diagramme de Dynkin ne contient aucune composante de type $E_7$ ou $E_8$; soit $R$ un système de représentants des classes de conjugaison de tores non ramifiés maximaux de $G$, et pour tout $T\in R$, soit $g_T$ un élément non ramifié de réduction régulière de $T$, $K_T$ l'unique sous-groupe parahorique de $T$ et ${\mth{T}}=K_T/K_T^1$. Soit $D\in\mathcal{D}_{c,1}$; on va montrer que la restriction à $\mathcal{H}'_0$ de $D$ est combinaison linéaire de celles des $J\left(g_T,\cdot\right)$.

Soit $K$ un sous-groupe parahorique de $G$ contenant $I$, et ${\mth{G}}=K/K^1$. D'après la proposition \ref{intorbf}, si $R_K$ est un ensemble de représentants des classes de conjugaison de tores maximaux de ${\mth{G}}$, il existe des constantes $\lambda_{\mth{T}}$, ${\mth{T}}\in R_K$, telles que pour toute fonction $f$ sur $G$ à support dans $K$ et biinvariante par $I$, on a:
\[D\left(f\right)=\sum_{{\mth{T}}\in R_K}\lambda_{\mth{T}}J_K\left(g_{\mth{T}},f\right),\]
où pour tout ${\mth{T}}$, $g_{\mth{T}}$ est un élément régulier de ${\mth{T}}$; de plus, si $K'$ est un sous-groupe parahorique de $G$ contenant $K$, et si ${\mth{G}}'=K'/K'{}^1$, ${\mth{G}}$ s'identifie à un sous-groupe de Levi du sous-groupe parabolique ${\mth{P}}$ image de $K$ dans ${\mth{G}}'$, et on a le lemme suivant:

\begin{lemme}
Pour tout ${\mth{T}}'\in R_{K'}$ et pour toute fonction $f$ sur $G$ à support dans $K$ et biinvariante par $I$, on a:
\[J_{K'}\left(g_{{\mth{T}}'},f\right)=\card\left(N_{{\mth{G}}'}\left({\mth{T}}'\right)\right)\sum_{\mth{T}}\dfrac{1}{\card\left(N_{\mth{G}}\left({\mth{T}}\right)\right)}J_K\left(g_{\mth{T}},f\right),\]
la somme portant sur l'ensemble des ${\mth{T}}\in R_K$ conjugués à ${\mth{T}}'$ dans ${\mth{G}}'$.
\end{lemme}

\begin{proof}
En effet, soit $R_{K,{\mth{T}}'}$ cet ensemble. Puisque $g_{{\mth{T}}'}$ est régulier dans ${\mth{T}}'$, il l'est aussi dans ${\mth{T}}$, donc pour tout $h\in{\mth{G}}'$, $h^{-1}g_{{\mth{T}}'}h$ appartient à au plus un élément de $R_{K,{\mth{T}}'}$; on en déduit:
\begin{equation*}\begin{split}J_{K'}\left(g_{{\mth{T}}'},f\right)&=\sum_{h\in{\mth{G}}'}f\left(h^{-1}g_{{\mth{T}}'}h\right)\\&=\sum_{{\mth{T}}\in R_{K,{\mth{T}}'}}\sum_{h\in{\mth{G}}',h^{-1}g_{{\mth{T}}'}h\in{\mth{T}}^{\mth{G}}}f\left(h^{-1}g_{{\mth{T}}'}h\right)\\&=\dfrac 1{\card\left({\mth{G}}\right)}\sum_{{\mth{T}}\in R_{K,{\mth{T}}'}}\sum_{h\in{\mth{G}}',h^{-1}g_{{\mth{T}}'}h\in{\mth{T}}^{\mth{G}}}J_K\left(h^{-1}g_{{\mth{T}}'}h,f\right)\\&=\dfrac 1{\card\left({\mth{G}}\right)}\sum_{{\mth{T}}\in R_{K,{\mth{T}}'}}\sum_{h\in{\mth{G}}',h^{-1}g_{{\mth{T}}'}h\in{\mth{T}}}\dfrac{\card\left({\mth{G}}\right)}{\card\left(N_{\mth{G}}\left({\mth{T}}\right)\right)}J_K\left(g_{\mth{T}},f\right)\\&=\card\left(N_{{\mth{G}}'}\left({\mth{T}}'\right)\right)\sum_{{\mth{T}}\in R_{K,{\mth{T}}'}}\dfrac{1}{\card\left(N_{\mth{G}}\left({\mth{T}}\right)\right)}J_K\left(g_{\mth{T}},f\right),\end{split}\end{equation*}
ce qui démontre le lemme.
\end{proof}

D'après ce lemme, on a pour tous $K,K',f$:
\begin{equation*}\begin{split}D\left(f\right)&=\sum_{{\mth{T}}\in R_K}\lambda_{\mth{T}}J_K\left(g_{\mth{T}},f\right)\\&=\sum_{{\mth{T}}'\in R_{K'}}\lambda_{{\mth{T}}'}\card\left(N_{{\mth{G}}'}\left({\mth{T}}'\right)\right)\sum_{{\mth{T}}}\dfrac{1}{\card\left(N_{\mth{G}}\left({\mth{T}}\right)\right)}J_K\left(g_{\mth{T}},f\right).\end{split}\end{equation*}
On en déduit que l'on peut supposer, pour tous $K,K'$ et tout ${\mth{T}}\in R_K$, si ${\mth{T}}'$ est l'unique élément de $R_{K'}$ tel que ${\mth{T}}\in R_{K,{\mth{T}}'}$:
\[\lambda_{\mth{T}}=\dfrac{\card\left(N_{{\mth{G}}'}\left({\mth{T}}'\right)\right)}{\card\left(N_{\mth{G}}\left({\mth{T}}\right)\right)}\lambda_{{\mth{T}}'}.\]
Pour tout tore non ramifié maximal $T$ de $G$, tel qu'il existe un sous-groupe parahorique $K$ de $G$ contenant $I$ et $K_T$, si ${\mth{T}}$ est l'image de $T$ dans ${\mth{G}}=K/K^1$, posons:
\[\mu_T=\card\left(N_{\mth{G}}\left({\mth{T}}\right)\right)\lambda_{\mth{T}};\]
on déduit de ce qui précède que $\mu_T$ ne dépend que de $T$ et pas du choix de $K$.

Soit $d$ le plus petit entier vérifiant la propriété suivante: il existe une face $B$ de $A$ de dimension $d$ et un tore non ramifié maximal $T$ de $G$ tels que $K_T\subset K_B$, où $K_B$ est le fixateur connexe de $B$, et que $\mu_T\neq 0$; si une telle face n'existe pas, on pose $d=-1$. On va montrer l'assertion cherchée par récurrence sur $d$.

Si $d=-1$, la restriction de $D$ à $\mathcal{H}'_0$ est nulle, donc $D$ est nulle sur $\mathcal{H}$ d'après le théorème \ref{dh} et l'assertion est triviale. Supposons donc $d\geq 0$, et soit $A_1,\dots,A_s$ les faces de $A$ de dimension $d$; pour tout $j$, soit $\mathcal{H}_{A_j}$ l'espace des fonctions sur $K_j=K_{A_j}$ biinvariantes par $I$.
D'après ce qui précède, pour tout $j\in\left\{1,\dots,s\right\}$ et tout $f\in\mathcal{H}_{A_j}$, on a:
\[D\left(f\right)=\sum_{{\mth{T}}\in R_{A_j,ell}}\lambda_{{\mth{T}}}J_{K_j}\left(g_{\mth{T}},f\right),\]
où $R_{A_j,ell}$ est un ensemble de représentants des classes de conjugaison de tores maximaux elliptiques de ${\mth{G}}_j=K_j/K_j^1$.

Si $T$ et $T'$ sont deux tores de $G$ conjugués entre eux, alors $K_T$ et $K_{T'}$ le sont également, de même que $K_T^1$ et $K_{T'}^1$; $G$ agit donc par conjugaison sur l'ensemble des tores finis $K_T/K_T^1$, où $T$ est un tore de $G$. Montrons d'abord le lemme suivant:

\begin{lemme}\label{tcj}
Soit $j,j'\in\left\{1,\dots,s\right\}$, soit ${\mth{T}}\in R_{A_j,ell}$ et ${\mth{T}}'\in R_{A_{j'},ell}$; supposons ${\mth{T}}$ et ${\mth{T}}'$ conjugués par $h\in G$. Alors $K_j$ et $K_{j'}$ sont associés; de plus, si $j=j'$, $h$ normalise ${\mth{G}}_j$.
\end{lemme}

\begin{proof}
Soit $h\in G$ tel que $h^{-1}{\mth{T}}h={\mth{T}}'$; on a alors, si $T$ est un tore de $G$ tel que ${\mth{T}}=K_T/K_T^1$:
\[K_T\subset K_j\cup hK_{j'}h^{-1}.\]
On en déduit que $K_T$ fixe à la fois $A_j$ et $hA_{j'}$.

Soit $\mathcal{A}'$ un appartement de $\mathcal{B}$ contenant $A_j$ et $hA_{j'}$ et $T'_0$ le tore déployé maximal de $G$ associé à $\mathcal{A}'$; posons $\left(M,K_M\right)=\kappa_0\left(K_j,T'_0\right)$, $\left(M',K_{M'}\right)=\kappa_0\left(K_{hA_{j'}},T'_0\right)$. Alors $T$ est un tore elliptique à la fois de $M$ et de $M'$; or puisque $M$ et $M'$ sont tous deux semi-standard relativement à $T'_0$, leur intersection est un sous-groupe de Levi à la fois de $M$ et de $M'$; comme elle contient $T$, on en déduit que $M=M'$. De plus, $K_M$ et $K_{M'}$ sont deux sous-groupes parahoriques maximaux de $M$ contenant $K_T$; puisque ${\mth{T}}$ est elliptique dans $K_M/K_M^1={\mth{G}}_j$, ils sont égaux, ce qui montre d'une part que $K_j$ et $K_{j'}$ sont associés, et d'autre part, si $j=j'$, que $h$ normalise $M$ et $K_M$, et le lemme est démontré.
\end{proof}

Revenons à la démonstration du théorème. Par le même raisonnement que dans le lemme \ref{reges}, si $K_j$ et $K_{j'}$ sont associés et si $\phi$ est un isomorphisme entre $\mathcal{H}_{A_j}$ et $\mathcal{H}_{A_{j'}}$ obtenu en les identifiant tous deux à l'espace des fonctions sur ${\mth{G}}_j={\mth{G}}_{j'}$ biinvariantes par le parabolique minimal image de $I$, on a $D\left(\phi\left(f\right)\right)=D\left(f\right)$ pour tout $f\in\mathcal{H}_{A_j}$. En particulier, pour tout $j$, si $h$ normalise ${\mth{G}}_j$, on a:
\[D\left(f\right)=\sum_{{\mth{T}}\in R_{A_j,ell}}\lambda_{{\mth{T}}}J_{K_j}\left(g_{\mth{T}},f\right)=\sum_{{\mth{T}}\in R_{A_j,ell}}\lambda_{h^{-1}{\mth{T}}h}J_{K_j}\left(g_{\mth{T}},f\right);\]
si ${\mth{T}}$ et ${\mth{T}}'$ sont deux éléments de $R_{A_j,ell}$ conjugués dans $G$, on peut donc supposer $\lambda_{{\mth{T}}}=\lambda_{{\mth{T}}'}$.

Si maintenant $T$ et $T'$ sont deux tores non ramifiés maximaux de $G$ conjugués entre eux et tels qu'il existe $j,j'$ tels que ${\mth{T}}=K_T/K_T^1$ et ${\mth{T}}'=K_{T'}/K_{T'}^1$ sont des tores elliptiques maximaux respectivement de ${\mth{G}}_j$ et ${\mth{G}}_{j'}$, d'après le lemme précédent, $K_j$ et $K_{j'}$ sont associés, et on déduit de ce qui précède que l'on a encore $\lambda_{{\mth{T}}}=\lambda_{{\mth{T}}'}$, d'où, puisque ${\mth{G}}_j$ et ${\mth{G}}_{j'}$ sont isomorphes, $\mu_T=\mu_{T'}$; $\mu_T$ ne dépend donc que de la classe de conjugaison de $T$ dans $G$, et on a, pour tout $j$ et pour tout $f\in\mathcal{H}_{A_j}$:
\[D\left(f\right)=\sum_{T\in R_d}\mu_T\sum_{{\mth{T}}''\in R_{A_j,ell,T}}\dfrac{1}{\card\left(N_{{\mth{G}}_j}\left({\mth{T}}''\right)\right)}J_{K_j}\left(g_{{\mth{T}}''},f\right),\]
où $R_{A_j,ell,T}$ est l'ensemble des éléments de $R_{A_j,ell}$ conjugués à ${\mth{T}}$, et où $R_d$ est l'ensemble des $T\in R$ tel que $R_{A_j,ell,T}$ est non vide pour au moins un $j$.

De plus, pour tout $T\in R$, en remplaçant $D$ par $J\left(g_T,\cdot\right)$ où $g_T$ est un élément non ramifié de réduction régulière de $T$, on obtient, pour tout $f\in\mathcal{H}_{A_j}$, en utilisant \ref{decj}:
\[J\left(g_T,f\right)=\nu_T\sum_{{\mth{T}}''\in R_{A_j,ell,T}}\dfrac{1}{\card\left(N_{{\mth{G}}_j}\left({\mth{T}}''\right)\right)}J_{K_j}\left(g_{{\mth{T}}''},f\right),\]
où $\nu_T$ est une constante non nulle. La distribution:
\[D'=D-\sum_{T\in R_d}\dfrac{\mu_T}{\nu_T}J\left(g_T,\cdot\right)\]
est alors nulle sur tous les $\mathcal{H}_{A_j}$. Par hypothèse de récurrence, $D'$ est combinaison linéaire d'intégrales orbitales de la forme $J\left(g_T,\cdot\right)$, $T\in R$; c'est donc également le cas de $D$ et l'assertion cherchée est démontrée.

Supposons maintenant $G$ déployé; on va montrer que les $J\left(g_T,\cdot\right)$, $T\in R$ sont linéairement indépendantes. Supposons que l'on a une relation de la forme:
\[\sum_{T\in R}c_TJ\left(g_T,\cdot\right)=0,\]
avec au moins un des $c_T$ non nul.
Soit $d$ le plus grand entier tel qu'il existe au moins un $T\in R_d$ tel que $c_T\neq 0$, soit $A_1,\dots,A_s$ les faces de $A$ de dimension $d$ et $K_1,\dots,K_s$ les parahoriques de $G$ correspondants. Alors pour tout $j$, en conservant les notations précédentes, on a, d'après ce qui précède:
\[\sum_{T\in R_d}c_T\nu_T\sum_{{\mth{T}}''\in R_{A_j,ell,T}}\dfrac{1}{\card\left(N_{{\mth{G}}_j}\left({\mth{T}}''\right)\right)}J_{K_j}\left(g_{{\mth{T}}''},\cdot\right)=0,\]
Or puisque $G$ est déployé, ${\mth{G}}_j=K_j/K_j^1$ l'est aussi, donc d'après la proposition \ref{intorbf}, les $J_{K_j}\left(g_{{\mth{T}}''},\cdot\right)$ sont linéairement indépendantes; l'égalité ci-dessus n'est alors possible que si pour tout $T$ tel que $R_{A_j,ell,T}$ est non vide, $c_T=0$. Comme ceci est vrai pour tout $j$, on en déduit que pour tout $T\in R_d$, $c_T=0$, ce qui contredit la définition de $d$; on en conclut que tous les $c_T$ sont nuls et que les $J\left(g_T,\cdot\right)$ sont linéairement indépendantes.
\end{proof}

\subsection{Démonstration de la proposition}

On va maintenant démontrer la proposition \ref{intorbf}. La première assertion se déduit immédiatement du lemme suivant:
\begin{lemme}\label{cddl}
Soit ${\mth{T}}$ un tore maximal de ${\mth{G}}$, $g$ un élément régulier de ${\mth{T}}$. Pour tout $f\in\mathcal{H}_{\mth{G}}$, on a:
\[1_g\left(f\right)=\dfrac 1{\card\left({\mth{T}}\right)}\langle R_{\underline{\mth{T}},1},f\rangle _{\mth{G}},\]
où $\underline{\mth{T}}$ est l'unique tore maximal de $\underline{\mth{G}}$ contenant ${\mth{T}}$.
\end{lemme}

\begin{proof}
En effet, on a, d'après \cite[proposition 7.5.5]{car}:
\[1_g\left(f\right)=\dfrac{\varepsilon_{Z_{\underline{\mth{G}}}^0\left(g\right)}}{\card\left(Z_{\mth{G}}\left(g\right)\right)\card\left(Z_{\underline{\mth{G}}}^0\left(g\right)\cap{\mth{G}}\right)_p}\sum_{\underline{\mth{T}}'\ni g\mid {\mathbf F}\left(\underline{\mth{T}}'\right)=\underline{\mth{T}}'}\varepsilon_{\underline{\mth{T}}'}\sum_\theta\theta\left(g\right)^{-1}\left(R_{\underline{\mth{T}}',\theta},f\right),\]
où pour tout $\underline{\mth{T}}'$, $\varepsilon_{\underline{\mth{T}}'}$ est défini comme dans \cite[6.5]{car} et vaut $\pm 1$, et où pour $a$ entier, $a_p$ est la plus grande puissance de $p$ divisant $a$.
Or puisque $g$ est régulier dans ${\mth{T}}$, $\underline{\mth{T}}$ est à la fois le centralisateur de $g$ dans $\underline{\mth{G}}$ et le seul tore maximal de $\underline{\mth{G}}$ le contenant; de plus, puisque ${\mth{T}}$ est un tore de ${\mth{G}}$, son cardinal n'est pas multiple de $p$. L'égalité ci-dessus se réduit donc à:
\[1_g\left(f\right)=\dfrac 1{\card\left({\mth{T}}\right)}\sum_\theta\theta\left(g\right)^{-1}\left(R_{\underline{\mth{T}},\theta},f\right).\]
Soit $\theta$ un caractère non trivial de ${\mth{T}}$; on va montrer que $\left(R_{\underline{\mth{T}},\theta},f\right)=0$ pour tout $f\in\mathcal{H}_{\mth{G}}$, ce qui montrera le lemme. Pour cela, soit ${\mth{T}}_0$ un tore maximal de ${\mth{G}}$ contenu dans ${\mth{P}}_0$, donc contenant un tore déployé maximal de ${\mth{G}}$, et soit $\underline{\mth{T}}_0$ le tore maximal de $\underline{\mth{G}}$ contenant ${\mth{T}}_0$.
Les couples $\left(\underline{\mth{T}}_0,1\right)$ et $\left(\underline{\mth{T}},\theta\right)$ ne sont pas géométriquement conjugués;
si $\chi$ est un caractère irréductible de ${\mth{G}}$ tel que $\langle R_{\underline{\mth{T}},\theta},\chi\rangle _{\mth{G}}\neq 0$, on a alors $\langle R_{\underline{\mth{T}}_0,1},\chi\rangle _{\mth{G}}=0$ d'après \cite[théorème 7.3.8]{car}, d'où, par \cite[proposition 7.2.4]{car}, $\langle \Ind_{{\mth{P}}_0}^{\mth{G}}1,\chi\rangle _{\mth{G}}=0$. Or d'après \cite[proposition 11.25]{cr}, cela implique $\langle \chi,f\rangle _{\mth{G}}=0$ pour tout $f\in\mathcal{H}_{\mth{G}}$; on en déduit l'assertion cherchée.
\end{proof}

Pour montrer la deuxième assertion de la proposition, il suffit, d'après le lemme \ref{cddl}, de montrer que les restrictions à $\mathcal{H}_{\mth{G}}$ des distributions $\langle R_{\underline{\mth{T}},1},\cdot\rangle _{\mth{G}}$, où ${\mth{T}}$ décrit un système de représentants des classes de conjugaison de tores maximaux de ${\mth{G}}$, sont linéairement indépendantes.

Montrons d'abord que l'on peut supposer $\underline{\mth{G}}$ simple et adjoint. Soit $\underline{\mth{Z}}$ le centre de $\underline{\mth{G}}$, et $\underline{\mth{Z}}^0$ sa composante neutre; en posant ${\mth{Z}}^0=\underline{{\mth{Z}}^0}^{\mathbf F}$, on a ${\mth{Z}}^0\subset{\mth{P}}_0$, donc $\mathcal{H}_{\mth{G}}$ est canoniquement isomorphe à l'algèbre de Hecke de ${\mth{G}}/{\mth{Z}}^0$, et pour tout $x\in{\mth{G}}$, si $\overline{x}$ est l'image de $x$ dans ${\mth{G}}/{\mth{Z}}^0$, on a pour tout $f\in\mathcal{H}_{\mth{G}}$:
\[\sum_{g\in{\mth{G}}/{\mth{Z}}^0}f\left(g^{-1}\overline{x}g\right)=\dfrac 1{\card\left({\mth{Z}}^0\right)}\sum_{g\in{\mth{G}}}f\left(g^{-1}xg\right).\]
De plus, l'application ${\mth{T}}\mapsto{\mth{T}}/{\mth{Z}}^0$ est une bijection canonique entre les tores maximaux de ${\mth{G}}$ et ceux de ${\mth{G}}/{\mth{Z}}^0$, et deux tores maximaux sont conjugués si et seulement si leurs images sont conjuguées. Enfin on a le lemme suivant:

\begin{lemme}
Le sous-groupe des points fixés par $\mathbf{F}$ de $\underline{\mth{G}}/\underline{\mth{Z}}^0$ est canoniquement isomorphe à ${\mth{G}}/{\mth{Z}}^0$.
\end{lemme}

\begin{proof}
Ce lemme est une conséquence immédiate de \cite[1.17]{car}.
\end{proof}

D'après ce lemme et ce qui précède, la proposition est donc vraie pour ${\mth{G}}$ si et seulement si elle est vraie pour ${\mth{G}}/{\mth{Z}}^0$, et on peut supposer que $\underline{\mth{Z}}^0$ est trivial; $\underline{\mth{G}}$ est alors semi-simple.

Montrons maintenant que l'on peut supposer $\underline{\mth{Z}}$ trivial. Par le même raisonnement que précédemment, d'une part la proposition est vraie pour ${\mth{G}}$ si et seulement si elle est vraie pour ${\mth{G}}/{\mth{Z}}$, et d'autre part, on a une injection canonique de ${\mth{G}}/{\mth{Z}}$ dans $\left(\underline{\mth{G}}/\underline{\mth{Z}}\right)^{\mathbf F}$, qui cette fois n'est plus surjective, mais dont l'image est un sous-groupe normal de $\left(\underline{\mth{G}}/\underline{\mth{Z}}\right)^{\mathbf F}$; de plus, l'ensemble des classes de $\left(\underline{\mth{G}}/\underline{\mth{Z}}\right)^{\mathbf F}$ modulo ${\mth{G}}/{\mth{Z}}$ est indexé par $\underline{\mth{Z}}/L\left(\underline{\mth{Z}}\right)$ ($L\left(\underline{\mth{Z}}\right)$ est un groupe car $\underline{\mth{Z}}$ est abélien). Ceci reste vrai en remplaçant $\underline{\mth{G}}$ par n'importe lequel de ses sous-groupes fermés ${\mathbf F}$-stables contenant ${\mth{Z}}$, et en particulier par un tore maximal $\underline{\mth{T}}$; on en déduit en particulier d'une part une bijection canonique entre les tores maximaux de ${\mth{G}}/{\mth{Z}}$ et ceux de $\left(\underline{\mth{G}}/\underline{\mth{Z}}\right)^{\mathbf F}$, et d'autre part, puisque $\left(\underline{\mth{T}}/\underline{\mth{Z}}\right)^{\mathbf F}$ rencontre toutes les classes de $\left(\underline{\mth{G}}/\underline{\mth{Z}}\right)^{\mathbf F}$ modulo ${\mth{G}}/{\mth{Z}}$, que deux tores maximaux de ${\mth{G}}/{\mth{Z}}$ sont conjugués dans $\left(\underline{\mth{G}}/\underline{\mth{Z}}\right)^{\mathbf F}$ si et seulement si ils le sont dans ${\mth{G}}/{\mth{Z}}$.
On a donc une bijection canonique entre les classes de conjugaison de tores maximaux de ${\mth{G}}/{\mth{Z}}$ et celles de $\left(\underline{\mth{G}}/\underline{\mth{Z}}\right)^{\mathbf F}$.

Soit maintenant ${\mth{P}}_0$ un sous-groupe parabolique minimal de $\left(\underline{\mth{G}}/\underline{\mth{Z}}\right)^{\mathbf F}$. C'est le groupe des points fixes par ${\mathbf F}$ d'un sous-groupe de Borel ${\mathbf F}$-stable de $\underline{\mth{G}}/\underline{\mth{Z}}$, donc d'après ce qui précède, ${\mth{P}}_0\cap\left({\mth{G}}/{\mth{Z}}\right)$ est un sous-groupe normal d'indice $\left[\underline{\mth{Z}}:L\left(\underline{\mth{Z}}\right)\right]$ de ${\mth{P}}_0$;
on en déduit que les algèbres de Hecke $\mathcal{H}\left(\left(\underline{\mth{G}}/\underline{\mth{Z}}\right)^{\mathbf F}\right)$ et $\mathcal{H}\left({\mth{G}}/{\mth{Z}}\right)$, cons\-ti\-tu\-ées des fonctions biinvariantes respectivement par ${\mth{P}}_0$ et ${\mth{P}}_0\cap\left({\mth{G}}/{\mth{Z}}\right)$, sont canoniquement isomorphes. D'autre part, les espaces des restrictions à ces deux algèbres des distributions invariantes sur les groupes respectifs sont é\-ga\-le\-ment canoniquement isomorphes: en effet, considérons la re\-pré\-sen\-ta\-tion $\pi=\Ind_{{\mth{P}}_0/{\mth{Z}}}^{{\mth{G}}/{\mth{Z}}}1$. D'après \cite[11.25]{cr}, l'espace des restrictions à $\mathcal{H}\left({\mth{G}}/{\mth{Z}}\right)$ des distributions invariantes sur ${\mth{G}}/{\mth{Z}}$ est engendré par les traces des sous-représentations de $\pi$; on a également une assertion similaire pour $\left(\underline{\mth{G}}/\underline{\mth{Z}}\right)^{\mathbf F}$ et $\pi'$ définie de façon similaire à $\pi$. D'autre part, d'après \cite[1.17]{car} appliqué à $\underline{\mth{G}}$ et $\underline{{\mth{P}}_0}$, les éléments de l'espace $V$ de $\pi$ sont les restrictions à ${\mth{G}}/{\mth{Z}}$ des éléments de l'espace $V'$ de $\pi'$, et l'application restriction est bijective;
de plus, si $W'$ est un sous-espace irréductible de $V'$ et $W$ son image dans $V$, $W$ est stable par $\pi$, et somme de sous-espaces irréductibles sur lesquels ${\mth{G}}/{\mth{Z}}$ agit par des représentations toutes conjuguées entre elles par des éléments de $\left(\underline{\mth{T}}/\underline{\mth{Z}}\right)^{\mathbf F}$; puisque $\mathcal{H}\left({\mth{G}}/{\mth{Z}}\right)$ est invariante par conjugaison par ce groupe, toutes ces représentations ont même trace sur $\mathcal{H}\left({\mth{G}}/{\mth{Z}}\right)$, ce qui démontre l'assertion cherchée.

Enfin, si $x$ est un élément semi-simple régulier de ${\mth{G}}/{\mth{Z}}$ et si $f\in\mathcal{H}\left(\left(\underline{\mth{G}}/\underline{\mth{Z}}\right)^{\mathbf F}\right)$, on a:
\[\sum_{g\in{\mth{G}}/{\mth{Z}}}f\left(g^{-1}xg\right)=\dfrac 1{\left[\underline{\mth{Z}}:L\left(\underline{\mth{Z}}\right)\right]}\sum_{g\in\left(\underline{\mth{G}}/\underline{\mth{Z}}\right)^{\mathbf F}}f\left(g^{-1}xg\right);\]
la proposition est donc vraie pour ${\mth{G}}/{\mth{Z}}$, et donc pour ${\mth{G}}$, si et seulement si elle est vraie pour $\left(\underline{\mth{G}}/\underline{\mth{Z}}\right)^{\mathbf F}$, et donc on peut, quitte à remplacer $\underline{\mth{G}}$ par $\underline{\mth{G}}/\underline{\mth{Z}}$, supposer $\underline{\mth{Z}}$ trivial, c'est-à-dire $\underline{\mth{G}}$ adjoint. Enfin, on peut supposer $\underline{\mth{G}}$ simple par le même raisonnement que dans \cite[12.1]{car}.

Soit $W$ le groupe de Weyl de $\underline{\mth{G}}$ relativement à $\underline{\mth{T}}_0$. Puisque $\underline{\mth{T}}_0$ est stable par $\mathbf{F}$, $\mathbf{F}$ agit sur $W$; on dit que deux éléments $w,w'$ de $W$ sont $\mathbf{F}$-conjugués s'il existe $x\in W$ tel que $w'=x^{-1}w{\mathbf F}\left(x\right)$. La $\mathbf{F}$-conjugaison est clairement une relation d'équivalence sur $W$; de plus, d'après \cite[proposition 3.3.3]{car}, il existe une bijection canonique entre les classes de conjugaison de tores maximaux de ${\mth{G}}$ et les classes de $\mathbf{F}$-conjugaison de $W$. Soit $W'$ un système de représentants de ces classes; pour tout $w\in W'$, on fixera un tore maximal $\underline{\mth{T}}_w$ de $\underline{\mth{G}}$ tel que ${\mth{T}}_w$ appartient à la classe correspondant à $w$.

Si ${\mth{G}}$ est déployé, l'action de $\mathbf{F}$ sur $W$ est triviale, et les classes de $\mathbf{F}$-conjugaison de $W$ sont les classes de conjugaison ordinaires. Pour tout caractère irréductible $\phi$ de $W$, posons alors:
\[R_\phi=\dfrac 1{\card\left(W\right)}\sum_{w\in W'}\phi\left(w\right)R_{\underline{\mth{T}}_w,1}.\]
Si ${\mth{G}}$ n'est pas déployé, soit $W_1$ le groupe fini d'endomorphismes de ${\mth{T}}_0$ engendré par $W$ et $\mathbf{F}$, et posons pour tout caractère irréductible $\phi$ de $W_1$:
\[R_\phi=\dfrac 1{\card\left(W\right)}\sum_{w\in W'}\phi\left({\mathbf F}w\right)R_{\underline{\mth{T}}_w,1}.\]
Puisque l'on a également pour tout $w\in W'$:
\[R_{\underline{\mth{T}}_w,1}=\sum_\phi\phi\left(w\right)R_\phi\]
si ${\mth{G}}$ est déployé, et:
\[R_{\underline{\mth{T}}_w,1}=\sum_\phi\phi\left({\mathbf F}w\right)R_\phi\]
sinon, l'espace engendré par les $R_\phi$ est le même que celui engendré par les $R_{\underline{\mth{T}}_w,1}$; on va donc montrer que les restrictions à $\mathcal{H}_{\mth{G}}$ des $\langle R_\phi,\cdot\rangle _{\mth{G}}$ engendrent l'espace des distributions sur $\mathcal{H}_{\mth{G}}$, et en constituent une base si ${\mth{G}}$ est déployé.

Pour cela, considérons les caractères irréductibles unipotents de ${\mth{G}}$, c'est-à-dire ceux qui interviennent comme composantes irréductibles des $R_{\underline{\mth{T}}_w,1}$. Lusztig (\cite{l3}) les classe en familles de la manière suivante: les familles sont les plus petites parties de l'ensemble des caractères unipotents de $G$ telles que si $\chi$ et $\chi '$ sont deux caractères intervenant comme composantes de $R_\phi$, pour $\phi$ donné, ils sont dans une même famille, ou encore : deux caractères unipotents $\chi$ et $\chi '$ de ${\mth{G}}$ sont dans une même famille si et seulement si il existe une suite $\chi=\chi_0,\chi_1,\dots,\chi_r=\chi '$ de caractères unipotents et des caractères irréductibles $\phi_1,\dots,\phi_r$ de $W$ tels que pour tout $i$, $\chi_{i-1}$ et $\chi_i$ interviennent tous deux dans $R_{\phi_i}$.

Supposons d'abord ${\mth{G}}$ déployé. D'après \cite[propositions 10.1.2 et 10.11.2]{car}, il existe une bijection $\phi\mapsto\chi_\phi$ entre les caractères irréductibles de $W$ et ceux de ${\mth{G}}$ intervenant dans la série principale. Considérons donc la matrice indexée par les caractères irréductibles de $W$ suivante:
\[M_{\mth{G}}=\left(\langle R_\phi,\chi_{\phi '}\rangle _{\mth{G}}\right)_{\phi,\phi '}.\]
Puisque les $\langle \chi_\phi,\cdot\rangle _{\mth{G}}$ forment une base de l'espace des distributions sur $\mathcal{H}_{\mth{G}}$, il suffit de montrer que cette matrice est inversible.

Soit $\mathcal{F}$ une famille de caractères unipotents. D'après \cite[12.3]{car}, on peut lui attacher un groupe fini $\Gamma=\Gamma\left(\mathcal{F}\right)$, qui est soit $\left({\mth{Z}}/2{\mth{Z}}\right)^n$, $n\geq 0$, soit le groupe symétrique $S_n$, $n=3,4,5$, et qui vérifie la propriété suivante: fixons un ensemble de représentants $R$ des classes de conjugaison de $\Gamma$. Soit $M\left(\Gamma\right)$ l'ensemble des couples $\left(x,\sigma\right)$, où $x\in R$ et $\sigma$ est un caractère irréductible de $Z_\Gamma\left(x\right)$; il existe alors une bijection $\left(x,\sigma\right)\mapsto\chi_{\left(x,\sigma\right)}$ entre $M\left(\Gamma\right)$ et $\mathcal{F}$ telle que pour tout $\left(x,\sigma\right)\in M\left(\Gamma\right)$ et tout caractère irréductible $\phi$ de $W$, si $\left(y,\tau\right)$ est l'élément de $M\left(\Gamma\right)$ tel que $\chi_\phi=\chi_{\left(y,\tau\right)}$, $\langle R_\phi,\chi_{\left(x,\sigma\right)}\rangle _{\mth{G}}$ vaut:
\begin{multline*}\left\{\left(x,\sigma\right),\left(y,\tau\right)\right\}\\=\dfrac 1{\card\left(Z_\Gamma\left(x\right)\right)\card\left(Z_\Gamma\left(y\right)\right)}\sum_{g\in\Gamma,xgyg^{-1}=gyg^{-1}x}\sigma\left(gyg^{-1}\right)\overline{\tau\left(g^{-1}xg\right)}.\end{multline*}
En particulier, les $\left\{\left(x,\sigma\right),\left(x,\sigma\right)\right\}$ sont non nuls, donc pour tout $\phi$, $\chi_\phi$ intervient dans $R_\phi$; on en déduit que si $\phi,\phi '$ sont tels que $\chi_\phi$ et $\chi_{\phi '}$ ne sont pas dans la même famille, $\langle R_\phi,\chi_{\phi '}\rangle _{\mth{G}}=0$. La matrice $M_{\mth{G}}$ est donc (en ordonnant convenablement les $\phi$) diagonale par blocs, chaque bloc étant constitué des termes d'indice $\left(\phi,\phi '\right)$ tels que $\chi_\phi$ et $\chi_{\phi '}$ appartiennent à une famille donnée;
le bloc de $M_{\mth{G}}$ correspondant à $\mathcal{F}$ est constitué des $\left\{\left(x,\sigma\right),\left(y,\tau\right)\right\}$, où $\left(x,\sigma\right)\in M\left(\Gamma\left(\mathcal{F}\right)\right)$ est tel qu'il existe un $\phi$ tel que $\chi_\phi=\chi_{\left(x,\sigma\right)}$, et de même pour $\left(y,\tau\right)$. Il suffit donc de vérifier que chacun de ces blocs est inversible; pour cela, il faut considérer séparément les différents cas possibles.

Si maintenant ${\mth{G}}$ n'est pas déployé, on obtient une matrice:
\[M_{\mth{G}}=\left(\langle R_\phi,\chi\rangle _{\mth{G}}\right)_{\phi,\chi}\]
qui n'est plus forcément une matrice carrée; ici, les blocs à considérer sont composés des termes d'indice ${\left(R_\phi,\chi\right)}$ tels que les $\chi$ sont les caractères de la série principale appartenant à une famille donnée et les $R_\phi$ sont ceux dans lesquels les $\chi$ interviennent; il suffit de vérifier pour chacun de ces blocs que, si $r$ est le cardinal de l'ensemble des $\chi$ de la série principale contenus dans la famille correspondante, le bloc est de rang $r$. Là aussi, il faut considérer séparément les différents cas possibles.

Supposons d'abord ${\mth{G}}$ de type $A_n$, $n\geq 1$. Alors d'après \cite[13.9]{car}, toutes les familles de caractères unipotents sont des singletons, et le groupe qui leur est attaché est $\left\{1\right\}$; la matrice $M_{\mth{G}}$ est alors une matrice diagonale dont les termes diagonaux valent $\left\{\left(1,1\right),\left(1,1\right)\right\}=1$; c'est donc la matrice identité, qui est inversible.

Supposons ensuite ${\mth{G}}$ de type ${}^2A_n$, $n\geq 1$. Alors d'après \cite[13.9]{car}, toutes les familles de caractères unipotents sont également des singletons, et si $\chi$ est un caractère unipotent, s'il appartient à la série principale, il existe un $\phi$ tel que $R_\phi=\chi$ et le bloc correspondant vaut $\left\{1\right\}$; si $\chi$ n'est pas dans la série principale, le bloc est vide. Ce bloc est donc toujours de même rang que le nombre de caractères de la série principale contenus dans la famille, ce qui montre l'assertion cherchée.

Supposons maintenant ${\mth{G}}$ de type $B_n$ ou $C_n$, $n\geq 2$. D'après \cite[13.9]{car}, les caractères unipotents de ${\mth{G}}$ sont alors en bijection avec les symboles de la forme suivante:
\[\left(\begin{array}{c}\lambda_1,\lambda_2,\dots,\lambda_a\\\mu_1,\dots,\mu_b\end{array}\right),\]
où $a-b$ est impair et positif, les suites $\left(\lambda_i\right)$ et $\left(\mu_i\right)$ sont des suites strictement croissantes d'entiers naturels, $\left(\lambda_1,\mu_1\right)\neq\left(0,0\right)$, et on a:
\[\sum_{i=1}^a\lambda_i+\sum_{i=1}^b\mu_i-\left[\left(\dfrac{a+b-1}2\right)^2\right]=n,\]
où pour $x\in{\mth{R}}$, $\left[x\right]$ représente la partie entière de $x$; le membre de gauche de l'égalité ci-dessus est appelé le {\em rang} du symbole. Parmi ces symboles, ceux correspondant à des caractères de la série principale sont ceux tels que $a-b=1$.

Soit $\left(\begin{array}{c}\lambda_1,\lambda_2,\dots,\lambda_a\\\mu_1,\dots,\mu_b\end{array}\right)$ et $\left(\begin{array}{c}\lambda'_1,\lambda'_2,\dots,\lambda'_{a'}\\\mu'_1,\dots,\mu'_{b'}\end{array}\right)$ deux symboles de rang $n$; ils correspondent à deux caractères de la même famille si et seulement si les multi-ensembles $\left\{\lambda_1,\dots,\lambda_a,\mu_1,\dots,\mu_b\right\}$ et $\left\{\lambda'_1,\dots,\lambda'_{a'},\mu'_1,\dots,\mu'_{b'}\right\}$ sont égaux. Si un élément apparaît dans les deux lignes d'un symbole, il apparaîtra également dans les deux lignes de tous les symboles de la même famille; pour une famille donnée, si $Z$ est l'ensemble des entiers apparaissant exactement une fois dans les symboles de la famille, et si son cardinal est $2z+1$, lesdits symboles sont alors en bijection avec les parties de $Z$ de cardinal congru à $z$ modulo $2$ (l'image d'un symbole étant alors l'ensemble des éléments de $Z$ se trouvant sur sa ligne du haut si $a-b$ est congru à $3$ modulo $4$ et sur sa ligne du bas sinon), et ceux correspondant à des caractères de la série principale sont en bijection avec les parties de $Z$ de cardinal exactement $z$.

Fixons une famille $\mathcal{F}$ de symboles. Elle contient un unique symbole spécial; c'est celui pour lequel on a $a-b=1$ et:
\[\lambda_1\leq\mu_1\leq\lambda_2\leq\dots\leq\mu_b\leq\lambda_a.\]
La partie $Z_*$ de $Z$ image de ce symbole par la bijection ci-dessus est l'ensemble des $\mu_i$ différents à la fois de $\lambda_i$ et de $\lambda_{i+1}$; posons également $Z^*=Z-Z_*$.

Soit un symbole de $\mathcal{F}$, et soit $M\subset Z$ son image par la bijection ci-dessus. Posons:
\[M^\#=\left(M\cup Z_*\right)-\left(M\cap Z_*\right).\]
Le cardinal de $M^\#$ est pair; de plus, on a:
\[M=\left(M^\#\cup Z_*\right)-\left(M^\#\cap Z_*\right).\]
L'application $M\mapsto M^\#$ est donc une bijection entre les parties de $Z$ de cardinal congru à $z$ modulo $2$ et les parties de $Z$ de cardinal pair; de plus, les parties $M$ de $Z$ de cardinal $z$ sont en bijection avec les parties $M^\#$ de $Z$ telles que $M^\#\cap Z_*$ et $M^\#\cap Z^*$ sont de même cardinal.

Soit donc $V$ l'ensemble des parties de $Z$ de cardinal pair. Si $M^\#$ et $N^\#$ sont deux éléments de $V$, leur différence symétrique est encore un élément de $V$; de plus, cette différence définit une loi de groupe sur $V$, pour laquelle $V$ est isomorphe à $\left({\mth{F}}_2\right)^{2z}$. Considérons la forme alternée sur $V$ définie par:
\[\langle M^\#,N^\#\rangle =\card\left(M^\#\cap N^\#\right)\mod 2.\]
Ecrivons $Z=\left\{z_0,\dots,z_{2z}\right\}$ et posons, pour $1<i<2z$:
\[f_i=\left\{z_0,\dots,z_i\right\}\]
si $i$ est impair, et:
\[f_i=\left\{z_0,\dots,z_{i-2},z_i\right\}\]
si $i$ est pair. Alors $\left(f_1,\dots,f_{2z}\right)$ est une base de $V$, pour laquelle la matrice de la forme alternée $\langle \cdot,\cdot\rangle $ est la matrice diagonale par blocs dont tous les blocs diagonaux sont $\left(\begin{array}{cc}0&1\\1&0\end{array}\right)$. Soit $\Gamma$ (resp. $\Gamma '$) le sous-groupe de $V$ engendré par les $f_i$, $i$ impair (resp. $i$ pair). Pour tout $g'\in\Gamma '$ non nul, l'application linéaire $\left(-1\right)^{\langle g,\cdot\rangle }$ de $\Gamma$ dans ${\mth{C}}$ est non nulle; $\Gamma '$ est donc canoniquement isomorphe à $\Gamma^*$. De plus, puisque $\Gamma$ est abélien, on a $M\left(\Gamma\right)=\Gamma\times\Gamma^*$; $M\left(\Gamma\right)$ est donc en bijection canonique avec $V$. Enfin, soit $M^\#,N^\#\in V$, et soit $\left(x,\sigma\right)$ et $\left(y,\tau\right)$ les éléments correspondants de $M\left(\Gamma\right)$; on a:
\[\left\{\left(x,\sigma\right),\left(y,\tau\right)\right\}=\dfrac 1{2^z}\sigma\left(y\right)\tau\left(x\right)=\dfrac 1{2^z}\left(-1\right)^{\langle x',y\rangle +\langle y',x\rangle },\]
où $x'$ et $y'$ sont les éléments de $\Gamma$ correspondant respectivement à $\sigma$ et $\tau$. De plus, on a, puisque $\Gamma=\Gamma^\perp$ et $\Gamma '=\Gamma '{}^\perp$ pour $\langle \cdot,\cdot\rangle $:
\[\langle x',y\rangle +\langle y',x\rangle =\langle x+x',y+y'\rangle =\langle M^\#,N^\#\rangle ,\]
d'où:
\[\left\{\left(x,\sigma\right),\left(y,\tau\right)\right\}=\dfrac 1{2^z}\left(-1\right)^{\card\left(M^\#\cap N^\#\right)}.\]
Il suffit donc de montrer que la matrice suivante:
\[\mathcal{M}=\left(\left(-1\right)^{\card\left(M^\#\cap N^\#\right)}\right)_{M^\#,N^\#}\]
où $M^\#$ et $N^\#$ décrivent l'ensemble $E$ des parties de $Z$ dont les intersections avec $Z_*$ et avec $Z^*$ sont de même cardinal, est inversible.
On va en fait montrer le lemme un peu plus général suivant que l'on appliquera à la transposée de $\mathcal{M}$:

\begin{lemme}\label{intorbf1}
Soit $z,z',d\in{\mth{N}}$, $d'\in{\mth{Z}}$ vérifiant $d\leq z$, $z-d=z'+d\pm 1$ et $-z'\leq d'\leq z$. Soit $Z^*$ (resp $Z_*$) un ensemble de cardinal $z$ (resp. $z'$), et soit $Z$ l'union disjointe de $Z^*$ et $Z_*$; soit $E$ (resp $F$) l'ensemble des parties $M$ de $Z$ vérifiant $\card\left(M\cap Z^*\right)=\card\left(M\cap Z_*\right)+d$ (resp. $\card\left(M\cap Z^*\right)=\card\left(M\cap Z_*\right)+d'$).
Alors la matrice:
\[\mathcal{M}'=\left(\left(-1\right)^{\card\left(N\cap M\right)}\right)_{N\in F, M\in E}\]
a pour rang le cardinal de $F$. En particulier, si $d=d'$, elle est inversible.
\end{lemme}

\begin{proof}
On va montrer le lemme par récurrence sur le cardinal $z+z'$ de $Z$ (en remarquant tout d'abord que la condition $z-d=z'+d\pm 1$ impose $z+z'$ impair). Si $z+z'=1$, on est dans un des cas suivants:
\begin{itemize}
\item $z=1$, $z'=0$, $d$ et $d'$ valent $0$ ou $1$;
\item $z=0$, $z'=1$, $d=d'=0$.
\end{itemize}
Dans tous les cas, on a $E=\left\{Z\right\}$ si $d=1$, et $E=\left\{\emptyset\right\}$ si $d=0$; il en est de même pour $F$ en fonction de $d'$. $\mathcal{M}'$ vaut alors $\left(1\right)$ ou $\left(-1\right)$ suivant les cas, et est donc de rang $\card\left(F\right)=1$.

Supposons maintenant que $z+z'>1$, et remarquons d'abord que l'on n'a ni $d=z$ ni $d=z'=0$: en effet, le premier cas impose $0=z'+d\pm 1=z'+z\pm 1$, d'où $z'+z=1$; quant au second, il impose $z=\pm 1$, d'où $z+z'=1$.

Appelons $\lambda_1,\dots,\lambda_z$ (resp. $\mu_1,\dots,\mu_{z'}$) les éléments de $Z^*$ (resp. $Z_*$); le groupe symétrique $S_z$ (resp. $S_{z'}$) agit sur $Z^*$ (resp. $Z_*$), et cette action induit une action de $S_z\times S_{z'}$ sur l'ensemble des parties de $Z$; $E$ et $F$ sont stables par cette action. Soit $r=\card\left(E\right)$, $r'=\card\left(F\right)$, et $M_1,\dots,M_r$ (resp. $N_1,\dots,N_s$) les éléments de $E$; si $\left(e_1,\dots,e_r\right)$ (resp. $\left(f_1,\dots,f_{r'}\right)$ est la base canonique de ${\mth{R}}^r$ (resp. ${\mth{R}}^{r'}$), on posera, pour tout $i$, $e_{M_i}=e_i$ (resp. $f_{N_i}=f_i$).

Pour tout $i\in\left\{0,\dots,z'\right\}$, soit $E_i$ (resp. $F_i$) l'ensemble des éléments $M$ de $E$ (resp. $F$) tels que $M\cap Z_*$ est de cardinal $i$; posons:
\[v_i=\sum_{M\in E_i}e_{M};\]
\[v'_i=\sum_{N\in F_i}f_{N}.\]
Pour toute transposition élémentaire $s$ de $S_z$ et tout élément $M$ de $E$ (resp. $N$ de $F$) non invariant par $\left(s,1\right)$, posons:
\[w_{M,s}=e_{M}-e_{\left(s,1\right)\left(M\right)};\]
\[w'_{N,s}=f_{N}-f_{\left(s,1\right)\left(N\right)}.\]
Pour toute transposition élémentaire $s'$ de $S_{z'}$ et tout élément $M$ de $E$ (resp. $N$ de $F$) non invariant par $\left(1,s'\right)$, posons enfin:
\[u_{M,s'}=e_{M}-e_{\left(1,s'\right)\left(M\right)};\]
\[u'_{N,s'}=f_{N}-f_{\left(1,s'\right)\left(N\right)}.\]
Considérons la famille d'éléments de ${\mth{R}}^r$ constituée des $v_i$, des $w_{M,s}$ et des $u_{M,s'}$. Cette famille engendre ${\mth{R}}^r$; en effet, on a pour tout $M$, si $i=\card\left(M\cap Z_*\right)$:
\[e_{M}=\dfrac 1{\card\left(E_i\right)}\left(v_i-\sum_{N\in E_i}\left(e_{N}-e_{M}\right)\right);\]
de plus, $S_z\times S_{z'}$ agit transitivement sur $E_i$, et il existe donc pour tout $N\in E_i$, une suite de transpositions élémentaires $s_1,\dots,s_l$ de $S_z$ et $s'_1,\dots,s'_{l+l'}$ de $S_{z'}$ telles que, si on pose $N_0=N$ et pour tout $i$, par récurrence, $N_i=\left(s_i,1\right)\left(N_{i-1}\right)$ si $1\leq i\leq l$ et $N_i=\left(1,s'_{i-l}\right)\left(N_{i-1}\right)$ si $l<i\leq l+l'$, on a $N_{l+l'}=M$; on en déduit:
\[e_{N}-e_{M}=\sum_{i=1}^{l+l'}\left(e_{N_i}-e_{N_{i-1}}\right)=\sum_{i=1}^lw_{N_{i-1},s_i}+\sum_{i=l+1}^{l'}u_{N_{l+i-1},s'_i}.\]
Il en est de même pour la famille d'éléments de ${\mth{R}}^s$ constituée des $v'_i$, des $w'_{N,s}$ et des $u'_{N,s'}$; il suffit donc de montrer que tous ces éléments sont dans l'image de $\mathcal{M}$. Considérons d'abord les $w'_{N,s}$. Si $d'=z$, il n'en existe pas et l'assertion cherchée est triviale; supposons donc $d'<z$. Pour tout $s$ et tout $M\in E$, si $s$ permute les éléments $\lambda_i$ et $\lambda_{i+1}$ de $Z^*$, avec par exemple $\lambda_i\in M$ et $\lambda_{i+1}\not\in M$, on a:
\begin{equation*}\begin{split}\mathcal{M}'w_{M,s}&=\sum_{N\in F}\left(\left(-1\right)^{\card\left(M\cap N\right)}-\left(-1\right)^{\card\left(\left(s,1\right)\left(M\right)\cap N\right)}\right)e_{N}\\&=\sum_{N\in F}\left(-1\right)^{\card\left(M\cap N\right)}\left(e_{\left(s,1\right)N}-e_{N}\right)\\&=-\sum_{N\in F,s\left(N\right)\neq N}\left(-1\right)^{\card\left(M\cap N\right)}w'_{N,s}\\&=-2\sum_{N\in F,\lambda_i\in N,\lambda_{i+1}\not\in N}\left(-1\right)^{\card\left(M\cap N\right)}w'_{N,s}.\end{split}\end{equation*}

Fixons un $s$; il existe une bijection entre l'ensemble des $M\in E$ tels que $\lambda_i\in M$ et $\lambda_{i+1}\not\in M$ et l'ensemble $E'$ de parties de $Z-\left\{\lambda_i,\lambda_{i+1}\right\}$ défini de manière similaire à $E$ en remplaçant $d$ par $d-1$ (si $d=0$, on inverse les rôles de  $Z^*-\left\{\lambda_i,\lambda_{i+1}\right\}$ et $Z_*$ et on remplace $d$ par $1$), donnée par $M\mapsto\left(Z-\left\{\lambda_{i+1}\right\}\right)-M$; on a de même une bijection entre l'ensemble des $N\in F$ contenant $\lambda_i$ et pas $\lambda_{i+1}$ et l'ensemble $F'$ défini de manière similaire à $E'$. De plus, on a pour tous $M,N$ contenant $\lambda_i$ et pas $\lambda_{i+1}$, puisque $Z-\left\{\lambda_{i+1}\right\}$ est de cardinal pair:
\begin{equation*}\begin{split}\left(-1\right)^{\card\left(M\cap N\right)}&=\left(-1\right)^{\card\left(M\cap\left(\left(Z-\left\{\lambda_{i+1}\right\}\right)-N\right)\right)}\\&=\left(-1\right)^{\card\left(\left(\left(Z-\left\{\lambda_{i+1}\right\}\right)-M\right)\cap\left(\left(Z-\left\{\lambda_{i+1}\right\}\right)-N\right)\right)}.\end{split}\end{equation*}
On vérifie aisément que si $d>0$, le septuplet:
\[\left(z-2,z',d-1,d'-1,Z^*-\left\{\lambda_i,\lambda_{i+1}\right\},Z_*,Z-\left\{\lambda_i,\lambda_{i+1}\right\}\right)\]
(resp. si $d=0$, le septuplet:
\[\left(z',z-2,1,1-d',Z_*,Z^*-\left\{\lambda_i,\lambda_{i+1}\right\},Z-\left\{\lambda_i,\lambda_{i+1}\right\}\right))\]
vérifie les conditions du lemme; en lui appliquant l'hypothèse de récurrence, on voit que la matrice constituée des termes ci-dessus est surjective; pour tout $M$, $w'_{N,s}$ est donc dans l'image de $\mathcal{M}'$, et ceci est vrai pour tout $s$. Par un raisonnement similaire, tous les $u'_{N,s'}$ sont également dans l'image de $\mathcal{M}$ (il n'en existe pas si $d'=-z'$, et si $d'>-z'$, on applique l'hypothèse de récurrence au septuplet $\left(z,z'-2,d+1,d'+1,Z^*,Z_*-\left\{\mu_i,\mu_{i+1}\right\},Z-\left\{\mu_i,\mu_{i+1}\right\}\right)$), où $\mu_i,\mu_{i+1}$ sont les éléments de $Z_*$ permutés par $s'$.

Considérons maintenant les $v_i$: on a, pour tout $i$:
\[\mathcal{M}'v_i=\sum_{M\in E_i}\sum_{N\in F}\left(-1\right)^{\card\left(M\cap N\right)}f_{N}=\sum_{j=0}^{z'}c_{j,i,z,z',d,d'}v'_j,\]
avec pour tous $i,j$, en fixant $N_j\in F_j$:
\[c_{j,i,z,z',d,d'}=\sum_{M\in E_i}\left(-1\right)^{\card\left(M\cap N_j\right)}.\]
Soit $i_1=\inf\left(z-d,z'\right)$, $j_1=\inf\left(z-d',z'\right)$; $E_i$ (resp $F_j$) est non vide si et seulement si $0\leq i\leq i_1$ (resp. $0\leq j\leq j_1$). On va montrer que pour tous ces $j$-là, $v'_j$ est dans l'image de $\mathcal{M}$, ce qui est le cas si et seulement si la matrice:
\[\left(c_{j,i,z,z',d,d'}\right)_{0\leq i\leq i_1,0\leq j\leq j_1}\]
est de rang $j_1+1$. Pour cela, fixons une suite croissante $M_0\subset M_1\subset\dots\subset M_{i_1}$, où pour tout $i$, $M_i\in E_i$; en posant $I=\left\{1,\dots,i_1\right\}$, pour toute partie $I'$ de $I$, posons:
\[M_{I'}=M_0\cup\bigcup_{i\in I'}\left(M_i-M_{i-1}\right).\]
C'est un élément de $E_{\card\left(I'\right)}$: en effet, pour tout $i$, $M_i-M_{i-1}$ contient exactement un élément de $Z^*$ et un élément de $Z_*$, et tous les $M_i-M_{i-1}$ sont disjoints.
Pour toute partie $I'$ de $I$, et tout $N\in F$, soit $\chi$ l'application de l'ensemble des parties de $I'$ dans ${\mth{C}}$ donnée par:
\[\chi:K\subset I'\mapsto\left(-1\right)^{\card\left(N\cap\left(M_K-M_0\right)\right)}:\]
si on munit l'ensemble des parties de $I'$ de la loi de groupe définie par la différence symétrique, $\chi$ est un caractère de ce groupe, et on a donc:
\[\sum_{K\subset I'}\chi\left(K\right)=0\]
sauf si $\chi=1$, ce qui est le cas si et seulement si tous les $M_i-M_{i-1}$, $i\in I$, sont soit inclus dans $N$ soit disjoints de $N$. Soit $F_{j,I'}$ l'ensemble des $N\in F_i$ vérifiant cette condition; posons, pour toute partie $I'$ de $I$ de cardinal $i$:
\begin{equation*}\begin{split}k_{j,i,z,z',d,d'}&=2^{-i}\sum_{k=0}^{i}\left(\begin{array}{c}i\\k\end{array}\right)c_{j,k,z,z',d,d'}\\&=2^{-i}\sum_{N\in F_i}\sum_{K\subset I'}\left(-1\right)^{\card\left(N\cap M_K\right)}\\&=\sum_{N\in F_{j,I'}}\left(-1\right)^{\card\left(N\cap M_0\right)}.\end{split}\end{equation*}
On a le produit matriciel suivant:
\[\left(k_{j,i,z,z',d,d'}\right)_{j,i}=\left(c_{j,i,z,z',d,d'}\right)_{j,i}\left(2^{-i}\left(\begin{array}{c}i\\i'\end{array}\right)\right)_{i,i'\in\left\{0,\dots,i_1\right\}}.\]
La seconde matrice du membre de droite est ne matrice triangulaire supérieure dont les coefficients diagonaux sont $1,2^{-1},\dots,2^{-i_1}$; elle est donc inversible, et on en déduit que la matrice des $c_{j,i,z,z',d,d'}$ est de rang $j_1-1$ si et seulement si celle des $k_{j,i,z,z',d,d'}$ l'est. On va donc monter que les lignes de cette matrice sont linéairement indépendantes.

Pour tout $i\in\left\{1,\dots,i_1\right\}$, en posant $I'=\left\{1,\dots,i\right\}$, écrivons $F_{j,I'}=F_{j,I',1}\cup F_{j,I',2}$, où $F_{j,I',1}$ (resp. $F_{j,I',2}$) est l'ensemble des éléments de $F_{j,I'}$ qui contiennent $M_i-M_{i-1}$ (resp. disjoints de $M_i-M_{i-1}$). Soit $E'_i$ (resp. $F'_i$) l'ensemble des parties de $Z-\left(M_i-M_{i-1}\right)$ défini de manière similaire à $E$ (resp. $F$); définissons également $\left(F'_i\right)_{j,K}$, pour $K\subset I'-\left\{i\right\}$, de manière similaire à $F_{j,K}$. Alors les éléments de $F_{j,I',2}$ sont exactement ceux de $\left(F'_i\right)_{j,I'-\left\{i\right\}}$, et l'application $N\mapsto N-\left(M_i-M_{i-1}\right)$ est une bijection canonique de $F_{j,I',1}$ dans $\left(F'_i\right)_{j-1,I'-\left\{i\right\}}$; on en déduit:
\begin{multline*}\sum_{N\in F_{j,I'}}\left(-1\right)^{\card\left(N\cap M_0\right)}\\=\sum_{N\in\left(F'_i\right)_{j,I'-\left\{i\right\}}}\left(-1\right)^{\card\left(N\cap M_0\right)}+\sum_{N\in\left(F'_i\right)_{j-1,I'-\left\{i\right\}}}\left(-1\right)^{\card\left(N\cap M_0\right)},\end{multline*}
soit:
\[k_{j,i,z,z',d,d'}=k_{j,i-1,z-1,z'-1,d,d'}+k_{j-1,i-1,z-1,z'-1,d,d'}.\]
Soit $\lambda_0,\dots,\lambda_{j_1}$ tels que pour tout $i$, $\sum_j\lambda_jk_{j,i,z,z',d,d'}=0$; on en déduit pour $i\geq 1$:
\[\sum_{j=0}^{j_1-1}\left(\lambda_j+\lambda_{j+1}\right)k_{j,i-1,z-1,z'-1,d,d'}=0.\]
Par hypothèse de récurrence, la matrice $\left(k_{j,i,z-1,z'-1,d,d'}\right)_{j,i}$ est de rang $j_1$, donc tous les $\lambda_j+\lambda_{j+1}$ sont nuls; on en déduit, pour tout $j$:
\[\lambda_j=\left(-1\right)^j\lambda_0.\]
Pour achever la démonstration, il suffit donc de vérifier que l'on a:
\[\sum_{j=0}^{j_1}\left(-1\right)^jk_{j,0,z,z',d,d'}\neq 0;\]
cela montrera que tous les $\lambda_j$ sont nuls et donc que les lignes de la matrice des $k_{j,i,z,z',d,d'}$ sont bien linéairement indépendantes. Pour tout $j$, on a:
\begin{equation*}\begin{split}k_{j,0,z,z',d,d'}&=c_{j,0,z,z',d,d'}\\&=\sum_{N\in F_j}\left(-1\right)^{\card\left(N\cap M_0\right)}\\&=\left(\begin{array}{c}z'\\j\end{array}\right)\sum_{k=0}^j\left(-1\right)^k\left(\begin{array}{c}d\\k\end{array}\right)\left(\begin{array}{c}z-d\\j+d'-k\end{array}\right).\end{split}\end{equation*}
On en déduit:
\begin{equation*}\begin{split}\sum_{j=0}^{j_1}\left(-1\right)^jk_{j,0,z,z',d,d'}&=\sum_{j=0}^{j_1}\left(-1\right)^j\left(\left(\begin{array}{c}z'\\j\end{array}\right)\sum_{k=0}^j\left(-1\right)^k\left(\begin{array}{c}d\\k\end{array}\right)\left(\begin{array}{c}z-d\\j+d'-k\end{array}\right)\right)\\&=\sum_{l=0}^{z-d}\left(-1\right)^l\left(\begin{array}{c}z-d\\l\end{array}\right)\left(\sum_{j=0}^{z'}\left(\begin{array}{c}z'\\j\end{array}\right)\left(\begin{array}{c}d\\j+d'-l\end{array}\right)\right).\end{split}\end{equation*}
Montrons donc le lemme suivant:

\begin{lemme}
Pour tout $l$, on a:
\[\sum_{j=0}^{z'}\left(\begin{array}{c}z'\\j\end{array}\right)\left(\begin{array}{c}d\\j+d'-l\end{array}\right)=\left(\begin{array}{c}z'+d\\l+d-d'\end{array}\right).\]
\end{lemme}

\begin{proof}
En effet, considérons le polynôme $\left(1+X\right)^{z'+d}$. Son terme de degré $l+d-d'$ est $\left(\begin{array}{c}z'+d\\l+d-d'\end{array}\right)$; or on a:
\[\left(1+X\right)^{z'+d}=\left(1+X\right)^{z'}\left(1+X\right)^{d},\]
et le terme de degré $l+d-d'$ du membre de droite de cette égalité vaut:
\[\sum_{j=0}^{z'}\left(\begin{array}{c}z'\\j\end{array}\right)\left(\begin{array}{c}d\\l+d-d'-j\end{array}\right).\]
Comme pour tout $j$, on a $\left(\begin{array}{c}d\\l+d-d'-j\end{array}\right)=\left(\begin{array}{c}d\\j+d'-l\end{array}\right)$, on trouve bien l'égalité de l'énoncé.
\end{proof}

On déduit de ce lemme que l'on a:
\[\sum_{j=0}^{j_1}\left(-1\right)^jk_{j,0,z,z',d,d'}=\sum_{l=0}^{z-d}\left(-1\right)^l\left(\begin{array}{c}z-d\\l\end{array}\right)\left(\begin{array}{c}z'+d\\l+d-d'\end{array}\right).\]
Quitte à remplacer $z$ par $z-d$, $z'$ par $z'+d$ et $d'$ par $d'-d$, on peut donc supposer $d=0$; on a alors $z=z'\pm 1$. Supposons par exemple $z'=z+1$; on a le lemme suivant:

\begin{lemme}
Pour tout $k\in{\mth{Z}}$, on a:
\[\sum_{l=0}^z\left(-1\right)^l\left(\begin{array}{c}z\\l\end{array}\right)\left(\begin{array}{c}z+1\\l+k\end{array}\right)=\left(-1\right)^{\left[\dfrac{z-k+1}2\right]}\left(\begin{array}{c}z\\{\left[\dfrac{z-k+1}2\right]}\end{array}\right).\]
\end{lemme}

\begin{proof}
En effet, le membre de gauche de cette égalité est égal à:
\[\sum_{l=0}^z\left(-1\right)^l\left(\begin{array}{c}z\\l\end{array}\right)\left(\begin{array}{c}z+1\\z+1-l-k\end{array}\right),\]
qui est le terme de degré $z+1-k$ du polynôme:
\[\left(X-1\right)^z\left(X+1\right)^{z+1}=\left(X^2-1\right)^z\left(X+1\right).\]
Or pour tout $l\in{\mth{Z}}$, les termes de degré $2l$ et $2l+1$ de ce polynôme sont égaux et valent $\left(-1\right)^l\left(\begin{array}{c}z\\l\end{array}\right)$. Le terme de degré $z-k+1$ est donc égal au membre de droite de l'égalité de l'énoncé, ce qui démontre le lemme.
\end{proof}

En appliquant ce lemme au cas $k=-d'$, on a:
\[\sum_{l=0}^z\left(-1\right)^l\left(\begin{array}{c}z\\l\end{array}\right)\left(\begin{array}{c}z+1\\l-d'\end{array}\right)=\left(-1\right)^{\left[\dfrac{z+d'+1}2\right]}\left(\begin{array}{c}z\\{\left[\dfrac{z+d'+1}2\right]}\end{array}\right).\]
Puisque $d'$ est compris entre $-z'=-z-1$ et $z$, on a $0\leq\left[\dfrac{z+d'+1}2\right]\leq z$, le membre de droite de l'égalité ci-dessus est non nul, ce qui achève la démonstration.
\end{proof}

Supposons maintenant ${\mth{G}}$ de type $D_n$, $n\geq 4$. D'après \cite[13.9]{car}, il existe une surjection de l'ensemble des caractères unipotents de ${\mth{G}}$ dans l'ensemble des symboles de la forme suivante (modulo inversion des deux lignes):
\[\left(\begin{array}{c}\lambda_1,\lambda_2,\dots,\lambda_a\\\mu_1,\dots,\mu_b\end{array}\right),\]
où $a-b$ est multiple de $4$, les suites $\left(\lambda_i\right)$ et $\left(\mu_i\right)$ sont des suites strictement croissantes d'entiers naturels, $\left(\lambda_1,\mu_1\right)\neq\left(0,0\right)$, et on a:
\[\sum_{i=1}^a\lambda_i+\sum_{i=1}^b\mu_i-\left[\left(\dfrac{a+b-1}2\right)^2\right]=n\]
(le membre de gauche de l'égalité ci-dessus est ici encore appelé le rang du symbole); un symbole $\left(\begin{array}{c}\lambda_1,\lambda_2,\dots,\lambda_a\\\mu_1,\dots,\mu_b\end{array}\right)$ a deux antécédents par cette surjection si $a=b$ et pour tout $i$, $\lambda_i=\mu_i$, et un seul dans tous les autres cas.
Parmi ces symboles, ceux correspondant à des caractères de la série principale sont ceux tels que $a=b$.

Soit $\left(\begin{array}{c}\lambda_1,\lambda_2,\dots,\lambda_a\\\mu_1,\dots,\mu_b\end{array}\right)$ un symbole de rang $n$. Si $a=b$ et $\lambda_i=\mu_i$ pour tout $i$, les deux caractères unipotents correspondants sont dans deux familles distinctes et réduites à eux-mêmes; sinon, si $\left(\begin{array}{c}\lambda'_1,\lambda'_2,\dots,\lambda'_{a'}\\\mu'_1,\dots,\mu'_{b'}\end{array}\right)$ est un autre symbole de rang $n$, ces deux symboles correspondent à deux caractères de la même famille si et seulement si les multi-ensembles $\left\{\lambda_1,\dots,\lambda_a,\mu_1,\dots,\mu_b\right\}$ et $\left\{\lambda'_1,\dots,\lambda'_{a'},\mu'_1,\dots,\mu'_{b'}\right\}$ sont égaux. On en déduit, de même que dans le cas $B_n/C_n$, que si $Z$ est l'ensemble des entiers apparaissant exactement une fois dans les symboles d'une famille, et si son cardinal est $2z$, lesdits symboles sont en bijection avec les couples non ordonnés de parties de $Z$ de cardinal congru à $z$ modulo $2$ de la forme $\left(M,Z-M\right)$ (l'image d'un symbole étant alors le couple formé de l'ensemble des éléments de $Z$ se trouvant sur chacune des deux lignes), et ceux correspondant à des caractères de la série principale sont en bijection avec les couples de parties de $Z$ de cardinal exactement $z$ de la forme ci-dessus.

Si $z=0$, on est dans le cas où $a=b$ et $\lambda_i=\mu_i$ pour tout $i$ et $\mathcal{M}=\left(1\right)$ est clairement inversible; supposons donc $z>0$.
Soit $M$ une partie de $Z$ de cardinal congru à $z$ modulo $2$; posons:
\[M^\#=\left(M\cup Z_*\right)-\left(M\cap Z_*\right).\]
D'après ce que l'on a déjà vu, l'application $M\mapsto M^\#$ est une bijection entre les parties de $Z$ de cardinal congru à $z$ modulo $2$ et les parties de $Z$ de cardinal pair; de plus, les parties $M$ de $Z$ de cardinal $z$ sont en bijection avec les parties $M^\#$ de $Z$ telles que $M^\#\cap Z_*$ et $M^\#\cap Z^*$ sont de même cardinal; enfin; on a pour tout $M$:
\[\left(Z-M\right)^\#=Z-M^\#.\]
Considérons le ${\mth{F}}_2$-espace vectoriel $V$ constitué de l'ensemble des parties de $Z$ de cardinal pair muni de la différence symétrique, et la forme alternée sur $V$ définie par:
\[\langle M^\#,N^\#\rangle =\card\left(M^\#\cap N^\#\right)\mod 2.\]
Puisque le cardinal de $Z$ est pair, on a pour tous $M^\#,N^\#$:
\[\langle M^\#,Z-N^\#\rangle =\langle M^\#,N^\#\rangle .\]
Considérons donc l'espace $V_1=V/\left\{0,Z\right\}$; les éléments de $V_1$ sont les couples non ordonnés de la forme $\left(M^\#,Z-M^\#\right)$, où $M^\#$ est de cardinal pair. $V_1$ est muni de la forme alternée induite par celle de $V$, et par un raisonnement similaire à celui du cas $B_n/C_n$, il suffit de montrer que la matrice:
\[\mathcal{M}=\left(\left(-1\right)^{\langle {\mathcal{C}}_1,{\mathcal{C}}_2\rangle }\right)_{\mathcal{C}_1,\mathcal{C}_2}\]
où $\mathcal{C}_1$ et $\mathcal{C}_2$ décrivent l'ensemble des couples non ordonnés $\left(M^\#,Z-M^\#\right)$ de parties de $Z$ tels que les intersections de $M^\#$ avec $Z_*$ et avec $Z^*$ sont de même cardinal, est inversible.

Soit $x$ un élément de $Z_*$ quelconque; l'ensemble des couples non ordonnés de parties de $Z$ de la forme $\left(M^\#,Z-M^\#\right)$ tels que $M^\#$ est de cardinal pair (resp. $M^\#$ est de cardinal pair et $M^\#\cap Z_*$ et $M^\#\cap Z^*$ ont même cardinal) est en bijection avec l'ensemble des parties $M_1^\#$ de $Z-\left\{x\right\}$ de cardinal pair (resp. de cardinal pair et dont les intersections avec $Z_*-\left\{x\right\}$ et $Z^*$ ont même cardinal); l'image de $\left(M^\#,Z-M^\#\right)$ est celle des deux qui ne contient pas $x$. On a donc:
\[\mathcal{M}=\left(\left(-1\right)^{\card\left(M_1^\#\cap N_1^\#\right)}\right)_{M_1^\#,N_1^\#}\]
où $M_1^\#$ et $N_1^\#$ décrivent l'ensemble des parties de $Z-\left\{x\right\}$ dont les intersections avec $Z_*-\left\{x\right\}$ et avec $Z^*$ sont de même cardinal, et une telle matrice est inversible d'après la démonstration du cas $B_n/C_n$.

Supposons maintenant ${\mth{G}}$ de type ${}^2D_n$, $n\geq 4$. D'après \cite[13.9]{car}, il existe une bijection de l'ensemble des caractères unipotents de ${\mth{G}}$ dans l'ensemble des symboles de la forme suivante:
\[\left(\begin{array}{c}\lambda_1,\lambda_2,\dots,\lambda_a\\\mu_1,\dots,\mu_b\end{array}\right),\]
où $a-b$ est positif et congru à $2$ modulo $4$, les suites $\left(\lambda_i\right)$ et $\left(\mu_i\right)$ sont des suites strictement croissantes d'entiers naturels, $\left(\lambda_1,\mu_1\right)\neq\left(0,0\right)$, et on a:
\[\sum_{i=1}^a\lambda_i+\sum_{i=1}^b\mu_i-\left[\left(\dfrac{a+b-1}2\right)^2\right]=n\]
(le membre de gauche de l'égalité ci-dessus est ici encore appelé le rang du symbole).
Parmi ces symboles, ceux correspondant à des caractères de la série principale sont ceux tels que $a-b=2$.

Soit $\left(\begin{array}{c}\lambda_1,\lambda_2,\dots,\lambda_a\\\mu_1,\dots,\mu_b\end{array}\right)$ et $\left(\begin{array}{c}\lambda'_1,\lambda'_2,\dots,\lambda'_{a'}\\\mu'_1,\dots,\mu'_{b'}\end{array}\right)$ deux symboles de rang $n$; ces deux symboles correspondent à deux caractères de la même famille si et seulement si les multi-ensembles $\left\{\lambda_1,\dots,\lambda_a,\mu_1,\dots,\mu_b\right\}$ et $\left\{\lambda'_1,\dots,\lambda'_{a'},\mu'_1,\dots,\mu'_{b'}\right\}$ sont égaux. Encore une fois, on en déduit de même que dans le cas $B_n/C_n$ que si $Z$ est l'ensemble des entiers apparaissant exactement une fois dans les symboles d'une famille, et si son cardinal est $2z$, lesdits symboles sont en bijection avec les parties de $Z$ de cardinal inférieur à $z$ et congru à $z-1$ modulo $2$, et ceux correspondant à des caractères de la série principale sont en bijection avec les parties de $Z$ de cardinal exactement $z-1$.

Soit $Z^*$ (resp; $Z_*$) les parties de $Z$ définies comme dans le cas $B_n/C_n$; si  $M$ est une partie de $Z$, posons ici aussi:
\[M^\#=\left(M\cup Z_*\right)-\left(M\cap Z_*\right).\]
L'application $M\mapsto M^\#$ induit une bijection de l'ensemble des parties de $Z$ de cardinal $z$ (resp. $z-1$) dans l'ensemble $Z_0$ (resp. $Z_1$) des parties $M^\#$ de $Z$ telles que $\card\left(M^\#\cap Z^*\right)=\card\left(M^\#\cap Z_*\right)$ (resp. $\card\left(M^\#\cap Z_*\right)-1$).

D'après \cite[3.15.2]{l1}, il suffit de vérifier que la matrice:
\[\mathcal{M}=\left(\left(-1\right)^{\card\left(M^\#\cap N^\#\right)}\right)_{M^\#\in Z_0,N^\#\in Z_1}\]
est de rang $r=\card\left(Z_1\right)$; c'est un cas particulier du lemme \ref{intorbf1}, que l'on applique à la transposée de $\mathcal{M}$ en posant $z'=z-1$, $d=0$ et $d'=-1$.

Supposons maintenant ${\mth{G}}$ de type ${}^3D_4$; d'après \cite[1.17]{l2},
la série principale de ${\mth{G}}$ possède six caractères unipotents, quatre constituant quatre familles à eux seuls et les deux restants dans une même famille $\mathcal{F}$ de cardinal $4$; d'autre part, $W$ possède sept caractères irréductibles, et les éléments de $\mathcal{F}$ interviennent dans les $R_\phi$ correspondant à trois d'entre eux; d'après \cite[1.18]{l2}, la matrice bloc correspondante est:
\[M_1=\dfrac 12\left(\begin{array}{cc}1&1\\1&-1\\1&1\end{array}\right),\]
qui est clairement de rang $2$.

Supposons maintenant ${\mth{G}}$ de type $G_2$. D'après \cite[13.9]{car}, la série principale de ${\mth{G}}$ possède six caractères unipotents, deux constituant deux familles à eux seuls et les quatre restants dans une même famille $\mathcal{F}$ de cardinal $8$; d'autre part, $W$ possède six caractères irréductibles, et les éléments de $\mathcal{F}$ interviennent dans les $R_\phi$ correspondant à quatre d'entre eux; d'après \cite[13.6]{car}, la matrice bloc correspondante est:
\[M_2=\left(\begin{array}{cccc}\dfrac 16&\dfrac 12&\dfrac 13&\dfrac 13\\\dfrac 12&\dfrac 12&0&0\\\dfrac 13&0&\dfrac 23&-\dfrac 13\\\dfrac 13&0&-\dfrac 13&\dfrac 23\end{array}\right),\]
dont le déterminant vaut $-\dfrac 16$ et qui est donc inversible.

Supposons maintenant ${\mth{G}}$ de type $F_4$. D'après \cite[13.6, 13.9]{car}, il suffit de considérer les deux matrices suivantes:
\[M_3=\dfrac 12\left(\begin{array}{ccc}1&1&1\\1&1&-1\\1&-1&1\end{array}\right),\]
\[M_4=\left(\begin{array}{ccccccccccc}\dfrac 1{24}&\dfrac 18&\dfrac 18&\dfrac 1{12}&\dfrac 14&\dfrac 14&\dfrac 18&\dfrac 18&\dfrac 18&\dfrac 13&\dfrac 14\\\dfrac 18&\dfrac 38&\dfrac 38&\dfrac 14&\dfrac 14&\dfrac 14&-\dfrac 18&-\dfrac 18&-\dfrac 18&0&-\dfrac 14\\\dfrac 18&\dfrac 38&\dfrac 38&\dfrac 14&-\dfrac 14&-\dfrac 14&-\dfrac 18&-\dfrac 18&-\dfrac 18&0&\dfrac 14\\\dfrac 1{12}&\dfrac 14&\dfrac 14&\dfrac 16&0&0&\dfrac 14&\dfrac 14&\dfrac 14&-\dfrac 13&0\\\dfrac 14&\dfrac 14&-\dfrac 14&0&\dfrac 12&0&\dfrac 14&\dfrac 14&-\dfrac 14&0&0\\\dfrac 14&\dfrac 14&-\dfrac 14&0&0&\dfrac 12&-\dfrac 14&-\dfrac 14&\dfrac 14&0&0\\\dfrac 18&-\dfrac 18&-\dfrac 18&\dfrac 14&\dfrac 14&-\dfrac 14&\dfrac 38&-\dfrac 18&\dfrac 38&0&\dfrac 14\\\dfrac 18&-\dfrac 18&-\dfrac 18&\dfrac 14&\dfrac 14&-\dfrac 14&-\dfrac 18&\dfrac 38&-\dfrac 18&0&-\dfrac 14\\\dfrac 18&-\dfrac 18&-\dfrac 18&\dfrac 14&-\dfrac 14&\dfrac 14&\dfrac 38&-\dfrac 18&\dfrac 38&0&-\dfrac 14\\\dfrac 13&0&0&-\dfrac 13&0&0&0&0&0&\dfrac 23&0\\\dfrac 14&-\dfrac 14&\dfrac 14&0&0&0&\dfrac 14&-\dfrac 14&-\dfrac 14&0&\dfrac 12\end{array}\right).\]
Le déterminant de $M_3$ vaut $-\dfrac 12$, celui de $M_4$ vaut $\dfrac 1{192}$; elles sont donc toutes deux inversibles.

Supposons maintenant ${\mth{G}}$ de type $E_6$. D'après \cite[13.6, 13.9]{car}, il suffit de considérer la matrice $M_3$ et la matrice suivante:
\[M_5=\left(\begin{array}{ccccc}\dfrac 16&\dfrac 12&\dfrac 13&\dfrac 13&\dfrac 16\\\dfrac 12&\dfrac 12&0&0&-\dfrac 12\\\dfrac 13&0&\dfrac 23&-\dfrac 13&\dfrac 13\\\dfrac 13&0&-\dfrac 13&\dfrac 23&\dfrac 13\\\dfrac 16&-\dfrac 12&\dfrac 13&\dfrac 13&\dfrac 16\end{array}\right),\]
dont le déterminant vaut $\dfrac 16$ et qui est donc inversible.

Supposons enfin ${\mth{G}}$ de type ${}^2E_6$. D'après \cite[1.15]{l2} et \cite[13.6, 13.9]{car}, il suffit de considérer les matrices suivantes:
\[M_6=\dfrac 12\left(\begin{array}{cccc}1&1&1&1\\1&1&-1&-1\\1&-1&1&-1\\1&-1&-1&1\end{array}\right),\]
\[M_7=\left(\begin{array}{ccccc}\dfrac 23&\dfrac 13&0&0&-\dfrac 13\\\dfrac 13&\dfrac 16&-\dfrac 12&-\dfrac 12&\dfrac 13\\0&-\dfrac 12&\dfrac 12&-\dfrac 12&0\\0&-\dfrac 12&-\dfrac 12&\dfrac 12&0\\-\dfrac 13&\dfrac 13&0&0&\dfrac 23\end{array}\right).\]
La première est hermitienne par \cite[13.9]{car}, donc inversible; le déterminant de la seconde est $-\dfrac 16$, donc elle est également inversible. $\Box$

{\bf Remarque:} si ${\mth{G}}$ est de type $E_7$ ou $E_8$, d'après \cite[13.9]{car}, il existe des familles (une pour $E_7$ et deux pour $E_8$) de caractères unipotents de $G$ qui sont de cardinal $4$ et qui contiennent exactement deux caractères de la série principale, et les matrices blocs correspondantes sont égales à:
\[M_8=\dfrac 12\left(\begin{array}{cc}1&1\\1&1\end{array}\right)\]
qui n'est pas inversible; les $R_\phi$ ne sont alors pas linéairement indépendants, et ne peuvent pas non plus engendrer l'espace des distributions sur $\mathcal{H}_{\mth{G}}$ puisqu'ils sont en nombre égal à la dimension de cet espace; ${\mth{G}}$ ne vérifie donc aucune des deux assertions de la proposition.

\subsection{Calcul de $q_\mathcal{T}$}\label{cqt}

D'après la démonstration du lemme \ref{srg}, dont on reprendra les notations, il suffit de considérer le cas où $G$ est simple et admet un tore non ramifié maximal $T$ déployé, ce que l'on supposera par la suite. Les systèmes de racines de $G$ et de $G_{nr}$ relativement à un tore déployé maximal sont alors isomorphes; on notera $\left\{\alpha_1,\dots,\alpha_n\right\}$ un système de racines simples, et ${\alpha_1}^\vee,\dots,{\alpha_i}^\vee$ le système de coracines simples associé.
Pour chaque type, la nomérotation des racines est la même que dans la proposition \ref{gcj} (ici, on ne considère plus que le diagramme de Dynkin non étendu). Toujours d'après la démonstration du lemme \ref{srg}, il suffit de considérer le système de racines réduit de $G$ relativement à $T$.

\begin{itemize}
\item Supposons d'abord $G$ de type $A_n$. D'après \cite{bou}, le groupe $P/X$ est cyclique d'ordre $r$ divisant $n+1$, et on a:
\[\phi\left({\alpha_1}^\vee\right)=\left(2,-1,0,0,\dots,0\right);\]
\[\phi\left({\alpha_2}^\vee\right)=\left(-1,2,-1,0,\dots,0\right);\]
\ldots
\[\phi\left({\alpha_{n-1}}^\vee\right)=\left(0,\dots,0,-1,2,-1\right);\]
\[\phi\left({\alpha_n}^\vee\right)=\left(0,\dots,0,0,-1,2\right).\]
On voit donc que pour tout $i$, $\phi\left({\alpha_i}^\vee\right)$ appartient au sous-groupe des é\-lé\-ments $\left(a_1,\dots,a_n\right)$ de ${\mth{Z}}^n$ vérifiant:
\[a_1+2a_2+\dots+na_n\in(n+1){\mth{Z}}.\]
Ce sous-groupe étant un sous-groupe d'ordre $n+1$ de ${\mth{Z}}^n$, ils l'engendrent; on en déduit que $\phi\left(X_*\left(T\right)\right)$ est un sous-groupe d'ordre $r$ de ${\mth{Z}}^n$ le contenant, et le seul sous-groupe vérifiant ces propriétés est le sous-groupe des éléments $\left(a_1,\dots,a_n\right)$ de ${\mth{Z}}^n$ vérifiant:
\[a_1+2a_2+\dots+na_n\in r{\mth{Z}}.\]
On a $\left(1,\dots,1\right)\in\phi\left(X_*\left(T\right)\right)$ si et seulement si $r$ divise $\dfrac{n\left(n+1\right)}2$, ce qui est toujours le cas si $r$ est impair; si $r$ est pair, c'est vrai si et seulement si $r$ divise $\dfrac{n+1}2$. Dans ce cas, puisque d'après \cite{bou}, pour tout $i$, $c_i=1$, on a $q_\mathcal{T}=n+1$.

Supposons donc que $r$ est pair et ne divise pas $\dfrac{n+1}2$. Alors $n$ est forcément impair, et l'élément $\left(1,\dots,1,2,1,\dots,1\right)$, où le terme d'indice $\dfrac{n+1}2$ vaut $2$, appartient à $X_*\left(T\right)$: en effet, on a:
\begin{equation*}\begin{split}1+2+\dots+\dfrac{n-1}2+\left(n+1\right)+\dfrac{n+3}2+\dots+n&=\dfrac{n\left(n+1\right)}2+\dfrac{n+1}2\\&=\dfrac{n+1}2\left(n+1\right),\end{split}\end{equation*}
qui est un multiple entier de $n+1$, donc de $r$. On a alors $q_\mathcal{T}=n+2$.

\item Supposons maintenant $G$ de type $B_n$, $n\geq 2$. D'après \cite{bou}, le groupe $P/X$ est d'ordre $1$ ou $2$ ($2$ si $G$ est simplement connexe, $1$ si $G$ est adjoint), et on a:
\[\phi\left({\alpha_1}^\vee\right)=\left(2,-1,0,0,\dots,0\right);\]
\[\phi\left({\alpha_2}^\vee\right)=\left(-1,2,-1,0,\dots,0\right);\]
\ldots
\[\phi\left({\alpha_{n-1}}^\vee\right)=\left(0,\dots,0,-1,2,-1\right);\]
\[\phi\left({\alpha_n}^\vee\right)=\left(0,\dots,0,0,-2,2\right).\]
Les $\phi\left({\alpha_i}^\vee\right)$ engendrent donc le sous-groupe des éléments $\left(a_1,\dots,a_n\right)$ de ${\mth{Z}}^n$ vérifiant:
\[a_1+a_3+a_5+\dots+a_m\in 2{\mth{Z}},\]
où $m=n-1$ si $n$ est impair et $m=n$ si $n$ est pair. Si $G$ est simplement connexe, $\phi\left(X_*\left(T\right)\right)$ est égal à ce groupe.

D'après \cite{bou}, on a $c_1=1$ et $c_i=2$ pour tout $i>1$; on en déduit que si $G$ est adjoint,
ou si $G$ est simplement connexe et $n$ est congru à $3$ ou $4$ modulo $4$, on a $\left(1,\dots,1\right)\in\phi\left(X_*\left(T\right)\right)$, d'où $q_\mathcal{T}=2n$. Si maintenant $G$ est simplement connexe et $n$ est congru à $1$ ou $2$ modulo $4$, $\phi\left(X_*\left(T\right)\right)$ ne contient pas $\left(1,\dots,1\right)$, mais contient $\left(2,1,\dots,1\right)$, d'où $q_\mathcal{T}=2n+1$.

\item Supposons maintenant $G$ de type $C_n$, $n\geq 2$. D'après \cite{bou}, le groupe $P/X$ est d'ordre $1$ ou $2$ ($2$ si $G$ est simplement connexe, $1$ si $G$ est adjoint), et on a:
\[\phi\left({\alpha_1}^\vee\right)=\left(2,-1,0,0,\dots,0\right);\]
\[\phi\left({\alpha_2}^\vee\right)=\left(-1,2,-1,0,\dots,0\right);\]
\ldots
\[\phi\left({\alpha_{n-1}}^\vee\right)=\left(0,\dots,0,-1,2,-2\right);\]
\[\phi\left({\alpha_n}^\vee\right)=\left(0,\dots,0,0,-1,2\right).\]
Les $\phi\left({\alpha_i}^\vee\right)$ engendrent donc le sous-groupe des éléments $\left(a_1,\dots,a_n\right)$ de ${\mth{Z}}^n$ vérifiant $a_n\in 2{\mth{Z}}$. Ce groupe ne contient jamais $\left(1,\dots,1\right)$, mais contient toujours $\left(1,\dots,1,2\right)$; comme d'après \cite{bou}, on a $c_i=2$ pour $i<n$ et $c_n=1$, on en déduit $q_\mathcal{T}=2n+1$ si $G$ est simplement connexe et $q_\mathcal{T}=2n$ si $G$ est adjoint.

\item Supposons maintenant $G$ de type $D_n$, $n\geq 4$. D'après \cite{bou}, le groupe $P/X$ est isomorphe à un sous-groupe de ${\mth{Z}}/4{\mth{Z}}$ si $n$ est impair (resp. $\left({\mth{Z}}/2{\mth{Z}}\right)^2$ si $n$ est pair), et on a:
\[\phi\left({\alpha_1}^\vee\right)=\left(2,-1,0,0,\dots,0\right);\]
\[\phi\left({\alpha_2}^\vee\right)=\left(-1,2,-1,0,\dots,0\right);\]
\ldots
\[\phi\left({\alpha_{n-2}}^\vee\right)=\left(0,\dots,0,-1,2,-1,-1\right);\]
\[\phi\left({\alpha_{n-1}}^\vee\right)=\left(0,\dots,0,0,-1,2,0\right);\]
\[\phi\left({\alpha_n}^\vee\right)=\left(0,\dots,0,0,-1,0,2\right).\]
Supposons d'abord $n$ impair. Les $\phi\left({\alpha_i}^\vee\right)$ engendrent le sous-groupe des é\-lé\-ments $\left(a_1,\dots,a_n\right)$ de ${\mth{Z}}^n$ vérifiant:
\[2a_1+2a_3+2a_5+\dots+2a_{n-2}+a_{n-1}-a_n\in 4{\mth{Z}}.\]
En posant $r=\left[P:X\right]$, le groupe $\phi\left(X_*\left(T\right)\right)$ est alors le sous-groupe des éléments $\left(a_1,\dots,a_n\right)$ de ${\mth{Z}}^n$ vérifiant:
\[2a_1+2a_3+2a_5+\dots+2a_{n-2}+a_{n-1}-a_n\in r{\mth{Z}}.\]
Si $r$ vaut $1$ ou $2$, ou si $r=4$ et $n$ est congru à $1$ modulo $4$, on a $\left(1,\dots,1\right)\in\phi\left(X_*\left(T\right)\right)$. Comme on a, d'après \cite{bou}, $c_1=c_{n-1}=c_n=1$ et $c_i=2$ pour tout $i$ différent des précédents, on en déduit que l'on a $q_\mathcal{T}=2n-2$. Si maintenant $r=4$ et $n$ est congru à $3$ modulo $4$, $\phi\left(X_*\left(T\right)\right)$ ne contient pas $\left(1,\dots,1\right)$ mais contient $\left(2,1,\dots,1\right)$, d'où $q_\mathcal{T}=2n-1$.

Supposons maintenant $n$ pair. Les $\phi\left({\alpha_i}^\vee\right)$ engendrent le sous-groupe des éléments $\left(a_1,\dots,a_n\right)$ de ${\mth{Z}}^n$ vérifiant les conditions suivantes:
\[(E_1):a_1+a_3+a_5+\dots+a_{n-1}\in 2{\mth{Z}};\]
\[(E_2):a_{n-1}+a_n\in 2{\mth{Z}}.\]
Le groupe $\phi\left(X_*\left(T\right)\right)$ est alors le groupe des éléments de ${\mth{Z}}^n$ vérifiant une partie de $\{(E_1),(E_2)\}$ de cardinal $m$, avec $r=2^m$. Si $r=\left[P:X\right]=1$, ou si $r$ vaut $2$ ou $4$ et si $n$ est multiple de $4$, $\left(1,\dots,1\right)$ vérifie les mêmes conditions et on a $q_\mathcal{T}=2n-2$. Si maintenant $r$ vaut $2$ ou $4$ et $n$ est congru à $2$ modulo $4$, $\left(1,\dots,1\right)$ vérifie toujours $(E_2)$, mais plus $(E_1)$, et $\left(2,1,\dots,1\right)$ vérifie les deux; on a donc $q_\mathcal{T}=2n-1$ sauf dans le cas où $r=2$ et $\phi\left(X_*\left(T\right)\right)$ est le sous-groupe des éléments de ${\mth{Z}}^n$ vérifiant $(E_2)$, auquel cas on a $q_\mathcal{T}=2n-2$.

\item Supposons maintenant $G$ de type $E_6$. D'après \cite{bou}, le groupe $P/X$ est d'ordre $1$ ou $3$ ($3$ si $G$ est simplement connexe, $1$ si $G$ est adjoint) et on a:
\[\phi\left({\alpha_1}^\vee\right)=\left(2,0,-1,0,0,0\right);\]
\[\phi\left({\alpha_2}^\vee\right)=\left(0,2,0,-1,0,0\right);\]
\[\phi\left({\alpha_3}^\vee\right)=\left(-1,0,2,-1,0,0\right);\]
\[\phi\left({\alpha_4}^\vee\right)=\left(0,-1,-1,2,-1,0\right);\]
\[\phi\left({\alpha_5}^\vee\right)=\left(0,0,0,-1,2,-1\right);\]
\[\phi\left({\alpha_6}^\vee\right)=\left(0,0,0,0,-1,2\right).\]
Les $\phi\left({\alpha_i}^\vee\right)$ engendrent le sous-groupe des éléments $\left(a_1,\dots,a_6\right)$ de ${\mth{Z}}^6$ vérifiant:
\[a_1-a_3+a_5-a_6\in 3{\mth{Z}}.\]
L'élément $\left(1,\dots,1\right)$ vérifie toujours cette condition, et on déduit de \cite{bou} que l'on a $q_\mathcal{T}=12$.

\item Supposons maintenant $G$ de type $E_7$. D'après \cite{bou}, le groupe $P/X$ est d'ordre $1$ ou $2$ ($2$ si $G$ est simplement connexe, $1$ si $G$ est adjoint) et on a:
\[\phi\left({\alpha_1}^\vee\right)=\left(2,0,-1,0,0,0,0\right);\]
\[\phi\left({\alpha_2}^\vee\right)=\left(0,2,0,-1,0,0,0\right);\]
\[\phi\left({\alpha_3}^\vee\right)=\left(-1,0,2,-1,0,0,0\right);\]
\[\phi\left({\alpha_4}^\vee\right)=\left(0,-1,-1,2,-1,0,0\right);\]
\[\phi\left({\alpha_5}^\vee\right)=\left(0,0,0,-1,2,-1,0\right);\]
\[\phi\left({\alpha_6}^\vee\right)=\left(0,0,0,0,-1,2,-1\right);\]
\[\phi\left({\alpha_7}^\vee\right)=\left(0,0,0,0,0,-1,2\right).\]
Les $\phi\left({\alpha_i}^\vee\right)$ engendrent donc le sous-groupe des éléments $\left(a_1,\dots,a_7\right)$ de ${\mth{Z}}^7$ vérifiant:
\[a_2+a_5+a_7\in 2{\mth{Z}}.\]
L'élément $\left(1,\dots,1\right)$ de ${\mth{Z}}^7$ ne vérifie pas cette condition, mais $\left(1,\dots,1,2\right)$ la vérifie. On déduit donc de \cite{bou} que l'on a $q_\mathcal{T}=18$ si $G$ est adjoint et $q_\mathcal{T}=19$ si $G$ est simplement connexe.

\item Supposons maintenant $G$ de type $E_8$, $F_4$ ou $G_2$. Dans tous ces cas, $P/X$ est trivial et on a donc toujours $\left(1,\dots,1\right)\in\phi\left(X_*\left(T\right)\right)$. On déduit donc de \cite{bou} que l'on a:
\begin{itemize}
\item $q_\mathcal{T}=30$ si $G$ est de type $E_8$;
\item $q_\mathcal{T}=12$ si $G$ est de type $F_4$;
\item $q_\mathcal{T}=6$ si $G$ est de type $G_2$;
\end{itemize}
\end{itemize}

On a ainsi, dans tous les cas, montré le résultat suivant:

\begin{prop}
Supposons que $q\geq h+1$, où $h$ est le nombre de Coxeter de $G_{nr}$. Alors tout tore non ramifié de $G$ admet des éléments de réduction régulière.

Si $G$ est adjoint, c'est vrai également si $q=h$.
\end{prop}

\section{Intégrales orbitales unipotentes}

\subsection{Orbites unipotentes et sous-groupes parahoriques}

Considérons maintenant les distributions intégrales orbitales unipotentes sur $G$. Pour tout $u\in G$ unipotent et tout $f\in C_c^\infty\left(G\right)$, posons, si $U$ est l'orbite de $u$ dans $G$:
\[J_U\left(f\right)=\int_{u\in Z_G^0\left(u\right)\backslash G}f\left(h^{-1}uh\right)dh;\]
cette intégrale converge, et la distribution $J\left(u,\cdot\right)$ est invariante et à support dans l'ensemble des éléments unipotents de $G$, qui est contenu dans l'ensemble des éléments compacts de $G^1$; c'est donc un élément de $\mathcal{D}_{c,1}$. Le but de cette partie est de montrer que, moyennant certaines conditions sur $G$, les restrictions à $\mathcal{H}$ des distributions intégrales orbitales unipotentes engendrent $\mathcal{D}_{c,1}|_\mathcal{H}$ tout entier.

Plus précisément, soit $F_{nr}$ l'extension non ramifiée maximale de $F$, et soit $G_{nr}=\underline{G}\left(F_{nr}\right)$; on supposera que $p$ est un "bon" nombre premier pour $G_{nr}$, c'est-à-dire que si $\Phi_{nr}$ est le système de racines de $G_{nr}$ relativement à un tore déployé maximal et $\Delta=\left\{\alpha_1,\dots,\alpha_n\right\}$ un système de racines simples de $\Phi_{nr}$, pour tout $\beta=\sum_ic_i\alpha_i\in\Phi_{nr}$, si $c_i\neq 0$, $c_i$ n'est pas multiple de $p$. Plus précisément, on aura:
\begin{itemize}
\item si $G_{nr}$ est de type $A_n$, $p$ quelconque;
\item si $G_{nr}$ est de type $B_n$, $C_n$, $D_n$ ou $BC_n$, $p\neq 2$;
\item si $G_{nr}$ est de type $E_6$, $E_7$, $F_4$ ou $G_2$, $p\neq 2,3$;
\item si $G_{nr}$ est de type $E_8$, $p\neq 2,3,5$;
\item si le diagramme de Dynkin de $G_{nr}$ est non connexe, $p$ est bon pour chacune des composantes.
\end{itemize}

Si $p$ est bon pour $G_{nr}$, il l'est également pour $\underline{G}$; il existe alors, d'après \cite{beal}, une application de Springer sur $\underline{G}$, c'est-à-dire un morphisme bijectif $\psi$ de variétés entre la variété des éléments unipotents de $\underline{G}$ et celle des éléments nilpotents de $Lie\left(\underline{G}\right)$, vérifiant les propriétés suivantes:
\begin{itemize}
\item pour tout $u\in\underline{G}$ unipotent et tout $g\in\underline{G}$, $\psi\left(g^{-1}ug\right)=Ad\left(g^{-1}\right)\psi\left(u\right)$;
\item l'image par $\psi$ de chaque orbite unipotente de $\underline{G}$ est une orbite nilpotente de $Lie\left(\underline{G}\right)$ toute entière.
\end{itemize}
On voit aisément que $\psi$ est entièrement déterminée par sa valeur en un élément $u$ de l'unique orbite unipotente de dimension maximale de $\underline{G}$, et que si on a à la fois $u\in G$ et $\psi\left(u\right)\in Lie\left(G\right)$, $\psi$ est fixée par $\Gal\left(\overline{F}/F\right)$; dans ce dernier cas, $\psi$ induit un morphisme bijectif de variétés entre la variété des éléments unipotents de $G$ et celle des éléments nilpotents de $Lie\left(G\right)$, que l'on appellera également application de Springer. On a une assertion similaire en remplaçant $F$ et $G$ par $F_{nr}$ et $G_{nr}$.

Soit $u$ un élément unipotent de $G$, et soit $K$ un sous-groupe parahorique de $G$ contenant $u$; l'image $\overline{u}$ de $u$ dans ${\mth{G}}=K/K^1$ est un élément unipotent de ${\mth{G}}$. Si $f$ est une fonction sur $G$ à support dans $K$ et biinvariante par $K^1$, on peut donc écrire, en identifiant $f$ à une fonction sur ${\mth{G}}$:
\begin{equation}\label{intu}
J\left(u,f\right)=\sum_{h\in G_u}\dfrac{vol\left(Z_G^0\left(u\right)\backslash Z_G^0\left(u\right)hK\right)}{vol\left(K\right)}1_{\overline{h^{-1}uh}}\left(f\right),
\end{equation}
où $G_u$ est un système de représentants des doubles classes de $G$ modulo $Z_G^0\left(u\right)$ à gauche et $K$ à droite, et où, si $u'$ est un élément unipotent de ${\mth{G}}$, $1_{u'}\left(f\right)$ est la somme des valeurs de $f$ sur l'orbite de $u'$. Grâce au théorème \ref{dh}, l'étude des distributions intégrales orbitales unipotentes sur $G$ peut donc se faire à l'aide des intégrales orbitales unipotentes sur les groupes finis.

Soit $\mu\in\mathcal{M}$, $\left(M,K_\mu\right)$ un représentant de $\mu$, ${\mth{M}}=K_\mu/K_\mu^1$ et $m$ un élément de ${\mth{M}}$ n'appartenant à aucun sous-groupe de Levi de ${\mth{M}}$; on appellera {\em $G$-orbite} (anisotrope) de $m$ dans ${\mth{M}}$ l'ensemble des conjugués de $m$ par un élément de $N_G\left(K_\mu\right)$; une telle orbite ne dépend pas du choix de $\left(M,K_\mu\right)$.
Si $m$ appartient à un sous-groupe de Levi ${\mth{M}}'$ de ${\mth{M}}$,
on appellera {\em $G$-orbite} de $m$ dans ${\mth{M}}$ la réunion des orbites de ${\mth{M}}$ rencontrant la $G$-orbite de $m$ dans ${\mth{M}}'$.

Dane ce qui suit, on va associer à chaque couple $\left(\mu,O\right)$, où $\mu\in\mathcal{M}$ et $O$ est une $G$-orbite unipotente anisotrope du groupe ${\mth{M}}$ associé à $\mu$, une orbite unipotente de $G$. On notera $\Lambda$ l'ensemble des couples $\left(\mu,O\right)$.

Considérons l'immeuble $\mathcal{B}_{nr}$ du groupe $G_{nr}=\underline{G}\left(F_{nr}\right)$. On a une action du groupe de Galois $\Gal\left(F_{nr}/F\right)$ sur cet immeuble, et $\mathcal{B}$ s'identifie à l'ensemble des points de $\mathcal{B}_{nr}$ fixes pour cette action. De plus, l'injection $\mathcal{B}\mapsto\mathcal{B}_{nr}$ préserve les facettes, c'est-à-dire que pour toute facette $A$ de $\mathcal{B}$, il existe une facette $A'$ de $\mathcal{B}_{nr}$ telle que $A=A'\cap\mathcal{B}$; $A'$ est alors stable par $\Gal\left(F_{nr}/F\right)$. Réciproquement, par le théorème du point fixe de Bruhat-Tits (\cite[I. 3.2.4]{bt}), toute facette de $\mathcal{B}_{nr}$ stable par $\Gal\left(F_{nr}/F\right)$ contient un élément de $\mathcal{B}$; on a donc une bijection canonique entre l'ensemble des facettes de $\mathcal{B}$ et l'ensemble des facettes de $\mathcal{B}_{nr}$ stables par $\Gal\left(F_{nr}/F\right)$. On en déduit une bijection canonique entre l'ensemble des parahoriques de $G$ et l'ensemble des parahoriques de $G_{nr}$ stables par $\Gal\left(F_{nr}/F\right)$: si $K$ est un parahorique de $G$ d'image $K_{nr}$, on a $K=K_{nr}\cap G$.

Soit $K_{nr}$ un sous-groupe parahorique de $G_{nr}$ stable par $\Gal\left(F_{nr}/F\right)$; son premier sous-groupe de congruence $K_{nr}^1$ est également stable par $\Gal\left(F_{nr}/F\right)$, et le quotient $K_{nr}/K_{nr}^1$ est un groupe réductif $\underline{\mth{G}}$ sur $\overline{{\mth{F}}_q}$ stable par $\Gal\left(F_{nr}/F\right)\simeq \Gal\left(\overline{{\mth{F}}_q}/{\mth{F}}_q\right)$, donc défini sur ${\mth{F}}_q$; si $K=K_{nr}\cap G$, le quotient $K/K^1$ s'identifie au groupe ${\mth{G}}$ des ${\mth{F}}_q$-points de $\underline{\mth{G}}$.

Soit ${\mathbf F}$ le Frobenius de $\Gal\left(F_{nr}/F\right)$, qui s'identifie à celui de $\Gal\left(\overline{{\mth{F}}_q}/{\mth{F}}_q\right)$; on peut définir sur $G_{nr}$ l'application de Lang de la même façon que sur un groupe défini sur un corps fini, par:
\[L: g\longmapsto g^{-1}{\mathbf F}\left(g\right).\]
On a le lemme suivant:

\begin{lemme}\label{apls}
Soit $H_{nr}$ le groupe des $F_{nr}$-points d'un groupe algébrique défini sur $F$, et $K_{nr}$ un sous-groupe ouvert de $H_{nr}$ stable par $\Gal\left(F_{nr}/F\right)$, contenu dans un sous-groupe parahorique $K_{0,nr}$ de $H_{nr}$ et tel que le quotient $K_{nr}/\left(K_{0,nr}^1\cap K_{nr}\right)$, en tant que groupe algébrique sur $\overline{\mth{F}}_q$, est connexe; la restriction à $K_{nr}$ de l'application de Lang est surjective.
\end{lemme}

\begin{proof}
Pour tout entier $s>0$, définissons le $s$-ème groupe de con\-gru\-ence $K_{0,nr}^s$ de $K_{0,nr}$ de la même façon que pour les parahoriques de $G$. Fixons un entier $s$ tel que $K_{0,nr}^s\subset K_{nr}$; un tel $s$ existe car les $K_{0,nr}^s$ constituent une base de voisinage de l'unité; de plus, $K_{0,nr}^s$ est normal dans $K_{nr}$. Considérons le groupe $\underline{\mth{H}}=K_{nr}/K_{0,nr}^s$; c'est un groupe algébrique connexe de dimension finie défini sur ${\mth{F}}_q$. L'application de Lang sur un tel groupe étant surjective, on en déduit que pour tout $g\in K_{nr}$, il existe $h_0\in K_{nr}$ tel que $g\in h_0^{-1}{\mathbf F}\left(h_0\right)K_{0,nr}^s=h_0^{-1}K_{0,nr}^s{\mathbf F}\left(h_0\right)$. Posons maintenant:
\[g_s=h_0g{\mathbf F}\left(h_0\right)^{-1}\in K_{0,nr}^s\]
et considérons le groupe $\underline{\mth{G}}_s=K_{0,nr}^s/K_{0,nr}^{s+1}$. C'est encore un groupe algébrique connexe de dimension finie définie sur ${\mth{F}}_q$; l'application de Lang y est donc également surjective, et on en déduit l'existence de $h_s\in K_{0,nr}^s$ tel que $g_{s+1}=h_sg_s{\mathbf F}\left(h_s\right)^{-1}\in K_{0,nr}^{s+1}$. Par une récurrence évidente, on obtient, pour tout $i\geq s$, un $h_i\in K_{0,nr}^i$ et un $g_{i+1}=h_ig_i{\mathbf F}\left(g_i\right)^{-1}\in K_{0,nr}^{i+1}$. Considérons maintenant l'élément:
\[h=\ldots h_i\ldots h_{s+2}h_{s+1}h_sh_0\in K_{nr};\]
le produit converge puisque pour tout $i\geq s$, $h_i\in K_{0,nr}^i$, et l'on a:
\[hg{\mathbf F}\left(h\right)^{-1}\in\bigcap_{i}K_{0,nr}^i=\left\{1\right\},\]
d'où $g=L\left(h\right)$, ce qui démontre le lemme.
\end{proof}

Soit $O$ une orbite unipotente anisotrope de ${\mth{G}}$, et soit $O'$ l'orbite de $\underline{\mth{G}}$ contenant $O$; $O'$ n'est pas toujours anisotrope, mais elle est toujours stable par $\Gal\left(F_{nr}/F\right)$.

Soit $\underline{\mth{P}}=\underline{\mth{M}}\underline{\mth{U}}$ un sous-groupe parabolique de $\underline{\mth{G}}$, et soit $\underline{\mth{M}}_{der}$ (resp. $\underline{\mth{U}}_{der}$) le groupe dérivé de $\underline{\mth{M}}$ (resp. $\underline{\mth{U}}$); on dit que $\underline{\mth{P}}$ est distingué si $dim\left(\underline{\mth{M}}_{der}\right)=dim\left(\underline{\mth{U}}\right)-dim\left(\underline{\mth{U}}_{der}\right)$.

Soit $\underline{\mth{T}}$ un tore maximal de $\underline{\mth{G}}$, et $\underline{\mth{B}}$ un sous-groupe de Borel de $\underline{\mth{G}}$ contenant $\underline{\mth{T}}$, que l'on supposera stables par $\Gal\left(F_{nr}/F\right)$. Soit $\Phi_{K_{nr}}$ le système de racines de $\underline{\mth{G}}$ relativement à $\underline{\mth{T}}$, et $\Phi_{K_{nr}}^+$ l'ensemble des racines de $\Phi_{K_{nr}}$ positives relativement à $\underline{\mth{B}}$. Puisque $p$ est bon pour $G_{nr}$, il l'est également pour $\underline{\mth{G}}$, et il existe alors un sous-groupe parabolique $\underline{\mth{P}}=\underline{\mth{M}}\underline{\mth{U}}$ de $\underline{\mth{G}}$ contenant $\underline{\mth{B}}$, un sous-groupe parabolique distingué $\underline{{\mth{P}}'}=\underline{{\mth{M}}'}\underline{{\mth{U}}'}$ de $\underline{\mth{M}}$ contenant $\underline{\mth{B}}\cap\underline{\mth{M}}$ et un élément $\overline{u}$ de $O'$ tels que $\overline{u}$ est un élément de $\underline{{\mth{U}}'}$ appartenant à l'orbite de Richardson de $\underline{{\mth{P}}'}$; c'est une conséquence de \cite[corollaire au théorème I.3.5 et théorème II.1.5]{pomm} lorsque $\underline{\mth{G}}$ est semi-simple, et dans le cas général, c'est vrai dans le groupe semi-simple $\underline{\mth{G}}_{der}$, donc a fortiori dans $\underline{\mth{G}}$.

Soit $\Phi_{\underline{\mth{G}}}=\Phi_{K_{nr}}$ le système de racines de $\underline{\mth{G}}$ relativement à $\underline{\mth{T}}$, et $\Phi_{K_{nr}}^+$ l'ensemble des racines positives de $\Phi_{K_{nr}}$ relativement à $\underline{\mth{B}}$.

Soit $\overline{\mathcal{U}}$ (resp. $\overline{\mathcal{N}}$) la variété des éléments unipotents (resp. nilpotents) de $\underline{\mth{G}}$ (resp. $Lie\left(\underline{\mth{G}}\right)$), et $\overline{\psi}$ une application de Springer de $\overline{\mathcal{U}}$ dans $\overline{\mathcal{N}}$, que l'on supposera invariante par $\Gal\left(\overline{\mth{F}}_q/{\mth{F}}_q\right)$; on a le lemme suivant:

\begin{lemme}\label{sli}
On peut supposer $\overline{u}$ tel que $\overline{n}=\overline{\psi}\left(\overline{u}\right)$ est de la forme:
\[\overline{n}=\sum_{\alpha\in\Phi_n}\overline{n}_\alpha,\]
où $\Phi_n$ est une partie de $\Phi_{K_{nr}}^+$ constituée d'éléments linéairement indépendants, et où pour tout $\alpha\in\Phi_n$, $n_\alpha$ est un élément de $Lie\left(\underline{\mth{U}}_\alpha\right)$, où $\underline{\mth{U}}_\alpha$ est le sous-groupe radiciel de $\underline{\mth{G}}$ correspondant à $\alpha$.
\end{lemme}

\begin{proof}
La démonstration est identique à celle de \cite[E.III.4.2.9]{beal}, sachant que la généralisation de \cite[E.III.4.2.8]{beal} au cas $p$ bon quelconque se déduit de \cite{pomm} et que l'on utilise la graduation suivante de $Lie\left(\underline{\mth{M}}\right)$: pour toute racine simple $\alpha_i$ de $\underline{\mth{M}}$ relativement à $\underline{\mth{T}}$ et $\underline{\mth{B}}\cap\underline{\mth{M}}$, on pose $c\left(\alpha_i\right)=0$ si $\underline{\mth{U}}_\alpha\subset\underline{{\mth{M}}'}$ et $c\left(\alpha_i\right)=2$ sinon; pour tout $\alpha=\sum_i\lambda_i\alpha_i$, on pose $c\left(\alpha\right)=\sum_i\lambda_ic\left(\alpha_i\right)$. On a alors:
\[Lie\left(\underline{\mth{M}}\right)=\bigoplus_{j\in 2{\mth{Z}}}Lie\left(\underline{\mth{M}}\right)_j,\]
avec:
\[Lie\left(\underline{\mth{M}}\right)_0=Lie\left(\underline{\mth{T}}\right)\oplus\bigoplus_{\alpha,c\left(\alpha\right)=0}Lie\left(\underline{\mth{U}}_\alpha\right)=Lie\left(\underline{{\mth{M}}'}\right),\]
et, pour tout $j\neq 0$:
\[Lie\left(\underline{\mth{M}}\right)_j=\bigoplus_{\alpha,c\left(\alpha\right)=j}Lie\left(\underline{\mth{U}}_\alpha\right).\]
Remarquons que si $p$ est assez grand, on retrouve ainsi la graduation associée à un ${\mathfrak{sl}}_2$-triplet contenant $\overline{n}$.
\end{proof}

D'autre part, on a le lemme suivant:

\begin{lemme}\label{ptor}
Soit $\overline{v}\in O'\cap{\mth{G}}$; les orbites unipotentes de ${\mth{G}}$ contenues dans $O'$ sont en bijection avec les éléments de ${\mth{G}}\backslash L^{-1}\left(Z_{\underline{\mth{G}}}\left(\overline{v}\right)\right)/Z_{\underline{\mth{G}}}\left(\overline{v}\right)$.
\end{lemme}

\begin{proof}
En effet, soit $\overline{v}'\in O'\cap{\mth{G}}$, et soit $g\in\underline{\mth{G}}$ tel que $g^{-1}\overline{v}'g=\overline{v}$. On a alors également:
\[{\mathbf F}\left(g\right)^{-1}\overline{v}'{\mathbf F}\left(g\right)=\overline{v},\]
d'où:
\[L\left(g\right)=g^{-1}{\mathbf F}\left(g\right)\in Z_{\underline{\mth{G}}}\left(\overline{v}\right).\]
Réciproquement, si $g'$ est un élément de $\underline{\mth{G}}$ qui vérifie l'assertion ci-dessus, on a $g'\overline{v}g'{}^{-1}={\mathbf F}\left(g'\right)\overline{v}{\mathbf F}\left(g'{}^{-1}\right)$, d'où $g'\overline{v}g'{}^{-1}\in{\mth{G}}$.

D'autre part, il est clair que si $g'\in\underline{\mth{G}}$ est tel que $\overline{v}''=g'\overline{v}g'{}^{-1}\in{\mth{G}}$, $\overline{v}''$ appartient à la même ${\mth{G}}$-orbite que $\overline{v}'$ si et seulement si on a $g'\in{\mth{G}}gZ_{\underline{\mth{G}}}\left(\overline{v}\right)$, ce qui montre le lemme.
\end{proof}

On va maintenant associer à la $G_{nr}$-orbite contenant $O'$ une orbite unipotente de $G_{nr}$; pour cela, on aura besoin de quelques préliminaires. Soit $A_{nr}$ la facette de $\mathcal{B}_{nr}$ fixée par $K_{nr}$, $\mathcal{A}_{nr}$ un appartement de $\mathcal{B}_{nr}$ stable par $\Gal\left(F_{nr}\right)$ et contenant $A_{nr}$, et $T_{nr}$ le tore déployé maximal de $G_{nr}$ correspondant à $\mathcal{A}_{nr}$; $\Phi_{K_{nr}}$ s'identifie à un sous-système de racines réduit et stable par $\Gal\left(F_{nr}/F\right)$ du système de racines $\Phi_{nr}$ de $G_{nr}$ relativement à $T_{nr}$. Soit $H'_{nr}=\underline{H'}\left(F_{nr}\right)$ le sous-groupe de $G_{nr}$ engendré par $T_{nr}$ et les sous-groupes radiciels $U_{\alpha,nr}$, $\alpha\in\Phi_{K_{nr}}$; c'est un sous-groupe réductif fermé de $G_{nr}$ défini sur $F$.

Montrons qu'il existe un sous-groupe $\underline{H}$ de $\underline{H}'$ défini sur $F$, fermé, réductif et déployé sur $F_{nr}$ dont le système de racines relativement à un tore maximal s'identifie canoniquement à $\Phi_{K_{nr}}$:
pour tout $\alpha\in\Phi_{K_{nr}}$, considérons le sous-groupe radiciel $U_{\alpha,nr}$ de $H_{nr}$. D'après \cite[II. 4.1.9 et 4.1.16]{bt}, il existe une extension séparable finie $E_{nr}$ de $F_{nr}$ telle que $U_{\alpha,nr}$ est de l'une des deux formes suivantes:
\begin{itemize}
\item $U_{\alpha,nr}$ est isomorphe à $E_{nr}$;
\item il existe une extension quadratique séparable $E'_{nr}$ de $E_{nr}$ telle que $U_{\alpha,nr}$ est isomorphe au groupe des éléments unipotents triangulaires supérieurs du la forme déployée de $SU_3\left(E'_{nr}/E_{nr}\right)$ (définie comme l'ensemble des éléments $g$ de $GL_3\left(E'_{nr}\right)$ vérifiant ${}^\tau g\sigma(g)=1$, où ${}^\tau g$ est la transposée de $g$ par rapport à la diagonale non principale et $\sigma$ l'élément non trivial de $\Gal\left(E'_{nr}/E_{nr}\right)$.
\end{itemize}
Dans les deux cas, soit $U_{0,\alpha,nr}$ le groupe des points de $U_{\alpha,nr}$ dont l'image est à coordonnées dans $F_{nr}$; c'est un sous-groupe fermé de dimension $1$ de $U_{\alpha,nr}$, et c'est le groupe des $F_{nr}$-points d'un sous-groupe fermé $\underline{U_{0,\alpha}}$ de $\underline{H'}$ de dimension 1; de plus, on déduit des relations de commutation entre éléments unipotents de $G_{nr}$ (cf. \cite[annexe]{bt}) que pour tous $\alpha,\beta\in\Phi_{K_{nr}}$, $[U_{0,\alpha,nr},U_{0,\beta,nr}]$ est inclus dans le produit des $U_{0,\gamma,nr}$, $\gamma\in\Phi_{K_{nr}}$. Considérons donc, en posant $T_{nr}=\underline{T}(F_{nr}$, le groupe $\underline{H}$ engendré par $\underline{T}$ et les $\underline{U_{0,\alpha}}$, $\alpha\in\Phi_{K_{nr}}$: d'après ce qui précède, il vérifie les conditions voulues.

On notera $H$ (resp. $H_{nr}$) l'ensemble de ses $F$-points (resp $F_{nr}$-points).
Il est clair que si $\underline{H'}$ est lui-même déployé sur $F_{nr}$, on a $\underline{H}=\underline{H'}$.

Soit $\mathcal{U}_{nr}$ (resp. $\mathcal{N}_{nr}$) la variété des éléments unipotents (resp. nilpotents) de $H_{nr}$ (resp. $Lie\left(H_{nr}\right)$); on a le lemme suivant:

\begin{lemme}
Il existe une application de Springer $\psi_{nr}$ de $\mathcal{U}_{nr}$ dans $\mathcal{N}_{nr}$ telle que pour tout $v'\in K_{nr}\cap\mathcal{U}_{nr}$ d'image $\overline{v'}$ dans $\underline{\mth{G}}$, l'image de $\psi_{nr}\left(v'\right)$ dans $Lie\left(\underline{\mth{G}}\right)$ est $\overline{\psi}\left(\overline{v'}\right)$.
\end{lemme}

\begin{proof}
En effet, soit $v$ un élément de $K_{nr}$ et $\overline{v}$ son image, que l'on supposera appartenir à l'orbite unipotente de dimension maximale de $\underline{\mth{G}}$, et soit $\psi'_{nr}$ une application de Springer quelconque de $\mathcal{U}_{nr}$ dans $\mathcal{N}_{nr}$. Soit $X=\psi_{nr}\left(v\right)$; on montre comme dans le lemme \ref{sli} qu'il existe $h\in\underline{H}$ tel que $Ad\left(h^{-1}\right)X$ est un élément de $Lie\left(\underline{H}\right)$ de la forme:
\[Ad\left(h^{-1}\right)X=\sum_{\alpha\in\Delta_{nr}}n_\alpha,\]
où $\Delta_{nr}$ est un ensemble de racines simples de $\Phi_{K_{nr}}$ et où pour tout $\alpha$, $n_\alpha$ appartient à $Lie\left(\underline{U}_\alpha\right)$, où $\underline{U}_\alpha$ est le sous-groupe radiciel de $\underline{H}$ correspondant à $\alpha$. Quitte à multiplier $h$ à droite par un élément de $\underline{T}$, on peut même supposer que l'on a $Ad\left(h^{-1}\right)X\in Lie\left(H_{nr}\right)$ et que pour tout $\alpha$, $n_\alpha$ appartient à $Lie\left(U_{\alpha,nr}\cap K_{nr}\right)$ et est d'image non nulle dans $\underline{\mth{G}}$. L'image de $Ad\left(h^{-1}\right)X$ dans $Lie\left(\underline{\mth{G}}\right)$ est alors un élément $\overline{X'}$ de l'orbite nilpotente de dimension maximale de $Lie\left(\underline{\mth{G}}\right)$, et il existe alors un élément $\overline{h'}\in\overline{\mth{G}}$ tel que $Ad\left(\overline{h'}^{-1}\right)\overline{X'}=\overline{v}$. Soit $h'$ un élément de $K_{nr}$ d'image $\overline{h'}$ dans $\underline{\mth{G}}$, et soit $\psi_{nr}$ l'unique application de Springer de $\mathcal{U}_{nr}$ dans $\mathcal{N}_{nr}$ telle que l'on a:
\[\psi_{nr}:v\longmapsto Ad\left(h'{}^{-1}h^{-1}\right)X;\]
l'image de $\psi_{nr}\left(v\right)$ dans $\underline{\mth{G}}$ est alors $\overline{\psi}\left(\overline{v}\right)$, et pour tout $v'\in K_{nr}$ dont l'image $\overline{v'}$ appartient à l'orbite unipotente de dimension maximale de $\underline{\mth{G}}$, puisque $v'$ est conjugué à $v$ dans $\underline{H}$, l'image de $\psi_{nr}\left(v'\right)$ dans $\underline{\mth{G}}$ vaut $\overline{\psi}\left(\overline{v'}\right)$.

Supposons maintenant que $v'\in\mathcal{U}_{nr}\cap K_{nr}$ est tel que son image $\overline{v'}$ dans $\underline{\mth{G}}$ n'appartient pas à l'orbite unipotente de dimension maximale de $\underline{\mth{G}}$. Quitte à le conjuguer, on peut supposer que $v'$ appartient à $U_{0,nr}=\prod_{\alpha\in\Phi_{K_{nr}}^+}U_{\alpha,nr}$; on trouve alors facilement une droite $D$ de $Lie\left(U_{0,nr}\right)$ contenant $\psi_{nr}\left(v'\right)$, telle qu'il existe un isomorphisme de $D$ dans $F_{nr}$ qui envoie $D\cap Lie\left(K_{nr}\right)$ dans l'anneau des entiers $\mathcal{O}_{nr}$ de $F_{nr}$ et que les images dans $Lie\left(\underline{\mth{G}}\right)$ de tous les éléments de $D\cap Lie\left(K_{nr}\right)$ correspondant à des éléments de $\mathcal{O}_{nr}^*$ (resp. l'idéal maximal $\mathfrak{p}_{nr}$ de $\mathcal{O}_{nr}$) sont contenues dans l'orbite nilpotente de dimension maximale de $Lie\left(\underline{\mth{G}}\right)$ (resp. sont égales à $\overline{\psi}\left(\overline{v}\right)$): si $\alpha_1,\dots,\alpha_k$ sont les racines simples de $\Phi_{K_{nr}}^+$ n'intervenant pas dans l'écriture de $\overline{\psi}\left(\overline{v}\right)$ comme somme d'éléments des sous-groupes radiciels associés aux éléments de $\Phi_{K_{nr}}^+$, et si pour tout $i$, $n_i$ est un élément de $Lie\left(U_{\alpha_i}\cap H_{nr}\cap K_{nr}\right)$ d'image non nulle dane $\underline{\mth{G}}$, la droite:
\[D=\{\psi_{nr}(v)+z\sum_{i=1}^kn_i\mid z\in F_{nr}\}\]
convient. Soit $V=\psi_{nr}^{-1}\left(D\right)$: c'est une sous-variété fermée de dimension $1$ de $U_{0,nr}$, et d'après ce qui précède, si $\overline{V}$ (resp $\overline{v'}$) est l'image de $V\cap K_{nr}$ (resp. $v'$) dans $\underline{\mth{G}}$, l'image de $D\cap Lie\left(K_{nr}\right)$ dans $Lie\left(\underline{\mth{G}}\right)$ est la réunion de $\overline{\psi}\left(\overline{V}\right)-\left\{\overline{\psi}\left(\overline{v'}\right)\right\}$ et d'un singleton; or elle est fermée, donc elle contient également $\overline{\psi}\left(\overline{v'}\right)$; puisque l'image de $\psi_{nr}\left(v'\right)$ n'appartient pas à l'orbite nilpotente de dimension maximale de $Lie\left(\underline{\mth{G}}\right)$, on en déduit qu'elle appartient au singleton, donc qu'elle est égale à $\overline{\psi}\left(\overline{v'}\right)$, ce qui achève la démonstration.
\end{proof}

Soit $n$ un élément de $Lie\left(H_{nr}\right)$ de la forme:
\[n=\prod_{\alpha\in\Phi_n}n_\alpha,\]
où pour tout $\alpha\in\Phi_n$, $n_\alpha$ est un élément de $Lie\left(U_{\alpha,nr}\cap K_{nr}\right)$ dont l'image dans $Lie\left(\underline{\mth{U}}\right)$ est $\overline{n}_\alpha$. Un tel $n$ est clairement nilpotent.
Soit $u=\psi_{nr}^{-1}\left(n\right)$, et soit $O_{nr}$ l'orbite dans le groupe $H_{nr}\cap K_{nr}$.

Soit $\underline{P}_0$ le sous-groupe parabolique minimal de $\underline{H}$ engendré par $\underline{T}$ et les sous-groupes radiciels $\underline{U}_\alpha$, $\alpha\in\Phi_{K_{nr}}^+$;
soit $\underline{P}=\underline{M}\underline{U}$ l'unique sous-groupe parabolique de $\underline{H}$ standard (relativement à $\underline{P}_0$ et $\underline{T}$) tel que l'image de $\underline{P}\cap K_{nr}$ dans $\underline{\mth{G}}$ est $\underline{\mth{P}}=\underline{\mth{M}}\underline{\mth{U}}$,
et $\underline{P'}=\underline{M'}\underline{U'}$ le sous-groupe parabolique de $\underline{M}$ standard (relativement à $\underline{P}_0\cap \underline{M}$ et $\underline{T}$) tel que l'image de $\underline{P'}\cap K_{nr}$ dans $\underline{\mth{G}}$ est $\underline{{\mth{P}}'}=\underline{{\mth{M}}'}\underline{{\mth{U}}'}$; il est clair d'après la classification de \cite{bc} (valable pour tout $p$ bon par \cite{pomm}) que $\underline{P'}$ est un sous-groupe parabolique distingué de $\underline{M}$; d'autre part, l'image de $O_{nr}\cap U'_{nr}$ dans $\underline{\mth{G}}$ est l'intersection de $\underline{\mth{U}}'$ et de l'orbite de Richardson de $\underline{{\mth{P}}'}$; on en déduit que $O_{nr}\cap U'_{nr}$ est Zariski-dense dans $U'_{nr}$, et donc qu'elle est contenue dans l'orbite de Richardson de $\underline{P'}$. On notera $P_{0,nr}$ (resp. $P_{nr}$, $M_{nr}$, $U_{nr}$, $P'_{nr}$, $M'_{nr}$, $U'_{nr}$) le groupe des $F_{nr}$-points de $\underline{P}_0$ (resp. $\underline{P}$, $\underline{M}$, $\underline{U}$, $\underline{P'}$, $\underline{M'}$, $\underline{U'}$).

L'orbite $O_{nr}$ ainsi définie n'est pas unique et dépend du choix de $n$: on notera $\mathcal{U}_{nr}$ l'ensemble des $O_{nr}$ possibles, et $\mathcal{U}_{nr}^0$ l'ensemble des éléments de $\mathcal{U}_{nr}$ stables par $\Gal\left(F_{nr}/F\right)$.
Ce dernier ensemble n'est pas vide: en effet, il est clair que $P_{0,nr}$ et $T_{nr}$ sont stables par $\Gal\left(F_{nr}/F\right)$. Soit $\sigma\in \Gal\left(F_{nr}/F\right)$. Puisque $O'$ rencontre ${\mth{G}}$, $\overline{u}$ et $\sigma\left(\overline{u}\right)$ sont conjugués;
d'après \cite[théorème II.1.7 (ii)]{pomm}, on en déduit que $\underline{\mth{M}}$ et $\sigma\left(\underline{\mth{M}}\right)$ sont conjugués. De plus, si $w$ est un élément du groupe de Weyl $W\left(\underline{\mth{G}}/\underline{\mth{T}}\right)$ tel que $w^{-1}\underline{\mth{M}}w=\sigma\left(\underline{\mth{M}}\right)$ (il en existe puisque $\underline{\mth{M}}$ et $\sigma\left(\underline{\mth{M}}\right)$ sont semi-standard relativement à $\underline{\mth{T}}$ et conjugués entre eux) et $w^{-1}\left(\underline{\mth{M}}\cap\underline{\mth{B}}\right)w=\sigma\left(\underline{\mth{M}}\cap\underline{\mth{B}}\right)$, on a $w^{-1}\underline{{\mth{P}}'}w=\sigma\left(\underline{{\mth{P}}'}\right)$; en effet, les systèmes de racines de $\underline{{\mth{M}}'}$ et de $w\sigma\left(\underline{{\mth{M}}'}\right)w^{-1}$ sont les mêmes à un automorphisme du diagramme de Dynkin de $\underline{\mth{M}}$ près, et on vérifie immédiatement, en se ramenant au cas où le diagramme de Dynkin de $\underline{\mth{M}}$ est connexe et en considérant successivement chaque cas de la classification de \cite{bc}, que ${\overline{\mth{P}}'}$ est toujours invariant par un tel automorphisme. Dans ce qui suit, on supposera $O_{nr}\in\mathcal{U}_{nr}^0$.

On a alors le lemme suivant:

\begin{lemme}\label{ptor2}
Soit $v$ un élément de $O_{nr}\cap K\cap H$; les orbites unipotentes de $K\cap H$ contenues dans $O_{nr}$ sont en bijection avec les éléments de l'ensemble:
\[\left(K\cap H\right)\backslash L^{-1}\left(Z_{K_{nr}\cap H_{nr}}\left(v\right)\right)/Z_{K_{nr}\cap H_{nr}}\left(v\right).\]
\end{lemme}

\begin{proof}
La démonstration est identique à celle du lemme \ref{ptor}.
\end{proof}

\begin{cor}\label{bcorb}
Il existe une bijection canonique entre les orbites unipotentes de $K\cap H$ contenues dans $O_{nr}$ et les orbites unipotentes de ${\mth{G}}$ contenues dans $O'$, qui à l'orbite de $v\in O_{nr}\cap K\cap H$ associe celle de sa classe dans ${\mth{G}}$.
\end{cor}

\begin{proof}
En effet, soit $v\in O_{nr}\cap K\cap H$, et $\overline{v}$ sa classe dans ${\mth{H}}$; on déduit de la classification de \cite[13.1]{car} appliquée successivement à $\underline{H}/Z_{\underline{H}}$ et à $\underline{\mth{G}}/Z_{\underline{\mth{G}}}$ que l'ensemble des classes dans $\underline{\mth{G}}$ des éléments de $Z_{K_{nr}\cap H_{nr}}\left(v\right)$ est exactement $Z_{\underline{\mth{G}}}\left(\overline{v}\right)$. De plus, si $g\in L^{-1}\left(Z_{K_{nr}\cap H_{nr}}\left(v\right)\right)$, son image dans $\underline{\mth{G}}$ appartient à $L^{-1}\left(Z_{\underline{\mth{G}}}\left(\overline{v}\right)\right)$; d'autre part, posons:
\[C_v=Z_{K_{nr}\cap H_{nr}}\left(v\right)/Z_{K_{nr}\cap H_{nr}}\left(v\right)^0.\]
Puisque $v\in H$, $\Gal\left(F_{nr}/F\right)$ agit sur $C_v$, et l'application de Lang sur $K_{nr}$ induit une bijection de $K\cap H\backslash L^{-1}\left(Z_{K_{nr}\cap H_{nr}}\right)\left(v\right)/Z_{K_{nr}\cap H_{nr}}\left(v\right)$ dans l'ensemble des classes de $C_v$ pour la relation d'équivalence suivante: $x\equiv y$ si et seulement si il existe $z\in C_v$ tel que $x=z^{-1}y{\mathbf F}\left(z\right)$.
De même, en définissant $C_{\overline{v}}$ de manière analogue à $C_v$, l'application de Lang sur $\underline{\mth{G}}$ induit une bijection de ${\mth{G}}\backslash L^{-1}\left(Z_{\underline{\mth{G}}}\right)\left(\overline{v}\right)/Z_{\underline{\mth{G}}}\left(\overline{v}\right)$ dans l'ensemble des classes de $C_{\overline{v}}$ pour la relation d'équivalence similaire. Or on déduit de \cite[13.1]{car} que dans tous les cas, $C_v$ et $C_{\overline{v}}$ sont canoniquement isomorphes et que $\Gal\left(F_{nr}/F\right)$ agit de la même manière sur l'un et sur l'autre; on en déduit une bijection canonique:
\[K\cap H\backslash L^{-1}\left(Z_{K_{nr}\cap H_{nr}}\left(v\right)\right)/Z_{K_{nr}\cap H_{nr}}\left(v\right)\longrightarrow {\mth{G}}\backslash L^{-1}\left(Z_{\underline{\mth{G}}}\right)\left(\overline{v}\right)/Z_{\underline{\mth{G}}}\left(\overline{v}\right),\]
qui à chaque double classe associe l'ensemble de ses images dans $\underline{\mth{G}}$. L'assertion du lemme se déduit alors immédiatement des lemmes \ref{ptor} et \ref{ptor2}.
\end{proof}

D'après ce corollaire, on peut associer à chaque orbite unipotente $O$ de ${\mth{G}}$ contenue dans $O'$ une orbite unipotente $U_O$ de $G$, qui est celle contenant l'orbite unipotente de $K\cap H$ antécédent de $O$ par la bijection en question.

Considérons la relation de préordre sur les orbites unipotentes de $G$ (resp. $G_{nr}$) définie de la façon suivante: $U\leq U'$ si $U$ est contenue dans l'adhérence de $U'$. On dira qu'un élément d'un ensemble d'orbites unipotentes de $G$ (resp. $G_{nr}$) est minimal s'il l'est pour cette relation.

On a le lemme suivant:

\begin{lemme}\label{ormin}
Les orbites de $G_{nr}$ contenant les éléments de $\mathcal{U}_{nr}$ sont minimales parmi celles rencontrant $uK_{nr}^1$.
\end{lemme}

\begin{proof}
En effet, supposons d'abord $G_{nr}$ déployé, et $K_{nr}$ maximal spé\-cial; on a alors $H_{nr}=G_{nr}$. Soit $v\in uK_{nr}^1$ unipotent; si $w\in Z_{K_{nr}}\left(v\right)$, son image dans $\underline{\mth{G}}$ est clairement dans $Z_{\underline{\mth{G}}}\left(u\right)$. Considérons le morphisme canonique $Z_{K_{nr}}\left(v\right)\rightarrow Z_{\underline{\mth{G}}}\left(\overline{u}\right)$. Il induit un morphisme canonique $Lie\left(Z_{K_{nr}}\left(v\right)^0\right)\rightarrow Lie\left(Z_{\underline{\mth{G}}}\left(\overline{u}\right)^0\right)$, qui se factorise en:
\[Lie\left(Z_{K_{nr}}\left(v\right)^0\right)\rightarrow\overline{Lie\left(Z_{K_{nr}}\left(v\right)^0\right)}\rightarrow Lie\left(Z_{\underline{\mth{G}}}\left(\overline{u}\right)^0\right).\]
Puisque $Lie\left(Z_{K_{nr}}\left(v\right)^0\right)$ est un ouvert de $Lie\left(Z_{G_{nr}}\left(v\right)^0\right)$ et puisque le second morphisme est injectif, on en déduit:
\[dim\left(Lie\left(Z_{G_{nr}}\left(v\right)^0\right)\right)=dim\left(Lie\left(Z_{K_{nr}}\left(v\right)^0\right)\right)\leq dim\left(Lie\left(Z_{\underline{\mth{G}}}\left(\overline{u}\right)^0\right)\right).\]
Or d'après \cite[13.1]{car}, on a l'égalité pour $v=u$; on en déduit:
\[dim\left(Z_{G_{nr}}\left(v\right)\right)\leq dim\left(Z_{G_{nr}}\left(u\right)\right),\]
ce qui démontre que l'orbite de $u$ est minimale parmi les orbites unipotentes rencontrant $uK_{nr}^1$.

Supposons ensuite que $K_{nr}$ n'est plus nécessairement maximal, mais est contenu dans un sous-groupe parahorique maximal hyperspécial $K_{0,nr}$ de $G_{nr}$, toujours avec $G_{nr}$ déployé. Pour montrer l'assertion du lemme, il faut alors montrer que les orbites de $G_{nr}$ contenant les éléments de $\mathcal{U}_{nr}$ sont exactement les orbites unipotentes minimales parmi celles rencontrant $u\underline{{\mth{U}}_0}K_{0,nr}^1$, où $\underline{{\mth{U}}_0}$ est le radical unipotent du sous-groupe parabolique de $\underline{{\mth{G}}_0}=K_{0,nr}/K_{0,nr}^1$ image de $K_{nr}$.
Cette assertion se déduit, par une récurrence évidente, du cas précédent et du lemme suivant:

\begin{lemme}\label{uzero}
L'orbite de $\overline{u}$ est l'unique orbite unipotente minimale parmi celles qui rencontrent $\overline{u}\underline{{\mth{U}}_0}$.
\end{lemme}

\begin{proof}
Par la correspondance de Springer, on voit que l'assertion du lemme est équivalente à l'assertion correspondante concernant $\overline{n}$. Soit donc $\overline{n'}\in\overline{n}+Lie\left(\underline{{\mth{U}}_0}\right)$; écrivons-le:
\[\overline{n'}=\sum_{\alpha\in\Phi_{K_{0,nr}}^+}n'_\alpha,\]
où $\Phi_{K_{0,nr}}^+$ est défini de manière analogue à $\Phi_{K_{nr}}^+$.
Soit $\Delta_0=\left\{\alpha_1,\dots,\alpha_n\right\}$ l'ensemble des racines simples de $\Phi_{K_{0,nr}}^+$, les indices étant choisis de telle sorte que $\Delta_{nr}=\{\alpha_{r+1},\dots,\alpha_n\}$. Pour tout $\alpha\in\Phi_{K_{0,nr}}^+$ de la forme:
\[\alpha=\sum_{i=1}^nc_i\alpha_i,\]
posons:
\[c\left(\alpha\right)=\sum_{i=1}^rc_i.\]
Pour tout $x\in\overline{{\mth{F}}_q}$, posons:
\[\overline{n'_x}=\sum_{\alpha\in\Phi_{K_{0,nr}}^+}x^{c\left(\alpha\right)}n'_\alpha.\]
Il est clair que l'ensemble $\left\{\overline{n'_x}\mid x\in\overline{{\mth{F}}_q}\right\}$ est une sous-variété fermée de $\underline{\mth{G}}$, que $\overline{n'_0}=\overline{n}$ et que pour tout $x\neq 0$, $\overline{n'_x}$ appartient à l'orbite de $\overline{n'}$; on en déduit que $\overline{n}$ appartient à l'adhérence de l'orbite de $\overline{n'}$; comme ceci est vrai pour tout $\overline{n'}$, l'orbite de $\overline{n}$ est l'unique orbite minimale parmi celles rencontrant $\overline{n}+Lie\left(\underline{{\mth{U}}_0}\right)$, ce qui démontre le lemme.
\end{proof}

Revenons à la démonstration du lemme \ref{ormin}.
Supposons maintenant $K_{nr}$ quelconque, toujours avec $G_{nr}$ déployé, et soit $K_{0,nr}$ un sous-groupe parahorique maximal de $G_{nr}$ contenant $K_{nr}$; on a le lemme suivant:

\begin{lemme}\label{extsp}
Il existe une extension finie $E_{nr}$ de $F_{nr}$ telle que si $K_{0,E_{nr}}$ est le sous-groupe parahorique maximal de $\underline{G}\left(E_{nr}\right)$ contenant $K_{0,nr}$, $K_{0,E_{nr}}$ est spécial.
\end{lemme}

\begin{proof}
En effet, soit $x_0$ un point spécial quelconque de l'appartement $\mathcal{A}_{nr}$, et soit $\left\{\alpha_1,\dots,\alpha_n\right\}$ un ensemble de racines simples de $G_{nr}$ relativement à $T_{nr}$. Pour tout $i$, considérons une forme affine $f_i$ sur $\mathcal{A}_{nr}$ définie de la façon suivante:
\begin{itemize}
\item si $V_{nr}$ est l'espace vectoriel associé à $\mathcal{A}_{nr}$ et $H_i$ l'hyperplan de $V_{nr}$ associé à $\alpha_i$, le noyau de $f_i$ est $x_0+H_i$;
\item les points où $f_i$ prend des valeurs entières sont exactement les hyperplans de la forme $x'+H_i$, où $x'$ est un sommet de $\mathcal{A}_{nr}$, qui sont réunion de facettes de $\mathcal{A}_{nr}$.
\end{itemize}
Une telle forme affine est unique au signe près; de plus, l'ensemble des points spéciaux de $\mathcal{A}_{nr}$ est exactement l'ensemble des points où toutes les $f_i$ prennent des valeurs entières. D'autre part, si $E_{nr}$ est une extension finie galoisienne de $F_{nr}$, soit $\mathcal{B}_{E_{nr}}$ l'immeuble de $\underline{G}\left(E_{nr}\right)$, et $\mathcal{A}_{E_{nr}}$ l'appartement de $\mathcal{B}_{E_{nr}}$ correspondant à l'unique tore déployé maximal de $\underline{G}\left(E_{nr}\right)$ contenant $T_{nr}$; l'appartement $\mathcal{A}_{nr}$ s'identifie à l'ensemble des points de $\mathcal{A}_{E_{nr}}$ fixés par $\Gal\left(E_{nr}/F_{nr}\right)$, et l'ensemble des points de $\mathcal{A}_{nr}$ correspondant à un sous-groupe parabolique maximal spécial de $\underline{G}_{E_{nr}}$ est exactement l'ensemble des points $x'$ tels que pour tout $i$, $f_i\left(x'\right)\in\dfrac 1{\left[E_{nr}:F_{nr}\right]}{\mth{Z}}$.
Pour montrer le lemme, il suffit donc de vérifier qu'il existe $E_{nr}$ vérifiant ces conditions, ce qui est toujours le cas si, en notant $x$ le sommet de $\mathcal{A}_{nr}$ fixé par $K_{0,nr}$, pour tout $i$, $f_i\left(x\right)$ est un rationnel dont le dénominateur n'est pas multiple de $p$.

Il est clair que l'on peut supposer que le diagramme de Dynkin de $G_{nr}$ est connexe. Si $G_{nr}$ est de type $A_n$, tout sous-groupe parahorique maximal de $G_{nr}$ est spécial et il suffit de poser $E_{nr}=F_{nr}$; on supposera donc que $G_{nr}$ n'est pas de type $A_n$. Considérons, pour tout $i$, l'hyperplan de $\mathcal{A}_{nr}$ symétrique de $x_0+H_i$ par rapport à $x$; une telle symétrie préserve la structure simpliciale de $\mathcal{A}_{nr}$, donc il est réunion de facettes; de plus, il contient le sommet $x''$ de $\mathcal{A}_{nr}$ symétrique de $x_0$ par rapport à $x$. Comme ceci est vrai pour tout $i$, $x''$ est spécial, donc pour tout $i$, $f_i\left(x''\right)\in{\mth{Z}}$; comme $x$ est le milieu du segment $\left[x'',x_0\right]$, on en déduit que pour tout $i$, $f_i\left(x\right)\in \dfrac 12{\mth{Z}}$. Or puisque $p$ est bon pour $G_{nr}$, il est différent de $2$; on en déduit l'assertion cherchée.
\end{proof}

Revenons à la démonstration du lemme \ref{ormin}.
Soit $K_{0,E_{nr}}^1$ le premier sous-groupe de congruence de $K_{0,E_{nr}}$;
on a:
\[K_{0,nr}^1=K_{0,E_{nr}}^1\cap G_{nr};\]
en effet, si $x$ est un point de l'immeuble $\mathcal{B}_{E_{nr}}$ de $G_{E,nr}$ dont $K_{0,E_{nr}}$ est le fixateur connexe, les deux membres de l'égalité ci-dessus sont égaux au groupe des éléments de $G_{nr}$ fixant un voisinage de $x$ dans $\mathcal{B}_{E_{nr}}$. L'assertion du lemme se déduit alors immédiatement du cas précédent.

Supposons enfin $G_{nr}$ quelconque, et soit $E'_{nr}$ une extension de $F_{nr}$ sur laquelle $\underline{G}$ est déployé; l'assertion est vraie pour $\underline{G}\left(E'_{nr}\right)$ d'après ce qui précède, et si $K_{E'_{nr}}$ est un sous-groupe parahorique de $\underline{G}\left(E'_{nr}\right)$ contenant $K_{nr}$, par un argument similaire à celui du cas précédent, on a:
\[K_{nr}^1=K_{E'_{nr}}^1\cap G_{nr};\]
elle est donc également vraie pour $G_{nr}$ par restriction.
\end{proof}

On en déduit immédiatement, par restriction, le corollaire suivant:

\begin{cor}\label{orminc}
Parmi les orbites unipotentes de $G$ rencontrant $uK^1$, les orbites rencontrant un élément de $\mathcal{U}_{nr}^0$ sont minimales.
\end{cor}

Fixons maintenant une orbite $U_O$ de $G$ rencontrant un élément de $\mathcal{U}_{nr}^0$; on a le lemme suivant:

\begin{lemme}\label{unica}
Soit $K_1$ un autre sous-groupe parahorique de $G$, ${\mth{G}}_1=K_1/K_1^1$, $O_1$ une orbite anisotrope de ${\mth{G}}_1$, et $U_{O_1}$ une orbite de $G$ associée à $O_1$ de la même façon que $U_O$ l'est à $O$. Si $U_O$ rencontre l'image réciproque de $O_1$ dans $K_1$ et est minimale dans celle-ci, alors $K$ et $K_1$ sont associés et $O$ et $O_1$ sont contenues dans la même $G$-orbite de ${\mth{G}}\equiv{\mth{G}}_1$.
\end{lemme}

\begin{proof}

Soit $v\in U_0\cap O_1K_1^1$; quitte à conjuguer $v$ et $K_1$ par un élément de $G$ convenable, on supposera $v\in H$.

Tout d'abord, montrons que l'on peut supposer que $K_1$ contient $K_{T_0}$. Soit $A_1$ la facette de $\mathcal{B}$ dont le fixateur connexe est $K_1$ et $\mathcal{A}'$ un appartement de $\mathcal{B}$ contenant à la fois $A$ et $A_1$; d'après \cite[2.5.8]{bt}, il existe $k\in G^1$ tel que $k\mathcal{A}=\mathcal{A}'$ et que $k$ fixe $\mathcal{A}\cap\mathcal{A'}$; en particulier, $k$ fixe $A$, donc $k\in K$. De plus, il est clair que si dans ce qui précède, on remplace $T_0$ par $kT_0k^{-1}$, on conjugue tous les sous-groupes radiciels, et donc également $H$, par $k$; on peut donc choisir $kuk^{-1}$ à la place de $u$, et l'orbite associée à $O$ est alors $kU_Ok^{-1}=U_O$; on peut donc, quitte à remplacer $T_0$ par $kT_0k^{-1}$, $H$ pat $kHk^{-1}$ et $v$ par $kvk^{-1}$, supposer $K_{T_0}\subset K_1$. Le groupe $K_{1,H}=K_1\cap H$ est alors un sous-groupe parahorique de $H$; de plus, ${\mth{H}}_1=K_{1,H}/K_{1,H}^1$ est un sous-groupe réductif fermé de ${\mth{G}}_1$ qui contient l'image de $v$, donc qui rencontre $O_1$, et $O_1\cap{\mth{H}}_1$ est réunion d'orbites unipotentes anisotropes de ${\mth{H}}_1$. De plus, $K_H=K\cap H$ est un sous-groupe parahorique maximal spécial de $H$.

Supposons $K_H$ et $K_{1,H}$ associés dans $H$; ils sont alors conjugués par un élément $h$ de $H$, et en particulier, $K_{1,H}$ est également maximal spécial. Soit $(L,K_L)=\kappa_0(K,T_0)$, $(L_1,K_{L_1})=\kappa_0(K_1,T_0)$; $L_1$ est alors un sous-groupe de Levi de $G$ contenant $H$, et de même rang semi-simple que $H$; en effet, si ce n'était pas le cas, $H$ serait contenu dans un sous-groupe de Levi propre de $L_1$, et $O_1$ rencontrerait alors un sous-groupe de Levi propre de ${\mth{G}}_1=L_1/L_1^1$, ce qui est exclu par hypothèse. D'autre part, par définition de $H$, $L$ est également un sous-groupe de Levi de $G$ contenant $H$ et de même rang semi-simple que $H$; on a donc $L=L_1$. De plus, $K_L$ (resp. $K_{L_1}$) est l'unique sous-groupe parahorique de $L$ contenant $K_H$ (resp. $K_{1,H}$); on en déduit que $K_L$ et $K_{L_1}$ sont conjugués par $h$, et donc que $K$ et $K_1$ sont associés. Enfin, il est clair que si $O$ et $O_1$ rencontrent une même $H$-orbite de ${\mth{G}}$, elles sont contenues dans la même $G$-orbite; il suffit donc de montrer le lemme dans le cas où $G=H$.

Supposons donc $G=H$; $K$ est alors un sous-groupe parahorique maximal spécial de $G$, et $U_O$ est une orbite unipotente anisotrope. On a le lemme suivant:

\begin{lemme}
Soit $U'$ une orbite unipotente de $G$ minimale parmi celles rencontrant $O_1K_1^1$; alors $U'$ rencontre $L_1\cap O_1K_1^1$.
\end{lemme}

\begin{proof}
Par une récurrence évidente, il suffit de considérer le cas où $L_1$ est propre maximal, sachant que si $L_1=G$ il n'y a rien à démontrer. Soit $V_1$ et $V_1^-$ les radicaux unipotents des deux sous-groupes paraboliques de $G$ semi-standard (relativement à $T_0$) admettant $L_1$ pour composante de Levi; on a la décomposition d'Iwahori:
\[K_1=\left(V_1^-\cap K_1\right)\left(L_1\cap K_1\right)\left(V_1\cap K_1\right).\]
Soit $u'\in U'\cap O_1K_1^1$; on va montrer qu'il existe $v\in V_1\cap K_1$ tel que $v^{-1}u'v\in L_1V_1^-$. En inversant les rôles de $V_1$ et $V_1^-$, on obtiendra ensuite un $v^-\in V_1^-\cap K_1$ tel que ${v^-}^{-1}v^{-1}u'vv^-\in L_1V_1\cap L_1V_1^-=L_1$, ce qui démontrera le lemme.

Soit $A_1$ l'arête de $\mathfrak{A}$ fixée par $K_1$, $x'_1$ le sommet de $A_1$ tel que si $K'_1$ est le sous-groupe parahorique maximal de $G$ fixant $x'_1$, on a $K'_1\cap L_1V_1=K_1\cap L_1V_1$, $A_2$ l'arête de $\mathcal{A}$ symétrique de $A_1$ par rapport à $x'_1$, et $K_2$ le fixateur connexe de $A_2$; puisque $A_1$ et $A_2$ engendrent la même droite $D$ de $\mathcal{A}$, $K_1$ et $K_2$ sont associés. On va montrer qu'il existe $v_1\in K_1\cap V_1$ tel que $v_1^{-1}u_1v_1\in O_1K_1^1\cap O_1K_2^1$; en considérant ensuite les arêtes $A_3,A_4,\cdots$ de $D$ telles que pour tout $i$, $A_i$ est le symétrique de $A_{i-1}$ par rapport au sommet $x'_{i-1}$ de $A_{i-1}$ distinct de $x'_{i-2}$, on obtient ainsi par récurrence pour tout $i$ un élément $v_i$ de $K_i\cap V_1$ tel que $v_i^{-1}\cdots v_1^{-1}u_1v_1\cdots v_i\in\bigcap_{j=1}^iO_1K_j^1$. De plus, la suite des $K_i\cap V_i$ est une base de voisinage de l'unité dans $V_i$; on en déduit que le produit infini des $v_i$ converge vers un élément $v$ de $K_1\cap V_1$ tel que $v^{-1}u_1v$ appartient à $\bigcap_{i=1}^\infty O_1K_i^1$, qui n'est autre que $O_1\left(K_1^1\cap L_1V_1^-\right)\subset L_1V_1^-$.

Montrons donc l'existence de $v_1$. Soit $\overline{u'}$ (resp ${\mth{P}}'_1$) l'image de $u'$ (resp. $K_1$) dans ${\mth{G}}'_1=K'_1/{K'_1}^1$; ${\mth{G}}_1$ s'identifie à une composante de Levi de ${\mth{P}}'_1$, et $\overline{u'}$ appartient à une orbite unipotente de ${\mth{G}}'_1$ minimale parmi celles rencontrant $O_1{\mth{P}}'_1$. Pour montrer l'existence de $v_1$, il suffit de montrer l'existence d'un élément $\overline{v_1}$ du radical unipotent ${\mth{U}}'_1$ de ${\mth{P}}'_1$ tel que $\overline{v_1}^{-1}\overline{u'}\overline{v_1}\in O_1$; puisque $K'_1\cap L_1V_1=K_1\cap L_1V_1$, l'image de $V_1$ dans ${\mth{G}}'_1$ est égale à ${\mth{U}}'_1$, et il existe donc un représentant $v_1$ de $\overline{v_1}$ appartenant à $V_1$, qui vérifie la condition cherchée.

Posons ${\mth{G}}'_1=\underline{\mth{G}}'_1({\mth{F}}_q)$, ${\mth{G}}_1=\underline{\mth{G}}_1({\mth{F}}_q)$, ${\mth{U}}'_1=\underline{\mth{U}}'_1({\mth{F}}_q)$. D'après le lemme \ref{uzero} appliqué à $\underline{\mth{P}}'_1$, par minimalité de l'orbite de $\overline{u'}$, il existe $v'_1\in\underline{\mth{P}}'_1$ tel que ${v'_1}^{-1}\overline{u'}v'_1\in\underline{\mth{G}}_1$, et il est clair que l'on peut supposer $v'_1\in\underline{\mth{U}}'_1$; puisque la composante dans $\underline{\mth{G}}_1$ de ${v'_1}^{-1}\overline{u'}v'_1$ est la même que celle de $\overline{u'}$, on a ${v'_1}^{-1}\overline{u'}v'_1\in{\mth{G}}_1$; si ${\mathbf F}$ est le Frobenius de $\Gal\left(\overline{\mth{F}}_q/{\mth{F}}_q\right)$, on a donc ${\mathbf F}(v'_1)^{-1}\overline{u'}{\mathbf F}(v'_1)={v'_1}^{-1}\overline{u'}v'_1$. On en déduit:
\[v'_1{\mathbf F}(v'_1)^{-1}\in Z_{\underline{\mth{G}}'_1}(\overline{u'})\cap\underline{\mth{U}}'_1.\]
Or d'après \cite[13.1]{car}, le groupe $Z_{\underline{\mth{G}}'_1}(\overline{u'})\cap\underline{\mth{U}}'_1$ est toujours connexe; l'application de Lang y est surjective, et on obtient $v''_1\in Z_{\underline{\mth{G}}'_1}(\overline{u'})\cap\underline{\mth{U}}'_1$ tel que l'on a:
\[v'_1{\mathbf F}(v'_1)^{-1}={v''_1}^{-1}{\mathbf F}(v''_1),\]
c'est-à-dire $v''_1v'_1\in{\mth{G}}'_1$. En posant $\overline{v_1}=v''_1v'_1$, on obtient $\overline{v_1}^{-1}\overline{u'}\overline{v_1}\in O_1$, ce qui est l'assertion cherchée.
\end{proof}

Revenons à la démonstration du lemme \ref{unica}. D'après le lemme précédent, $U_O$ rencontre $L_1$; or $U_O$ est anisotrope, donc $L_1=G$ et $K_1$ est maximal.
Soit $x$ (resp. $x_1$) le sommet de $\mathcal{B}$ fixé par $K$ (resp. $K_1$); on va montrer qu'il existe $g\in G$ tel que $gx_1=x$, c'est-à-dire que $K$ et $K_1$ sont conjugués.

D'après la définition de $U_O$, $v$ est conjugué à $u$ par un élément de $K_{nr}$; on déduit donc de la définition de $u$ qu'il existe un appartement $\mathcal{A}'_{nr}$ de $\mathcal{B}_{nr}$ stable par $\Gal\left(F_{nr}/F\right)$ et tel que l'ensemble des points de $\mathcal{A'}_{nr}$ fixés par $v$ est l'intersection de $t=\card\left(\Phi_n\right)$ racines affines $\alpha'_1,\dots,\alpha'_t$ telles que pour tout $i$, la facette $A_{1,nr}$ de $\mathcal{B}$ contenant $x_1$ appartient à la frontière de $\alpha'_i$.
Soit $\mathfrak{C}_{nr}$ un quartier de $\mathcal{A}'_{nr}$ stable par $\Gal\left(F_{nr}/F\right)$ contenu dans cette intersection, et soit $\mathcal{A}''_{nr}$ un appartement de $\mathcal{B}_{nr}$ contenant à la fois $A_{1,nr}$ et un sous-quartier $\mathfrak{C}'_{nr}$ de $\mathfrak{C}_{nr}$ que l'on pourra, quitte à le réduire, supposer également stable par $\Gal\left(F_{nr}/F\right)$; si $P'_{0,nr}=M'_{0,nr}U'_{0,nr}$ est le sous-groupe parabolique minimal de $G_{nr}$ associé à la direction de $\mathfrak{C}_{nr}$, il existe $v\in U'_{0,nr}$ tel que $v\mathcal{A}''_{nr}=\mathcal{A}'_{nr}$; de plus, on déduit des relations de commutation entre éléments unipotents de $G_{nr}$ (cf. \cite{chev} et \cite[annexe]{bt}) que $vA_{1,nr}$ est situé sur l'intersection des $\alpha'_i$; enfin, on a le lemme suivant:

\begin{lemme}
On peut choisir $\mathcal{A}''_{nr}$ et $v$ tels que $v\in G$ et $\mathcal{A''}_{nr}$ est stable par $\Gal\left(F_{nr}/F\right)$.
\end{lemme}

\begin{proof}
En effet, puisque $A_{1,nr}$ et $\mathfrak{C}'_{nr}$ sont stables par $\Gal\left(F_{nr}/F\right)$, l'appartement ${\mathbf F}(\mathcal{A}''_{nr})$, où $\mathbf{F}$ est le Frobenius de $\Gal\left(F_{nr}/F\right)$, les contient également. De plus, on a le lemme suivant:

\begin{lemme}
L'application de Lang sur $U'_{0,nr}$ est surjective.
\end{lemme}

\begin{proof}
En effet, pour tout $y\in\mathfrak{C'}_{nr}$ stable par $\Gal\left(F_{nr}/F\right)$, elle est surjective sur $K_y\cap U'_{0,nr}$ par le lemme \ref{apls}, et la réunion de ces sous-groupes est $U'_{0,nr}$ tout entier.
\end{proof}

On en déduit qu'il existe $v'\in U'_{0,nr}$ tel que ${v'}^{-1}{\mathbf F}(v'){\mathbf F}(\mathcal{A}''_{nr})=\mathcal{A}''_{nr}$; l'appartement $v'\mathcal{A}''_{nr}$ est alors stable par $\Gal\left(F_{nr}/F\right)$. De plus, on a:
\[\mathbf{F}\left(v'v\right)\mathcal{A}''_{nr}=v'v\mathcal{A}''_{nr}=\mathcal{A}'_{nr},\]
donc $v^{-1}{v'}^{-1}\mathbf{F}\left(v'v\right)$ appartient à l'intersection de $U'_{0,nr}$ et de $N_{G_{nr}}\left(T'_{0,nr}\right)$, où $T'_{0,nr}$ est le tore maximal associé à $\mathcal{A}'_{nr}$. Or une telle intersection est forcément réduite à l'identité: en effet ses éléments sont des éléments unipotents de $N_{G_{nr}}\left(T'_{0,nr}\right)$ qui fixent un quartier de $\mathcal{A}'_{nr}$, donc $\mathcal{A}'_{nr}$ tout entier; ce sont donc des éléments de $Z_{G_{nr}}\left(T'_{0,nr}\right)$, qui est un tore et ne contient pas d'éléments unipotents autres que $1$. On en déduit $v'v\in G$. L'assertion du lemme est donc vraie en prenant $v'\mathcal{A}''_{nr}$ et $v'v$ à la place de $\mathcal{A}''_{nr}$ et $v$.
\end{proof}

On en déduit que $vx_1$ est un point de $\mathcal{A}'_{nr}$ fixé par $\Gal\left(F_{nr}/F\right)$ et situé sur le sous-espace affine $E$ de $\mathcal{A}'_{nr}$ engendré par $A_{nr}$. Or le seul point vérifiant ces propriétés est $x$: en effet, s'il existe un autre point $y$ de $E$ fixé par $\Gal\left(F_{nr}/F\right)$, ce groupe fixe point par point la droite de $E$ engendrée par $x$ et $y$; comme l'intersection de $A_{nr}$ et d'une telle droite est un segment ouvert, elle n'est pas réduite à $x$ et on obtient une contradiction. On en déduit que l'on a $vx_1=x$ et donc que $K$ et $K_1$ sont conjugués. Quitte à remplacer $K$ par $vK$, on peut donc supposer $K=K_1$; le fait que $O$ et $O_1$ sont contenues dans une même $G$-orbite de ${\mth{G}}$ se déduit immédiatement du corollaire \ref{bcorb}.

\end{proof}

Soit $\Lambda$ l'ensemble des couples $\left(\mu,O\right)$, où $\mu\in\mathcal{M}$ et $O$ est une $G$-orbite unipotente anisotrope du groupe fini ${\mth{M}}$ associé à $\mu$; on a ainsi montré le résultat suivant:

\begin{prop}\label{unipsj}
Il existe une injection $\phi$ de $\Lambda$ dans l'ensemble des orbites unipotentes de $G$, qui possède la propriété suivante: pour tout couple $\left(\mu,O\right)$, si $K$ est un sous-groupe parahorique de $G$ dont l'élément de $\mathcal{M}$ associé est $\mu$, $U_O=\phi\left(\mu,O\right)$ est une orbite unipotente de $G$ minimale parmi celles rencontrant l'image réciproque de $O$ dans $K$; de plus, si $\left(M,K_\mu\right)$ est un représentant de $\mu$, $U_O$ rencontre $K_\mu$.
\end{prop}

Soit maintenant $\mu\in\mathcal{M}$, $\left(M,K_M\right)$ un représentant de $\mu$, ${\mth{G}}=K_M/K_M^1$, $O$ une $G$-orbite unipotente de ${\mth{G}}$ et $U=U_O=\phi\left(\mu,O\right)$. Si $O$ est anisotrope, on posera $r\left(O\right)=dim\left(U_O\right)$; si $O$ et quelconque, si ${\mth{M}}$ est un Levi de ${\mth{G}}$ et $O_{\mth{M}}$ une $G$-orbite unipotente anisotrope de ${\mth{M}}$ tels que $O_{\mth{M}}\subset O$, on posera $r\left(O\right)=r\left(O_{\mth{M}}\right)$. On a le résultat suivant:

\begin{lemme}\label{co}
Soit $\mu'\in\mathcal{M}$, $K'$ un sous-groupe parahorique de $G$ tel que $\kappa\left(K'\right)=\mu'$, ${\mth{G}}'=K'/K'{}^1$, et $\mathcal{U}'$ l'ensemble des $G$-orbites unipotentes de ${\mth{G}}'$. On a une égalité de la forme:
\[J_U|_{\mathcal{C}_{K',0}}=\sum_{O'\in\mathcal{U}'}c_{U,O'}1_{O'},\]
avec pour tout $O'$:
\begin{itemize}
\item si $K'$ est associé à un sous-groupe de $K$ et si $O'\subset O$, on a $c_{U,O'}\neq 0$;
\item si $O'$ ne vérifie pas la condition ci-dessus et si $r\left(O'\right)\geq r\left(O\right)$, on a $c_{U,O'}=0$.
\end{itemize}
\end{lemme}

\begin{proof}
En effet, on a par (\ref{intu}), pour tout $f\in\mathcal{C}_{K',0}$, si $u$ est un élément quelconque de $U$:
\[J_U\left(f\right)=\sum_{O'\in\mathcal{U}'}\left(\sum_{h\in G_u,\overline{h^{-1}uh}\in O'}\dfrac{vol\left(Z_G^0\left(u\right)\backslash Z_G^0\left(u\right)hK'\right)}{vol\left(K'\right)}\right)1_O\left(f\right),\]
d'où pour tout $O'\in\mathcal{U}$:
\[c_{U,O'}=\sum_{h\in G_u,\overline{h^{-1}uh}\in O'}\dfrac{vol\left(Z_G^0\left(u\right)\backslash Z_G^0\left(u\right)hK'\right)}{vol\left(K'\right)}.\]
Soit $O'\in\mathcal{U}'$ telle que $r\left(O'\right)\geq r\left(O\right)$; si l'image réciproque de $O'$ dans $K'$ ne rencontre pas $U$, il n'existe aucun $h\in G_u$ tel que $h^{-1}uk\in K'$ et $\overline{h^{-1}uh}\in O'$, et donc $c_{U,O'}=0$. Supposons donc que l'image réciproque de $O'$ dans $K'$ rencontre $U$.
Soit ${\mth{M}}$ un Levi de ${\mth{G}}'$ et $O'_{\mth{M}}$ une orbite unipotente anisotrope de ${\mth{G}}'$ tels que $O'_{\mth{M}}\subset O'$. Soit ${\mth{P}}$ un sous-groupe parabolique de ${\mth{G}}'$ de Levi ${\mth{M}}$, et $K'_{\mth{P}}$ le sous-groupe parahorique de $G$ image réciproque de ${\mth{P}}$ dans $K'$; d'après la proposition \ref{unipsj}, $U_{O'_{\mth{M}}}$ est minimale parmi les orbites unipotentes de $G$ rencontrant l'image réciproque d'un élément $u'$ quelconque de $O'_{\mth{M}}$ dans $K'_{\mth{P}}$; comme $K'{}^1\subset K'_{\mth{P}}{}^1$, a fortiori, $U_{O'_{\mth{M}}}$ est minimale parmi les orbites unipotentes de $G$ rencontrant $u'K'{}^1$, et on a donc $r\left(O'\right)=r\left(O'_{\mth{M}}\right)\leq r\left(O\right)$. De plus, par le lemme \ref{unica}, si $U$ est également minimale parmi les orbites unipotentes de $G$ rencontrant $u'K'{}^1$, $K$ et $K'_{\mth{P}}$ sont associés et, en identivfiant ${\mth{G}}$ à ${\mth{M}}$, on a $O'_{\mth{M}}=O$, ce qui démontre la deuxième assertion du lemme.

Montrons maintenant la première: supposons que $K'$ est associé à un sous-groupe de $K$ et que $O'\subset O$.
Alors la somme définissant $c_{U,O'}$ comporte au moins un terme non nul; puisque ces termes sont tous positifs, on en déduit que $c_{U,O'}\neq 0$.
\end{proof}

\begin{cor}\label{ioli}
Les intégrales orbitales $J_{U_O}$ sont linéairement in\-dé\-pen\-dan\-tes sur $\mathcal{C}_0$.
\end{cor}

\begin{proof}
En effet, supposons qu'il existe des constantes $\lambda_{O'}$, $O'\in\Lambda$, telles que $\sum_{O'}\lambda_{O'}J_{U,O'}|_{\mathcal{C}_0}=0$; on va montrer par récurrence sur $r\left(O'\right)$ qu'elles sont toutes nulles. Fixons donc $\left(\mu',O'\right)\in\Lambda$, et supposons que pour tout $\left(\mu,O\right)\in\Lambda$ telle que $r\left(O\right)>r\left(O'\right)$, on a $\lambda_O=0$. Soit $K$ un sous-groupe parahorique de $G$ contenant $I$, et ${\mth{G}}=K/K^1$; écrivons, comme dans le lemme précédent, pour tout $\left(\mu,O\right)\in\Lambda$:
\[J_{U_O}|_{\mathcal{C}_{K,0}}=\sum_{O''}c_{U_O,O''}1_{O''},\]
où $O''$ décrit l'ensemble des orbites unipotentes de ${\mth{G}}$. On en déduit, en identifiant $1_{O'}$ à un élément de $\mathcal{C}_0$:
\[J_{U_O}\left(1_{O'}\right)=c_{U_O,O'},\]
d'où:
\[0=\sum_{O,r\left(O\right)\leq r\left(O'\right)}\lambda_OJ_{U_O}\left(1_{O'}\right)=\sum_{O,r\left(O\right)\leq r\left(O'\right)}\lambda_Oc_{U_O,O'}.\]
Or d'après le lemme précédent, si $r\left(O\right)\leq r\left(O'\right)$, $c_{U_O,O'}=0$ si $O\neq O'$; l'égalité précédente se réduit donc à:
\[\lambda_{O'}c_{U_{O'},O'}=0.\]
Comme, toujours par le lemme précédent, $c_{U_{O'},O'}\neq 0$, on a $\lambda_{O'}=0$, ce qui démontre le corollaire.
\end{proof}

\subsection{Quelques résultats sur les groupes finis}

\begin{lemme}\label{ssuf}
Soit ${\mth{G}}$ un groupe fini réductif et connexe, ${\mth{P}}_0$ un sous-groupe parabolique minimal de ${\mth{G}}$, $\mathcal{H}_{\mth{G}}$ l'algèbre de Hecke des fonctions sur ${\mth{G}}$ biinvariantes par ${\mth{P}}_0$. Soit ${\mth{T}}$ un tore maximal de ${\mth{G}}$, $g$ un élément régulier de $T$; la restriction à $\mathcal{H}_{\mth{G}}$ de l'intégrale orbitale $1_g$ est combinaison linéaire des restrictions des intégrales orbitales unipotentes sur ${\mth{G}}$.
\end{lemme}

\begin{proof}
Soit $\underline{\mth{G}}$ un groupe réductif connexe défini sur ${\mth{F}}_q$ et $\mathbf F$ une application de Frobenius de $\underline{\mth{G}}$ dans lui-même tels que ${\mth{G}}=\underline{\mth{G}}^{\mathbf F}$; pour tout tore maximal ${\mth{T}}$ de $G$, soit $\underline{\mth{T}}$ le tore maximal ${\mathbf F}$-stable de ${\mth{G}}$ tel que ${\mth{T}}=\underline{\mth{T}}^{\mathbf F}$. Si $g_{\mth{T}}$ est un élément régulier de ${\mth{T}}$, on a, d'après le lemme \ref{cddl}, pour tout $f\in\mathcal{H}_{\mth{G}}$:
\[1_{g_{\mth{T}}}\left(f\right)=\dfrac 1{\card\left({\mth{T}}\right)}\langle R_{\underline{\mth{T}},1},f\rangle _{\mth{G}},\]
où $R_{\underline{\mth{T}},1}$ est l'application de Deligne-Lusztig de ${\mth{G}}$ dans ${\mth{C}}$ définie dans la dé\-mons\-tra\-tion du théorème \ref{intorb}. Or d'après la démonstration du lemme \ref{cddl}, tous les $\langle R_{\underline{\mth{T}},\theta},f\rangle _{\mth{G}}$, $\theta\neq 1$, sont nuls, et on peut donc écrire:
\begin{equation*}\begin{split}\langle R_{\underline{\mth{T}},1},f\rangle _{\mth{G}}&=\sum_\theta\langle R_{\underline{\mth{T}},\theta},f\rangle _{\mth{G}}\\&=\dfrac 1{\card\left({\mth{G}}\right)}\sum_{g\in{\mth{G}}}f\left(g\right)\sum_\theta R_{\underline{\mth{T}},\theta}\left(g\right).\end{split}\end{equation*}
D'autre part, on a pour tout $g\in{\mth{G}}$, d'après la démonstration de \cite[théorème 7.5.1]{car}:
\[\sum_\theta R_{\underline{\mth{T}},\theta}\left(g\right)=\mathcal{L}\left(\left(g,1\right),\tilde{X}\right),\]
où, comme dans la démonstration du théorème \ref{intorb}, $\tilde{X}=L^{-1}\left(\underline{\mth{U}}\right)$, $\underline{\mth{U}}$ étant le radical unipotent d'un sous-groupe de Borel $\underline{\mth{B}}$ de $\underline{\mth{G}}$ contenant $\underline{\mth{T}}$, et où $\mathcal{L}\left(\left(g,1\right),\tilde{X}\right)$ est le nombre de Lefschetz; si $g=su$ est la décomposition de Jordan de $g$, ce nombre vérifie, d'après \cite[7.1.10]{car}:
\[\mathcal{L}\left(\left(g,1\right),\tilde{X}\right)=\mathcal{L}\left(\left(u,1\right),\tilde{X}^{\left(s,1\right)}\right),\]
où $\tilde{X}^{\left(s,1\right)}$ est l'ensemble des points de $\tilde{X}$ fixés par $\left(s,1\right)$, c'est-à-dire l'ensemble des éléments $h$ de $\tilde{X}\subset\underline{\mth{G}}$ vérifiant $sh=h$. Or cet ensemble est vide si $s\neq 1$, c'est-à-dire si $g$ n'est pas unipotent; on en déduit:
\[\langle R_{\underline{\mth{T}},1},f\rangle _{\mth{G}}=\sum_{g\in{\mth{G}},g\mbox{ unipotent}}f\left(g\right)\sum_\theta R_{\underline{\mth{T}},\theta}\left(g\right).\]
Puisque les $R_{\underline{\mth{T}},\theta}$ sont des fonctions sur ${\mth{G}}$ invariantes par conjugaison, si $\mathcal{U}$ est l'ensemble des orbites unipotentes de ${\mth{G}}$, on peut reformuler l'égalité ci-dessus de la manière suivante:
\[\langle R_{\underline{\mth{T}},1},f\rangle _{\mth{G}}=\sum_{u\in\mathcal{U}}\left(\sum_\theta R_{\underline{\mth{T}},\theta}\left(u\right)\right)1_u\left(f\right),\]
où pour tout $u\in\mathcal{U}$, $1_u$ est l'intégrale orbitale sur $u$, ce qui démontre le lemme.
\end{proof}

\begin{cor}\label{ssufc}
Si ${\mth{G}}$ est tel que tout tore maximal ${\mth{T}}$ de ${\mth{G}}$ admet au moins un élément régulier, et si aucune composante connexe du diagramme de Dynkin de ${\mth{G}}$ n'est de type $E_7$ ou $E_8$, alors les intégrales orbitales unipotentes engendrent l'espace des restrictions à $\mathcal{H}_{\mth{G}}$ des distributions sur ${\mth{G}}$.
\end{cor}

\begin{proof}
Ce corollaire se déduit immédiatement du lemme précédent et de la proposition \ref{intorbf}.
\end{proof}

Soit maintenant ${\mth{P}}={\mth{M}}{\mth{U}}$ un sous-groupe parabolique de ${\mth{G}}$ et $O'$ une $G$-orbite unipotente anisotrope de ${\mth{M}}$. On peut assimiler $1_{O'}$ à la fonction caractéristique de $O'$ dans ${\mth{M}}$; on peut alors définir la distribution $\Ind_{\mth{M}}^{\mth{G}}1_{O'}$ sur ${\mth{G}}$. On déduit de la définition de l'induite que l'on a:
\[\Ind_{\mth{M}}^{\mth{G}}1_{O'}=\sum_O\dfrac{\card\left(Z_{\mth{G}}\left(u\right)\right)}{\card\left({\mth{P}}\right)}\card\left(O'{\mth{U}}\cap O\right)1_O,\]
la somme portant sur les $G$-orbites unipotentes de ${\mth{G}}$; $u$ désigne un élément quelconque de $O$. De plus, on a le lemme suivant:

\begin{lemme}\label{infgf}
Avec les notations précédentes, pour tout $O$, $O'{\mth{U}}\cap O$ ne peut être non vide que si on a soit $O'\subset O$, soit $r\left(O'\right)>r\left(O\right)$.
\end{lemme}

\begin{proof}
En effet, soit $K_{\mth{P}}$ le sous-groupe parahorique de $G$ image ré\-ci\-pro\-que de ${\mth{P}}$ dans $K$, $\left(M,K_M\right)=\kappa_0\left(K_{\mth{P}}\right)$, $\mu'$ la classe de $\left(M,K\right)$ dans $\mathcal{M}'$ et $u$ un représentant dans $K_M$ de $U_{O'}$. Soit $P=MU$ l'unique sous-groupe parabolique de $G$ de Levi $M$ tel que l'image de $P\cap K$ dans ${\mth{G}}$ est ${\mth{P}}$; si $O'{\mth{U}}\cap O$ est non vide, $O$ admet un représentant dans $G$ de la forme $uv$, $v\in U$; un tel représentant est clairement unipotent, et si $O'$ n'est pas contenue dans $O$, l'orbite $V$ de $uv$ dans $G$ ne contient pas $u$. Or il existe $t\in N_G\left(K_M\right)$ tel que $t^{-1}\left(U\cap K\right)t\subset U\cap K^1$: si $T_0$ est un tore déployé maximal de $G$ contenu dans $M$, il suffit de choisir $t$ tel que pour toute racine $\alpha$ de $G$ relativement à $T_0$ telle que $U_\alpha\subset U$, $\alpha\left(t\right)$ est de valuation strictement positive. L'élément $t^{-1}uvt$ est alors un élément de $K$ dont l'image dans ${\mth{G}}$ est dans $O'$, et on déduit de la proposition \ref{unipsj} que $\phi\left(\mu',O'\right)\subset\overline{V}$. Puisque $O'$ n'est pas contenue dans $O$, on déduit du lemme \ref{unica} que $\phi\left(\mu',O'\right)$ n'est pas minimale dans $\overline{V}$, d'où $r\left(O'\right)>r\left(O\right)$.
\end{proof}

Soit maintenant $f$ une fonction sur ${\mth{G}}$; on dira que $f$ est {\em cuspidale} si la fonction centrale $f^{\mth{G}}=\dfrac 1{\card\left({\mth{G}}\right)}\sum_{g\in{\mth{G}}}Ad\left(g\right)f$ est cuspidale. Pour une telle fonction $f$, on a en particulier, pour tout sous-groupe de Levi ${\mth{M}}$ propre de ${\mth{G}}$ et toute orbite unipotente $O'$ de ${\mth{M}}$:
\begin{equation*}\begin{split}\left(\Ind_{\mth{M}}^{\mth{G}}1_{O'}\right)\left(f\right)&=\left(\Ind_{\mth{M}}^{\mth{G}}1_{O'}\right)\left(f^{\mth{G}}\right)\\&=\langle 1_{O'},r_{\mth{M}}^{\mth{G}}f^{\mth{G}}\rangle _{\mth{M}}\\&=0.\end{split}\end{equation*}
Fixons un sous-groupe parahorique minimal ${\mth{P}}_0={\mth{M}}_0{\mth{U}}_0$ de ${\mth{G}}$ et un tore déployé maximal ${\mth{T}}_0$ de ${\mth{G}}$. Soit $\mathcal{H}_{\mth{G}}$ l'algèbre de Hecke des fonctions sur ${\mth{G}}$ biinvariantes par ${\mth{P}}_0$; on notera $\mathcal{H}_{{\mth{G}},cusp}$ le sous-espace des éléments cuspidaux de $\mathcal{H}_{\mth{G}}$. 

Soit ${\mth{P}}={\mth{M}}{\mth{U}}$ un sous-groupe parabolique de ${\mth{G}}$ standard (relativement à ${\mth{P}}_0$ et ${\mth{T}}_0$); par restriction, $\mathcal{H}_{\mth{M}}$ s'identifie canoniquement à une sous-algèbre de $\mathcal{H}_{\mth{G}}$. On a le lemme suivant:

\begin{lemme}\label{dchk}
Soit $R$ un système de représentants des classes de conjugaison de sous-groupes de Levi standard de ${\mth{G}}$; on a la décomposition suivante:
\[\mathcal{H}_{\mth{G}}=\sum_{{\mth{M}}\in R}\mathcal{H}_{{\mth{M}},cusp},\]
et la somme est directe modulo $C_0\left({\mth{G}}\right)\cap\mathcal{H}_{\mth{G}}$.
\end{lemme}

\begin{proof}
En effet, soit $f\in\mathcal{H}_{\mth{G}}$. D'après la proposition \ref{dcc}, on peut écrire de manière unique, en posant pour tout ${\mth{M}}\in R$:
\[f^{\mth{G}}=\sum_{{\mth{M}}\in R}\Ind_{\mth{M}}^{\mth{G}}f_{\mth{M}},\]
où $f_{\mth{M}}$ est une fonction centrale et cuspidale sur ${\mth{M}}$; de plus, pour tout ${\mth{M}}$, si ${\mth{P}}={\mth{M}}{\mth{U}}$ est un parabolique de ${\mth{G}}$ contenant ${\mth{M}}$, d'après le lemme \ref{find}, $\Ind_{\mth{M}}^{\mth{G}}f_{\mth{M}}$ et $f_{\mth{M}}$ ne diffèrent que par un élément de $C_0\left({\mth{G}}\right)$. Pour montrer le lemme, il suffit donc de montrer que pour tout ${\mth{M}}$, $f_{\mth{M}}\in\mathcal{H}_{\mth{M}}+C_0\left({\mth{G}}\right)$.

Pour toute fonction $\phi$ sur ${\mth{G}}$ et tout ${\mth{M}}\in R$, on posera $r_{{\mth{M}}}^{\mth{G}}\phi=r_{{\mth{M}}}^{\mth{G}}\left(\phi^{\mth{G}}\right)$. Montrons le lemme suivant:

\begin{lemme}
Pour tout ${\mth{M}}\in R$, $r_{{\mth{M}}}^{\mth{G}}f\in\mathcal{H}_{\mth{M}}+C_0\left({\mth{M}}\right)$.
\end{lemme}

\begin{proof}
En effet, soit ${\mth{P}}_0={\mth{M}}_0{\mth{U}}_0$ un sous-groupe parabolique minimal de ${\mth{G}}$ tel que ${\mth{M}}_0\subset{\mth{M}}$, ${\mth{T}}_0$ le tore déployé maximal de ${\mth{G}}$ contenu dans ${\mth{M}}_0$, $W$ le groupe de Weyl de ${\mth{G}}$ relativement à ${\mth{T}}_0$, et $W'$ un système de représentants de $W$ modulo $W\cap{\mth{M}}$ à gauche; on a:
\begin{equation*}\begin{split}r_{{\mth{M}}}^{\mth{G}}f&=\dfrac 1{\card\left({\mth{G}}\right)\card\left({\mth{U}}\right)}\sum_{g\in{\mth{G}}}\sum_{u\in{\mth{U}}}f\left(g^{-1}\cdot ug\right)\\&=\sum_{w\in W}\dfrac{\card\left({\mth{P}}w{\mth{P}}_0\right)}{\card\left({\mth{M}}\right)\card\left({\mth{G}}\right)\card\left({\mth{U}}\right)}\sum_{m\in{\mth{M}}}\sum_{u\in{\mth{U}}}f\left(w^{-1}m^{-1}\cdot muw\right).\end{split}\end{equation*}
Cette dernière fonction est égale à un élément de $C_0\left({\mth{M}}\right)$ près à:
\[\sum_{w\in W}\dfrac{\card\left({\mth{P}}w{\mth{P}}_0\right)}{\card\left({\mth{G}}\right)\card\left({\mth{U}}\right)}\sum_{u\in{\mth{U}}}f\left(w^{-1}\cdot uw\right).\]
Considérons, pour tout $w$, la fonction $\sum_uf\left(w^{-1}\cdot uw\right)$. Elle est biinvariante par $w{\mth{P}}_0 w^{-1}\cap{\mth{M}}$; or ce groupe est un sous-groupe parabolique minimal de ${\mth{M}}$, et est donc conjugué dans ${\mth{M}}$ à ${\mth{P}}_0\cap{\mth{M}}$, et on en déduit l'assertion du lemme.
\end{proof}

Revenons à la démonstration du lemme \ref{dchk}; on va montrer, par récurrence sur le rang semi-simple de ${\mth{M}}$, que pour tout ${\mth{M}}$, $f_{\mth{M}}\in\mathcal{H}_{\mth{M}}+C_0\left({\mth{G}}\right)$. Fixons donc un ${\mth{M}}$, et supposons que pour tout ${\mth{M}}'\in R$ de rang strictement inférieur, on a $f_{{\mth{M}}'}\in\mathcal{H}_{{\mth{M}}'}\subset\mathcal{H}_{\mth{G}}$. D'après le lemme précédent, on a:
\[\sum_{{\mth{M}}'\in R}r_{{\mth{M}}}^{\mth{G}}\Ind_{{\mth{M}}'}^{\mth{G}}f_{{\mth{M}}'}\in\mathcal{H}_{{\mth{M}}}+C_0\left({\mth{M}}\right).\]
Or si $\phi\in C_0\left({\mth{G}}\right)$, $\phi^{\mth{G}}=0$, donc $r_{{\mth{M}}}^{\mth{G}}\phi=0$; on en déduit, en utilisant le lemme \ref{find}:
\[\sum_{{\mth{M}}'\in R}r_{\mth{M}}^{\mth{G}}f_{{\mth{M}}'}\in\mathcal{H}_{{\mth{M}}}+C_0\left({\mth{M}}\right);\]
d'autre part, si ${\mth{M}}'$ n'est pas contenu dans un conjugué de ${\mth{M}}$, le lemme \ref{vdx} montre que $r_{{\mth{M}}}^{\mth{G}}\Ind_{{\mth{M}}'}^{\mth{G}}f_{{\mth{M}}'}=0$.
Enfin, pour tout ${\mth{M}}'$ distinct de ${\mth{M}}$ et contenu dans un conjugué de ${\mth{M}}$, par hypothèse de récurrence, $f_{{\mth{M}}'}\in\mathcal{H}_{{\mth{M}}'}+C_0\left({\mth{G}}\right)\subset\mathcal{H}_{{\mth{M}}}+C_0\left({\mth{G}}\right)$ (cette inclusion est évidente si ${\mth{M}}'\subset {\mth{M}}$, et le lemme \ref{find} permet de ce ramener à ce cas), d'où $r_{\mth{M}}^{\mth{G}}f_{{\mth{M}}'}\in\mathcal{H}_{{\mth{M}}}+C_0\left({\mth{G}}\right)$ d'après le lemme précédent; on en déduit que $r_{\mth{M}}^{\mth{G}}\Ind_{\mth{M}}^{\mth{G}}f_{\mth{M}}\in\mathcal{H}_{\mth{M}}+C_0\left({\mth{G}}\right)$; enfin, on déduit immédiatement du lemme \ref{fdm} que l'on a:
\[r_{\mth{M}}^{\mth{G}}\Ind_{\mth{M}}^{\mth{G}}f_{\mth{M}}\in f_{\mth{M}}+C_0\left({\mth{G}}\right),\]
ce qui permet de conclure.
\end{proof}

\subsection{Le résultat principal}

On va maintenant démontrer le résultat principal de cette partie:

\begin{theo}\label{intun}
Supposons que $G$ vérifie (C1), que son diagramme de Dynkin ne contient aucune composante connexe de type $E_7$ ou $E_8$,et que $p$ est bon pour $G_{nr}$. Alors les restrictions à $\mathcal{H}$ des distributions intégrales orbitales unipotentes sur $G$ engendrent $\mathcal{D}_{c,1}|_\mathcal{H}$.
\end{theo}

\begin{proof}
Ici encore, grâce au théorème \ref{dh}, il suffit de considérer les restrictions des distributions au sous-espace $\mathcal{H}'_0$ des éléments de $\mathcal{H}$ à support dans la réunion des parahoriques de $G$ contenant $I$. On va montrer l'assertion suivante: si $f$ est un élément de $\mathcal{H}'_0$ tel que pour tout $\left(\mu,O\right)\in\Lambda$, $J_{U_O}\left(f\right)=0$, alors $f$ annule toute distribution invariante sur $G$; comme $\mathcal{H}'_0$ est de dimension finie, cela implique l'assertion du théorème.

Soit donc $D\in\mathcal{D}_{c,1}$. Pour tout sous-groupe parahorique $K$ de $G$ contenant $I$, si $\mathcal{H}_K$ est l'ensemble des éléments de $\mathcal{H}'_0$ à support dans $K$, on peut écrire, grâce au corollaire \ref{ssufc}:
\[D|_{\mathcal{H}_K}=\sum_{O\in\mathcal{U}'_K}\lambda'_O1_O,\]
où $\mathcal{U}'_K$ est l'ensemble des orbites unipotentes de ${\mth{G}}=K/K^1$. De plus, puisque $D$ est invariante par conjugaison, si $O$ et $O'$ sont deux orbites unipotentes contenues dans une même $G$-orbite, on peut supposer $\lambda'_O=\lambda'_{O'}$; on peut donc reformuler la somme ci-dessus de la façon suivante:
\[D|_{\mathcal{H}_K}=\sum_{O\in\mathcal{U}_K}\lambda_O1_O,\]
où $\mathcal{U}_K$ est l'ensemble des $G$-orbites unipotentes de ${\mth{G}}$.

Considérons l'ensemble des classes de $G$-conjugaison de couples $\left({\mth{M}},O'\right)$, où ${\mth{M}}$ est un sous-groupe de Levi de ${\mth{G}}$ et $O'$ une $G$-orbite unipotente anisotrope de ${\mth{M}}$. Cet ensemble est en bijection canonique avec l'ensemble des $G$-orbites unipotentes de ${\mth{G}}$: à une telle $G$-orbite $O$, on associe un sous-groupe de Levi ${\mth{M}}$ minimal rencontrant $O$, et $O'=O\cap{\mth{M}}$; de plus, on déduit du lemme \ref{infgf} que la matrice:
\[\left(\card\left(O'{\mth{U}}\cap O\right)\right)_{O,O'},\]
où $O'$ décrit l'ensemble des seconds membres des éléments de $\Lambda$, est triangulaire supérieure, à coefficients diagonaux non nuls puisqu'il est clair que si $O'\subset O$, $\card\left(O'{\mth{U}}\cap O\right)\neq 0$. Elle est donc inversible, et on peut alors écrire $D|_{\mathcal{H}_K}$ sous la forme:
\[D|_{\mathcal{H}_K}=\sum_{\left({\mth{M}},O'\right)\in R_K}d_{O',{\mth{G}},D}\Ind_{\mth{M}}^{\mth{G}}1_{O'}.\]
On en déduit, si $\mathcal{U}_{K,anis}$ est l'ensemble des $G$-orbites unipotentes anisotropes de ${\mth{G}}$:
\[D|_{\mathcal{H}_{K,cusp}}=\sum_{O'\in \mathcal{U}_{K,anis}}d_{O',{\mth{G}},D}1_{O'}.\]
Notons que si $D=J_U$, où $U$ est une orbite unipotente de $G$, l'égalité ci-dessus est valable pour la restriction de $D$ à l'espace $\mathcal{C}_{K,0,cusp}$ des éléments cuspidaux de $\mathcal{C}_{K,0}$ tout entier.

Pour toute $O'\in \mathcal{U}_{K,anis}$, on posera $d_{O',U}=d_{O',{\mth{G}},J_U}$; on obtient alors une matrice carrée $A=\left(d_{O',U_{O''}}\right)_{O',O''}$, où $O'$ et $O''$ décrivent l'ensemble des seconds membres des éléments de $\Lambda$, qui est inversible d'après le corollaire \ref{ioli}.

Soit maintenant $K_1,\dots,K_t$ les sous-groupes parahoriques maximaux de $G$ contenant $I$.
Soit $f\in\mathcal{H}_0$;
posons:
\[f=\sum_{i=1}^tf_i,\]
où pour tout $i$, $f_i$ est un élément de $\mathcal{H}_{K_i,0}$.
Pour tout $i$, soit $R_i$ un système de représentants des classes d'association de sous-groupes parahoriques $K'$ de $G$ tels que $I\subset K'\subset K$; d'après le lemme \ref{dchk}, on peut écrire, de façon unique modulo $C_0\left(G\right)$:
\[f_i=\sum_{K'\in R_i}f_{i,K'},\]
avec pour tout $K'$, $f_{i,K'}\in\mathcal{H}_{K',cusp}$. Pour tout $\mu\in\mathcal{M}$, posons:
\[f_\mu=\sum_{K'\in\kappa^{-1}\left(\mu\right)}\sum_{i,K'\in R_i}f_{i,K'}.\]
Il est clair que l'on a:
\[f=\sum_{\mu\in\mathcal{M}}f_\mu;\]
d'autre part, on déduit du lemme \ref{find} que pour tout sous-groupe parahorique $K'$ de $G$ contenant $I$, si $\mu=\kappa^{-1}\left(K'\right)$, $f_\mu\in\mathcal{H}_{K',cusp}+C_0\left(G\right)$. Fixons donc, pour tout $\mu$, un sous-groupe parahorique $K_\mu$ de $G$ contenant $I$ et tel que $\kappa\left(K_\mu\right)=\mu$; quitte à rajouter à $f$ un élément de $\mathcal{H}\cap C_0\left(G\right)$, on pourra supposer $f_\mu\in\mathcal{H}_{K_\mu,cusp}$ pour tout $\mu$.

Pour tout $\left(\mu,O'\right)\in\Lambda$, écrivons:
\[J_{U_{O'}}\left(f\right)=\sum_{\mu'\in\mathcal{M}}J_{U_{O'}}\left(f_{\mu'}\right)=\sum_{\mu'\in\mathcal{M}}\sum_{O''\in\mathcal{U}_{K_{\mu'},anis}}d_{O'',O'}1_{O''}\left(f_{\mu'}\right).\]
Pour tout $\left(\mu',O''\right)\in\Lambda$, considérons la distribution $D_{O''}$ sur $\mathcal{H}'_0$ définie par:
\[D_{O''}|_{\mathcal{H}_{K_{\mu''},cusp}}=1_{O''}|_{\mathcal{H}_{K_{\mu''},cusp}}\]
pour tout $\mu''$ tel que $K_{\mu''}$ est contenu dans $K_{\mu'}+C_0\left(G\right)$, et:
\[D_{O''}|_{\mathcal{H}_{K_{\mu''},cusp}}=0\]
pour tout $\mu''$ ne vérifiant pas cette condition; on obtient alors:
\[J_{U_{O'}}\left(f\right)=\sum_{\left(\mu',O''\right)\in\Lambda}d_{O'',O'}D_{O''}\left(f\right).\]
Or on sait que la matrice des $d_{O'',O'}$ est inversible;
on en déduit $D_{O''}\left(f\right)=0$ pour tout $O''$, d'où, par une récurrence évidente,
pour tout $(\mu',O'')\in\Lambda$, $1_{O''}\left(f_{\mu'}\right)=0$, d'où $D\left(f_{\mu''}\right)=0$ pour tout $\mu''$, donc $D\left(f\right)=0$, ce qui démontre le théorème.
\end{proof}

\begin{cor}\label{intuncor}
Si $G$ vérifie les conditions du théorème précédent, les espaces engendrés par les restrictions à $\mathcal{H}$ des distributions intégrales orbitales respectivement sur les éléments semi-simples réguliers non ramifiés de $G$ et sur les éléments unipotents de $G$ sont égaux.
\end{cor}

\begin{proof}
En effet, d'après le théorème précédent et le théorème \ref{intorb}, ces deux espaces sont égaux à $\mathcal{D}_{c,1}|_\mathcal{H}$.
\end{proof}

\end{document}